\documentclass[oneside]{book}

\usepackage{amsmath,amsthm,amssymb,amscd,amsfonts,amsopn}
\usepackage{latexsym}
\usepackage{mathabx}
\usepackage[crop=pdfcrop,process=all,cleanup={.tex,.dvi,.ps,.pdf,.log}]{pstool}
\usepackage{latexsym}
\usepackage{stmaryrd}
\usepackage{tikz-cd}
\usepackage{enumerate}
\usepackage[small,nohug,heads=vee]{diagrams}
\diagramstyle[labelstyle=\textstyle]
\usepackage{bbm}
\usepackage{centernot}
\usepackage{mathdots}
\usepackage{epstopdf}
\usepackage{tikz}
\usepackage{psfrag}
\usepackage{graphicx}
\usepackage{euscript}
\usepackage[all]{xy}
\usepackage{appendix}
\usepackage[ansinew]{inputenc}
\usepackage[T1]{fontenc}
\usepackage{rotating}
    \newenvironment{dedication}
        {\vspace{6ex}\begin{quotation}\begin{center}\begin{em}}
        {\par\end{em}\end{center}\end{quotation}}
\usetikzlibrary{decorations.markings}

\newcommand{\vertex}{\node[vertex]}

\newcommand{\FR}{\FF \mathcal{R}}
\newcommand{\CFR}{C\FF \mathcal{R}}
\newcommand{\FRtc}{\FR_{\text{tot-com.}}}
\newcommand{\FRone}{\FR_{\text{1-com.}}}
\newcommand{\FRx}{\FR_{\text{$\times$-com.}}}
\newcommand{\FRce}{\FR_{\text{cent'l}}}

\newcommand{\Op}{\mathcal{O}_{\kk,\eta}^{1/p}}
\newcommand{\Opsig}{\mathcal{O}_{\kk,\eta}^{\sigma}}
\newcommand{\Optwo}{\mathcal{O}_{\kk,\eta}^{1/2}}

\newcommand{\CG}{\mathcal{C}G}
\newcommand{\catrig}{\mathcal{R}ig}
\newcommand{\catrigs}{\mathcal{R}igs}
\newcommand{\catring}{\mathcal{R}ing}
\newcommand{\catrings}{\mathcal{R}ings}
\newcommand{\Mon}{\mathcal{M}on}
\newcommand{\CMon}{C\mathcal{M}on}
\newcommand{\slfrac}[2]{\left.#1\middle/#2\right.}

\newcommand{\FS}{\;\mathfrak{g}\FF\mathcal{S}c\;}
\newcommand{\FRS}{\;\FF\mathcal{RS}p\;}
\newcommand{\LFRS}{\;\mathcal{L}\FF\mathcal{RS}p\;}

\newcommand{\FSt}{\;\mathfrak{g}\FF\mathcal{S}c^t\;}
\newcommand{\LFRtS}{\;\mathcal{L}\FF\mathcal{R}^t\mathcal{S}p\;}
\newcommand{\frsp}{\;\FF\text{-}ringed\text{-}space\;}
\newcommand{\FRt}{\;\FF\mathcal{R}^t}
\newcommand{\frsps}{\;\FF\text{-}ringed\text{-}spaces\;}
\newcommand{\lfrsp}{\;\FF\text{-}locally\text{-}ringed\text{-}space\;}
\newcommand{\lfrsps}{\;\FF\text{-}locally\text{-}ringed\text{-}spaces\;}

\newcommand{\Amod}{A\text{-mod}}
\newcommand{\Der}{\mathcal{D}er}
\newcommand{\gFRst}{\;\mathfrak{g}\FF\mathcal{R}^{(t)}\mathcal{S}c\;}

\newcommand{\cG}{\mathcal{G}}
\newcommand{\cR}{\mathcal{R}}

\newcommand{\Amd}{A\text{-module}}

\newcommand{\be}{\begin{enumerate}}
\newcommand{\ee}{\end{enumerate}}
\newcommand{\FAT}[1]{\mbox{{$\mathbb{#1}$}}}
\newcommand{\Fat}[1]{\mbox{{$\scriptstyle\mathbb{#1}$}}}

\newcommand{\ZZ}{\FAT{Z}}
\newcommand{\NN}{\FAT{N}}

\newcommand{\QQ}{\FAT{Q}}

\newcommand{\RR}{\FAT{R}}

\newcommand{\A}{\mathbb{A}}

\newcommand{\FF}{\FAT{F}}
\newcommand{\CC}{\FAT{C}}

\newcommand{\zz}{\Fat{Z}}
\newcommand{\nn}{\Fat{N}}
\newcommand{\qq}{\Fat{Q}}

\newcommand{\kk}{\Bbbk}

\renewcommand{\bar}[1]{\overline{#1}}
\newcommand{\IFF}{\Leftrightarrow}
\newcommand{\THEN}{\Rightarrow}

\newcommand{\inj}{\hookrightarrow}
\newcommand{\sur}{\twoheadrightarrow}
\newcommand{\iso}{\xrightarrow{\sim}}
\newcommand{\fe}{\varphi}
\newcommand{\fr}{\FF\text{-}\mathcal{R}\text{ing}}

\newcommand{\frs}{\FF\text{-}\mathcal{R}\text{ings}}
\newcommand{\Cfrs}{C\FF\text{-}\mathcal{R}\text{ings}}
\newcommand{\catfr}{\;\FF\text{-}\mathcal{R}ings\;}

\newcommand{\Set}{$Set$_{\bullet}}

\newcommand{\KER}{\mathcal{KER}}
\newcommand{\loc}[1]{S^{-1}#1}
\newcommand{\g}[1]{\mathfrak{#1}}

\newcommand{\Def}{\stackrel{\text{def}}{=}}

\newcommand{\en}[1]{\eqno (#1)}
\newcommand{\minus}{\smallsetminus}

\newcommand{\lem}[1]{\subsection*{Lemma #1}}
\newcommand{\thm}[1]{\subsection*{Theorem #1}}
\newcommand{\cor}[1]{\subsection*{Corollary #1}}
\newcommand{\cl}[1]{\subsection*{Claim #1}}
\newcommand{\prop}[1]{\subsection*{Proposition #1}}
\newcommand{\pro}[1]{\subsection*{Proposition}}

\newcommand{\defin}[1]{\subsection*{Definition #1}}
\newcommand{\exmpl}[1]{\subsection*{Example #1}}
\newcommand{\rem}[1]{\subsection*{{\bf Remark #1}}}
\newcommand{\property}[1]{\subsection*{Property #1}}

\def\Ker{\operatorname{Ker}}
\def\Coker{\operatorname{Coker}}
\def\Spec{\operatorname{Spec}}
\def\Hom{\operatorname{Hom}}
\def\Sup{\operatorname{Sup}}
\def\Ht{\operatorname{ht}}
\def\Vec{\operatorname{Vec}}
\def\Val{\operatorname{Val}}

\def\Max{\operatorname{Max}}

\def\cK{\mathcal{K}\mathcal{E}\mathcal{R}}

\newcommand{\sqa}{\mathfrak{a}}
\newcommand{\sqb}{\mathfrak{b}}
\newcommand{\sqp}{\mathfrak{p}}
\newcommand{\sqq}{\mathfrak{q}}
\newcommand{\sqm}{\mathfrak{m}}
\newcommand{\cO}{\mathcal{O}}
\newcommand{\cX}{\mathcal{X}}
\newcommand{\cZ}{\mathcal{Z}}
\newcommand{\cY}{\mathcal{Y}}
\newcommand{\cC}{\mathcal{C}}
\newcommand{\cL}{\mathcal{L}}
\newcommand{\cN}{\mathcal{N}}
\newcommand{\cS}{\mathcal{S}}
\newcommand{\Fbb}{\mathbb{F}}
\newcommand{\smallsqp}{\tiny \mathfrak{p}}
\newcommand{\smallsqq}{\tiny \mathfrak{q}}
\newcommand{\smallsqm}{\tiny \mathfrak{m}}

\def\vsp{\vspace{6mm}}
\def\vsm{\vspace{3mm}}
\def\vp{\vspace{1.5mm}}

\def\hfar{\hspace*{45mm}}

\def\hsf{\hspace*{15mm}}
\def\hsp{\hspace*{6mm}}
\def\hsm{\hspace*{3mm}}
\def\sp{\hspace*{1.5mm}}
\addtocounter{chapter}{-1}

\begin{document}

\title{\textbf{New Foundations for Geometry} \\ Two non- additive languages for arithmetical geometry}
\let\cleardoublepage\clearpage
\author{M.J.\;\;Shai Haran}
\maketitle
\let\cleardoublepage\clearpage
\tableofcontents
\bigskip
\let\cleardoublepage\clearpage
\chapter{Introduction}
\bigskip

We spend our first years in the world of mathematics doing
addition (using our fingers), and then we learn of multiplication
(of natural numbers) as a kind of generalized addition. There is
no wonder that the vast majority of structures in mathematics
begin with addition, either explicitly (as an abelian group), or
abstractly (as additive/abelian categories). The language of
Grothendieck's algebraic geometry is based on commutative rings
having addition and  multiplication. But when we compare
arithmetic and geometry, we see that it is precisely the presence
of addition in our language that causes all the problems. The ring
of integers $\mathbb Z$ is similar to the polynomial ring in one
variable $F[x]$ over a field $F$. Taking for simplicity $F$
algebraically closed, we have analogous diagrams of embeddings of
rings in arithmetic and geometry:

$$\begin{tikzcd}
\underline{Arithmetic}&&\\
&\ZZ_p\arrow[hook]{r}&\QQ_p\\
\ZZ\arrow[hook]{ru}\arrow[hook]{rd}&& \\
 &\QQ\arrow[hook]{r}\arrow[hook]{uur}&\QQ_{\eta}=\RR\\
&&\ZZ_{\eta}=[-1,1] \arrow[hook]{u}
\end{tikzcd}
\begin{tikzcd}
\underline{Geometry}&&\\
&F[[x-\alpha]]\arrow[hook]{r}&F((x-\alpha)) \\
F[x]\arrow[hook]{ru}\arrow[hook]{rd}&& \\
&F(x)\arrow[hook]{r}\arrow[hook]{uur}&F((\frac{1}{x}))\\
&&F[[\frac{1}{x}]]\arrow[hook]{u}
\end{tikzcd} $$

\vspace{10pt}

Here the rational numbers $\mathbb Q$
are analogous to the field of rational functions $F(x)$.
The field of $p$-adic numbers ${\mathbb Q}_p$
is analogous to the field of Laurent series $F((x-\alpha))$.
The $p$-adic numbers contain the (one dimensional local) ring of $p$-adic integers
${\mathbb Z}_p = \lim\limits_{\stackrel{\longleftarrow}{n}}{\mathbb Z}/(p^n)$,
which is analogous to the (one dimensional local) ring of power series
$F[|x-\alpha|] = \lim\limits_{\stackrel{\longleftarrow}{n}}
   F[x]/(x-\alpha)^n$.\\
The embedding $\mathbb Q \hookrightarrow \mathbb Q_p$
is analogous to the embedding $F(x) \hookrightarrow F((x-\alpha))$,
i.e. of expanding a rational function as a Laurent series around the point $\alpha$. \\
There are three basic problems where this analogy breaks down. \vspace{3mm}

{\bf The problem of the \emph{arithmetical plane}:}
In geometry when we have two objects we have their product.
In particular, the product of the (affine) line with itself gives us the (affine) plane. This translates in the language of commutative rings into the fact that
\[
F[X] \bigotimes\limits_F F[X] = F[X_1, X_2]
\]
the polynomial ring in two variables. When we try to find the
analogous \emph{arithmetical plane}, we find
\[
\mathbb Z \otimes \mathbb Z = \mathbb Z
\]
Having addition as part of the structures of a commutative ring forces the integers $\mathbb Z$
to be the initial object of the category of commutative rings, hence its categorical sum with itself reduces to $\mathbb Z$,
and the \emph{arithmetical plane} reduces to its diagonal.

\vspace{10pt}

{\bf The problem of the \emph{absolute point}:}
The category of $F$-algebras has $F$ as an initial object,
hence in geometry (over $F$) we have the point $spec(F)$
as a final object. Addition forces the integers $\mathbb Z$
to
be the first object of commutative rings,
hence $spec(\mathbb Z)$ is the final object of Grothendieck's
algebraic geometry, and we are missing the absolute point
\emph{$spec(\mathbb F)$},
where $\mathbb F$ is the
\emph{field with one element} - the non-existing common field of all finite fields $\mathbb F_{p}$, $p$ prime.

\vspace{10pt}

{\bf The problem of the real prime:}
In geometry over $F$ we realize that if we want to have global
theorems we need to pass from affine to projective geometry.
In particular, we have to add the point at infinity $\infty$
to the affine line to obtain the projective line.
In our language of commutative rings this translates into the
extra embedding $F(x) \hookrightarrow F((\frac{1}{x}))$,
i.e. expanding a rational function as a Laurent series at infinity.
There is nothing special about the point $\infty$,
all the points of the projective line are the same,
and the field $F((\frac{1}{x}))$ is isomorphic to each of the fields $F((x-\alpha))$,
in particular: $F((\frac{1}{x}))$ contains the
(one-dimensional local) ring of power series $F[|\frac{1}{x}|]$.
The analog of the infinite point $\infty$ in arithmetic is the
\emph{real prime},
which we denote by $\eta$,
and the analog of the extra embedding $F(x) \hookrightarrow F((\frac{1}{x}))$
is the embedding $\mathbb Q \hookrightarrow \mathbb Q_{\eta} = \mathbb R$
of the rational numbers in the reals. The analog of the power-series subring $F[|\frac{1}{x}|] \subseteq F((\frac{1}{x}))$,
is the \emph{ring of real integers}
$\mathbb Z_{\eta}=[-1, 1] \subseteq \mathbb R$.
But the real interval $[-1, 1]$
is not closed under addition,
and is not a commutative ring.
The language of Grothendieck's algebraic geometry cannot
see the real prime, hence cannot produce global theorems in
arithmetic (this is the point of Arakelov's geometry where
the real prime $\eta$ is added to $spec(\mathbb Z)$
in an ad hoc way).\vspace{3mm}

All this brings us to the inevitable conclusion,
that if we want to see arithmetic as a true geometry,
we have to change the language of geometry,
and moreover, we have to give up addition as part of this language.\\
Kurokawa, Ochiai, and Wakayama \cite{KOW} were the first to
suggest simply abandoning addition, and work instead with the
language of multiplicative monoids.
This approach of using (multiplicative) monoids
was further developed by Deitmar \cite{De},
and indeed is the \emph{minimal} concept included in all other approaches, as it will be in our approach, cf. \S 2.4.
 But this approach creates many
new problems: the spectra of monoids always looks like the
spectra of a local ring
(the non-invertible elements of a monoid are the unique maximal
ideal),
and the primes of the (multiplicative) monoid $\mathbb Z$
are arbitrary subsets of the (usual) primes.
What we need is another operation that will replace addition. \vspace{3mm} \\

The very same problem of the inadequacy of addition appears in physics in
the theory of relativity:
the interval of speeds $(-c,c)$,
$c$ being the speed of light, is not closed under addition.
Einstein's solution, from which all of the theory of special relativity can be deduced, is to change addition into $c$-addition given by
\[
x+_c y = \frac{x+y}{1-\frac{x \cdot y}{c^2}}
\]
Like addition, this operation is associative, commutative, $0$ is the unit, and $-x$ is the inverse:
\[
\begin{array}{l}
(x+_c y) +_c z = x+_c (y+_c z) \ ,\\
x+_c y = y +_c x \ , \\
x+_c 0 = x \ , \\
x+_c (-x) = 0
\end{array}
\]
There is also a kind of distributive law:
\[
z \cdot (x+_c y) = (z \cdot x) +_{|z|\cdot c} (z \cdot y)
\]
This approach of Einstein's is not useful in arithmetic where we have also \underline{complex} primes
(e.g. the unique prime of ${\mathbb Z}[i]$ over $\eta$).
For complex primes, the real interval, ${\mathbb Z}_{\eta} = [-1, 1]$,
is replaced by the complex unit ball
${\mathbb D}_{\eta} = \left\{z \in {\mathbb C} \ , \ |z| \leq 1 \right\}$,
and the complex-$c$-addition, under which
$c \cdot {\mathbb D}_{\eta}^{\circ} =
\left\{ z \in {\mathbb C} \ , \ |z| < c \right\}$
is closed, is given by (the non-associative non-commutative) operation:
\[
x+_c y = \frac{x+y}{1+\frac{x \cdot \bar{y}}{|c|^2}} \ , \
\mbox{($\bar{y}$ is the complex conjugate of $y$).}
\]
Perhaps also Nature was trying to tell us one of her secrets, when Heisenberg found matrix multiplication as the basic language of the microscopic world...  \vspace{3mm} \\
The hint comes from the mysteries of the real prime.
We can change the underlying additive group we use to represent
a ring, $\mathbb G_a (B)= Hom(\mathbb Z[x], B)= B$,
to the general linear groups $GL_n(B)$.
Then the analog for the real prime of the (maximal compact) subgroup $GL_n(\mathbb Z_p)\subseteq GL_n(\mathbb Q_p)$
is the orthogonal group - the (maximal compact) subgroup
$"GL_n(\mathbb Z_{\eta})" = O_n \subseteq GL_n(\mathbb Q_{\eta})=
GL_n(\mathbb R)$.
And for a complex prime of a number field it is the unitary group
$U_n \subseteq GL_n(\mathbb C)$.
Indeed, Macdonald \cite{Mac} gives a q-analog interpolation
between the zonal-spherical-functions on
$GL_n(\mathbb Q_p)/GL_n(\mathbb Z_p)$,
and the zonal-spherical-functions on
$GL_n(\mathbb R)/O_n$
 and $GL_n(\mathbb C)/U_n$.
Similarly, there is a q-analog interpolation between the
 zonal-spherical-function on the p-adic Grassmanian
$GL_n(\mathbb Z_p)/B_{n_1, n_2}$,\\
$B_{n_1,n_2}= \left\{\left(\begin{array}{c|c}*&*\\ \hline  0&* \\ \end{array}\right) \in GL_n(\mathbb Z_p)\right\}$,
and the real and complex Grassmanians
$O_n/O_{n_1}\times O_{n_2}$
and $U_n/U_{n_1}\times U_{n_2}$,
$n_1 +n_2 = n$ (see \cite{H08}, \cite{O}).\vspace{5mm}\\

This point of view of the general-linear-group suggests also that
for the \emph{field with one element} $\mathbb F$ we have
$"GL_n(\mathbb F)"= S_n$, the symmetric group, which embeds as a
common subgroup of all
the finite group $GL_n(\mathbb F_p)$, $p$ prime (or the "field" $\FF\{\pm 1\}$, with "$GL_n(\FF\{\pm 1\})"=\{\pm 1\}^n\rtimes S_n$).\vspace{5mm}\\
But we need zero as a part of the language for geometry,
and so we represented a commutative ring $B$ by the collection
of all $m$ by $n$ matrices over $B$, $Mat_{m,n}(B)$,
for all $m,n$.\\
For the \emph{real integers} $\mathbb Z_{\eta}$,
$Mat_{m,n}(\mathbb Z_{\eta})$ are the matrices in
$Mat_{m,n}(\mathbb R)$ that carry the $n$ dimensional unit ball
into the $m$ dimensional unit ball.

$$\ZZ_{\eta}^n=\{ f \in \RR^n,\;\;\; \sum_{i=1}^n\;|f(i)|^2 \leq 1 \}$$

$$(\ZZ_{\eta})_{n,m}= Mat_{n,m}(\ZZ_{\eta})="\Hom(\ZZ_{\eta}^m,\ZZ_{\eta}^n)"=\{f\in Mat_{n,m}(\RR),\; f(\ZZ_{\eta}^m)\subseteq \ZZ_{\eta}^n\}.$$
There is a natural involution:
$$(,)^t:(\ZZ_{\eta}){_{n,m}\rightarrow (\ZZ_{\eta})}_{m,n}$$
which only exists when we considered the $l_2$ -norm. \\
The initial object $\mathbb F$ is represented by
$"Mat_{m,n}(\mathbb F)"$, the $m$ by $n$ matrices with entries
$0,1$ having at most one 1 in every row and column.
As a category $\mathbb F$ is equivalent to the category with objects
the finite sets, and with morphisms the partial bijections.\\
We also keep as part of our language, the operations of matrix
multiplication, as well as the operations of direct-sum and transposition
of matrices. \\
This language for geometry sees the real prime,
there is a natural compactification
$\overline{spec(\mathbb Z)}= spec(\mathbb Z) \cup \{\eta\}$
as a pro-object of the associated category of schemes.
The arithmetical plane does not reduce to its diagonal,
and yet one can do
\emph{algebraic-geometry-Grothendieck-style} over it.\vspace{3mm} \\

Recently, there have been a few approaches to "geometry over $\FF_1$", such as Borger $\cite{Bo09}$, Connes Consani $\cite{CC09}$,$\cite{CC14}$, Durov $\cite{Du}$, Lorscheid $\cite{Lo12}$, Soule $\cite{S}$, Takagi  $\cite{Tak12}$, T\"{o}en- Vaquie $\cite{TV}$ and Haran $\cite{H07}$ and $\cite{H09}$. \\
For relations between these see $\cite{PL}$. There are no inclusion relations between Haran's and other approaches, except that Durov's generalized rings are a subset of the $\frs$ of $\cite{H07}$.But Durov's use of monads forces him to use the $\ell_1$ -metric at the real primes, so that the unit ball $\ZZ_{\eta}^n$ is replaced by the polytop
\begin{equation*}
(\ZZ_{\eta}^n)_{\ell_1}=\{f\in \RR^n, \; \sum_{i=1}^{n}|f(i)|\leq 1 \}
\end{equation*}
and $O_n=GL_n(\ZZ_{\eta})$ is replaced by the (finite) subgroup of $GL_n(\RR)$ preserving this polytop. \\

The present language of $\frs$ is the same as $\cite{H07}$, \emph{except} that we omit the tensor- product of matrices from the structure: we use only matrix multiplication and direct- sum, and we add the involution to the structure. \\

An important observation of the present approach is that we do not need the tensor product to do geometry, and that the addition of an involution to the structure means that we have to work with the \emph{symmetric} spectrum. We also analyze the notion of commutativity more carefully, and define the "commutative-$\frs$" over which we can do geometry. The $\frs$ of $\cite{H07}$ are the subset of "totally- commutative- $\frs$" of the present approach. We show the arithmetical surface (the categorical- sum of $\ZZ$ with itself in the category of commutative- $\frs$) does not reduce to its diagonal, while in categories of totally- commutative $\frs$ \cite{H07}, or Durov's \cite{Du}, as in ordinary commutative rings it does reduce to its diagonal: $\ZZ\otimes \ZZ= \ZZ$. 
\vspace{3mm} \\

We observe though that the geometric object $\Spec(A)$, only depend on the operads $\{A_{1,X}\}$ and $\{A_{X,1}\}$, (and these can be identified in the presence of an involution). We therefore axiomatize the properties of the "self-adjoint operad" $\{A_{1,X}\}$, for an $\fr$ with involution $A$. This gives us the "Generalized- Rings" of \cite{H10}, the geometry of which was developed in \cite{H09}. But in $\cite{H09}$ we assumed our generalized rings to be totally commutative, and self- adjoint, which are unnatural and limiting assumptions. Here we avoid these assumptions, and show that using the language of commutative- generalized- Rings, one can do "algebraic- geometry- Grothendieck- style". It includes "classical" algebraic geometry (fully- faithfully), and yet it solves the three basic problems of the analogy of arithmetic and geometry, as there are: \\
\begin{enumerate}[-]
\item A final object for Geometry, the "absolute point": $\Spec \FF$. \\
\item natural compactifications $\bar{\Spec(\ZZ)}=\Spec(\ZZ)\cup \{\eta\}$, and similarly for a number field $K$,
 $$\bar{\Spec(\mathcal{O}_K)}=\Spec(\mathcal{O}_K)\cup \{\eta_i\}_{\eta:K\inj {\scriptsize \CC}\text{ mod conj'n}},$$
as pro- object of the associated category of Grothendieck generalized schemes.
\item The arithmetical plane does not reduce to the diagonal. \vspace{3mm}\\
Unfortunately, it seems that to understand the geometric theory of generalized- Rings, one has first to go through the same theory using $\frs$ with involution. Especially, the role of the symmetric spectrum $\Spec^t (A)$, (as opposed to the usual spectrum $\Spec (A)$) as the basic building blocks for schemes of $\frs$ with involution, and for schemes of generalized rings (non self- adjoint as in $\cite{H09}$). \vspace{3mm}\\
\end{enumerate}
The contents of the chapters are as follows: \vspace{3mm}\\
\S1.1 We define "the field with one element $\FF$" (cf. \cite{S},\cite{H07}) to be the category of finite sets and partial bijections.
$$\FF(X,Y)=\{f:D(f)\xrightarrow{\sim} I(f)\sp\text{ bijection}\sp, D(f)\subseteq X,I(f)\subseteq Y\} $$
and let $\oplus$ be the disjoint union on this category.
$$\bigoplus:\FF\times \FF\rightarrow \FF$$
$$X_0\oplus X_1=\{(i,x), i\in \{0,1\}, x\in X_i\}\sp ,$$ \vp
$$(f_0\oplus f_1)(i,x)=(i,f_i(x)), f_i\in \FF_{Y_i,X_i}\sp $$
We have the associativity, commutativity, unit isomorphisms:
$$a=a_{X_0,X_1,X_2}:(X_0\oplus X_1)\oplus X_2\xrightarrow{\sim} X_0\oplus (X_1\oplus X_2)$$
$$c=c_{X_0,X_1}: X_0\oplus X_1\xrightarrow{\sim} X_1\oplus X_0$$
$$u=u_X: X\oplus [0]\xrightarrow{\sim} X\vspace{3mm}$$
We have an involution on $\FF$,
$$(\;)^t:\FF\iso \FF^{op}$$
$$(f:D(f)\iso I(f))^t=(f^{-1}:I(f)\iso D(f)).$$  \vspace{1.5mm}\\
We shall assume the objects of $\FF$ form a (countable) set, containing $[n]=\{1,...,n\}, n\geq 0$. \vspace{5mm}

\S1.2 We define an $\fr$ (cf. \cite{H07}) to be a symmetric
monoidal category over $\FF$: An $\fr ~A$ is a category, together
with a faithful functor $\epsilon:\FF\rightarrow A$, which is the
\underline{identity on objects}, and a symmetric monoidal
structure
$$\bigoplus:A\times A\rightarrow A$$
such that $\epsilon$ is strict monoidal; thus $A$ has objects the
finite sets, for $X,Y\in A$, $X\bigoplus Y$ is the disjoint union,
and $\bigoplus$ on arrows has associativity(resp.
commutativity,unit) isomorphism given by $\epsilon(a)$
(resp. $\epsilon(c)$, $\epsilon (u)$). We also demand that $[0]$ is the initial and final object of $A$.\vspace{3mm} \\
We write the arrows in $A$ from $X$ to $Y$ in "matrix" form $A(X,Y):=A_{Y,X}$.\vspace{3mm} \\
An $\fr$ with involution is an $\fr$ $A$ with a functor

$$\begin{tikzcd}
A\arrow{r}[swap]{\sim}{(\;)^t}&A^{op} \\
\FF\arrow{r}{\sim}\arrow{u}& \FF^{op}\arrow{u}
\end{tikzcd}$$
such that $(a^t)^t=a$, and $(a_0\oplus a_1)^t=a_0^t\oplus a_1^t$. \vspace{1.5mm} \\

\S1.3 We discuss the notion of commutativity for an $\fr$. An $\fr$ is commutative if the following condition holds:
$$\forall a\in A_{Y,X} ,  b\in A_{1,J} ,  d\in A_{J,1}:$$
$$ a\circ \underset{X}{\bigoplus}( b\circ d)= \underset{Y}{\bigoplus}(b\circ d)\circ
a=(\underset{Y}{\bigoplus} b)\circ (\underset{J}{\bigoplus}
a)\circ (\underset{X}{\bigoplus} d)\in A_{Y,X}.$$
It is totally commutative if
$$\forall a\in A_{Y,X} , \forall b\in A_{J,I}\;, \;(\underset{J}{\bigoplus} a)\circ (\underset{X}{\bigoplus} b)=
(\underset{Y}{\bigoplus} b)\circ (\underset{I}{\bigoplus} a)\in A_{Y\otimes J,X\otimes I}.$$  \vspace{3mm}\\

\S1.4 For an $\fr$ $A$ we have a mapping
$$A(X,Y)=A_{Y,X}\rightarrow (A_{1,1})^{Y\times X}$$
$$a\mapsto a_{y,x}=j^t_y\circ a\circ j_x$$
where $(j_x:\{1\}\iso \{x\})\in \FF_{X,[1]}$. Although for most of the examples this mapping is an injection (we say $A$ is a "matrix - $\frs$"), it need not be in general: for the "residue field" at the real prime it is not an injection! a replacement is the notion of tame $\fr$.  \vspace{3mm}\\

\S A.1 We show that the category of $\frs$ has push-outs.  \vspace{1.5mm}\\

\S A.2 We discuss equivalence - ideals and quotients.  \vspace{3mm}\\

\S2.1 We show the category of rings is (fully) embedded in $\frs$\\ ($R\mapsto \FF(R)$). \vspace{1.5mm} \\

\S2.2 We show the category of monoids is (fully) embedded in $\frs$ ($M\mapsto \FF\{M\}$).  \vspace{1.5mm}\\

\S2.3 The category of finite sets, $Set$; its opposite, $Set^{op}$; the category of finite sets and relations, $\mathcal{R}el (\supseteq\; Set,Set^{op})$, are all examples of $\frs$. \vspace{1.5mm} \\

\S2.4 For every $p\geq 1$ the sub -$\fr$ of real matrices $\FF(\RR)$ with operator $l_p$- norm $\leq 1$ is an $\fr$. For $p=2$ (and only $p=2$) it has involution.  \vspace{1.5mm}\\

\S2.5 We discuss valuation- $\frs$, and prove Ostrowski theorem that for a number field $K$ the valuation $\frs$ with involution correspond bijectively with the finite and infinite primes of $K$. (the proof itself is in the appendix B.1.).  \vspace{1.5mm}\\

\S2.6 We show the finite directed graphs with no loops form an $\fr$. \vspace{1.5mm} \\

\S2.7 We construct the free - $\fr$ on one generator of "degree" $(Y,X)$.  \vspace{1.5mm}\\

\S2.8 We construct the $\fr$ representing the functor $A\mapsto GL_n (A)$ .
\vspace{3mm}\\

\S2.9 We consider the arithmetical surface $\FF(\NN) \underset{\FF}{\otimes} \FF(\NN)$. Its totally commutative quotient is reduced to the diagonal, but we prove that its commutative quotient is Not reduced to the diagonal!  \vspace{1.5mm}\\

\S2.10 The $\fr$ $\FF(\NN)$ (respectively, $\FF(R)$ for a ring $R$) is generated by the matrices $(1,1)$ and $\left( \begin{matrix} 1 \\ 1 \end{matrix}\right)$, (resp. and $(r),r\in R$). We give the precises relations satisfied by these matrices.  \vspace{3mm}\\

$\S$3 We show one can do algebraic geometry Grothendieck style (cf. \cite{EGA}) over any commutative $\fr$.  \vspace{3mm}\\

$\S$4 The point of this chapter can be explained with ordinary commutative rings. If we want to develop algebraic geometry for commutative rings $A$ with (a possibly non trivial) involution, we can forget the involution and consider $\Spec A$  - a topological space with an involution. But when we localize, or glue, such objects we lose the involution. The "right" way is to consider $\Spec A^+$ , where $A^+=\{a\in A\; |\; a^t=a\}$, over this space we have a sheaf of rings with involution.  \vspace{3mm}\\

$\S$5 The point of this chapter can be explained with ordinary schemes (= the locally ringed spaces which are locally affine). While the locally - ringed - spaces (and the affine schemes) are closed under inverse limits, schemes are not. Given a point $x=\{x_j\}\in \varprojlim X_j$, $X_j$ schemes, while each $x_j\in X_j$ has an affine neighborhood , these neighborhoods can shrink so that $x$ will not have an affine neighborhood. The real and complex primes of a number field are such points. We show that in the pro - category of schemes there exists the compactification of $\Spec \ZZ$, and $\Spec \mathcal{O}_K$, $K$ a number field. (this is the compactification of \cite{H07}, reproduced also in \cite{Du}).  \vspace{3mm}\\

$\S$6 We show that the process of real completion creating the continuum $\RR^+=GL_1(\RR)/\{\pm 1\}$ (and similarly $GL_n(\RR)/O_n$, $GL_n(\CC)/U_n$) can naturally be embedded in the language of pro- schemes. We define rank- $n$ vector bundle for an arbitrary pro- $\FF$ scheme in such a way that for the compactified $\Spec \mathcal{O}_K$, $K$ a number field, we obtain $GL_n(\mathbb{A}_K)/\underset{\mathfrak{p}}{\Pi} GL_n(\mathcal{O}_{K,\mathfrak{p}})\times O_n^{r_{\tiny \RR}}\times U_n^{r_{\tiny \CC}}.$ \vspace{3mm}\\

\S7.1 For any $\fr$ $A$, we define $\Amod$ as the full subcategory of the functor category $(Ab)^{A\times A^{op}}$ of $M$'s such that $M_{0,X}=M_{Y,0}=\{0\}$. The category of $\Amod$ is complete and co- complete abelian category with enough projective and injectives. (and similarly the category $\Amod^t$, of $A$- modules with involution, when $A$ has an involution).
 \vspace{1.5mm}\\

\S7.2 We discuss the notion of commutativity for modules.
\vspace{1.5mm} \\

\S7.3 We define the category of $\mathcal{O}_X\text{- mod}^{(t)}$ (possibly with involution) and its full subcategory $q.c.\; \mathcal{O}_X\text{- mod}^{(t)}$ of quasi- coherent $\mathcal{O}_X$- modules, so that for affine $X=\Spec A$ localization gives an equivalence $\Amod\iso q.c.\; \mathcal{O}_X\text{- mod}$. \vspace{1.5mm} \\

\S7.4 For every homomorphism of $\frs^{(t)}$, $\varphi:B\to A$,
there is an adjunction (analogous with the commutative algebra's
extension and restriction of scalars):
$$\Amod^{(t)}(M^A,N)\cong B\text{-mod}^{(t)}(M,N_B),$$
$$M\in B\text{- mod}^{(t)},\;\;\; N\in \Amod^{(t)}.\vspace{1mm}$$

\S7.5 We define the infinitesimal extension of  $A\in \FR$ with an $M\in \Amod$, $A\prod M$, which is an abelian group object in $\FR/A$.  \vspace{1.5mm} \\

\S7.6 We define the concept of derivations with values in $A$- modules and the module of K\"{a}hler differentials representing them. For $\varphi \in \FR^{(t)}(C,A), \; B\in C \diagdown \FR^{(t)}\diagup A$ we have the adjunction,

$$ \bigg(C \diagdown \FR^{(t)}\diagup A\bigg)(B,A\Pi M)\equiv \Der^{(t)}_C(B,M_B )\equiv \Amod^{(t)}(\Omega(B/C)^A,M).$$
 \vspace{1.5mm} \\

\S7.7 We list the properties of differentials (such as the first and second exact sequences). \vspace{1.5mm} \\

\S7.8 We give an explicit description of the modules
$\Omega(\slfrac{\FF(\ZZ)}{\FF\{\pm 1\}})$ and
$\Omega(\slfrac{\FF(\NN)}{\FF})$. We show there is an exact
sequence of $\FF(\ZZ)$- modules,
$$\Omega(\slfrac{\FF(\ZZ)}{\FF\{\pm 1\}})\xrightarrow{\partial} N \xrightarrow{\pi} \FF(\ZZ) \rightarrow 0,$$
where $N$ is also an $\fr$ and $\pi$ is a homomorphism of $\frs$ (with involution). \vspace{1.5mm} \\

\S7.9 We generalize the previous description to have an explicit description for  $\Omega (\slfrac{\FF(R)}{\FF\{S\}})$, $R$ a commutative $\mathcal{R}$ing and $S\subseteq R$, a multiplicative set. \vspace{1.5mm} \\

\S7.10 Here we sketch the application of Quillen's "non- additive homological algebra", or "homotopical algebra", to our context. We give Quillen model structures for modules and algebras and define the Quillen cotangent complex.\vspace{3mm}\\

$\S 8$ Here we begin our "de-je-vu", with $\frs$ (with involution) replaced by generalized-rings.
We give more emphasis to the definitions and principal examples, than to the proofs - which are very similar, and usually even simpler, than the corresponding proofs for $\frs$.\vspace{1.5mm}\\
\S 8.1 is devoted to the definition of a generalized- ring, and in \S8.2 we give important remarks. \vspace{1.5mm}\\
In \S 8.3 we associate with every ($\times$- commutative) $\fr$ with involution, a (commutative) generalized ring. In particular, we describe the generalized rings:\\

$\FF$- the initial object of the category of generalized rings $\cG\cR$, \S8.3.1 .\\

$\cG(B)$- the generalized ring associated with a commutative rig $B$, \S8.3.2  .\\

$\cO_{\tiny K,\eta}$- the generalized ring associated to an embedding $\eta: K \inj \CC$, \S8.3.3 . \\

Valuation- generalized- rings of a generalized- field are given in \S8.3.4, and we have Ostrowski's theorem that for a number field $K$, the valuations of $\cG(K)$ correspond to $\cG(\cO_{K,\mathfrak{p}})$, $\mathfrak{p}$ a finite prime, or to $\cO_{K,\eta}$, $\eta:K\inj \CC$ (mod conjugation). \vspace{3mm}\\

In \S8.3.5 we give the generalized ring $\FF\{M\}$ associated to a commutative monoid $M$. \vspace{1.5mm}\\

In \S8.3.6 we describe the free- commutative generalized- ring $\Delta^W$, such that for any commutative generalized- ring $A$ we have $\cG\cR(\Delta^W,A)\equiv A_W$.\vspace{3mm}\\

In $\S 9$ we give the various notions of ideals for a commutative generalized ring $A$, and the relations between them. \vspace{1.5mm} \\

In \S9.1 we give the equivalence ideals $eq(A)$; in \S9.2 the functorial ideals, $fun \cdot il(A)$; in \S9.3 the operations on $fun \cdot il(A)$; in \S9.4 we give the homogeneous functorial ideals $[1]\mbox{-}il(A)$; and finally in \S9.5 we give the useful notations of ideals $il(A)$, and symmetric ideals $il^t(A)$; generally, we have $il^t(A)\subseteq [1]\mbox{-}il(A)\subseteq il(A)$. For $A=\mathcal{G}(B)$, the generalized ring associated with a commutative ring $B$ with an involution $b\mapsto b^t, il(A)$ corresponds bijectively with the ideals of $B$, $[1]-il(A)$ corresponds to the ideals $\mathfrak{b}$ of $B$ that are invariant under the involution: $\mathfrak{b}=\mathfrak{b}^t$, and $il^t(A)$ corresponds to the ideals of $B$ that are generated as ideals by elements that are invariant under the involution $\mathfrak{b}=(b_j)_{j\in J}, b_j^t=b_j$.   
\vspace{3mm}\\

In \S10 we give the geometry of commutative- generalized- rings (with no restrictions such as self- adjointness or total- commutativity). We discuss maximal (symmetric) ideals, and (symmetric) primes, and we have contravariant functors associating to a commutative- generalized- ring $A$ its space of (symmetric) primes, $\Spec(A)$, (resp. $\Spec^t(A)$), with its compact sober Zariski topology. There is a continuous map $\Spec(A)\rightarrow \Spec^t(A)$. \vspace{3mm}\\

In \S11 we give the localizations of a commutative- generalized- ring $A$, and the associated sheaf $\cO_A$ of generalized rings over $\Spec^t(A)$. \vspace{3mm}\\

In \S12 we briefly describe the category $\cL\cG\cR\cS$ of locally-generalized- ringed- spaces, and its full subcategory of Grothendieck- generalized- schemes $\cG\cG\cS$. The category of generalized- schemes $\cG\cS$ is the pro- category of $\cG\cG\cS$: $\cG\cS=\text{pro-}\cG\cG\cS$. It contains the compactified $\bar{\Spec \cO_K}$, $K$ a number field.\vspace{3mm}\\

In \S13.1 we briefly describe the tensor- products of commutative- generalized- rings, and we write down the precise relations for the presentations
\begin{equation*}
\begin{split}
\Phi:&\Delta^{[2]}\sur \cG(\NN), \\
&\delta^{[2]}\mapsto (1,1),\\
\Phi_{\small \ZZ}:&\FF\{\pm 1\}\underset{\FF}{\otimes}\Delta^{[2]}\sur \cG(\ZZ), \\
\Phi_B:&\FF\{B^!\}\underset{\FF}{\otimes}\Delta^{[2]}\sur \cG(B), \; B\text{ a commutative ring}.
\end{split}
\end{equation*}\vspace{3mm}\\

In \S13.2 we give a description of $\cG(\NN)\underset{\FF}{\otimes}\cG(\NN)$. \\
The tensor product gives the products in $\cG\cG\cS$, \S13.3, and $\cG\cS$, \S13.4, and we give the basic special case of the \emph{compactified arithmetical plane} 
$$\bar{\Spec \ZZ}\underset{\FF\{\pm 1\}}{\otimes}\bar{\Spec \ZZ}.$$
To quote \cite{CC14}:
\begin{center}
"This note provides the algebraic geometric space underlying the non- commutative approach to RH. It gives a geometric framework reasonably suitable to transpose the conceptual understanding of the Weil proof in finite characteristic. This translation would require in particular an adequate version of the Riemann- Roch theorem in characteristic $1$". see \cite{H09} for more hints.\end{center}\vspace{3mm}
In \S14.1 we briefly describe the theory of $A$- modules, for $A$ a generalized ring, and its localization - the $\mathcal{O}_X$- modules, $X\in \mathcal{GGS}$.\vspace{3mm}\\

In \S14.2 we give the derivations and differentials (but only the
even or odd ones, where the involution on $M$ is the identity or
minus the identity). We give explicitly the basic
examples of the even derivations
\begin{equation*}
\begin{split}
d=\frac{1}{2}d^+_{[1]}:\ZZ\rightarrow \bar{\Omega}_{[1]}^{\ZZ} \\
\text{and } d=\frac{1}{2}d_{[1]}^+: \NN \rightarrow \bar{\Omega}_{[1]}^{\NN}
\end{split}
\end{equation*}
Here $\bar{\Omega}_{[1]}^{\NN}$ is the free abelian group with generators $\left\{\begin{matrix} a\\ a'\end{matrix} \right\}$, $a,a'\in \NN_{>0}$, modulo the relations: \\
\begin{equation*}
\begin{matrix}
\text{cocycle} & \left\{\begin{matrix} a+a'\\ a''\end{matrix} \right\}+\left\{\begin{matrix} a\\ a'\end{matrix} \right\}=
\left\{\begin{matrix} a\\ a'+a''\end{matrix} \right\}+\left\{\begin{matrix} a'\\ a''\end{matrix} \right\}\\
\text{symmetric} & \left\{\begin{matrix} a\\ a'\end{matrix} \right\}=\left\{\begin{matrix} a'\\ a\end{matrix} \right\}\\
\text{linear} & \left\{\begin{matrix} k\cdot a\\ k\cdot a'\end{matrix} \right\}=k\cdot \left\{\begin{matrix} a\\ a'\end{matrix} \right\}
\end{matrix}
\end{equation*}
and $d(n)= \left\{\begin{matrix} 1\\ 1\end{matrix} \right\}+ \left\{\begin{matrix} 2\\ 1\end{matrix} \right\}+\dots + \left\{\begin{matrix} n-1\\ 1\end{matrix} \right\}$ satisfies
\begin{equation*}
d(n+m)=d(n)+d(m)+ \left\{\begin{matrix} n\\ m\end{matrix} \right\}
\end{equation*}
and Leibnitz:
\begin{equation*}
d(n\cdot m)=n\cdot d(m)+m\cdot d(n) .
\end{equation*}
hence $d(n)=\sum_p v_p(n)\cdot \frac{n}{p}\cdot d(p)$ and $\bar{\Omega}_{[1]}^{\NN}$ is the free abelian group on $d(p)$, (or $\left\{\begin{matrix} p-1\\ 1\end{matrix} \right\}$), $p$ prime. \vspace{7mm} \\

In a final appendix, we make contact with the analytic theory, and we give yet another explanation to the fact that the Gamma function gives the real analogue of the Euler factor $(1-p^{-s})^{-1}$ (or equivalently, via Mellin transform, that the Gaussian gives the real analogue of the charateristic function of $\ZZ_p$. For one explanation, which goes via the "quantum" $q$- analogue interpolation between the real and $p$- adic worlds, see \cite{H01}, \cite{H08}). We show that the Haar- Maak $O(N)$ or $GL_N(\ZZ_p)$- invariant- probability measure $\sigma_N$ on the real or $p$- adic sphere $S_p^N\subseteq \QQ_p^N$, behaves naturally with respect to the operations of multiplication and contraction, giving rise to the Beta function, and in particular, contracting $\sigma_N$ with the vector $(1,\dots,1)\in \QQ_p^N$ we get in the limit $N\rightarrow \infty$,
\begin{equation*}
\lim_{N\rightarrow \infty} \int_{S_p^N}|x_1+\dots +x_N|^{s-1}\sigma_N(dx) = \frac{\zeta_p(s)}{\zeta_p(1)}
\end{equation*} 
with $\zeta_p(s)=(1-p^{-s})^{-1}$ for the $p$- adic numbers, $\zeta_{\eta}(s)=2^{\frac{s}{2}}\Gamma (\frac{s}{2})$ for the reals.\vspace{7mm}\\

Alas, our idea and message are very simple: If one changes the
operations of addition and multiplication to the more fundamental
operations of multiplication and contraction of vectors
(respectively, multiplication, direct- sum, and transposition of
matrices), then one obtains a language in which arithmetic can be
viewed as a geometry, and one proceeds exactly as in Grothendieck
(with Quillen's non- additive homotopical algebra).
\newpage
\begin{dedication}
\begingroup
\centering
\emph{This book would have not existed without the continuous efforts of Itai Cohen, who typed them into latex, and edited them.\\ It is dedicated to my father,\\ Prof. Menachem Haran,\\ and to my mother,\\ Dr. Raaya Haran- Twerski, \\ who taught us at an early stage that \\ "the trick is finding the simple explanation behind the complex phenomena" \\ and to the memory of\\ Daniel Quillen, \\ a teacher and a mentor, who taught us that mathematics is a language - \\ "when you have the right language to speak about a problem, you solved the problem".}
\endgroup
\vspace*{\fill}
\end{dedication}
\let\cleardoublepage\clearpage
\numberwithin{equation}{section}
\let\cleardoublepage\clearpage
\part{$\frs$}
\let\cleardoublepage\clearpage
\chapter{Definition of $\frs$.}
\let\cleardoublepage\clearpage
\bigskip
\setcounter{equation}{0}
\section{$\FF$ the field with one element}
\smallbreak
For a category $ C $ we write $ X \in C $ for ``$ X $ is an object of
$ C $'', and we let $ C(X,Y) $ denote the set of maps in $ C $ from $
X $ to $ Y $. We denote by $Set_0$ the category with objects sets
$X$ with a distinguished element $O_X\in X$, and with maps
preserving the distinguished elements
\begin{equation}    
 Set_0(X,Y) = \Big\{f\in Set(X,Y) \;,\; f(O_X)=O_Y\Big\}.
\end{equation}
The category $Set_0$ has direct and inverse limits. The set
$[o]=\big\{o\big\}$ is the initial and final object of $Set_0$. For
$f\in Set_0(X,Y)$ we have
\begin{equation}    
   ker f=f^{-1}(O_Y);
\end{equation}
\[
coker f=Y/f(X),
\]
the set obtained from $Y$ by collapsing $f(X)$ to a point.

There is a canonical map
\begin{equation}    
   coker \ ker f= X/f^{-1}(O_Y) \quad \longrightarrow \quad ker \ coker  f=f(X)  .
\end{equation}\vspace{3mm}\\
Our first instinct is to take for $\FF_0$, "the field with one element", or rather, the "finite dimensional $\FF_0$- vector spaces", the full subcategory of $Set_0$ consisting of the finite sets. But note that the map that identifies two points to one (non-zero) point, has $\RR$- linear extension the map $\RR\times \RR\rightarrow \RR$, $(x_1,x_2)\mapsto x_1+x_2$, and this map does not takes the unit $L_2$- disc into the unit interval ($(\frac{1}{\sqrt{2}},\frac{1}{\sqrt{2}})\mapsto \sqrt{2}>1$).  \vspace{3mm}\\
Thus we denote by $\mathbb{F}_0$ the subcategory of $Set_0$ with objects
the \underline{finite} pointed sets, and with maps
$$
\hspace{10mm} \FF_0(X,Y)  =
$$
\[=
\big\{ f \in Set_0(X,Y),\ coker \ ker f
                    \xrightarrow{\sim} ker \ coker  f
       \mbox{ is  an  isomorphism} \big\}=
\]
\[=
\big\{ f \in Set_0(X,Y),\; f|_{X \setminus
                  f^{-1}(O_Y)}\, \mbox{ is  an  injection} \big\}
\]

We let $ Set_{\bullet} $ denote the category with objects sets and with
partially defined maps
\begin{equation}    
   Set_{\bullet}(X,Y)= \underset{X' \subseteq X}{\perp\mkern-7mu\perp} Set(X',Y).
\end{equation}
Thus to $ f \in Set_{\bullet}(X,Y) $ there is associates its domain $ D(f)
\subseteq X $, and $ f \in Set(D(f),Y) $.\\
We have an isomorphism of categories

\[
\begin{array}{lrcl}
& \hspace{-0.7in}Set_0& \textstyle \stackrel{\thicksim }{ \longleftrightarrow }& Set_{\bullet}\\
\mbox{given by} & & & \\
&  X &\longmapsto & X_+ = X \smallsetminus \big\{ O_X \big\} \\
  &  f& \longmapsto& f_+ \;,\; D(f_+) = X \smallsetminus f^{-1}( O_Y );\\
\mbox{and inversly} & & & \\
 & \big\{ O_X \big\}  \coprod X  = X_0 &\mapsfrom& X\\
 \quad \begin{array}{r}x \in D(f):\\x \not \in D(f):
\end{array}&
       \left.\begin{array}{r}
f(x)\\
 O_Y
\end{array}\right\}
 = f_0 (x)& \mapsfrom& f
\end{array}
\]

We let $\FF$ denote the subcategory of $Set_{\bullet}$
corresponding to $\FF_0 $ under this isomorphism, it has objects the \underline{finite} sets, and maps are the
partial bijections
\begin{equation}    
   \FF(X,Y) = \big\{ f:D(f) \xrightarrow{\thicksim} I(f)
\big.  \ bijections,\  \big. D(f) \subseteq X,\ I(f) \subseteq Y
\big\} .
\end{equation}

It is crucial that the objects of $\FF$
are finite sets, but we do not need $\FF$
to contain all finite sets.

To avoid problems with set theory we shall work with a countable
set-model of $\FF$ that contains $ [0] = \varnothing $ (the
empty set, the initial and final object), $ [1] = \big\{ 1 \big\}, \
\ldots, \ [n] = \big\{ 1,\ldots,n \big\} \ ,\ldots $ and is closed
under the operations of pull-back and push-out in $Set_{\bullet}$.

The operation "$\circ$" will denote composition of partial - bijections, but note that if $g\circ f$ is defined,
than $D(g\circ f)=f^{-1}(I(f)\cap D(g))=D(f)\cap f^{-1}(D(g))$.

Note that we can identify $\FF(X,Y)$ with $Y\times X$ - matrices with value in $\{0,1\}$, having at most one $1$ in every row or column,
(and than $\circ$ is matrix multiplication), and we will denote this set by $\FF_{Y,X}$.

Note that the category $\FF$ has \underline{no} sums or products, but we do have two symmetric monoidal structures on $\FF$, the disjoint union:
\begin{equation}
\bigoplus:\FF\times \FF\rightarrow \FF
\end{equation}

\begin{equation}
X_0\oplus X_1=\{(i,x), i\in \{0,1\}, x\in X_i\}\sp,
\end{equation}

\begin{equation}
(f_0\oplus f_1)(i,x)=(i,f_i(x)), f_i\in \FF_{Y_i,X_i}
\end{equation}
(this is the categorical sum in $Set_{\bullet}$). \\
and the cartesian product:

\begin{equation}
\bigotimes:\FF\times \FF\rightarrow \FF
\end{equation}
\begin{equation}
X_0\otimes X_1=\{(x_0,x_1), x_i\in X_i\}.
\end{equation}
\begin{equation}
(f_0\otimes f_1)(x_0,x_1) = (f_0(x_0),f_1(x_1)), f_i\in \FF_{Y_i,X_i}.
\end{equation}
Thus in $Set_{\bullet}$
\[
X \bigotimes Y  = coker \big\{ X \perp\mkern-7mu\perp  Y \rightarrow X \top\mkern-7mu\top Y \big\}
\]
or $\displaystyle{X \top\mkern-7mu\top Y =
X \perp\mkern-7mu\perp Y \perp\mkern-7mu\perp ( X\bigotimes Y)}$. \\
We have associativity, commutativity, and unit isomorphisms:

\begin{equation}
a=a_{X_0,X_1,X_2}:(X_0\oplus X_1)\oplus X_2\xrightarrow{\sim} X_0\oplus (X_1\oplus X_2)
\end{equation}
\begin{equation}
a(1,x_2)=(1,(1,x_2))
\end{equation}
\begin{equation}a(0,(i,x_i))=\begin{cases}
   (1,(0,x_1)),&  i= 1\\
    (0,x_0)       & i=0
\end{cases}
\end{equation}
\begin{equation}
c=c_{X_0,X_1}: X_0\oplus X_1\xrightarrow{\sim} X_1\oplus X_0
\end{equation}
\begin{equation}
c(i,x)=(1-i,x)
\end{equation}
\begin{equation}
u=u_X: X\oplus [0]\xrightarrow{\sim} X
\end{equation}
\begin{equation}
u(0,x)=x.
\end{equation}
where $[0]=\text{the empty set}$, is the initial and final object of $\FF$. Similarly, there are associativity $a^*$, commutativity $c^*$,
and unit $u^*$ isomorphisms for the operation $\otimes$, (the unit object for $\otimes$ being $[1]=\text{the one point set}$).
Moreover, there is a distributivity isomorphism:
\begin{equation}
d=d_{X;Y_0,Y_1} :X\otimes (Y_0\oplus Y_1)\xrightarrow{\sim} (X\otimes Y_0)\oplus (X\otimes Y_1).
\end{equation}
Given a finite collection of finite sets $\{Y_x\}_{x\in X}$, we can form the disjoint union:
\begin{equation*}
\bigoplus_X Y_x=\{(x,y),\sp x\in X,\sp y\in Y_x\}
\end{equation*}
and we have canonical isomorphisms:
\begin{equation}
\bigoplus_X Y\xrightarrow{\sim} X\bigotimes Y\xrightarrow{\sim} \bigoplus_Y X.
\end{equation}

In order to keep our formulas simple, we shall abuse language and will not write these canonical isomorphisms. \\

Note that the category $\FF$ has involution:

\begin{equation}
(\sp)^t:\FF^{op}\xrightarrow{\sim} \FF
\end{equation}
\begin{equation}
(f:D(f)\iso I(f))^t=(f^{-1}:I(f)\iso D(f)),
\end{equation}
$$ (\text{or the transposed } \{0,1\} \text{ - matrix})$$ \vp
\begin{equation}
(f\circ g)^t=g^t\circ f^t,\sp\;\;\; (id_X)^t=id_X,
\end{equation}
\begin{equation}
 (f^t)^t=f,
\end{equation}
and this involution preserves the sum $\oplus$ (and the product
$\otimes$):
\begin{equation}
(f_0\oplus f_1)^t=f^t_0\oplus f^t_1.
\end{equation}

We usually let $ X,Y,Z,W $ denote objects of $\FF$,
without explicitly saying so, and when we consider "$Set_{\bullet}(X,Y) $"'
it is usually implicitly assumed that $ X,Y \in \FF$.

\section{$\frs$}
\smallbreak

In this section we define the category of $\frs$, $\FR$, and the category of $\frs$ with involution, $\FR^t$. We show these categories are bi - complete.

\defin{1.2.1} \textit{An $\fr ~A$ is a category, together with a faithfull functor $\epsilon:\FF\rightarrow A$,
which is the \underline{identity on objects}, and a symmetric monoidal structure
\begin{equation}
\bigoplus:A\times A\rightarrow A
\end{equation}
such that $\epsilon$ is strict monoidal; thus $A$ has objects the
finite sets, for $X,Y\in A$, $X\bigoplus Y$ is the disjoint union,
and $\bigoplus$ on arrows has associativity(resp.
commutativity,unit) isomorphism given by $\epsilon(a)$
(resp. $\epsilon(c)$, $\epsilon (u)$, which we abuse language and omit from our formulas!). We also demand that $[0]$ is the initial and final object of $A$.}\\

Thus an $\fr$ is a collection of pointed sets $A_{Y,X}=A(X,Y)$ for
all finite sets $X,Y$, together with the operation of composition:
\begin{equation}
A_{Z,Y}\times A_{Y,X}\rightarrow A_{Z,X}
\end{equation}
\begin{equation}
g,f\mapsto g\circ f
\end{equation}
which is associative:
\begin{equation}
(h\circ g)\circ f=h\circ(g\circ f),
\end{equation}
unital:
\begin{equation}
f\circ id_X=f=id_Y\circ f \text{ for } f\in A_{Y,X},
\end{equation} \\
and agree with composition on arrows of $\FF$, where $\FF_{Y,X}\subseteq A_{Y,X}$. \vp

We also have the operation of direct sum:
\begin{equation}
A_{Y_0,X_0}\times A_{Y_1,X_1}\rightarrow A_{Y_0\oplus Y_1,X_0\oplus X_1}
\end{equation}
\begin{equation}
f_0,f_1\mapsto f_0\oplus f_1
\end{equation}
such that we have: $\vsp$ \\
naturality:
\begin{equation}
\left(g_0\oplus g_1\right)\circ \left(f_0\oplus f_1\right)=(g_0\circ f_0)\oplus \left(g_1\circ f_1\right)
\end{equation}
\begin{equation}
id_{X_0}\oplus id_{X_1}=id_{X_0\oplus X_1}
\end{equation}
 \\  associativity:
\begin{equation}
\left(f_0\oplus f_1\right)\oplus f_2=f_0\oplus \left(f_1\oplus f_2\right)
\end{equation}
\begin{equation}
(\text{i.e. }\;\;\; a_{Y_0,Y_1,Y_2}\circ ((f_0\oplus f_1)\oplus f_2)= (f_0\oplus (f_1\oplus f_2))\circ a_{X_0,X_1,X_2}, \;\;\;f_i\in A_{Y_i,X_i})
\end{equation}
 \\  commutativity:
\begin{equation}
f_0\oplus f_1=f_1\oplus f_0
\end{equation}
\begin{equation}
 (\text{i.e. }\;\;\;  c_{Y_0,Y_1}\circ (f_0\oplus f_1)=(f_1\oplus f_0)\circ c_{X_0,X_1}\;,\; f_i\in A_{Y_i,X_i})
\end{equation}
 \\  unit:
\begin{equation}
f\oplus 0=f
\end{equation}
\begin{equation}
(i.e. \;\;\; (f\oplus id_0)\circ u_X=u_Y\circ f \;\;\;, \;\;\; f\in A_{Y,X})
\end{equation}
 and $f_0\oplus f_1$ agree with the sum in $\FF$ for $f_i\in \FF_{Y_i,X_i}\subseteq A_{Y_i,X_i}$.  \\

\defin{1.2.2} \textit{For $\frs \; A,A'$ we denote by $\FR(A,A')$ the collection of functors $\varphi: A\rightarrow A'$, over $\FF$,
and strict monoidal. i.e. $\varphi$ is a collection of set mappings $\varphi_{Y,X}:A_{Y,X}\rightarrow A'_{Y,X}$,
with $\varphi(g\circ f)=\varphi (g)\circ \varphi (f)$; $\varphi(f_0\oplus f_1)=\varphi (f_0)\oplus \varphi (f_1)$;
$\varphi(\epsilon (f))=\epsilon '(f)$ for $f\in \FF_{Y,X}$.
Thus we have a category of $\frs:\FR .$ }\vsp \\
We have for $X\in \FF$ a functor $(\sp)^X:\FR\rightarrow \FR$,
\begin{equation}
\sp\hsp  (A^X)_{Z,Y}:=A_{X\otimes Z,X\otimes Y}.
\end{equation}
We have a functor $(\sp)^{op}:\FR\rightarrow \FR$,
\begin{equation}
(A^{op})_{Y,X}:=A_{X,Y}.
\end{equation}
We denote by $\FR^t$ the $\frs$ with involution, i.e. $A\in \FR$ and $(~)^t:A^{op}\rightarrow A$ is an involution:
\begin{equation}
f^{tt}=f
\end{equation}
\begin{equation}
(g\circ f)^t=f^t\circ g^t
\end{equation}
\begin{equation}
\epsilon (f)^t=\epsilon (f^t) \hsp \text{for} \sp f\in \FF_{Y,X}
\end{equation}
\begin{equation}
\left(f_0 \oplus f_1\right)^t=f_0^t\oplus f_1^t .
\end{equation}
For $A,A'\in \FR^t$, we let
\begin{equation}
\FR^t(A,A')=\{\varphi \in \FR(A,A'), \hsp \varphi (a^t)=\varphi (a)^t \}.
\end{equation}
We have a category of $\frs$ with involutions: $\FR^t$.\vp \\
For $A\in \FR, X\in \FF$, the set $A_{X,X}$ is an associative monoid, and we let $GL_X(A)$ denote its invertible elements:
\begin{equation}
GL_X(A)=\{a\in A_{X,X}, \sp \text{there exists} \sp a^{-1}\in A_{X,X}, \sp a\circ a^{-1}=a^{-1}\circ a=id_X \}
\end{equation}
Clearly this constitute a functor form $\FR$ to the category of groups,
\begin{equation}\label{glx}
GL_X: \FR\rightarrow Grps.
\end{equation}
We have embeddings
$$GL_{X_0}(A)\times GL_{X_1}(A)\hookrightarrow GL_{X_0\oplus X_1}(A)$$
$$(a_0,a_1)\mapsto a_0\oplus a_1$$
In particular, we have the group $GL_{\infty}(A)=\underset{n}{\varinjlim} GL_{[n]}(A)$, the direct limit with respect to $a\mapsto a\oplus 1$ (hence Quillen's higher K-groups $K_i(A)=\pi_i(BGL_{\infty}(A)^+)$ , cf. \cite{Q73}).
\thm{1.1} \textit{The categories $\FR$ and $\FR^t$ are bi-complete.} \\
\begin{proof} Inverse limits exists and can be calculated in sets:
\begin{equation}
 ( \lim_{\underset{j}{\longleftarrow}} A^{(j)})_{Y,X}=\lim_{\underset{j}{\longleftarrow}} ( A_{Y,X}^{(j)})
\end{equation}
Also co-limits over a directed partially ordered set $J$, exists and can be calculated in Sets:
\begin{equation}
( \lim_{\underset{j}{\longrightarrow}} A^{(j)})_{Y,X}=\lim_{\underset{j}{\longrightarrow}} ( A_{Y,X}^{(j)})
\end{equation}
The pushout $B^0\underset{A}{\coprod} B^1$ is denoted (as usual) by $B^0\underset{A}{\otimes} B^1$, it is constructed in the appendix A.1 .
 Equalizers also exist, and in fact we can factorize an $\fr$ by any equivalence - ideal, this is described in appendix A.2 .
 This suffices for the existence of arbitrary co-limits.
\end{proof}

\section{Commutativity}
\smallbreak
Note that we usually do not write canonical-isomorphisms of $\FF$, especially, $\underset{Y}{\bigoplus} X\overset{\sim}{\leftrightarrow} X\bigotimes Y$.

\defin{1.3.1}{\it Let $A\in \fr$. We say $A$ is:\vp \\

\underline{Totally - commutative:}
\begin{equation}
\forall a\in A_{Y,X} , \forall b\in A_{J,I}\;, \;(\underset{J}{\bigoplus} a)\circ (\underset{X}{\bigoplus} b)=
(\underset{Y}{\bigoplus} b)\circ (\underset{I}{\bigoplus} a)\in A_{Y\otimes J,X\otimes I}.
\end{equation}

\underline{Left - commutative:}
\begin{equation}
\forall a\in A_{[1],X} , b\in A_{[1],I}\;, a\circ (\underset{X}{\bigoplus} b)=b\circ (\underset{I}{\bigoplus} a)\in A_{[1],X\otimes I}.
\end{equation}

\underline{Right - commutative:}
\begin{equation}
\forall a\in A_{Y,[1]} , b\in A_{J,[1]}\;, (\underset{J}{\bigoplus} a)\circ b= (\underset{Y}{\bigoplus} b)\circ a\in A_{Y\otimes J,[1]}.
\end{equation}

\underline{1 - commutative:} If $A$ is both Left-  and Right- commutative. \vsm  \\

\underline{$\times$ - commutative:}
\begin{equation} \forall a\in A_{Y,[1]},\forall b\in A_{[1],I}\;,\;\;\;, a\circ b=(\underset{Y}{\bigoplus} b)\circ
(\underset{I}{\bigoplus} a)\in A_{Y,I}.
\end{equation}

\underline{Central:}
\begin{equation}
\forall a\in A_{Y,X} , b\in A_{1,1} : a\circ (\underset{X}{\bigoplus} b)= (\underset{Y}{\bigoplus} b)\circ a=:b\cdot a\in A_{Y,X}.
\end{equation}
i.e.  $A_{1,1}$ is a commutative monoid and it acts centrally on
$A_{Y,X}$, \\
and we shall denote this action by $b\cdot a$.

\underline{Commutative:}
$$\forall a\in A_{Y,X} ,  b\in A_{1,J} ,  d\in A_{J,1}:$$
\begin{equation}
a\circ (\underset{X}{\bigoplus} b\circ d)= (\underset{Y}{\bigoplus}b\circ d)\circ
a=(\underset{Y}{\bigoplus} b)\circ (\underset{J}{\bigoplus}
a)\circ (\underset{X}{\bigoplus} d)\in A_{Y,X}.
\end{equation}}

We let $C\FR$ ( resp. $\FR_{\text{tot-com.}}, \FR_{1-\text{com.}},
\FR_{\times-\text{com.}},\FR_{\text{cent'l}}$ ) denote the full
subcategories of $\FR$ consisting of the commutative (resp.
totally-commutative, 1-commutative, $\times$-commutative, central)
$\frs$. Noting that, $A$ commutative $\implies$ $A$ central, we
have embeddings of categories

\begin{equation}\label{embeddingfr}
\begin{tikzcd}
& \CFR\arrow[hook]{r}&\FRce\arrow[hook]{rd}&\\
\FRtc\arrow[hook]{ru}\arrow[hook]{rd}\arrow[hook]{r}& \FRx\arrow[hook]{rr}&&\FR \\
&\FRone\arrow[hook]{rru}
\end{tikzcd}
\end{equation}

\section{Matrix coefficients and tameness.}
\smallbreak

For a set $X\in \FF$, and an element $x\in X$, we denote by $j_x=j_x^X$, the morphisms of $\FF$ given by
\begin{equation}
j_x:[1]\rightarrow X\;,\; j_x(1)=x\in X,
\end{equation}
and where
\begin{equation}
j_x^t:X\rightarrow [1] \text{   is the partial bijection } \{x\}\rightarrow \{1\}.
\end{equation}

\defin{1.4.1} \textit {Define the martix coefficients $J_{Y,X}:A_{Y,X}\rightarrow (A_{1,1})^{Y\times X}$  via
\begin{equation}
a\mapsto \{j_y^t\circ a\circ j_x\},
\end{equation}
where $j_x\in \FF_{X,1}\;,\; j_x(1)=x\;,\; \text{and } j_y^t\in \FF_{1,Y}$.}

\defin{1.4.2}\textit{ We say $A\in \FR$ is a "matrix $\FF$ - ring" if $J_{Y,X}$ is injective for any $X,Y\in \FF$.}

\defin{1.4.3}\textit{ We say $A\in \FR$ is \underline{tame:}\quad $\forall a,a'\in A_{Y,X},$}
\begin{equation}
\forall b\in A_{1,Y},\forall d\in A_{X,1}: b\circ a\circ d=b\circ a'\circ d\in A_{1,1} \implies a=a'.
\end{equation}
We have the implication (taking $b=j_y^t,d=j_x$):
\begin{equation}
A\text{ matrix} \implies A\text{ tame.}
\end{equation}
We also have the implication:
\begin{equation}
A\text{ commutative} +\text{tame}\implies \text{ A} \times\text{-commutative.}
\end{equation}
(indeed, for $a\in A_{Y,1},b\in A_{1,J}$, and any $d\in A_{J,1},
d'\in A_{1,Y}$, commutativity gives $d'\circ a\circ b\circ
d=d'\circ\underset{Y}{\bigoplus} b\circ \underset{J}{\bigoplus} a
\circ d$, hence by tameness $a\circ b=\underset{Y}{\bigoplus}
b\circ \underset{J}{\bigoplus} a$).

\let\cleardoublepage\clearpage
\begin{appendices}
\chapter{}
\section{Proof of existence of $B^0\coprod_{A} B^1=B^0\otimes_A B^1$.}
\thm{A.1}{ The category $\FR$ has pushouts $\otimes_A=\coprod_A$.
\begin{equation}\label{pushoutdiagram}
\begin{tikzcd}
&B^0\arrow[dotted]{rd}{\psi^0}\arrow[bend left]{rrd}{f^0}&& \\
A\arrow{ru}{\varphi^0}\arrow[swap]{rd}{\varphi^1}&&B\arrow[dotted]{r}{f^0\otimes f^1}&C \\
&B^1\arrow[dotted]{ru}{\psi^1}\arrow[swap,bend right]{rru}{f^1}&&
\end{tikzcd}
\end{equation}
}
\begin{proof}
Define the sets of chains of arrows,

\begin{equation}
\mathcal{B}_{Y,X}=\{(b_l,...,b_{\delta})\;,\; l\geqq \delta=0,1\;,\; b_j\in B_{X_{j+1,X_j}}^{j \text{ mod} 2},\;Y=X_{l+1}\;,\;X=X_{\delta} \}
\end{equation}



Let $\sim$ be the equivalence relation on chains generated by:

$$\en{1} \hsm \big(\dots,b_{j+1}\circ\varphi^{j+1}(a),b_j,\dots\big)\sim \big(\dots,b_{j+1},\varphi^j(a)\circ b_j,\dots \big),$$
$$\en{2} \hsm \big(\dots,b_{j+1},f,b_j,\dots\big)\sim \big(\dots,b_{j+1}\circ f\circ b_j,\dots \big), f\in \FF, \hsm \hsm$$
and the boundary cases:
$$  \big(f,b_j,\dots\big)\sim \big(f \circ b_j,\dots\big) \;\;,\;\;  \big(\dots ,b_{\delta} ,f\big) \sim \big(\dots ,b_{\delta}\circ f\big)$$
\begin{equation}
 f\in\FF \;\;\;(\text{or  } f=\varphi^j(a))
\end{equation}
$$\en{3}\hsm \big(\dots,b_{j+2}\;,\;id_{X'_{j+2}}\oplus \bar{b}_{j+1}\;,\;\bar{b}_j\oplus id_{X''_j}\;,\;b_{j-1}\;,\;\dots\big)\sim  \hsf\hsm\;\;\;$$
$$\hsf \big(\dots,b_{j+2}\circ (\bar{b}_j\oplus id_{X''_{j+2}})\;, \;(id_{X'_j}\oplus \bar{b}_{j+1})\circ b_{j-1}\;,\;\dots\big)\;\;$$
$$ \bar{b}_{j+1}\in B_{X''_{j+2},X''_j}^{j+1\text{ (mod~2})}\;,\;\bar{b}_{j}\in B_{X'_{j+2},X'_j}^{j\text{ (mod~2)}}\;,
\; X'_j\oplus X''_j=X_j\;,\;X''_j\oplus
X'_{j+2}=X_{j+1}\;,\;X'_{j+2}\oplus X''_{j+2}=X_{j+2} $$ or in
diagram:
$$ (\dots, \xleftarrow{b_{j+2}} \underset{X''_{j+2}}{\overset{X'_{j+2}}{\bigoplus}} \xleftarrow[id]{\bar{b}_{j}}
\underset{X''_{j+2}}{\overset{X'_j}{\bigoplus}}
\xleftarrow[\bar{b}_{j+1}]{id}
\underset{X''_j}{\overset{X'_j}{\bigoplus}} \xleftarrow{b_{j-1}}
\dots) \sim (\dots, \xleftarrow{b_{j+2}}
\underset{X''_{j+2}}{\overset{X'_{j+2}}{\bigoplus}}
\xleftarrow[\bar{b}_{j+1}]{id}
\underset{X''_j}{\overset{X'_{j+2}}{\bigoplus}}
\xleftarrow[id]{\bar{b}_{j}}
\underset{X''_j}{\overset{X'_j}{\bigoplus}} \xleftarrow{b_{j-1}}
\dots)$$
$$ \text{and the boundary cases:   } l=j+1, \delta=j.$$
We let $B_{Y,X}=\mathcal{B}_{Y,X}/\sim$.\\
The composition $B_{Z,Y}\times B_{Y,X}\rightarrow B_{Z,X}$ is
given by,
\begin{equation}
\big(b'_{l'},\dots,b'_{\delta'}\big)/\sim \circ \big(b_{l},\dots,b_{\delta}\big)/\sim=\begin{cases}
\big(b'_{l'},\dots,b'_{\delta'},b_{l},\dots,b_{\delta}\big)/\sim & \delta'\not\equiv_2 l\\
\big(b'_{l'},\dots,b'_{\delta'}\circ b_{l},\dots,b_{\delta}\big)/\sim & \delta'\equiv_2 l
\end{cases}
\end{equation}
It is well defined, independent of representatives (since $\circ$
is associative). Furthermore it is associative and unital,
\begin{equation}
(\text{with identity: }\; (id_X^0)/\sim\;=\;(id_X^1)/\sim\;=\;id_X)
\end{equation}
and therefore $B$ is a category. It has the natural maps
$\psi^j:B^j\rightarrow B\;,\text{and } \; \psi^0\varphi^0=
\psi^1\varphi^1$
by $(1)$.\\
Define the map,
\begin{equation}
\bigoplus:B_{Y_0,X_0}\times B_{Y_1,X_1}\rightarrow B_{Y_0\oplus Y_1,X_0\oplus X_1}
\end{equation}
\begin{equation}
\text{by:    }  \big(b'_{l},\dots,b'_{\delta}\big)/\sim \oplus \big(b_{l},\dots,b_{\delta}\big)/\sim\; = \big(b'_{l}\oplus b_{l}
,\dots,b'_{\delta}\oplus b_{\delta} \big)/\sim
\end{equation}
Note that without loss of generality we can assume that
$l'=l,\delta'=\delta$, (otherwise add identities). The map
$\bigoplus$ is well defined and independent of representatives.
Indeed, it is invariant by the 3 possible moves:

move $(1):$
$$ (\dots, b_{j+1}',b_j',\dots)\bigoplus (\dots, b_{j+1}\circ \varphi^{j+1}(a),b_j,\dots)=$$
$$(\dots, b_{j+1}'\oplus (b_{j+1}\circ \varphi^{j+1}(a)),b_j'\oplus b_j, \dots )=(\dots, (b_{j+1}'\oplus b_{j+1})\circ \underbrace{(id_{X_{j+1}'}\oplus \varphi^{j+1}(a))}_{\varphi^{j+1}(id_{X_{j+1}'}\oplus a)}, b_j'\oplus b_j,\dots ) $$
$$= (\dots, b_{j+1}' \oplus b_{j+1}, b_j'\oplus (\varphi^j (a)\circ b_j),\dots)$$
$$\sim(\dots, b_{j+1}'\oplus b_{j+1},\underbrace{(id_{X_{j+1}'}\oplus \varphi^j(a))}_{\varphi^j (id_{X_{j+1}'}\oplus a)}\circ (b_j'\oplus b_j),\dots ) $$
\begin{equation}
 =(\dots, b_{j+1}',b_j',\dots)\bigoplus (\dots, b_{j+1},\varphi^{j}(a)\circ b_j,\dots)
\end{equation}
move $(2):$
$$ (\dots, b_{j+1}',b_j',b_{j-1}',\dots)\bigoplus (\dots,b_{j+1},f,b_{j-1}, \dots)=$$
$$(\dots,b_{j+1}'\oplus b_{j+1}, b_j'\oplus f, b_{j-1}'\oplus b_{j-1},\dots)=(\dots,b_{j+1}'\oplus b_{j+1}, (b_j'\oplus id_{X_{j+1}})\circ (id_{X'_{j}}\oplus f), b_{j-1}'\oplus b_{j-1},\dots)$$
$$\sim (\dots b_{j+1}'\oplus b_{j+1},b_{j}'\oplus id_{X_{j+1}}, b_{j-1}'\oplus (f\circ b_{j-1}),\dots )$$
$$=(\dots, (b_{j+1}'\oplus id_{X_{j+2}})\circ (id_{X'_{j+1}}\oplus b_{j+1}),(b_j'\oplus id_{X_{j+1}})\circ id_{X_j' \oplus X_{j+1}}, b_{j-1}'\oplus (f\circ b_{j-1}),\dots )$$
$$\sim (\dots, (b_{j+1}'\oplus id_{X_{j+2}}),(b_j'\oplus id_{X_{j+2}}), (id_{X'_j}\oplus b_{j+1}),id_{X_j' \oplus X_{j+1}},b_{j-1}'\oplus (f\circ b_{j-1}),\dots)$$
$$ \sim (\dots,b_{j+1}'\oplus id_{X_{j+2}}, b_j'\oplus id_{X_{j+2}},b_{j-1}'\oplus (b_{j+1}\circ f\circ b_{j-1}), \dots)=$$
\begin{equation}
=(\dots, b_{j+1}',b_j',b_{j-1}',\dots)\bigoplus (\dots,id_{X_{j+2}},id_{X_{j+2}},(b_{j+1}\circ f\circ b_{j-1}), \dots)
\end{equation}

move $(3):$
$$(\dots, b_{j+2}',b_{j+1}',b_j',b_{j-1}',\dots)\bigoplus (b_{j+2},(id\oplus b_{j+1}),(b_j\oplus id),b_{j-1})=$$
$$(\dots, b_{j+2}'\oplus b_{j+2}, b_{j+1}'\oplus (id\oplus b_{j+1}), b_j'\oplus (b_j\oplus id), b_{j-1}'\oplus b_{j-1},\dots )=$$
$$= (\dots, b_{j+2}'\oplus b_{j+2}, (b_{j+1}'\oplus id)\circ (id \oplus (id\oplus b_{j+1})),(id\oplus (b_j\oplus id))\circ (b_j'\oplus id), b_{j-1}'\oplus b_{j-1},\dots) $$
$$\sim (\dots, b_{j+2}'\oplus b_{j+2},b_{j+1}'\oplus id, id\oplus (b_j\oplus id), id\oplus (id\oplus b_{j+1}),b_j'\oplus id,b_{j-1}'\oplus b_{j-1},\dots )$$
$$\sim (\dots,b_{j+2}'\oplus (b_{j+2}\circ (b_j\oplus id)),(b_{j+1}'\oplus id), (b_j'\oplus id), b_{j-1}'\oplus ((id\oplus b_{j+1})\circ b_{j-1}),\dots)$$
\begin{equation}
=(\dots, b_{j+2}',b_{j+1}',b_j',b_{j-1}',\dots )\bigoplus (\dots, (b_{j+2}\circ (b_j\oplus id)), id,id, ((id\oplus b_{j+1})\circ b_{j-1})
\end{equation}
It is than straightforward to check that $B$ satisfy all the
axioms of an $\fr$.

For $\FR^t$ define the induced involution on $B$:
\begin{equation}
(\;)^t:B_{Y,X}\rightarrow B_{X,Y}
\end{equation}
$$ ((b_l,\dots,b_{\delta})/\sim)^t=(b_{\delta}^t,\dots, b_l^t)/\sim.$$
This operation is well defined since all moves are self dual, and
$B$ is an $\fr$ with involution.

The map,
\begin{equation}
f^0\otimes f^1 ((b_l,...,b_{\delta})/\sim)= f^l(b_l)\circ \dots \circ f^\delta(b_\delta), \hsm f^j:= f^{j\;mod\;2}
\end{equation}
is well defined since it remains invariant under each move
$(1),(2),(3)$, and it is (the unique) homomorphism of $\fr$
solving the universal property ($\ref{pushoutdiagram}$).

\end{proof}

\section{Equivalence ideals and quotients.}
\defin{A.2.1} \textit{For $A\in \FR$,( resp. $A\in \FR^{(t)}$) an \underline{equivalence ideal} (resp. t-equivalence ideal) is a sub-$\fr$ (resp. with involution) $\mathcal{E}\subseteq A\prod A$ s.t.
\begin{equation}
\mathcal{E}_{Y,X}\subseteq A_{Y,X}\prod A_{Y,X}\hsm \text{is an equivalence  relation on}\sp A_{Y,X},
\end{equation}
equivalently, we have a collection of equivalence relations $\sim$
on the $A_{Y,X}$'s, compatible with the operations:
\begin{equation}\label{A.2.2}
a\sim a'\implies \;\;\; b\circ a\circ d\sim b\circ a'\circ d
\end{equation}
}
\begin{equation}
a\sim a'\implies\;\;\; a\oplus id\sim a'\oplus id
\end{equation}
$$(\textit{and therefore}\hsm a_0\sim a_0',a_1\sim a_1'
\implies a_0\oplus a_1\sim a'_0\oplus a'_1)$$
\begin{equation}\label{A.2.4}
\textit{and in the presence of an involution:} \hsm a\sim a'\implies\;\;\; a^t\sim (a')^t
\end{equation}\\
Given an equivalence ideal $\mathcal{E}$ of $A$, let
\begin{equation}
A/\mathcal{E} =
\coprod_{Y,X\in\FF}A_{Y,X}/\mathcal{E}_{Y,X},
\end{equation}
and let $\pi:A\to
A/\mathcal{E}$ denote the canonical map which associates with $a\in
A_{Y,X}$ its equivalence class $\pi(a)\in A_{Y,X}/\mathcal{E}_{Y,X}$.
It follows from (\ref{A.2.2}-\ref{A.2.4}) that we have well defined operations on $A/\mathcal{E}$,
$$\pi(f)\circ\pi(g) = \pi(f\circ g),$$
$$\pi(f)\oplus\pi(g) = \pi(f\oplus g)$$
\begin{equation}
\text{resp. } \hsm \pi(f)^t= \pi (f^t)
\end{equation}
making $A/\mathcal{E}$ into an $\fr$ such that $\pi:A\to A/\mathcal{E}$
is a homomorphism of $\frs$ (resp. with involutions). \\
Given a homomorphism of $\frs$ $\fe:A\to B$ denote by
$$\cK(\fe)=A\prod_B A = \coprod_{Y,X\in\FF}\cK_{Y,X}(\fe)\;,$$
\begin{equation}
\cK_{Y,X}(\fe)=\{(a,a')\in A_{Y,X}\times A_{Y,X}|\fe(a)=\fe(a')\}.
\end{equation}

It is clear that $\cK(\fe)$ is an equivalence-ideal of $A$, and that
$\fe$ induces an injection of $\frs$
$\bar{\fe}:A/\cK(\fe)\;\inj\;B$, such that $\fe = \bar{\fe}\circ
\pi$, i.e.

\begin{equation}
\begin{tikzcd}
A\arrow[two heads]{d}\arrow{r}{\fe}\arrow[two heads]{dr}&B \\
A/\cK(\fe)\arrow{r}{\sim}&\fe(A)\arrow[hook]{u}
\end{tikzcd}
\end{equation}
is a commutative diagram. Thus every map $\fe$ of $\frs$ factors as
epimorphism $(\pi)$ followed by an injection $(\bar{\fe})$.\vsp\\

For a family $\{(a_i,a_i')\}\in(A\times A)^I,\; a_i,a'_i\in
A_{Y_i,X_i}$ let $\mathcal{E}$ be the equivalence ideal generated by
$\{(a_i,a_i')\}$. \\
We have $(b,b')\in \mathcal{E}$ iff $\exists$ path
$b=b_0,b_1,...,b_l=b'$, s.t. $\{b_{j-1},b_j\}$ has the form

\begin{equation}\label{A.2.9}
\{c_j\circ (a_{i(j)}\oplus id_{V_j})\circ d_j,\;\; c_j\circ (a'_{i(j)}\oplus id_{V_j})\circ d_j\}\hsf
\end{equation}
$$(\text{or the form } \{c_j\circ (a^t_{i(j)}\oplus id_{V_j})\circ d_j,\;\;\;c_j\circ (a'^t_{i(j)}\oplus id_{V_j})\circ d_j\},
\text{ in the presence of involution}). $$ Indeed, if there is
such a path $b=b_0,\dots, b_l=b'$, than $(b_{j-1},b_j)\in \mathcal{E}$,
and $\mathcal{E}$ is an equivalence relation so $(b,b')\in \mathcal{E}$.
On the other hand, the set of $(b,b')$ s.t. there is a path
$b=b_0,\dots,b_l=b'$, with $\{b_{j-1},b_j\}$ of the form
$(\ref{A.2.9})$, is an equivalence relation $b\sim b'$, satisfying
\begin{equation}
b\sim b'\implies c\circ b\circ d\sim c\circ b'\circ d\;\;\;,\;\;\;b\oplus
id\sim b'\oplus id, \;\;\;\; (\text{and } b^t\sim b'^t),
\end{equation}
 and hence it is precisely the equivalence ideal
$\mathcal{E}$ generated by $\{(a_i,a'_i)\}$. \vspace{3mm}

For example, for any $\fr \; A$, we have the equivalence ideals
$\mathcal{E_?}$ generated by the "?- commutativity" relation, and the
associated quotient $A^?=A/\mathcal{E}_?$, giving rise to the diagram
of surjections,

\begin{equation}\label{A.2.11}
\begin{tikzcd}
& A^{cent'l}=A/\mathcal{E}_{cent'l}\arrow[two heads]{r}&A^{com}=A/\mathcal{E}_{com}\arrow[two heads]{rd}&\\
A\arrow[two heads ]{ru}\arrow[two heads]{rd}&&&A^{tot.com}=A/\mathcal{E}_{tot-com.} \\
&A^{1-com}=A/\mathcal{E}_{1-com}\arrow[two heads]{rru}{\sim}
\end{tikzcd}
\end{equation}
The functor $A\mapsto A^?$ is left adjoint to the embeddings (\ref{embeddingfr}). $\FR_? \hookrightarrow \FR$ \\

\defin{A.2.2}\textit{For an $\mathcal{E} \subseteq A\prod A$ equivalence ideal, we say $\mathcal{E}$ is \underline{tame} if $A/\mathcal{E}$ is tame.
Thus a tame equivalence ideal $\mathcal{E}$ is completely determined by the equivalence relation $\mathcal{E}_{1,1}$ on $A_{1,1}$:
$$(a,a')\in \mathcal{E}_{Y,X}\iff \:\forall \; b\in A_{1,Y\oplus Z},\;\; d\in A_{X\oplus
Z,1},\;(b\circ(a\oplus id_Z)\circ d,\:\;\; b\circ (a'\oplus
id_Z)\circ d)\in \mathcal{E}_{1,1}. $$ }

We have the following bijection:
$$\{\text{tame equivalence ideals}\}$$
$$\overset{1:1}{\longleftrightarrow}$$
$$\{\mathcal{E}\subseteq A_{1,1}\times A_{1,1}\hsm \text{equivalence relation such that }$$
$$\text{for } (a_i,a'_i)\in \mathcal{E}^{I},b\in A_{1,I},d\in A_{I,1}\implies
(b\circ (\oplus a_i)\circ d,b\circ (\oplus a'_i)\circ d)\in
\mathcal{E}\}$$
\vspace{3mm}\\
For $A\in \FR^t$ take also $(a,a')\in \mathcal{E}\implies
(a^t,(a')^t)\in \mathcal{E}$ .
\defin{A.2.3}\label{defA.2.3} \textit{ $\mathfrak{a}\subseteq A_{1,1}$ is an ideal if the following property hold:
$$a_j\in \mathfrak{a},b\in A_{1,J}, d\in A_{J,1}\implies b\circ (\underset{J}{\oplus}a_j)\circ d\in \mathfrak{a} $$
 and is a t-ideal if also:
$$ a\in \mathfrak{a}\implies a^t\in \mathfrak{a}. $$}
\prop{A.2.4}\textit{There exists a Galois connection of,
$$ \{\mathfrak{a}\subseteq A_{1,1} \; \text{ideal} \}\underset{E}{\overset{Z}{\leftrightarrows}} \{\mathcal{E}\subseteq A_1\times A_1 \; \text{is tame equivalence ideal}\}$$
$$ \{\mathfrak{a}\subseteq A_{1,1} \; \text{t-ideal} \}\leftrightarrows \{\mathcal{E}\subseteq A_1\times A_1 \; \text{is tame t-equivalence ideal}\}$$ where
$$ Z(\mathcal{E})=\{a\in A_{1,1},\hsm (a,0)\in \mathcal{E} \} $$
$$ E(\mathfrak{a})= \underset{(\mathfrak{a},0)\subseteq \mathcal{E}}{\bigcap} \mathcal{E} $$ }

\end{appendices}

\chapter{Examples of $\FF\text{-} \mathcal{R}$ings.}
\bigskip

\section{$\mathcal{R}$ings}
\smallbreak

\defin{2.1.1} \textit{A "$\catrig$" is a ring without negatives. Thus a $\catrig$ is a set with two associative operations $(+,\bullet)$, with units $0$ and $1$, addition $+$ being commutative, and multiplication distributive over addition. A morphism between $\catrigs$ $A\rightarrow B$ is a map $\varphi:A\rightarrow B$ which perserves operations and units. Thus we have a category: $\catrigs$. We let $C\catrigs \subseteq \catrigs$ denote the full subcategory of commutative $\catrigs$, i.e. where $x\cdot y=y\cdot x$. We let $\catrigs^t$ denote the category of $\catrigs$ with involution, its objects are rigs $R$ with involution:
$$ (\;)^t:R\rightarrow R$$
$$x^{tt}=x$$
$$(x\cdot y)^t=y^t\cdot x^t$$
$$(x+y)^t=x^t+y^t$$
$$1^t=1$$
\begin{equation}
0^t=0.
\end{equation}
and the morphisms are morphisms of $\catrigs$ that preserve the involution.\\ Note that the identity is always an involution on a commutative $\catrig$, and so we have a diagram of categories and functors}
\begin{equation}
\begin{tikzcd}
&\catrig^t\arrow{d} \\
C\catrig\arrow[hook]{r}\arrow[dotted]{ru}&\catrig
\end{tikzcd}
\end{equation}

We shall write $\NN$ (resp. $[0,\infty)^{(1)}$ ) for the rigs of
natural numbers (resp. non-negative reals) with the usual
operation of multiplication $\bullet$ and addition $+$. We shall
write $\NN^0$ (resp. $[0,\infty)^{(0)}$) for the "frozen" rigs
where $x+y:= \Max\{x,y\}$. Moreover, for $p\geq 1$ (or $1/p\in
(0,1])$, we write $[0,\infty)^{(1/p)}$ for the rig of non-
negative reals with $x+y:= (x^p+y^p)^{(1/p)}$, and with the usual
multiplication $\bullet$. Note that the rigs
$[0,\infty)^{(\sigma)}$ interpolate continuously between the
frozen $(\sigma =0)$ and the usual $(\sigma =1)$ reals, and that
for $\sigma \in (0,1]$ they are all isomorphic via
$\Phi_{\sigma_2}^{\sigma_1}: [0,\infty)^{(\sigma_1)}\iso
[0,\infty)^{(\sigma_2)}\;,\; \Phi_{\sigma_2}^{\sigma_1}(x)=
x^{\sigma_2/ \sigma_1}$. The multiplicative group of positive reals act as automorphisms of the frozen rig $[0,\infty)^{(0)}$ via $x\mapsto \Phi_{\sigma}(x)=x^{\sigma}$. For any $q\in (0,1)$, we have the sub-rigs $M_q=\{0\}\cup q^{\NN}\subseteq [0,1]^{(0)}\subseteq [0,\infty)^{(0)}$, and $N_q=\{0\}\cup q^{\ZZ}\subseteq [0,\infty)^{(0)}$, and the multiplicative monoid pf positive natural numbers act as endomorphisms of $M_q$ and $N_q$ via $q^j\mapsto \Phi_n(q^j)=q^{jn}$.

\defin{2.1.2} \textit{Let $R\in \mathcal{R}$ig. Define $R^{+/-}$ by
$$R^{+/-}\equiv R\times R,$$
$$(n_+,n_- ) +(m_+,m_-) = (n_++m_+,n_-+m_-), $$
$$(n_+,n_- ) \cdot (m_+,m_-)= (n_+\cdot m_++ n_-\cdot m_-, n_+\cdot m_- +n_-\cdot m_+).$$
}
There is an equivalence relation $\sim$ on $R^{+/-}$:
$$(n_+,n_-)\sim (m_+,m_-) \iff n_++m_-+r=m_++n_-+r,\; \text{some } r\in R.$$
such that we have,
$$a\sim a' , b\sim b' \implies \begin{matrix}
  a+b \sim a'+b' \\ a\cdot b\sim a'\cdot b'
\end{matrix}.$$
We have the sequence of Rig homomorphisms,

$$\xymatrix@R-2pc{
R\;\ar@{^{(}->}[r] &  R^{+/-}\ar@{->>}[r]^-{\pi}& \slfrac{R^{+/-}}{\sim}\equiv K(R) \\
r\ar@{|->}[r] & (r,0) & \\
&(n_+,n_-)\ar@{|->}[r]^-{\pi} &n_+-n_-:= (n_+,n_-)/ \sim.&}$$

The rig $K(R)$ is a ring, and for any ring $A$,
$$\mathcal{R}ing (K(R),A)= \mathcal{R}ig (R,A).$$
If $R$ has involution (resp. commutative), than $R^{+/-}$ and $K(R)$ have involution (resp. are commutative).
We get adjoint functors,
$$\xymatrix{C\mathcal{R}ing\ar[r]\ar@/^/[d]^{U} & \mathcal{R}ing^t \ar[r]\ar@/^/[d]^{U}& \mathcal{R}ing\ar@/^/[d]^{U} \\
C\mathcal{R}ig\ar[r]\ar@/^/[u]^{K} & \mathcal{R}ig^t \ar[r]\ar@/^/[u]^{K} & \mathcal{R}ig\ar@/^/[u]^{K}}$$

\defin{2.1.3}\textit{For any Rig $A$, let $A\cdot X=A^X$ be the free $A$ - module with basis $X$, and define:}
$$\FF(A)_{Y,X} =\Hom_A (A\cdot X,A\cdot Y)=Y\times X~ \text{ - matrices with values in} ~A.$$
The composition $\circ$ is matrix multiplication, and $\bigoplus$ is the direct sum of matrices. Clearly $\FF (A)$ is an $\fr$.
Note that a morphism of Rigs $\varphi : A\to B$, induces a map
of $\frs$
$$\FF(\varphi) :\FF(A)\to \FF(B),$$
hence we have a functor
$$\FF(\;) :\catrigs\to \FR .$$

Note that if $A\in \catrigs^{t}$ has involution, than $\FF(A)\in
\FR^t$ also has involution: for $a=(a_{y,x})\in \FF(A)_{Y,X}$ we
have $a^t\in \FF(A)_{X,Y}\;,\; (a^t)_{x,y}=(a_{y,x})^t\;,$ and it
satisfies $(a\circ b)^t=b^t\circ a^t$. Thus we have a functor
$\FF(\;):\catrigs^t\rightarrow \FR^t$. Note that if $A\in
C\catrig$ is commutative, than $\FF(A)$ is totally-commutative.
Thus we have the diagram of categories and functors

\begin{equation}
\begin{tikzcd}
&\catrigs^t\arrow{r}{\FF(\;)}\arrow{d}&\FR^t\arrow{d} \\
&\catrigs\arrow{r}{\FF(\;)}&\FR \\
C\catrigs\arrow{r}{\FF(\;)}\arrow[hook]{ru}\arrow[bend left]{ruu}&\FRtc\arrow[hook]{ru}&
\end{tikzcd}
\end{equation}

Note that $\FF(A)$ is always a matrix-$\;\fr$.\vp \\

Moreover let $\varphi :\FF(A)\to \FF(B)$ be a map of $\frs$. For
$a\in\FF(A)_{Y,X}$ write
\begin{equation}
a_{y,x} = j^t_y\circ a\circ j_x\;\in A = \FF(A)_{[1],[1]},
\end{equation}
for its matrix coefficients. \vp

Since $\fe$ is a functor over $\FF$, and $j^t_y,j_x\;\in\FF$, we
have $\fe(a)_{y,x} = \fe(a_{y,x})$ and $\fe$ is determined by $\fe:
A = \FF(A)_{[1],[1]}\to B = \FF(B)_{[1],[1]}$. This map is
multiplicative: $\fe(a_1\cdot a_2) =
\fe(a_1)\cdot\fe(a_2)\;,\;\fe(1) = 1$, and moreover it is additive:
\begin{equation}
\fe(a_1 + a_2) = \fe\left[(a_1,a_2)\circ
 \begin{pmatrix}
1\\
1
\end{pmatrix}\right] = (\fe(a_1),\fe(a_2))\circ
\begin{pmatrix}
1\\
1
\end{pmatrix} = \fe(a_1) + \fe(a_2)
\end{equation}
Thus the functor $\FF(\;)$ is fully faithful.\;\;\\
\defin{2.1.4} \textit{Define,
$$(\FR^{(t)})^{Add}\equiv \text{commutative monoid objects in } \FR^{(t)}, $$
That is, $A\in (\FR^{(t)})^{Add}$ is the same as an addition map $A_{Y,X}\times A_{Y,X}\xrightarrow{+} A_{Y,X}$ which is associative, commutative, has a unit element $0_{Y,X}$, and satisfies,
$$a\circ (b+b')=a\circ b+a\circ b' \;, (a+a')\circ b=a\circ b+a'\circ b,\; , (\text{resp.} (a+b)^t=a^t+b^t).  $$
The category of abelian group objects is the full subcategory $(\FR^{(t)})^{Ab}\subseteq (\FR^{(t)})^{Add}$. We have $A\in (\FR^{(t)})^{Ab} \iff A_{Y,X}\in Ab, \forall X,Y\in \FF$.
}\\
For a rig (resp. ring) $R$, the $\fr$ $\FF(R)$ is in
$(\FR)^{Add}$, (resp. $(\FR)^{Ab}$), and we have a similar diagram
as before,
\begin{equation}
\begin{tikzcd}
&\catrigs^t\arrow{r}{\FF(\;)}\arrow{d}&(\FR^t)^{Add}\arrow{d} \\
&\catrigs\arrow{r}{\FF(\;)}&(\FR)^{Add} \\
C\catrigs\arrow{r}{\FF(\;)}\arrow[hook]{ru}\arrow[bend left]{ruu}&(\FRtc)^{Add}\arrow[hook]{ru}&
\end{tikzcd}\begin{tikzcd}
&\catrings^t\arrow{r}{\FF(\;)}\arrow{d}&(\FR^t)^{Ab}\arrow{d} \\
&\catrings\arrow{r}{\FF(\;)}&(\FR)^{Ab} \\
C\catrings\arrow{r}{\FF(\;)}\arrow[hook]{ru}\arrow[bend left]{ruu}&(\FRtc)^{Ab}\arrow[hook]{ru}&
\end{tikzcd}
\end{equation}

\section{Monoids}
\smallbreak

\defin{2.2.1}\textit{ Let $M$ be a monoid with a unit $1$ and a zero
element $0$. Thus we have an associative operation
$$M\times M\to M\;\;,\;\;(a,b)\;\mapsto\;a\cdot b,$$
\begin{equation}
a\cdot(b\cdot c) = (a\cdot b)\cdot c
\end{equation}
and $1\in M$ is the (unique) element such that
\begin{equation}
1\cdot a= a\cdot 1 = a\;\;,\;\;a\in M,
\end{equation}
and $0\in M$ is the (unique) element such that
\begin{equation}
0\cdot a= a\cdot 0 = 0 \;\;\;,\;\;\;a\in M .
\end{equation}
Let $\FF\{ M\}$ denote the $\fr$ with $\FF\{
M\}_{Y,X}$ the $Y\times X$ matrices with values in $M$ with at
most one non-zero entry in every row and column. Note that this is
indeed an $\fr$ with the usual "multiplication" of matrices $\circ$
(there is no addition involved  - only multiplication in $M$) and
direct sum $\oplus$.}\vsm \par Denoting by
$\mathcal{M}on$ the category of monoids with unit
and zero elements, and with maps respecting the operation and the
elements $0,1$, the above construction yields a functor
\begin{equation}
\mathcal{M}on\to \FR\;\;,\;\; M\;\mapsto\;\FF\{ M\}.
\end{equation}
This is the functor left-adjoint to the functor
\begin{equation}
\FR\to \mathcal{M}on\;\;,\;\;A\;\mapsto\;A_{[1],[1]} :
\end{equation}
\begin{equation}
\FR(\FF\{ M\},A) = \mathcal{M}on(M,A_{[1],[1]}).
\end{equation}
As a particular example, take $M=M_q$ to be the free monoid (with
zero) generated by one element $q$,
\begin{equation}
M_q=q^{\nn}\cup\{0\}.
\end{equation}
Then
\begin{equation}
\FR(\FF\{ M_q\},A)=A_{[1],[1]}.
\end{equation}
Denote by $\Mon^t$ the category of monoids with involution (i.e.
the objects are monoids $M\in\Mon$, with involution
$(\;)^t:M\rightarrow M,\;x^{tt}=x,\; (x\cdot y)^t=y^t\cdot x^t,\;
1^t=1,\;0^t=0$, and the morphisms respect the involutions).We
have an involution on $\FF\{M\}$ for $M\in \Mon^t$, so that we
have a functor, $\FF\{\}:\Mon^t\rightarrow\FR^t$. Denote by
$\CMon\subseteq \Mon$ the full subcategory of commutative monoids
(i.e. where $x\cdot y=y\cdot x$). For $M\in \CMon,\FF\{M\}$ is
totally commutative. Thus we have the diagram
\begin{equation}
\begin{tikzcd}
\Mon^t\arrow{r}{\FF\{\}}\arrow{d}&\FR^t\arrow{d} \\
\Mon\arrow{r}{\FF\{\}}& \FR \\
\CMon\arrow{r}{\FF\{\}}\arrow[hook]{u}&\FRtc\arrow[hook]{u}
\end{tikzcd}
\end{equation}
Note that $\FF\{M\}$ is always a matrix $\fr$, and the functors
$\FF\{\}$ are full and faithful.

\section{$ \mathcal{S}et , \mathcal{S}et^{op} \subseteq \mathcal{R}el \subseteq \FF (\NN^0)$. }
\smallbreak

\defin{2.3.1}\textit{Let $\mathcal{S}et$ denote the $\fr$ of sets. The objects of $\mathcal{S}et$
are the finite sets of $\FF$, and we let $\mathcal{S}et_{Y,X}$ be the
partially defined maps of sets from $X$ to $Y$
\begin{equation}
\mathcal{S}et_{Y,X}=\cS et_{\bullet}(X,Y)=\mathcal{S}et_0(X\cup \{0\},Y\cup \{0\} ).
\end{equation}
We can view the elements of $\mathcal{S}et_{Y,X}$ as $Y\times X$ -matrices with values in $\{0,1\}$, such that every column has at most one $1$:
\begin{equation}
f\leftrightarrow (f)_{y,x} \;\;\;\text{with}\;\;\; (f)_{y,x}=
\begin{cases}
  1  & y=f(x) \\
  0 & \mbox{otherwise}
\end{cases}
\end{equation}
Than composition $\circ$ corresponds to matrix multiplication; The disjoint union $\bigoplus$ correspond to direct sum of matrices. These make $\mathcal{S}et$ into an $\fr$ (with no involution), which is matrix and totally-commutative.} \vp \\
We have the opposite $\fr$ $\mathcal{S}et^{op}$
with
\begin{equation}
\mathcal{S}et^{op}_{Y,X}=\mathcal{S}et_{X,Y}.
\end{equation}
We have the $\fr$ of relations $\mathcal{R}el$ that contains both $\mathcal{S}et$
and $\mathcal{S}et^{op}$, with
\begin{equation}
\mathcal{R}el_{Y,X}=\{F\subseteq Y\times X\; \text{a subset}\}:=\{0,1\}^{Y\times X}.
\end{equation}
The composition of
$F\in\mathcal{R}el_{Y,X}\;\text{and}\;G\in\mathcal{R}el_{Z,Y}$ is given by
\begin{equation}
G\circ F=\{(z,x)\in Z\times X|\exists\;y\in Y\;\text{with}\;(z,y)\in G, (y,x)\in F\},
\end{equation}
and $G\circ F\in\mathcal{R}el_{Z,X}$.\vp \\
When we view $G,F$ as $\{0,1\}$ -matrices, this composition correspond to "matrix - multiplication" where we replace addition by $Max\{x,y\}$. \vp \\
The sum $F_0\oplus F_1\in\mathcal{R}el_{Y_0\oplus Y_1,X_0\oplus X_1}$ of
$F_i\in\mathcal{R}el_{Y_i,X_i}$ is given by the disjoint union of $F_0$
and $F_1$,
\begin{equation}
F_0\oplus F_1=\{((x,i),(y,i))|(x,y)\in F_i\}\text{ or by direct sum of matrices}.
\end{equation}
Thus $\mathcal{R}el=\FF(\{0,1\})$ is the $\fr$ with involution associated to the rig $\{0,1\}$ with usual multiplication $\bullet$, and $i+j=Max\{i,j\}$ as addition. \vp \\
Thus $\mathcal{R}el$ is totally -commutative,matrix $\fr$ with involution. The embedding of rigs $\{0,1\}\hookrightarrow \NN^0$
gives an embedding of $\frs\;\;\; \mathcal{R}el=\FF\bigg(\{0,1\}\bigg)\hookrightarrow \FF(\NN^0).$

\section{Real primes}
\smallbreak

Let $\eta:\kk\inj \CC$ be an embedding of the rig $\kk$ into the
complex numbers.  We have an injection of $\frs$ $\FF(\kk)\inj
\FF(\CC)$ . For $X\in\FF$,  let $\kk^X=\FF(\kk)_{X,[1]}$ denote
the free $\kk$-module over $X$. Thus for $a=(a_x)\in\kk^X$, and
$p\in [1,\infty]$, we have the vector $p$ - norm:
\begin{equation}
|a|_p = \left(\sum_{x\in X}|\eta(a_x)|^p\right)^{1/p} \;,\; p\in [1,\infty)
\end{equation}
\begin{equation}
|a|_{\infty}=\underset{x\in X}{Max\;}|\eta(a_x)|
\end{equation}
and for $a=(a_{y,x})\in\kk^{Y\times X}$ we have its operator $p$ - norm:
\begin{equation}
\|a\|_p =\underset{b\in \kk^X, |b|_p\leq 1}{\Sup} \Big\{ |a\circ b|_p\Big\}
\end{equation}
e.g.,
\begin{equation}
\|\underbrace{(1,\dots,1)}_{n}\|_p=\underset{|\xi_1|^p+\dots |\xi_n|^p\leq 1}{\Sup} | \xi_1+\dots+\xi_n |=n^{1-\frac{1}{p}}\equiv n^{1/p'}.
\end{equation}
\begin{equation}
\Bigg\| \left( \begin{array}{c}
1 \\
\vdots \\
1 \end{array} \right)\Bigg\|_p=\underset{|\xi|\leq 1}{\sup} \left(
\overset{n}{\underset{i=1}{\Sigma}} |\xi|^p\right)^{1/p}=n^{1/p}.
\end{equation}
\begin{equation}
\|(1,\dots,1)\|_p\cdot \Bigg\| \left( \begin{array}{c}
1 \\
\vdots \\
1 \end{array} \right)\Bigg\|_p= n^{1/p'}\cdot n^{1/p}=n.
\end{equation}

\defin{2.4.1}\textit{ Define the sub-\;$\fr\; \Op\;\subseteq \FF(\kk)$ as follows:
\begin{equation}
\left(\Op\right)_{Y,X}=\{a\in \FF(\kk)_{Y,X}=\kk^{Y\times X}\;\; ,\;\; \|a\|_p\leq 1\}.
\end{equation}
}
As a sub-$\fr$ of $\FF(\kk)$, the $\frs$ $\Opsig,\sigma\in [0,1]$, are matrix and totally - commutative.\\  Note that in general, $\Opsig$ has no involution, in fact we have(by H\"{o}lder's inequality: the dual of the $p$- norm is the $p'$- norm, where $\frac{1}{p'}=1-\frac{1}{p}$):

\begin{equation}
(\Opsig)^{op}\simeq \mathcal{O}_{\kk,\eta}^{1- \sigma}\;\; \text{for}\;\; 0 \leq\sigma\leq 1.
\end{equation}

\defin{2.4.2}\textit{Let $\mathcal{O}_{\kk,\eta}:=\Optwo\subseteq \FF(\kk)$. It is a sub $\fr$ of $\FF (\kk)$, hence totally - commutative, and matrix, and at $\sigma =\frac{1}{2}$ (i.e. using the $L_2$ -norm) it has involution! }

\defin{2.4.3}\textit{Define $\FF_{\kk,\eta}\in \FR^t$,
$$\left(\FF_{\kk,\eta}\right)_{Y,X}= \{ f:D(f)\iso I(f)\;, \; D(f)\subseteq \kk^X\;,\; I(f)\subseteq \kk^Y\; \text{are $\kk$ -subspaces, and, } $$
\begin{equation}
f\; is \;\text{$\kk$-linear and an isometry}: |f(v)|_2=|v|_2 \}
\end{equation}
}
Note that when the composition $g\circ f$ is defined we have:
\begin{equation}
D(g\circ f)=f^{-1}(D(g)\cap I(f))\;\; , \; \; I(g\circ f)=g(D(g)\cap I(f)).
\end{equation}
There is a surjective homomorphism of $\frs$ with involution, $\phi:\mathcal{O}_{\kk,\eta}\sur \FF_{\kk,\eta}$.  \\
For $a\in (\mathcal{O}_{\kk,\eta})_{Y,X}\subseteq \CC^{Y\times X}$, define $\Delta_X =\bar{a}^t\circ a, \Delta_Y =a\circ \bar{a}^t$, and let $V_X[\lambda]$ (resp. $V_Y[\lambda]$) denote the $\lambda$ - eigenspace of $\Delta_X$ (resp. $\Delta_Y$).

The operator $\Delta_X$ (resp. $\Delta_Y$) is non-negative and we have the spectral decomposition $\CC^X=\underset{i}{\bigoplus} V_X[\lambda_i]\;,\; \Delta_X=\underset{i}{\bigoplus}\lambda_i\cdot id_{V_X[\lambda_i]}$\hsm (resp. $\CC^Y=\underset{i}{\bigoplus} V_Y[\lambda_i],\Delta_Y=\underset{i}{\bigoplus}\lambda_i\cdot id_{V_Y[\lambda_i]}$ ), with eigenvalues $\lambda_i\in [0,1]$, and the non-zero eigenvalues are the same for $\Delta_X$ and $\Delta_Y$ including multiplicities.\\
For $\lambda>0$ we have isomorphisms:
\begin{equation}
\begin{tikzcd}
V_Y[\lambda] \arrow[bend left]{r}{\sim} &V_X[\lambda]\arrow[bend left]{l}{\sim}
\end{tikzcd}$$ \\
For $\lambda=1$ we have an isometry:
$$\phi (a) = \{ a: V_X[1]\iso V_Y[1]\}\in (\FF_{\kk,\eta})_{Y,X} .
\end{equation}
This defines the homomorphism $\phi: \mathcal{O}_{\kk,\eta}\sur \FF_{\kk,\eta}$.
Note that $\FF_{\kk,\eta}$ is (as quotient of $\mathcal{O}_{\kk,\eta} $) totally - commutative.
It is our first example of an $\fr$ which is \underline{Not} matrix:\\
Indeed, any vector $x=(x_1,\dots, x_n)\in (\FF_{\kk,\eta})_{1,n}$, with $\sum_{i=1}^n |\eta (x_i)|^2=1$, but with $|\eta (x_i)|< 1$ for $i=1,\dots, n$, is non- zero, but all its matrix coefficients are zero. \\
But note that $\FF_{\kk,\eta}$ is tame.

\section{Valuation $\frs$: Ostrowski theorem}\label{2.5}
\smallbreak

\defin{2.5.1}\textit{A commutative $\fr$ with involution $K\in \CFR^t$ is called an \underline{$\FF$ - field} if every non-zero $a\in K_{[1],[1]} \setminus \{0\}$ is invertible (have $a^{-1}\in K_{[1],[1]}$ with $a\circ a^{-1}=a^{-1}\circ a=1=id_{[1]}).$  \vsm \\
A sub-$\fr$ with involution $B\subseteq K$ is called
\underline{full} if
$$\en{1} \hsm \text{for every } a\in K_{Y,X},X,Y\in
\FF,\hfar$$
$$\text{there exists a non - zero element } d\in
B_{[1],[1]}\setminus \{0\} \sp \text{with } $$
$$d\cdot
a=(\underset{Y}{\oplus} d)\circ a=a\circ (\underset{X}{\oplus}
d)\in B_{Y,X}$$ (This means that $K$ is the fraction-field of the
domain $B$, i.e. $K=B_{(0)}=(B_{[1],[1]}\setminus \{0\})^{-1}\cdot
B$ the localization of $B$ at the prime $(0)$, cf. \S \ref{3.3}).\\ We will
say that $B$ is \underline{tame in $K$} if for $X,Y\in \FF$ we
have an equality:
$$\en{2} \hsm  B_{Y,X}=\{a\in K_{Y,X} \;|\; \text{for all} \; b\in B_{1,Y},d\in B_{X,1}, \; b\circ a\circ d\in B_{1,1} \}. \vsm $$
A sub-$\fr$ with involution $B$, full and tame in $K$, will be called a \underline{valuation-\FF-subring} of $K$
if for every non-zero $a\in K_{1,1}\setminus \{0\}$,
$$\en{3}\;\;\;  \text{either} \;\; a\in B_{1,1}\;\;\; \text{or}\;\;\; a^{-1}\in B_{1,1}. \;\;\hfar$$
Given $\FF$-fields $k\subseteq K$, we denote by $Val(K/k)$ the set of all valuation-$\FF$-subrings $B\subseteq K$, such that $B\supseteq k$.  }\vsm\\
Let $B$ be a valuation - $\FF$- subring of an $\FF$-field $K$. The group of units:
$$B^*=GL_1(B)=\{a\in B_{1,1}\;|\; \exists a^{-1}\in B_{1,1}\;,\; a\circ a^{-1}=1\} $$
is a subgroup of
$$K^*=GL_1(K)=K_{1,1}\setminus {\{0\}}.$$
The quotient group $\Gamma=K^*/B^*$ is ordered: $|x|\leq |y| \iff x\cdot y^{-1}\in B_{1,1}$, where $|x|=x\cdot B^*$
is the quotient map $|\;|:K^*\rightarrow \Gamma$. We extend this quotient map by $|0|=0$, to the map
$$|\;|:K_{1,1}\rightarrow K_{1,1}/B^*=\Gamma\cup \{0\}$$
satisfying
$$(I)\hspace{100mm}$$
$$\;\;\;\; |x|=0\iff x=0$$
$$\;\;|x_1\cdot x_2|=|x_1|\cdot |x_2| $$
$$\;\;\;\;\;\;|1|=1\;\;\; (=\text{unit of} \; \Gamma).$$
We can embed $\Gamma$ in a \underline{complete} ordered abelian group $\tilde{\Gamma}$, (e.g. $\tilde{\Gamma}=$
all dedekind subsets $D\subseteq \Gamma$), so that for every subset $Y\subseteq \tilde{\Gamma}$ which is bounded above
(resp. bellow) there is a unique least upper bound $\sup Y\in \tilde{\Gamma}$ (resp. maximal lower bound $\inf Y\in \tilde{\Gamma}$).
We can than define for $X,Y\in \FF$ the two maps
$$ |\;|_{Y,X}\;,\;|\;|_{Y,X}':K_{Y,X}\rightarrow \tilde{\Gamma}\cup\{0\} \hsm$$
$$(i) \hsm |y|_{Y,X}=\sup \{|b\circ y\circ b'|,b\in B_{1,Y},b'\in B_{X,1}\}$$
$$(ii) \hsm |y|_{Y,X}'=\inf \{|d^{-1}|,d\in K^* \;,\;d\cdot y\in B_{Y,X}\;\;\;\; \}$$

(II)\hsm \underline{\bf{Claim:}}  $\;\;|y|_{Y,X}= |y|_{Y,X}'$.
\begin{proof} Note that for $y\in K_{Y,X}, $ there exists $d\in K^*$ with $d\cdot y\in B_{Y,X}$, and for $b\in B_{1,Y},b'\in B_{X,1}$ we have $b\circ (d\cdot y)\circ b'\in B_{1,1}$, and so
$$ |b\circ y\circ b'|=|d^{-1}\cdot b\circ (d \cdot y)\circ b'|=|d|^{-1}\cdot |b\circ (d\cdot y)\circ b'|\leq |d|^{-1}.$$
This shows that the set in (i) (resp. (ii)) is bounded above (resp. below), and we have the inequality: $|y|_{Y,X}\leq |y|_{Y,X}'$.\\
Conversely, given $y\in K_{Y,X}$, if $d^{-1}\in K^*$ is such that $|d|_{1,1}^{-1}\geq |y|_{Y,X}$, that is: $|d|_{1,1}^{-1}\geq |b\circ y\circ b'|_{1,1}$
for all $b\in B_{1,Y},b'\in B_{X,1}$, than $b\circ (d\cdot y)\circ b'\in B_{1,1}$ for all $b\in B_{1,Y},b'\in B_{X,1}$, and (since B is tame in $K$)
this imply $d\cdot y \in B_{Y,X}$, hence $|d|_{1,1}^{-1}\geq |y|'_{Y,X}$, giving the reverse inequality: $|y|_{Y,X}\geq |y|'_{Y,X}$.
\end{proof}

(III)\hsm \underline{\bf{Claim:}}
$$(i)\;\;\; |a\circ a'|\leq |a|\cdot |a'| $$
$$(ii)\;\;\; |a_0\oplus a_1|=\max\{|a_0|,|a_1|\}$$
$$(iii)\;\;\; |a^t|=|a|$$

\begin{proof} (i): If $d,d'\in K^*$ are such that $|d|^{-1}\geq |a|, |d'|^{-1}\geq |a'|$, than $d\cdot a,d'\cdot a'$ are in $B$, so $(d\circ d')\cdot a\circ a'=(d\cdot a)\circ (d'\circ a')$ is in $B$, and $|d|^{-1}|d'|^{-1}\geq |a\circ a'|$. \vsp \\
(ii): If $d_0,d_1\in K^*$ are such that $d_0\cdot a_0,d_1\cdot a_1$ are in $B$, and if $|d_j|\leq |d_{1-j}|$ than $d_j\cdot a_0,d_j\cdot a_1$ are in $B$,
 so $d_j\cdot (a_0\oplus a_1)=(d_j\cdot a_0)\oplus (d_j\cdot a_1)$ is in $B$, so $|a_0\oplus a_1|\leq |d_j|^{-1}=\max\{|d_0|^{-1},|d_1|^{-1}\}$
 Taking the infimum over all such $d_0,d_1$ we get $|a_0\oplus a_1|\leq \max\{|a_0|,|a_1|\}=|a_{j_0}|$, say. The inverse inequality follows from (i)
  since $a_{j_0}=f'\circ (a_0\oplus a_1)\circ f$, with $f',f$ arrows of $\FF\subseteq B$, so $|a_{j_0}|\leq |f'|\cdot |a_0\oplus a_1|\cdot |f|\leq |a_0\oplus a_1|$.\vsp \\
(iii) This follows since we are assuming $B\subseteq K$ to be stable under the involution, so $a\in B_{Y,X}$ if and only if $a^t\in B_{X,Y}$.
\end{proof}
Let $\Gamma$ be a complete ordered abelian group, written multiplicatively, and form the ordered abelian monoid $\Gamma\cup\{0\}$, with $0\cdot x=0, 0<x$
for all $x\in \Gamma$. \\
Given a collection of mappings
$$ |\;|_{Y,X}:K_{Y,X}\rightarrow \Gamma\cup \{0\}$$
satisfying (III), with $|\;|_{1,1}$ satisfying (I), the subsets
$$B_{Y,X}=\{a\in K_{Y,X},\; |a|\leq 1\}$$
form a sub-$\FF$- ring with involution $B\subseteq K$. If we have
the equalities for $y\in K_{Y,X}$
$$(II)\hsp \hsp \;\;\;\;\; (i)\;\;\; |y|_{Y,X}=\sup\{|b\circ y\circ b'|_{1,1}\; ; \; |b|_{1,Y},|b'|_{X,1}\leq 1\}\hsf \;\;\;\hsm $$
$$ (ii)\hsp =\inf\{|d|_{1,1}^{-1};\;\; d\in K^* \;\;,|d\cdot y|_{Y,X}\leq 1 \} $$
than $B$ is full (by II ii), and tame in $K$ (by II i), and it follows that $B$ is a valuation-$\FF$- subring of $K$. \\

\thm{Ostrowski I} \textit{
$$Val(\FF(\QQ)/\FF\{\pm 1\})=\{\FF(\QQ),\;\FF(\ZZ_{(p)}),\; \text{p a finite prime},\; \mathcal{O}_{\QQ,\eta}\}.$$
where $\mathcal{O}_{\QQ,\eta}$ is the real prime.}

\begin{proof}
cf. Appendix \ \ref{appB}.
\end{proof}

\thm{Ostrowski II} \textit{For a number field $K$,
$$Val(\slfrac{\FF(K)}{\FF\{\mu_K\}})= \{\FF(K)\; ; \; \FF(\mathcal{O}_{K,\mathfrak{p}})\; ; \;\; \mathcal{O}_{K,\eta}\}$$
with $\mathcal{O}_{K,\mathfrak{p}}$ the localization of the ring of
integers $\mathcal{O}_K$ at prime ideals $\mathfrak{p}\subseteq
\mathcal{O}_K$, and with $\mathcal{O}_{K,\eta}$, the "real primes" of $K$,
$\eta$ varies over the embeddings $\eta:K\hookrightarrow \CC$
modulo conjugation.}
\begin{proof}
cf. Appendix \ \ref{appB}.
\end{proof}

\rem{2.5.1} Note that for any $\sigma\in [0,1]$, the sub-$\FF$- ring $\mathcal{O}_{K\eta}^{(\sigma)}\subseteq \FF(K)$ satisfies (1),(2),(3) of definition (2.5.1), i.e. it is full and tame valuation $\FF$- subring. Alternatively, the operator $p=1/\sigma$- norm satisfies (I),(II) and (IIIi),(IIIii). But only at $\sigma=\frac{1}{2},p=2$, we have an involution on $\mathcal{O}_{K,\eta}^{(\sigma)}$.\\ Thus it is the presence of the involution that singles out the $L_2$- metric at the real primes.

\section{Graphs}
\smallbreak

\defin{2.6.1} \textit{A graph $G$ is a pair of finite sets $(G_0,G_1)$ with two maps:}
\begin{equation}
G=\Big\{ \begin{tikzcd}
G^1\arrow[bend left]{r}{\pi^0}\arrow[bend right]{r}{\pi^1} & G^0
\end{tikzcd} \Big\}
\end{equation}
\textit{where $G^0$ - 'vertices' ,$G^1$ - 'edges'.}

Given such a finite graph we get a category $\CG$: the objects of $\CG$ are the elements of $G_0$,
and the arrows of $\CG$ are given by "paths":
$$Ob(\CG)\equiv G^0,$$
\begin{equation}
\CG (x,y)= \{e\equiv (e_l,\dots,e_1)\;|e_i\in G^1, \pi^0(e_{j+1})=\pi^1(e_j), \pi^1(e_l)=y,\pi^0(e_1)=x\}.
\end{equation}

\defin{2.6.2}\textit{Given such a path $e=(e_l,...,e_1)\in \CG (x,y)$ we shall say that $e$ "begins" at $x$, "ends" at $y$,
and for a vertice $z \in G^0$, we say $e$ "goes through z" and we write:
\begin{equation}
z \in e \iff  \exists k, z = \pi^i (e_k), i\in \{0,1\},
\end{equation}
for an edge $e_0\in G^1$, write:
\begin{equation}
e_0\in e \iff  \exists k, e_0 = e_k.
\end{equation}}

Note that for a vertice $x\in G^0$, we have
\begin{equation}
id_x=\text{"empty path" at } x\in \CG(x,x)
\end{equation}

Assume $G$ has \underline{no loops}: $\CG(x,x)=\{id_x\}, \forall x\in G^0.$ \vp

$\triangleright$ Every path can be extended to a maximal path. \vp

$\triangleright$ Every maximal path begins at $In(G)=G^0\setminus \pi^1(G^1)$. \vp

$\triangleright$ Every maximal path ends at  $Out(G)=G^0\setminus \pi^0(G^1)$. \vp

We denote by $m(G)$ the set of all maximal paths.
\begin{equation}
\begin{tikzcd}
& m(G)\arrow{ld}[swap]{\pi^1}\arrow{rd}{\pi^0}&  \\
\underset{\pi^1 (e_l,...,e_1)=\pi^1(e_l)}{Out(G)}&&\underset{\pi^0
(e_l,...,e_1)=\pi^0(e_1)}{In(G)}
\end{tikzcd}
\end{equation}
\vspace{5mm}\\
We define an $\fr$ with involution: $Graph\in \FR^t$.
\begin{equation}
(Graph)_{Y,X} =\begin{tikzcd}
\Big\{ G= G^1\arrow[bend left]{r}{\pi^0}\arrow[bend right]{r}{\pi^1} & G^0
\end{tikzcd}, \text{no loops}+\begin{array}{lcr} i: In\; G\hookrightarrow X \\
 o: Out\; G\hookrightarrow Y \end{array} \Big\}/ \text{isom.}
\end{equation}
(i.e. modulo isomorphisms of graphs that respects the embeddings $i$ and $o$.)\vp \\
For $G\in (Graph)_{Y,X}$ and $X_0\subseteq X$, define $G[X_0]$:
$$\text{vertices:    \;\;   } x\in G^0,\;\;\exists e\in m(G), \;\;\underline{x\in e},\;\;\pi^0(e)\in X_0.$$
$$\text{edges: \;\;   } e_0\in G^1,\;\;\exists e\in m(G), \;\;\underline{e_0\in e},\;\;\pi^0(e)\in X_0.$$
and for $Y_0\subseteq Y$, define $[Y_0]G$:
 $$\text{vertices:    \;\;   } x\in G^0,\;\;\exists e\in m(G), \;\;\underline{x\in e},\;\;\pi^1(e)\in Y_0.$$
$$\text{edges: \;\;   } e_0\in G^1,\;\;\exists e\in m(G), \;\;\underline{e_0\in e},\;\;\pi^1(e)\in Y_0.$$
We also let,
\begin{equation}
[Y_0]G[X_0]=[Y_0]G\cap G[X_0].
\end{equation}
For $G\in (Graph)_{Y,X},\;\; G'\in (Graph)_{Z,Y}$, let
$Y_0=In(G)\cap Out(G')\subseteq Y$,

\begin{equation}
\begin{tikzcd}
&Y_0\arrow[hook]{rd}\arrow[hook]{ld}& \\
In \;G'\arrow[hook]{rd}&&Out \; G\arrow[hook]{ld} \\
&Y&
\end{tikzcd}
\end{equation}
The operation of composition is defined by gluing $G'[Y_0]$ and
$[Y_0]G$ along $Y_0$:
\begin{equation}\label{glue}
G'\circ G= G'[Y_0]\underset{Y_0}{\amalg}[Y_0]G
\end{equation}
The sum in $Graph$ is given by the disjoint union:
\begin{equation}
G_i\in (Graph)_{Y_i,X_i},\;\; G_0\oplus G_1=G_0\amalg G_1.
\end{equation}
The involution is given by reversing the directions of the edges
of the graph:

\begin{equation}
\begin{tikzcd}  (G= G^1\arrow[bend left]{r}{\pi^0}\arrow[bend right]{r}[swap]{\pi^1} & G^0)^t \end{tikzcd}=
\begin{tikzcd}  (G= G^1\arrow[bend left]{r}{\pi^1}\arrow[bend right]{r}[swap]{\pi^0} & G^0) \end{tikzcd}:(Graph)_{Y,X}\iso (Graph)_{X,Y}.
\end{equation}

E.g. the "discrete" graphs $G$ in $(Graph)_{Y,X}$, (i.e.
$G^{1}=\emptyset$ and $G$ is just a set
\begin{equation}
G^0=In \;G=Out \;G
\end{equation}
with embeddings into $X$ and into $Y$), give the elements of $\FF$:
\begin{equation}
\FF_{Y,X}=\{G^0\hookrightarrow X,\; G^0\hookrightarrow Y\}/isom\subseteq (Graph)_{Y,X}
\end{equation}

$Graph\in \FR^t$ is not even central, \underline{not matrix}, but
it is \underline{tame}.\\ \vp Note that we have a homomorphism of
$\frs$ :
\begin{equation}
\phi: Graph\sur\FF(\NN), \;\; \phi(G)_{y,x}=\sharp \{e\in m(G),\pi^0(e)=x,\pi^1(e)=y\}
\end{equation}

\section{Free $\frs\; \FF[\delta_{Y,X}]$ }
\smallbreak

We have an $\fr$ $\FF[\delta_{Y,X}]\in \FR$, such that for any
$A\in \FR$,
\begin{equation}
\FR(\FF[\delta_{Y,X}],A)\equiv A_{Y,X},\;\;\; \varphi\mapsto \varphi (\delta_{Y,X}),
\end{equation}
and similarly, we have an $\frs$ with involution
$\FF[\delta_{Y,X},\delta_{Y,X}^t]\in \FR^t$, such that for any
$A\in \FR^t $,
\begin{equation}
\FR^t (\FF[\delta_{Y,X},\delta_{Y,X}^t],A)\equiv A_{Y,X}.
\end{equation}
The elements of $\FF[\delta_{Y,X}]_{W,Z}$ can be written as
sequences of maps in $\FF$,
$$(f_l,...,f_j,...,f_0),\; \text{with}\; f_j\in \FF_{(I_{j+1}\otimes X)\oplus V_{j+1},(I_j\otimes Y)\oplus V_j},\;\; l>j>0,$$
\begin{equation}\label{2.7.3}
f_0\in \FF_{(I_1\otimes X)\oplus V_1,Z}\;\;,f_l\in \FF_{W,(I_l\otimes Y)\oplus V_l}\;\;,
\end{equation}
modulo certain identifications. Such a sequence represents the
element
\begin{equation}\label{2.7.4}
f_l\circ \dots \circ f_j\circ ((\underset{I_j}{\oplus} \delta_{Y,X})\oplus
id_{V_j})\circ f_{j-1}\circ \dots \circ f_1\circ
((\underset{I_1}{\oplus} \delta_{Y,X})\oplus id_{V_1})\circ f_0.
\end{equation}
These elements can also be described as "$(Y,X)$- marked
graphs". The full $(Y,X)$ - graph is given by
\begin{equation}
\delta_{Y,X}\equiv \left(Y\otimes X\underset{\pi_1}{\overset {\pi_0}{\rightrightarrows}} Y\oplus X
\right), \hsm \;\; \begin{matrix}
  \pi_0(y,x)=x  \\
  \pi_1(y,x)=y.
\end{matrix}
\end{equation}
$$\delta_{Y,X}\in Graph_{Y,X}$$
$$ In(\delta_{Y,X})\equiv X $$
$$ Out(\delta_{Y,X})\equiv Y $$
e.g.,

\begin{equation}
\begin{tikzpicture}[scale=0.8, transform shape]
\tikzstyle{node}=[circle,draw,fill=black!100,inner sep=0.2pt,
minimum width=4pt]

\node[node] (a) at (0,0)[label=left: $y_1$] {};

\node[node](b) at (0,2)[label=left: $y_2$] {};

\node[node] (c) at (0,4)[label=left: $y_3$] {};

\node[node] (d) at (5,1)[label=right: $x_1$] {};

\node[node] (e) at (5,3)[label=right: $x_2$] {};

\foreach \from/\to in {a/d, a/e, b/d,b/e,c/d,c/e} \draw [->]
(\from) -- (\to);

 \draw [white]  (-2,1.5) -- (-1,1.5)
    node [above=1mm,midway,black,text width=8cm,text centered]
      {\textbf{$\delta_{Y,X}\equiv$}};

\end{tikzpicture}\hfar
\end{equation}

A $(Y,X)$- marked graph from $Z$ to $W$ is given by a graph of the
form
\begin{equation}
G=(J\otimes Y\otimes X\oplus G^1\underset{\pi_1}{\overset {\pi_0}{\rightrightarrows}}(J\otimes
Y)\oplus (J\otimes X)\oplus W_0\oplus Z_0)
\end{equation}
with $Z_0\subseteq
Z, W_0\subseteq W, \pi^0(j,y,x)=(j,x), \pi^1(j,y,x)=(j,y)$ , and
$\pi^0,\pi^1$ are injections on $G^1$:
$$\pi^0: G^1\hookrightarrow (\underset{j\in J}{\bigoplus}\; Y)\oplus Z_0$$
\begin{equation}
\pi^1: G^1\hookrightarrow (\underset{j\in J}{\bigoplus}\; X)\oplus W_0
\end{equation}
we shall assume it has no loops and that for every $j\in J$, there
is $y\in Y$, (resp. $x\in X$), and a path going from $(j,y)$ to
$W_0$, (resp. from $Z_0$ to $(j,x)$). Thus a $(Y,X)$ - marked
graph is a graph that can be made out of a disjoint union of the
full $(Y,X)$ - marked graphs $\delta_{Y,X}$ (one copy for each
$j\in J$), and some partial bijections. An isomorphism of such
$(Y,X)$ - marked graph $G=(J,G^1,\pi^0,\pi^1)$ and
$H=(I,H^1,\pi^0,\pi^1)$ is an isomorphism of graphs
$\phi:G\rightarrow H$, that is compatible with the maps into $W$
and $Z$, but is such that for some bijection
\begin{equation}
\sigma:J\rightarrow I\;\;, \hsm \begin{matrix}
                                    \\
  \phi(j,y,x)=(\sigma(j),y,x)  \\
  \phi(j,y)=(\sigma(j),y)  \\
  \phi(j,x)=(\sigma(j),x).
\end{matrix}
\end{equation}

The marked $(Y,X)$- graph associated to a sequence $(f_l,\dots,
f_0)$ as above, (\ref{2.7.3}),(\ref{2.7.4}),  is obtained by taking
$J=\overset{l}{\underset{j=1}{\bigoplus}} I_j$, and adding an edge
from $(i_1,y)$ to $(i_2,x)$ if and only if $(i_1,y)\in
D(f_{i_1}),f_{i_1}(i_1,y)\in V_{i_1+1}\cap
D(f_{i_1+1}),f_{i_1+1}\circ f_{i_1}(i_1,y)\in V_{i_1+2}\cap
D(f_{i_1+2}),\dots, f_{i_2-1}\circ \dots \circ
f_{i_1}(i_1,y)=(i_2,x)$ (and similarly for an element of $Z$,
resp. $W$, instead of $(i_1,y)$, resp. $(i_2,x)$). We can now
describe $\FF[\delta_{Y,X}]_{W,Z}$ as $(Y,X)$ - marked - graphs
from $Z$ to $W$ modulo isomorphism.

There is actually a "canonical form" (in fact two "dual" canonical
forms) for $G\in \FF[\delta_{Y,X}]$. Thinking of $G$ as a marked
$(Y,X)$ - graph $G=(J,G^1,\pi^0,\pi^1)$, let

$$J_1=\{j\in J|\forall e\in G^1, \pi^0(e)=(j,y)\implies \pi^1(e)\in
W\}$$
$$J_2=\{j\in J|\forall e\in G^1, \pi^0(e)=(j,y)\implies \pi^1(e)\in
W\; \text{or}\; \pi^1(e)=(j',x),j'\in J_1\}$$
$$\vdots$$
\begin{equation}
J_k=\{j\in J|\forall e\in G^1,\pi^0(e)=(j,y)\implies \pi^1(e)\in W\; \text{or} \; \pi^1(e)=(j',x)\; \text{with }\; j'\in \bigcup_{i<k} J_i\}
\end{equation}
$$u_{0,0}=\{e\in G^1|\pi^0(e)\in Z,\pi^1(e)\in W\} $$
$$u_{i,j}=\{e\in G^1|\pi^1(e)=(i_1,x),\pi^0(e)=(i_0,y) ,\:\text{for }\; i_0\in J_i\;,\; i_1\in J_j\},\; 0<i<j<k $$
$$u_{0,j}=\{e\in G^1|\pi^1(e)\in W, \;\; \pi^0(e)=(i,y)\;,\; i\in J_j \}$$
\begin{equation}
u_{j,j}=\{e\in G^1|\pi^1(e)=(i,x),\;,\; \pi^0(e)\in Z\;,\; i\in J_j\}
\end{equation}
Then the canonical form of $G$ is given by the non-empty sets $J_1,\dots,J_k$, the finite sets $\{u_{i,j}\}_{0\leq i\leq j\leq k}$, and the embeddings
$$(i) \hsp \underset{j}{\bigoplus}u_{0,j}\hookrightarrow W\;\;,\;\;\underset{j}{\bigoplus}u_{j,j}\hookrightarrow Z $$
$$(ii)\hsp\underset{i_0\leq j}{\bigoplus}u_{i_0,j}\hookrightarrow J_{i_0}\otimes X $$
\begin{equation}
(iii)\hsp\underset{i\leq j_0}{\bigoplus}u_{i,j_0}\hookrightarrow J_{j_0}\otimes X
\end{equation}
the embeddings in $(ii)$ and $(iii)$ are "dense" in the sense that
for each $j\in J_{i_0}$ (resp. $j\in J_{j_0}$), there is an $x\in
X$ (resp. $y\in Y$), such that $(j,x)$ (resp. $(j,y)$) is in the
image. \rem{1} There is a dual canonical form obtained by taking
instead
$$J_1=\{j\in J|\forall e\in G^1,\pi^1(e)=(j,x)\implies \pi^0(e)\in Z\}$$
$$J_2=\{j\in J|\forall e\in G^1,\pi^1(e)=(j,x)\implies \pi^0(e)\in Z\; \text{or }\; \pi^0(e)=(j',y),j'\in J_1 \}$$
etc. \rem{2} Similarly, any element $G\in
\FF[\delta_{Y,X},\delta_{Y,X}^t]_{W,Z}$ can be described as a
graph made out of disjoint copies of $\delta_{Y,X}$ and of
$\delta_{X,Y}$, and has canonical form
\begin{equation}
(J_1,I_1,J_2,I_2,\dots,J_k,I_k,u_{i,j}\;,\; 0\leq i < j\leq k)
\end{equation}
where now $J_i\oplus I_i\neq \emptyset$ for $i=1,\dots,k$, and the embeddings
$$(i) \hsp \underset{j}{\bigoplus}u_{0,j}\hookrightarrow W\;\;,\;\;\underset{j}{\bigoplus}u_{j,j}\hookrightarrow Z$$
$$(ii)\hsp\underset{i_0\leq j}{\bigoplus}u_{i_0,j}\hookrightarrow (J_{i_0}\otimes X)\oplus (I_{i_0}\otimes Y)$$
\begin{equation}
(iii)\hsp\underset{i\leq j_0}{\bigoplus}u_{i,j_0}\hookrightarrow (J_{j_0}\otimes Y)\oplus (J_{j_0}\otimes X)
\end{equation}
the embedding in (ii) and (iii) being dense (the image meets every copy of $X$ and $Y$ ).
\defin{2.7.1} \textit{When $G$ has such a cannonical form we shall write $\deg G=k$.\\ Thus an element of degree $0$ is just given by the set $u_{0,0}$ and the embeddings $u_{0,0}\hookrightarrow Z,u_{0,0}\hookrightarrow W$, i.e. it is just an element of $\FF$:
\begin{equation}
\FF_{W,Z}=\{G\in\FF[\delta_{Y,X}]_{W,Z} \;,\; \deg G=0 \}\;.
\end{equation}
When $\deg G=1$, we have the embeddings
$$u_{0,0}\oplus u_{0,1}\hookrightarrow W\;, \; u_{0,0}\oplus u_{1,1}\hookrightarrow Z$$
\begin{equation}
u_{0,1}\hookrightarrow \underset{J_1}{\bigoplus}Y\;,\; u_{1,1}\hookrightarrow \underset{J_1}{\bigoplus}X
\end{equation}
giving the elements
$$f_{00}=(W\hookleftarrow u_{00}\hookrightarrow Z)\in \FF_{W,Z}$$
$$f_{01}=(W\hookleftarrow u_{01}\hookrightarrow \underset{J_1}{\bigoplus}Y)\in \FF_{W,\underset{J_1}{\bigoplus}Y}$$
\begin{equation}
f_{11}=(\underset{J_1}{\bigoplus}X \hookleftarrow u_{11}\hookrightarrow Z)\in \FF_{\underset{J_1}{\bigoplus}X,Z}
\end{equation}
and
\begin{equation}
G=[f_{01}\circ (\underset{J_1}{\bigoplus}\delta_{Y,X})\circ
f_{11}]\oplus f_{00}.
\end{equation}
In general we have:
\begin{equation}
\deg (G_0\oplus G_1)=\max\{\deg G_0,\deg G_1\}
\end{equation}
\begin{equation}
\deg (H\circ G)\leq \deg(H)+\deg (G)
\end{equation}
and for $\FF[\delta_{Y,X},\delta_{Y,X}^t]$:
\begin{equation}
\deg G^t=\deg G.
\end{equation}
}\vspace{1.5mm}\\

More generally, given $C\in \FR$, and given $I\in \slfrac{Set}{\FF\times \FF}$, i.e. a set $I$ with a mapping $I\rightarrow \FF\times \FF$, $i\mapsto (Y_i,X_i)$, we have the $\FF$- ring
\begin{equation}
C[\delta_I]\equiv C[\delta_{Y_i,X_i};\; i\in I]:= C\underset{\FF}{\otimes}\underset{i\in I}{\bigotimes}\FF[\delta_{Y_i,X_i}],
\end{equation}
and we have the adjunction (with $U$ the forgetful functor):
\begin{equation}
\xymatrix{C\setminus \FR \ar@/^/[d]^{U} \\
 \slfrac{Set}{\FF \times \FF}\ar@/^/[u]^{C[\;]}}\hspace{5mm} C\setminus \FR \bigg(C[\delta_I],A\bigg)\equiv \slfrac{Set}{\FF\times \FF}\bigg(I,UA\bigg).
\end{equation}

Similarly, given $C\in\FR^t$, and given $I\in
(\slfrac{Set}{\FF\times \FF})^t$, i.e. $I$ is a set with maps
$d_i:I\rightarrow \FF, i=0,1$, and an involution $I\iso I,i\mapsto
i^t,(i^t)^t=i, d_0\circ (i^t)=d_1(i),d_1(i^t)=d_0(i)$, we have the
$\fr$ with involution \begin{equation}
C[\delta_I,\delta_I^t]:=C\underset{\FF}{\otimes}\underset{\slfrac{i\in
I}{i\sim i^t}}{\bigotimes} C[\delta_{Y_i,X_i},\delta_{Y_i,X_i}^t],
\end{equation}
and we have the adjunction
\begin{equation}
\xymatrix{C\setminus \FR^t \ar@/^/[d]^{U} \\
 (\slfrac{Set}{\FF \times \FF})^t\ar@/^/[u]^{C[\;]}}\hspace{5mm} C\setminus \FR^t \bigg(C[\delta_I,\delta_I^t],A\bigg)\equiv (\slfrac{Set}{\FF\times \FF})^t\bigg(I,UA\bigg).
\end{equation}

\section{$\FF[GL_X]$}
\smallbreak We have the functor $GL_X:\FR\rightarrow Grps\;
(\ref{glx})$,
\begin{equation}
GL_X(A)=\{a\in A_{X,X},\;\;\; \exists a^{-1}\in A_{X,X}\;\; a\circ a^{-1}=a^{-1}\circ a=id_{X} \}.
\end{equation}
It is representable:
\begin{equation}
GL_X(A)=\FR(\FF[GL_X],A)
\end{equation}
\begin{equation}
\FF[GL_X]=\slfrac{\FF[\delta_{X,X}]\underset{\FF}{\otimes} \FF[\delta_{X,X}']}{\{\delta_{X,X}\circ \delta_{X,X}'\sim
\delta_{X,X}'\circ \delta_{X,X}\sim id_X}\}
\end{equation}
The following structure exists on $\FF[GL_X]$:
\begin{equation}
\begin{split}
m^*:&\FF[GL_X]\rightarrow \FF[GL_X]\underset{\FF}{\otimes} \FF[GL_X]\hsf \text{(co-multiplication)} \\
&\delta_{X,X}\mapsto \delta_{X,X}^{(0)}\circ \delta_{X,X}^{(1)}\\
&\delta_{X,X}'\mapsto \delta_{X,X}'^{(1)}\circ \delta_{X,X}'^{(0)}
\end{split}
\end{equation}
\begin{equation}
\begin{split}
e^*:&\FF[GL_X]\rightarrow \FF \hsf\text{(co-unit)} \\
&\delta_{X,X}\mapsto id_X \\
&\delta_{X,X}'\mapsto id_X
\end{split}
\end{equation}
\begin{equation}
\begin{split}
S^*:&\FF[GL_X]\rightarrow \FF[GL_X]\;\;\;\;\text{(co-inverse (antipode))}\\
&\delta_{X,X}\mapsto \delta_{X,X}'\\
&\delta_{X,X}'\mapsto \delta_{X,X}
\end{split}
\end{equation}
These make $\FF[GL_X]$ into a group object in $(\FR)^{op}$.\\
Also if there is an involution then $GL_X:\FR^t\rightarrow Grps$
is represented by $\FF[GL_X^t]\in \FR^t$,
\begin{equation}
\slfrac{\FF[GL_X^t]=\FF[\delta_{X,X},\delta_{X,X}^t]\otimes \FF[\delta',(\delta_{X,X}')^t]}{\delta_{X,X}\circ \delta_{X,X}'\sim \delta_{X,X}'\circ \delta_{X,X}\sim id_X}\
\end{equation}

\section{The arithmetical surface: $\mathcal{L}=\FF(\NN)\otimes_{\FF} \FF(\NN)$}
\smallbreak

Consider $\mathcal{L}=\FF(\NN)\otimes_{\FF}\FF(\NN)$. As a particular example of $(\ref{A.2.11})$ we have:

\begin{equation}
\begin{tikzcd}
& \mathcal{L}^{cent'l}\arrow[two heads]{r}&\mathcal{L}^{com}\arrow[two heads]{rd}{\text{\underline{Thm}: not isomorphism}}[swap]{g}&\\
\mathcal{L}=\FF(\NN)\underset{\FF}{\otimes}\FF(\NN)\arrow[two heads ]{ru}\arrow[two heads]{rd}&&&\mathcal{L}^{tot.com}\arrow[two heads]{r}{d}&\FF(\NN) \\
&\mathcal{L}^{1-com}\arrow[two heads]{rru}{f}
\end{tikzcd}
\end{equation}

Note that the diagonal homomorphism $diag:
\mathcal{L}=\FF(\NN)\otimes_{\FF} \FF(\NN)\sur \FF(\NN)$, factors
through a surjection $d:\mathcal{L}^{tot-com}\sur \FF(\NN)$, since
$\FF(\NN)$ is totally commutative. We have the following

\thm{2.9.1} \textit{ (1) The composition $d\circ f:\mathcal{L}^{1-com}\sur
\FF(\NN)$ is an isomorphism (and hence both $d$ and $f$ are
isomorphisms).\\
(2) The composition $d\circ g:\mathcal{L}^{com}\rightarrow \FF(\NN)$
is \underline{not} an isomorphism. }

\begin{proof}
(1) Note that the $\fr \; \FF(\NN)$ is generated by the elements
$\sigma=(1,1)\in \FF(\NN)_{1,2}$ and $\sigma^t=\left(
\begin{array}{c}
1 \\
1 \end{array} \right)\in \FF(\NN)_{2,1} $, with certain relations:
$\FF(\NN)=\FF<\sigma,\sigma^t>$ . It follows that
$\mathcal{L}=\FF<\sigma,\sigma^t,\sigma',(\sigma')^t>$ is generated by
$\sigma, \sigma^t$ and $\sigma',(\sigma')^t$ coming from the left
and right factors of $\mathcal{L}=\FF(\NN)\otimes_{\FF} \FF(\NN)$.\\
From $1-commutativity$ we have:
$$\sigma\circ(\sigma'\oplus \sigma')=\sigma'\circ (\sigma\oplus \sigma)$$
$$\implies \sigma\circ (\sigma'\oplus \sigma')\circ \left(\left( \begin{array}{c}
0 \\
1 \end{array} \right)\oplus \left( \begin{array}{c}
1 \\
0 \end{array} \right)\right)=\sigma'\circ (\sigma \oplus
\sigma)\circ \left(\left( \begin{array}{c}
0 \\
1 \end{array} \right)\oplus \left( \begin{array}{c}
1 \\
0 \end{array} \right)\right) $$
$$\implies \sigma \circ \left(\sigma'\circ \left( \begin{array}{c}
0 \\
1 \end{array} \right)\oplus \sigma'\circ \left( \begin{array}{c}
1 \\
0 \end{array} \right)\right)=\sigma'\circ\left(\sigma\circ \left(
\begin{array}{c}
0 \\
1 \end{array} \right)\oplus \sigma'\circ \left( \begin{array}{c}
1 \\
0 \end{array} \right)\right) $$
$$\implies \sigma\circ (id_1\oplus id_1)=\sigma'\circ (id_1\oplus id_1)$$
\begin{equation}
\implies \sigma=\sigma\circ id_2 = \sigma'\oplus id_2=\sigma'.
\end{equation}
Thus $1-$commutativity imply $\sigma=\sigma'$, and similarly
$\sigma^t=(\sigma')^t$, and $d\circ
f:\mathcal{L}^{1\text{-com}}\rightarrow \FF(\NN)$ is an isomorphism.
\\
In a pictorial way we see this as follows:
$$\begin{tikzpicture}[scale=1,auto=left]
  \tikzstyle{no}=[circle,draw,fill=black!50,inner sep=0pt, minimum width=4pt]
  \begin{scope}
  \node[no] (n1) at (0,0) {};
  \node[no] (n2) at (1,1)  {};
  \node[no] (n3) at (1,-1)  {};
  \node[no] (n4) at (2,1.5) {};
  \node[no] (n5) at (2,0.5)  {};
  \node[no] (n6) at (2,-0.5){};
  \node[no] (n7) at (2,-1.5) {};

 \foreach \from/\to in {n1/n3,n1/n2}
    \draw (\from) -- (\to);

    \foreach \from/\to in {n2/n4,n2/n5,n3/n6,n3/n7}
    \draw (\from) -- (\to) [dashed];

 \end{scope}
 \begin{scope}[xshift=3cm]
  \node[no] (n8) at (0,0) {};
  \node[no] (n9) at (1,1)  {};
  \node[no] (n10) at (1,-1)  {};
  \node[no] (n11) at (2,1.5) {};
  \node[no] (n12) at (2,0.5)  {};
  \node[no] (n13) at (2,-0.5)  {};
  \node[no] (n14) at (2,-1.5) {};

    \foreach \from/\to in {n8/n10,n8/n9}
    \draw (\from) -- (\to) [dashed];

    \foreach \from/\to in {n9/n11,n9/n12,n10/n13,n10/n14}
    \draw (\from) -- (\to);
  \end{scope}

  \draw [white]  (2.3,-.4) -- (2.7,-.4)
    node [above=1mm,midway,black,text width=3cm,text centered]
      {\textbf{$=$}};
\end{tikzpicture} $$

$$ \Downarrow $$

$$\begin{tikzpicture}[scale=1,auto=left]
  \tikzstyle{no}=[circle,draw,fill=black!50,inner sep=0pt, minimum width=4pt]
  \begin{scope}
  \node[no] (n1) at (0,0) {};
  \node[no] (n2) at (1,1)  {};
  \node[no] (n3) at (1,-1)  {};
  \node[no] (n4) at (2,0.5)  {};
  \node[no] (n5) at (2,-0.5){};

 \foreach \from/\to in {n1/n2,n1/n3}
    \draw (\from) -- (\to);

    \foreach \from/\to in {n2/n4,n3/n5}
    \draw (\from) -- (\to) [dashed];

 \end{scope}
 \begin{scope}[xshift=3cm]
  \node[no] (n8) at (0,0) {};
  \node[no] (n9) at (1,1)  {};
  \node[no] (n10) at (1,-1)  {};
  \node[no] (n12) at (2,0.5)  {};
  \node[no] (n13) at (2,-0.5)  {};

    \foreach \from/\to in {n8/n9,n8/n10}
    \draw (\from) -- (\to) [dashed];

    \foreach \from/\to in {n9/n12,n10/n13}
    \draw (\from) -- (\to);
  \end{scope}

  \draw [white]  (2.3,-.4) -- (2.7,-.4)
    node [above=1mm,midway,black,text width=3cm,text centered]
      {$=$};
\end{tikzpicture} $$

$$ \Downarrow $$

\begin{equation}
\begin{tikzpicture}[scale=1,auto=left]
  \tikzstyle{no}=[circle,draw,fill=black!50,inner sep=0pt, minimum width=4pt]
  \begin{scope}
  \node[no] (n1) at (0,0) {};
  \node[no] (n2) at (1,1)  {};
  \node[no] (n3) at (1,-1)  {};

 \foreach \from/\to in {n1/n2,n1/n3}
    \draw (\from) -- (\to);
 \end{scope}

 \begin{scope}[xshift=3cm]
  \node[no] (n8) at (0,0) {};
  \node[no] (n9) at (1,1)  {};
  \node[no] (n10) at (1,-1)  {};

    \foreach \from/\to in {n8/n9,n8/n10}
    \draw (\from) -- (\to) [dashed];
  \end{scope}

  \draw [white]  (1.7,-.3) -- (2.3,-.3)
    node [above=1mm,midway,black,text width=3cm,text centered]
      {$=$};
\end{tikzpicture}
\end{equation}
(2) Recall that elements of $\mathcal{L}_{Y,X}=(\FF(\NN)\otimes_{\FF} \FF(\NN))_{Y,X}$, are given by a sequence of sets $X=X_0,X_1,\dots,X_l=Y$, and $X_{j+1}\times X_j$ -matrices over $\NN$, coming from the left $"l"$, or right $"r"$ copies of $\FF(\NN)$ (depending on the parity of $j$). Thinking of an $X_{j+1}\times X_j$ -matrix with values in $\NN$, $B_j$, as a set with maps
$$Z_j\overset{\pi_i}{\rightrightarrows} X_{j+i},  \hspace{10mm} i=0,1.$$
(via $(B_j)_{x',x}= \# \pi_1^{-1}(x')\cap \pi_0^{-1}(x) )$ , we obtain a graph with no loops
\begin{equation}
G=\{G^1=\underset{j}{\amalg} Z_j \underset{\pi_0}{\overset{\pi_1}{\rightrightarrows}} G^0=\underset{j}{\amalg} X_j  \},
\end{equation}
with a mapping $\mu: G^1\rightarrow \{l,r\}$.

Eliminating from our matrices rows $(i=1)$ , or columns $(i=0)$ , which are zero we may assume there are no such rows or columns.  \\
\underline{zero reduction}:
$$\text{If }\;\;\;\xi\in G^0\setminus (Y\amalg X)\hsm\text{is such that}\hsm \pi_{i}^{-1}(\xi)=\emptyset,\hsm i=0\;\text{or}\; i=1$$
$$ \implies G\sim (G^1\setminus \pi_{i-1}^{-1}(\xi),G^0\setminus \{\xi\}). $$
After a finite number of zero reductions we may assume with out loss of generality that $G$ is \underline{zero -reduced:}
every path in $CG$ extends to a maximal path
beginning in $In\;G\hookrightarrow X$ and ending in
$Out\;G\hookrightarrow Y$.\\
Set,
\begin{equation}
A_{Y,X}=\slfrac{\bigg\{\begin{matrix}
G=\{G^1\underset{\pi_1}{\overset{\pi_0}{\rightrightarrows}} G^0\},\;\; \mu:G^1\rightarrow
\{l,r\}\;,\; \mu(\pi_i^{-1}(x))\equiv \;l\;\;\;\text{or}\;\equiv \;\;r,\hsm\text{no loops \&} \\
\;\text{zero-reduced,}\hspace{2mm} In \;G=G^0\setminus \pi^1(G^1)\hookrightarrow X\;,\;
  Out\;G=G^0\setminus \pi^0(G^1)\hookrightarrow Y
\end{matrix}\bigg\}}{\sim}
\end{equation}\vspace{3mm} \\
If $(G,\mu)\in A_{Y,X},\; x',x\in G^0,\; e\in G^1$ are such that $\pi_0^{-1}(x')=\{e\}=\pi_1^{-1}(x)$, i.e. $e$ is the unique edge going out of $x'$, and also the unique edge going into $x$, than we can form the graph
\begin{equation}
G'= (G^1\setminus \{e\}\;,\; \slfrac{G^0}{\{x'\sim x\}}).
\end{equation}
by throwing out the edge $e$, and identifying the vertices $x'$ and $x$. We say $G'$ is obtained from $G$ by 1-reduction. If there are no such $x,x',e$ in $G$, we say $G$ is one-reduced. After a finite number of 1-reductions we obtain a one reduced graph.
If $(G,\mu)$ is both zero- $\&$ one- reduced, we say it is \underline{$\FF$ -reduced}. \vspace{2mm} \\
We can relax the condition on the matrices to alternate between $"l"$ and $"r"$, if we remember that consecutive matrices both $"l"$ or both $"r"$, are allowed to be multiplied, and in the description of $\NN$ -valued matrices as sets, matrix multiplication corresponds to taking fiber products.

Define \underline{$l/r$- reduction:} For a graph $G$, and $x\in G^0$ such that
\begin{equation}
\begin{matrix} \mu(\pi_1^{-1}(x))\equiv \mu(\pi_0^{-1}(x))\;\;\; (\equiv l\;\; \text{or}\;\;\equiv r)\vspace{2.5mm} \\ \implies G\sim \bar{G}= \bigg(\bigg[G^1\setminus (\pi_1^{-1}(x)\amalg \pi_0^{-1}(x))\bigg]\amalg \bigg[\pi_1^{-1}(x)\Pi \pi_0^{-1}(x)\bigg],\;\; G^0\setminus \{x\}\bigg)\end{matrix}
\end{equation}
$$\text{with:}\;\;\;\; (e,e')\in \pi_1^{-1}(x)\Pi \pi_0^{-1}(x),\; \pi_1(e,e')=\pi_1(e'),\; \pi_0(e,e')=\pi_0(e)$$
The inverse passage from $\bar{G}$ to $G$ will be called $l/r$ -inflation.
A Graph $(G,\mu) \in A_{Y,X}$ is then \underline{$l/r$ -reduced} when for all $x\in G^0\setminus (In\; G\cup Out\;G)$
\begin{equation}
\{\mu(\pi_1^{-1}(x)),\mu(\pi_0^{-1}(x))\}=\{l,r\}
\end{equation}
We have a canonical form to any $(G,\mu)\in A_{Y,X}$ obtained
after a finite number of $l/r$ -reduction and 1-reduction, and it
is characterized by the fact that $\#G^0$ is minimal, and so
actually:
\begin{equation}
\mathcal{L}=\FF(\NN)\otimes_{\FF}\FF(\NN)\equiv \{(G,\mu)\;,\; \FF\text{ -reduced}\;\text{ and }\;l/r\text{ -reduced} \}
\end{equation}
Next let us look at the commutative quotient of $\mathcal{L}$,
\begin{equation}\label{Lcom}
\mathcal{L}^{com}=\slfrac{\FF(\NN)\otimes_{\FF}\FF(\NN)}{\approx} \equiv \slfrac{\{(G,\mu)\;,\; \FF\text{ -reduced} \}}{\approx},
\end{equation}
where $\approx$ is the equivalence relation generated by $l/r$ -reduction, $l/r$ -inflation, and six commutative relations, passing from any vertice of the "commutativity triangle" to any other vertice, once one of the three patterns are recognized within our graph.

\underline{Commutativity triangle:} for $a\in A_{\bar{Y},\bar{X}}, b\in A_{1,J},d\in A_{J,1}$,

\begin{tikzpicture}[scale=0.2,auto=left]
\tikzstyle{no}=[circle,draw,fill=black!100,inner sep=0pt, minimum width=4pt]



\node[no] (n1) at (-1,3) {};
\node[no] (n2) at (-1,7)  {};
\node[no] (n3) at (-1,11)  {};
\node[no] (n4) at (-1,15) {};
\node[no] (n5) at (3,2)  {};
\node[no] (n6) at (3,3){};
\node[no] (n7) at (3,4) {};
\node[no] (n8) at (3,6){};
\node[no] (n9) at (3,7) {};
\node[no] (n10) at (3,8) {};
\node[no] (n11) at (3,10)  {};
\node[no] (n12) at (3,11){};
\node[no] (n13) at (3,12) {};
\node[no] (n14) at (3,14) {};
\node[no] (n15) at (3,15)  {};
\node[no] (n16) at (3,16){};
\node[no] (n17) at (5,15) {};
\node[no] (n18) at (5,11) {};
\node[no] (n19) at (5,7)  {};
\node[no] (n20) at (5,3)  {};
\node[no] (n21) at (16,6)  {};
\node[no] (n22) at (16,12)  {};

\foreach \from/\to in {n4/n14,n4/n16,n3/n11,n3/n13,n2/n8,n2/n10,n1/n7,n1/n5,n16/n17,n14/n17,n13/n18,n11/n18,n10/n19,n8/n19,n5/n20,n7/n20,n14/n16,n11/n13,n8/n10,n5/n7,n17/n22,n17/n20,n20/n21,n21/n22}
\draw[thick] (\from) -- (\to);


\node at (3,18) {$(\underset{\bar{Y}}{\oplus}(b\circ d))$};
\node at (11,18.5) {$\circ a$};
\node at (-2,9) {$\cdots$};
\node at (18,9) {$\cdots$};
\node at (20,11) {$\overset{\sim}{\longleftrightarrow}$};
\node at (22,9) {$\cdots$};


\node[no] (n1) at (24,3) {};
\node[no] (n2) at (24,7)  {};
\node[no] (n3) at (24,11)  {};
\node[no] (n4) at (24,15) {};
\node[no] (n5) at (33,12)  {};
\node[no] (n6) at (33,6){};
\node[no] (n7) at (37,11) {};
\node[no] (n8) at (37,12){};
\node[no] (n9) at (37,13) {};
\node[no] (n10) at (37,5) {};
\node[no] (n11) at (37,6)  {};
\node[no] (n12) at (37,7){};
\node[no] (n13) at (39,12) {};
\node[no] (n14) at (39,6) {};

\foreach \from/\to in {n1/n4,n1/n6,n5/n6,n4/n5,n6/n10,n6/n10,n6/n12,n5/n7,n5/n9,n7/n9,n10/n12,n7/n13,n9/n13,n10/n14,n12/n14}
\draw[thick] (\from) -- (\to);

\node at (36,16) {$(\underset{\bar{X}}{\oplus}(b\circ d))$};
\node at (29,16.5) {$a\circ $};


\node at (12,-1){$\searrow$};
\node at (11,0){$\nwarrow$};
\node at (11.9,0.1){\begin{turn}{140}
$\sim$
\end{turn}};

\node at (29,0){$\nearrow$};
\node at (28,-1){$\swarrow$};
\node at (28.4,-0.2){\begin{turn}{45}
$\sim$
\end{turn}};

\node[no] (n1) at (9,-8) {};
\node[no] (n2) at (9,-12)  {};
\node[no] (n3) at (9,-16)  {};
\node[no] (n4) at (9,-20) {};
\node[no] (n5) at (13,-7)  {};
\node[no] (n6) at (13,-8){};
\node[no] (n7) at (13,-9) {};
\node[no] (n8) at (13,-11){};
\node[no] (n9) at (13,-12) {};
\node[no] (n10) at (13,-13) {};
\node[no] (n11) at (13,-15)  {};
\node[no] (n12) at (13,-16){};
\node[no] (n13) at (13,-17) {};
\node[no] (n14) at (13,-19) {};
\node[no] (n15) at (13,-20)  {};
\node[no] (n16) at (13,-21){};

\node[no] (n17) at (17,-9) {};
\node[no] (n18) at (17,-11)  {};
\node[no] (n19) at (17,-14)  {};
\node[no] (n20) at (17,-16) {};
\node[no] (n21) at (18,-11)  {};
\node[no] (n22) at (18,-13){};
\node[no] (n23) at (18,-16) {};
\node[no] (n24) at (18,-18){};
\node[no] (n25) at (19,-13) {};
\node[no] (n26) at (19,-15) {};
\node[no] (n27) at (19,-18)  {};
\node[no] (n28) at (19,-20){};
\node[no] (n29) at (23,-11) {};
\node[no] (n30) at (23,-14) {};
\node[no] (n31) at (24,-13)  {};
\node[no] (n32) at (24,-16){};
\node[no] (n33) at (25,-15)  {};
\node[no] (n34) at (25,-18){};

\node[no] (n35) at (28,-9) {};
\node[no] (n36) at (28,-10)  {};
\node[no] (n37) at (28,-11)  {};
\node[no] (n38) at (28,-14) {};
\node[no] (n39) at (28,-15)  {};
\node[no] (n40) at (28,-16){};
\node[no] (n41) at (30,-10) {};
\node[no] (n42) at (30,-15){};

\foreach \from/\to in {n17/n20,n21/n24,n25/n28,n29/n30,n31/n32,n33/n34,n17/n29,n20/n30,n21/n31,n25/n33,n24/n32,n28/n34,n1/n5,n1/n7,n2/n8,n2/n10,n3/n11,n3/n13,n4/n14,n4/n16,n5/n7,n8/n10,n11/n13,n14/n16,n35/n37,n38/n40,n35/n41,n37/n41,n38/n42,n40/n42}
\draw[thick] (\from) -- (\to);

\foreach \from/\to in {n29/n35,n31/n36,n33/n37,n5/n17,n6/n21,n7/n25}
\draw (\from) -- (\to) [dashed];

\node at (11,-4){$(\underset{\bar{Y}}{\oplus}b)\circ$};
\node at (20,-4){$(\underset{J}{\oplus}a)\circ$};
\node at (29,-4){$(\underset{\bar{X}}{\oplus}d)$};


\node[no] (n1) at (0,-26) {};
\node[no] (n2) at (0,-33)  {};
\node[no] (n3) at (0,-28)  {};
\node[no] (n4) at (0,-31) {};
\node[no] (n5) at (6,-28)  {};
\node[no] (n6) at (6,-31){};
\node at (-2,-29.5) {$a=$};
\node at (9,-29.5) {$d=$};
\node[no] (n7) at (11,-28.5){};
\node[no] (n8) at (11,-29.5) {};
\node[no] (n9) at (11,-30.5) {};
\node[no] (n10) at (13,-29.5) {};
\node at (16,-29.5) {$b=$};
\node[no] (n11) at (18,-29.5){};
\node[no] (n12) at (22,-29.5) {};
\node[no] (n13) at (22,-28.5) {};
\node[no] (n14) at (22,-30.5)  {};

\foreach \from/\to in {n1/n2,n1/n5,n5/n6,n2/n6,n7/n9,n7/n10,n9/n10,n11/n13,n11/n14,n13/n14}
\draw[thick] (\from) -- (\to);

\node at (38,-29.5) {$\bar{X}=2,\;\bar{Y}=4,\;J=3.$};
\node at (-2,-22.5) {here have:};
\end{tikzpicture}

More specifically, $\approx$ is the equivalence relation such that $G\approx G'$ if and only if there is a path $G=G_0,\dots,G_j,\dots,G_l=G'$, where $\{G_{j-1},G_j\}$ is one of the following forms:\\
(1) $G_j$ is obtained from $G_{j-1}$ via $l/r$ -reduction.  \\
(2) $G_j$ is obtained from $G_{j-1}$ via $l/r$ -inflation.  \\
(3) $G_j$ is obtained from $G_{j-1}$ via (one out of $6$ possible) commutativity moves (see example above).\vspace{3mm}\\
\underline{\bf{Claim:}} $\approx$ is an equivalence ideal. \\
indeed, it is an equivalence relation and,
\begin{equation}
\begin{matrix} (i) & G\approx G' \implies G\oplus F\approx G'\oplus F \\
(ii)& G\approx G' \implies H\circ G\circ F\approx H\circ G'\circ F\end{matrix}
\end{equation}
The first implication $(i)$ is trivial, while the latter implication $(ii)$ follows by the fact that (unlike total commutativity and 1 -commutativity), the commutativity relation is \underline{$\FF$ -linear}:
\begin{equation}
\hspace{11mm}G\approx G'\implies [Y_0]G[X_0]\approx [Y_0]G'[X_0].
\end{equation}
Indeed, $\FF$- linearity implies $(ii)$ as (cf. (\ref{glue})):
$$H\circ G\circ F=H[Y_0]\underset{Y_0}{\coprod}[Y_0]G[X_0]\underset{X_0}{\coprod}[X_0]F.$$
To show $\FF$- linearity we have to show that if $G'$ is obtained from $G\in A_{Y,X}$ by any of the steps $(1),(2),(3)$, then so is $[Y_0]G'[X_0]$ obtained from $[Y_0]G[X_0]$. This is clear for steps $(1)$ and $(2)$. For the 6 possible commutativity moves, we have to identify within $G$ one of the 3 possible patterns, and in particular subsets $\bar{X},\bar{Y}\subseteq G^0$, and a subgraph $a\subseteq G$, $a$ from $\bar{X}$ to $\bar{Y}$.
If $G'$ is obtained from $G$ via commutativity relation with pattern $(\bar{Y},\bar{X},J,a,b,d)$,
then $[Y_0]G'[X_0]$ is obtained from $[Y_0]G[X_0]$ via the same commutativity relations with pattern $(\bar{Y}_0,\bar{X}_0,J,[\bar{Y}_0]a[\bar{X}_0],b,d)$, where
$$\bar{X}_0=\{x\in \bar{X}\;|\; \exists \text{ maximal path }e\in [Y_0]G[X_0]\; \text{with}\; x\in e \}\hspace{11mm}$$
\begin{equation}
\bar{Y}_0=\{y\in \bar{Y}\;|\; \exists \text{ maximal path }e\in [Y_0]G[X_0]\; \text{with}\; y\in e \}
\end{equation}
This proves that $\cong$ is $\FF$- linear, and hence an equivalence ideal.\vspace{1.5mm}\\

Note that commutative relations does not preserve $l/r$ -reduction, it is also needed to have $l/r$ -inflation to get full commutative relations, e.g.,

\vspace{4mm}

\begin{tikzpicture}[scale=0.3,auto=left]

\tikzstyle{no}=[circle,draw,fill=black!100,inner sep=0pt, minimum width=4pt]



\node[no] (n1) at (-2,7) {};
\node[no] (n2) at (-4,7.5)  {};
\node[no] (n3) at (-4,6.5)  {};
\node[no] (n4) at (1,9) {};
\node[no] (n5) at (1,5)  {};
\node[no] (n6) at (2.5,9.5){};
\node[no] (n7) at (2.5,8.5) {};
\node[no] (n8) at (2.5,5.5){};
\node[no] (n9) at (2.5,4.5) {};
\node[no] (n10) at (4,9) {};
\node[no] (n11) at (4,5)  {};
\node[no] (n12) at (6,9.5){};
\node[no] (n13) at (6,8.5) {};
\node[no] (n14) at (6,5.5) {};
\node[no] (n15) at (6,4.5)  {};

\foreach \from/\to in {n2/n1,n2/n3,n2/n1,n3/n1,n1/n4,n1/n5,n6/n10,n7/n10,n10/n12,n10/n13,n11/n14,n11/n15,n12/n13,n14/n15,n8/n11,n9/n11}
\draw[thick] (\from) -- (\to);

    \foreach \from/\to in {n4/n6,n4/n7,n5/n8,n5/n9}
    \draw (\from) -- (\to) [dashed];

\node[no] (n16) at (3,11.5)  {};
\node[no] (n17) at (1,12)  {};
\node[no] (n18) at (1,11)  {};

    \foreach \from/\to in {n10/n16,n16/n17,n16/n18,n17/n18}
    \draw[thick] (\from) -- (\to);

\draw [black, fill=gray] (3,11.5)-- (1,12) --(1,11)-- (3,11.5);
\draw [black, fill=gray] (-2,7)-- (-4,6.5) --(-4,7.5)-- (-2,7);
\draw [black, fill=gray] (6,9.5)-- (6,8.5) --(4,9)-- (6,9.5);
\draw [black, fill=gray] (6,4.5)-- (6,5.5) --(4,5)-- (6,4.5);

\node at (3,2.5) {$l/r$ -reduced};

\draw [thick, ->] (8,7.5) -- (18,7.5) ;

\node at (13,9) {$l/r$ -inflation};


\node[no] (n1) at (23,7) {};
\node[no] (n2) at (21,7.5)  {};
\node[no] (n3) at (21,6.5)  {};
\node[no] (n4) at (26,9) {};
\node[no] (n5) at (26,5)  {};
\node[no] (n6) at (27.5,9.5){};
\node[no] (n7) at (27.5,8.5) {};
\node[no] (n8) at (27.5,5.5){};
\node[no] (n9) at (27.5,4.5) {};
\node[no] (n10) at (29,9) {};
\node[no] (n11) at (29,5)  {};
\node[no] (n12) at (33,9.5){};
\node[no] (n13) at (33,8.5) {};
\node[no] (n14) at (31,5.5) {};
\node[no] (n15) at (31,4.5)  {};
\node[no] (n16) at (28.5,11.5)  {};
\node[no] (n17) at (26.5,12)  {};
\node[no] (n18) at (26.5,11)  {};
\node[no] (n19) at (31,9)  {};

\node at (24.2,7) {$b$};
\node at (26,10) {$a$};
\node at (28.5,10) {$d$};
\node at (26,6) {$a$};
\node at (28.5,6) {$d$};

\foreach \from/\to in {n2/n1,n2/n3,n2/n1,n3/n1,n1/n4,n1/n5,n6/n10,n7/n10,n10/n19,n11/n14,n11/n15,n12/n13,n14/n15,n8/n11,n9/n11,n19/n12,n19/n13,n17/n18,n16/n17,n16/n18,n19/n16}
\draw[thick] (\from) -- (\to);

\foreach \from/\to in {n4/n6,n4/n7,n5/n8,n5/n9}
\draw (\from) -- (\to) [dashed];

\draw [black, fill=gray] (23,7)-- (21,7.5) --(21,6.5)-- (23,7);
\draw [black, fill=gray] (28.5,11.5)-- (26.5,12) --(26.5,11)-- (28.5,11.5);
\draw [black, fill=gray] (31,9)-- (33,9.5) --(33,8.5)-- (31,9);
\draw [black, fill=gray] (29,5)-- (31,5.5) --(31,4.5) -- (29,5);

\draw [thick, ->] (28,3) -- (28,-1) ;
\draw [thick, ->] (17,6) -- (6,0) ;
\node at (15,1.5) {comm.:};
\node at (15,0) {$(\underset{1}{\oplus} b)(\underset{2}{\oplus} a)(\underset{2}{\oplus} d)=b\circ d\circ a$};
\node at (19.6,-1.5) {$=a \circ (\underset{2}{\oplus}b\circ d)$};






\node[no] (n1) at (-2,-5) {};
\node[no] (n2) at (-4,-4.5)  {};
\node[no] (n3) at (-4,-5.5)  {};
\node[no] (n4) at (1,-3) {};
\node[no] (n5) at (1,-7)  {};
\node[no] (n6) at (2.5,-2.5){};
\node[no] (n7) at (2.5,-3.5) {};
\node[no] (n8) at (2.5,-6.5){};
\node[no] (n9) at (2.5,-7.5) {};
\node[no] (n10) at (4,-3) {};
\node[no] (n11) at (4,-7)  {};
\node[no] (n12) at (8,-2.5){};
\node[no] (n13) at (8,-3.5) {};
\node[no] (n14) at (6,-6.5) {};
\node[no] (n15) at (6,-7.5)  {};
\node[no] (n16) at (3.5,-0.5)  {};
\node[no] (n17) at (1.5,0)  {};
\node[no] (n18) at (1.5,-1)  {};
\node[no] (n19) at (6,-3)  {};

\foreach \from/\to in {n2/n3,n4/n6,n4/n7,n5/n8,n5/n9,n2/n1,n3/n1,n6/n10,n7/n10,n10/n19,n11/n14,n11/n15,n12/n13,n14/n15,n8/n11,n9/n11,n19/n12,n19/n13,n17/n18,n16/n17,n16/n18,n19/n16}
\draw[thick] (\from) -- (\to);

\foreach \from/\to in {n1/n4,n1/n5}
\draw (\from) -- (\to) [dashed];

\draw [black, fill=gray] (-2,-5)-- (-4,-4.5) --(-4,-5.5)-- (-2,-5);
\draw [black, fill=gray] (6,-3)-- (8,-2.5) --(8,-3.5)-- (6,-3);
\draw [black, fill=gray] (4,-7)-- (6,-6.5) --(6,-7.5) -- (4,-7);
\draw [black, fill=gray] (3.5,-0.5)-- (1.5,0) --(1.5,-1)-- (3.5,-0.5);

\node at (-0.5,-5) {$a$};
\node at (1.5,-2) {$b$};
\node at (3.5,-2) {$d$};
\node at (1.5,-6) {$b$};
\node at (3.5,-6) {$d$};


\node[no] (n1) at (22.5,-6) {};
\node[no] (n2) at (24.5,-5) {};
\node[no] (n3) at (24.5,-7) {};
\node[no] (n4) at (26.5,-6) {};
\node[no] (n5) at (28.5,-5) {};
\node[no] (n6) at (28.5,-7) {};
\node[no] (n7) at (30.8,-5) {};
\node[no] (n8) at (28.5,-3){};
\node[no] (n10) at (20.5,-5.5)  {};
\node[no] (n11) at (20.5,-6.5)  {};
\node[no] (n12) at (30.5,-6.5) {};
\node[no] (n13) at (30.5,-7.5) {};
\node[no] (n15) at (26.5,-2.5)  {};
\node[no] (n16) at (26.5,-3.5)  {};
\node[no] (n17) at (32.8,-4.5)  {};
\node[no] (n18) at (32.8,-5.5) {};

\draw [black, fill=gray] (20.5,-5.5)-- (20.5,-6.5) --(22.5,-6)-- (20.5,-5.5);
\draw [black, fill=gray] (30.5,-6.5)-- (30.5,-7.5) --(28.5,-7)-- (30.5,-6.5);
\draw [black, fill=gray] (26.5,-2.5)-- (26.5,-3.5) --(28.5,-3) -- (26.5,-2.5);
\draw [black, fill=gray] (32.8,-4.5)-- (32.8,-5.5) --(30.8,-5)-- (32.8,-4.5);

\foreach \from/\to in {n2/n4,n3/n4,n5/n7,n7/n8,n1/n2,n1/n3}
\draw[thick] (\from) -- (\to);

\foreach \from/\to in {n4/n5,n4/n6}
\draw (\from) -- (\to) [dashed];

\node at (23.7,-6) {$b$};
\node at (25.2,-6) {$d$};
\node at (28,-6) {$a$};

\draw [thick, ->] (28,-8) -- (28,-12) ;
\draw [thick, ->] (1,-8) -- (1,-12) ;
\node at (14,-10) {$l/r$ -reduction};


\node[no] (n1) at (20.5,-17) {};
\node[no] (n2) at (20.5,-18) {};
\node[no] (n3) at (22.5,-17.5) {};
\node[no] (n4) at (26.5,-17.5) {};
\node[no] (n5) at (29.5,-19.5) {};
\node[no] (n6) at (29.5,-15.5) {};
\node[no] (n7) at (31.5,-15.5) {};
\node[no] (n8) at (31.5,-19) {};
\node[no] (n9) at (31.5,-20) {};
\node[no] (n10) at (33.5,-15) {};
\node[no] (n11) at (33.5,-16) {};
\node[no] (n12) at (30,-13) {};
\node[no] (n13) at (28,-13.5) {};
\node[no] (n14) at (28,-12.5) {};

\foreach \from/\to in {n1/n2,n2/n3,n1/n3,n5/n8,n6/n7,n5/n9,n8/n9,n7/n10,n7/n11,n10/n11,n12/n13,n12/n14,n13/n14,n7/n12}
\draw[thick] (\from) -- (\to);

\foreach \from/\to in {n4/n5,n4/n6}
\draw (\from) -- (\to) [dashed];

\draw [black, fill=gray] (20.5,-17)-- (20.5,-18) --(22.5,-17.5)-- (20.5,-17);
\draw [black, fill=gray] (29.5,-19.5)-- (31.5,-19) --(31.5,-20)-- (29.5,-19.5);
\draw [black, fill=gray] (30,-13)-- (28,-13.5) --(28,-12.5) -- (30,-13);
\draw [black, fill=gray] (33.5,-16)-- (33.5,-15) --(31.5,-15.5)-- (33.5,-16);


\draw [-,thick, black] (22.5,-17.5) to [out=30,in=150] (26.5,-17.5) ;
\draw [-,thick, black] (22.5,-17.5) to [out=-30,in=-150] (26.5,-17.5) ;


\node[no] (n1) at (-4,-17) {};
\node[no] (n2) at (-4,-18) {};
\node[no] (n3) at (-2,-17.5) {};
\node[no] (n4) at (1,-15.5) {};
\node[no] (n5) at (1,-19.5) {};
\node[no] (n6) at (5,-15.5) {};
\node[no] (n7) at (5,-19.5) {};
\node[no] (n8) at (7,-15) {};
\node[no] (n9) at (7,-16) {};
\node[no] (n10) at (7,-20) {};
\node[no] (n11) at (7,-19) {};
\node[no] (n12) at (4,-13) {};
\node[no] (n13) at (2,-13.5) {};
\node[no] (n14) at (2,-12.5) {};

\foreach \from/\to in {n1/n2,n2/n3,n1/n3,n6/n8,n6/n9,n8/n9,n7/n10,n7/n11,n10/n11,n6/n12,n12/n13,n12/n14,n13/n14}
\draw[thick] (\from) -- (\to);

\foreach \from/\to in {n3/n4,n3/n5}
\draw (\from) -- (\to) [dashed];


\draw [-,thick, black] (1,-15.5) to [out=30,in=150] (5,-15.5) ;
\draw [-,thick, black] (1,-15.5) to [out=-30,in=-150] (5,-15.5) ;

\draw [-,thick, black] (1,-19.5) to [out=30,in=150] (5,-19.5) ;
\draw [-,thick, black] (1,-19.5) to [out=-30,in=-150] (5,-19.5) ;

\draw [black, fill=gray] (-4,-17)-- (-4,-18) --(-2,-17.5)-- (-4,-17);
\draw [black, fill=gray] (5,-15.5)-- (7,-15) --(7,-16)-- (5,-15.5);
\draw [black, fill=gray] (5,-19.5)-- (7,-20) --(7,-19) -- (5,-19.5);
\draw [black, fill=gray] (4,-13)-- (2,-13.5)--(2,-12.5)-- (4,-13);

\end{tikzpicture}

\vspace{4mm}
Now using the presentation (\ref{Lcom}) for $\mathcal{L}^{com}$ it follows that,
$$d\circ g:\mathcal{L}^{com}\sur \FF(\NN)\;\text{ is an isomorphism }\; $$
$$\iff$$
$$\sigma=((1,1),\;l)\approx \sigma'=((1,1),\;r)$$
$$\iff$$
$$\exists\;\; \text{path}\;\;(G_0,\mu_0),\;\dots\;,(G_k,\mu_k):\; \text{with } \{G_{j-1},G_j\}\; \text{of the form }(1),(2),\text{or } (3),$$

\begin{equation}
\begin{tikzpicture}[scale=0.3,auto=left]
\tikzstyle{no}=[circle,draw,fill=black!100,inner sep=0pt, minimum
width=4pt]

\node[no](n1) at (2,6){};
\node[no](n2) at (4,4){};
\node[no](n3) at (4,8){};
\node[no](n4) at (15,6){};
\node[no](n5) at (17,4){};
\node[no](n6) at (17,8){};

\foreach \from/\to in {n1/n2,n1/n3}
\draw[thick] (\from) -- (\to);

\foreach \from/\to in {n4/n5,n4/n6}
\draw (\from) -- (\to) [dashed];

\node [left] at (2,6) {1};
\node [right] at (4,8) {1};
\node [right] at (4,4) {2};
\node [left] at (15,6) {1};
\node [right] at (17,8) {1};
\node [right] at (17,4) {2};

\node [above] at (3,7) {$l$};
\node [below] at (2.7,5) {$l$};
\node [above] at (16,7) {$r$};
\node [below] at (15.8,5) {$r$};



\node at (7,6) {\large{$,\dots,$}};
\node at (11.2,6) {$\tiny G_k=\sigma'=$};
\node at (-3,6) {and $G_0=\sigma=$};
\node at (19,6) {\large{$,$}};
\end{tikzpicture}
\end{equation}

where for any $0 \leq j\leq k$, $\;Out\; G_j\equiv\{1\},In\;
G_j\equiv\{1,2\}$, and $\exists !$ path from $i$ to 1, $i=1,2$.
(its image under the diagonal homomorphism is $(1,1)$), and $G_j$
has no circuits and is zero -reduced. This implies that $G_j$ must
be of the following form:

\begin{equation}
\begin{tikzpicture}[scale=0.3,auto=left]
\tikzstyle{no}=[circle,draw,fill=black!100,inner sep=0pt, minimum
width=4pt]

\node at (0,5){$G_j\equiv$};
\node at (24,5){$i',j',k'\geq 0$};

\node[no](n1) at (2,5){};
\node[no](n2) at (4,5){};
\node[no](n3) at (6,5){};
\node[no](n4) at (8,5){};
\node[no](n5) at (10,5){};

\node[no](n6) at (12,6){};
\node[no](n7) at (14,6){};
\node[no](n8) at (16,6){};
\node[no](n9) at (18,6){};
\node[no](n10) at (20,6){};

\node[no](n11) at (12,4){};
\node[no](n12) at (14,4){};
\node[no](n13) at (16,4){};
\node[no](n14) at (18,4){};
\node[no](n15) at (20,4){};

\foreach \from/\to in {n1/n2,n2/n3,n3/n4,n4/n5,n5/n6,n5/n11}
\draw[thick] (\from) -- (\to);

\foreach \from/\to in {n6/n7,n7/n8,n8/n9,n9/n10,n11/n12,n12/n13,n13/n14,n14/n15}
\draw (\from) -- (\to) [dashed];

\draw [|-] (2,4) -- (5,4);
\draw [-|] (7,4) -- (10,4);
\node at(6,4) {$k'$};

\draw [|-] (12,3) -- (15,3);
\draw [-|] (17,3) -- (20,3);
\node at(16,3) {$i'$};

\draw [|-] (12,7) -- (15,7);
\draw [-|] (17,7) -- (20,7);
\node at(16,7) {$j'$};
\end{tikzpicture}
\end{equation}
Therefore such a path from $\sigma$ to $\sigma'$ is not possible.
\end{proof}

\section{Generators and relations for $\FF$-rings}
\label{genandrel}
\smallbreak

We have a surjection
$$\FF[\delta_{1,2},\delta_{1,2}^t]\sur \FF(\NN)$$
\begin{equation}
 \delta_{1,2}\mapsto (1,1)
\end{equation}
and a surjection
\begin{equation}
\FF\{\pm 1\}[\delta_{1,2},\delta_{1,2}^t]\sur \FF(\ZZ).
\end{equation}
Here $\FF\{\pm1\}[\delta_{1,2},\delta_{1,2}^t]=\FF\{\pm1\}\otimes_{\FF}\FF[\delta_{1,2},\delta_{1,2}^t].$\\
Indeed, if $A\subseteq\FF(\NN)$, (resp. $A\subseteq \FF(\ZZ)$) is a sub-$\FF$-ring and contains $(1,1)\in A_{1,2},\begin{pmatrix}
  1  \\ 1    \end{pmatrix}\in A_{2,1}$, than $A_{Y,X}$ is closed under addition: $a,a'\in A_{Y,X}$:
\begin{equation}
a+a'=\bigg( \underset{Y}{\bigoplus} (1,1)\bigg)\circ (a\oplus a')\circ \bigg( \underset{X}{\bigoplus}\begin{pmatrix}
  1  \\ 1    \end{pmatrix}\bigg) \in A_{Y,X}
\end{equation}
And any matrix in $\FF(\NN)_{Y,X}$ (resp. $\FF(\ZZ)_{Y,X}$) is a sum of matrices in $\FF$ (resp. $\FF\{\pm 1\}$).

Given $R\in \mathcal{R}ig^{(t)},\;A=\FF(R)\in \FR^{(t)}$. We have a surjective homomorphism:
$$\Phi:B=\FF[\delta_{1,2},\delta_{1,2}^t;\delta_{1,1}(r),r\in R]\sur A=\FF(R),$$
\begin{equation}
\Phi(\delta_{1,2})=(1,1), \Phi(\delta_{1,2}^t)=\begin{pmatrix}
  1 \\
  1
\end{pmatrix}, \Phi(\delta_{1,1}(r))=(r).
\end{equation}
Elements of $B_{Y,X}$ can be represented by graphs with no loops $G$, with $In\;G\inj X, Out\;G\inj Y$,with a map $\mu:G^1\rightarrow R$ and the graph $G$ is made of the basic graphs

\tikzstyle{vertex}=[circle, draw,fill=black!50, inner sep=0pt, minimum size=6pt]
\begin{equation}
\begin{tikzpicture}[x=1.3cm, y=1cm,
    every edge/.style={
        draw,
        postaction={decorate,
                    decoration={markings,mark=at position 0.5 with {\arrow{>}}}
                   }
        }
]
    \vertex (ul1) at (-0.2,0.5) {};
    \vertex (ur1) at (0.8,1) {};
    \vertex (lr1) at (0.8,0) {};
           \node (lll1) at (-1,0.5)  {$\delta\equiv$};
    \node (l1) at (-0.5,0.5)  {$\bigg \{$};
    \node (r1) at (1.2,0.5)  {$\bigg \}$};
    \path
        (ur1) edge  node[auto,swap] {$1$}(ul1)
        (lr1) edge   node[auto] {$1$}(ul1)
        ;
\vertex (ul2) at (4.1,0.5) {};
    \vertex (ur2) at (3.1,1) {};
    \vertex (lr2) at (3.1,0) {};
           \node (lll2) at (2.1,0.5)  {$\delta^t\equiv$};
    \node (l2) at (2.6,0.5)  {$\bigg \{$};
    \node (r2) at (4.6,0.5)  {$\bigg \}$};
    \path
        (ul2) edge  node[auto] {$1$}(lr2)
        (ul2) edge  node[auto,swap] {$1$}(ur2)
    ;
\vertex (r3) at (7.1,0.5) {};
    \vertex (l3) at (6.1,0.5) {};
           \node (lll3) at (5.3,0.5)  {$(r)\equiv$};
    \node (lb3) at (5.9,0.5)  {$\bigg \{$};
    \node (rb3) at (7.4,0.5)  {$\bigg \}$};
    \node (t3) at (6.6,0.7)  {$r$};
    \path
        (r3) edge (l3)
    ;
    \node (a) at (1.6,0.5)  {,};
    \node (b) at (4.8,0.5)  {,};
\end{tikzpicture}
\end{equation}

The homomorphism $\Phi$ takes $(G,\mu)\in B_{Y,X}$ into the $Y\times X$ -matrix with values in $R$,
\begin{equation}
\Phi (G,\mu)_{y,x}=\sum_{\begin{matrix} (e_k,\dots, e_1)\in m(G) \\ \pi_1(e_k)=y \\ \pi_0(e_1)=x \end{matrix}} \mu(e_k)\dots \mu(e_1).
\end{equation}

The equivalence ideal $\cK (\Phi)=B\underset{A}{\prod}B$ contains the following elements:\\
i. \underline{(1)}:  $(1)\equiv 1 \hspace{7.5mm},\hspace{7.5mm} (0)\equiv 0.$ \\
ii. $  \begin{cases}
  \underline{Zero}: \delta\circ \left( \begin{matrix} 1\\ 0 \end{matrix}\right)\equiv (1). \\
  \underline{Zero^t}: (1,0)\circ \delta^t\equiv (1).
\end{cases}$ \\
iii.  $\begin{cases}
  \underline{Ass}: \delta\circ( \delta\oplus (1))\equiv \delta\circ ((1)\oplus \delta). \\
  \underline{Ass^t}:  (\delta^t\oplus (1))\circ \delta^t\equiv ((1)\oplus \delta^t)\circ \delta^t  .
\end{cases}$ \\
iv.  $\begin{cases}
  \underline{Comm}: \delta\circ \left( \begin{matrix} 0 & 1\\ 1& 0 \end{matrix}\right)\equiv \delta. \\
  \underline{Comm^t}:   \left( \begin{matrix} 0 & 1\\ 1& 0 \end{matrix}\right)\circ\delta^t \equiv \delta^t  .
\end{cases}$ \\
v. \underline{Total - commutativity}: $\delta^t\circ \delta\equiv (\delta\oplus \delta)\circ (\delta^t\oplus \delta^t).$ \\
i.e. we have,
\begin{equation}
\begin{tikzpicture}[x=1.3cm, y=1cm,
    every edge/.style={
        draw,
        postaction={decorate,
                    decoration={markings,mark=at position 0.5 with {\arrow{>}}}
                   }
        }
]
           \vertex(lu) at (-1,1) [label=above:$1$]{};
     \vertex (ld) at (-1,-1) [label=above:$2$] {};
           \vertex (lm) at (0,0) {};
           \vertex (lru) at (1,1) [label=above:$1$] {};
    \vertex (lrd) at (1,-1) [label=above:$2$]{};
    \node(m) at (2,0){$\equiv$};
    \vertex (ru) at (3,1) [label=above:$1$] {};
    \vertex (rd) at (3,-1) [label=above:$2$]{};
           \vertex (rmu) at (4,1.5) [label=above:${(1,1)}$] {};
           \vertex (rmmu) at (4,0.5) [label=above:${(1,2)}$] {};
    \vertex (rmmd) at (4,-0.5) [label=above:${(2,1)}$]{};
    \vertex (rmd) at (4,-1.5) [label=above:${(2,2)}$] {};
    \vertex (rru) at (5,0.5) [label=above:$1$] {};
    \vertex (rrd) at (5,-0.5) [label=above:$2$]{};

\path
        (lm) edge (lu)
        (lm) edge (ld)
        (lru) edge (lm)
        (lrd) edge (lm)

        (rmu) edge (ru)
        (rmmu) edge (ru)
        (rmmd) edge (rd)
        (rmd) edge (rd)

        (rru) edge (rmu)
        (rru) edge (rmmd)
        (rrd) edge (rmmu)
        (rrd) edge (rmd)
        ;

\end{tikzpicture}
\end{equation}
the $\Phi$- image of this relation is
\begin{equation}
\left( \begin{matrix} 1\\ 1\end{matrix}\right)\circ (1,1)\equiv \left( \begin{matrix} 1&1 \\ 1&1 \end{matrix}\right)\equiv \left( \begin{matrix} 1 & 1& 0&0\\ 0& 0&1& 1 \end{matrix}\right)\circ \left( \begin{matrix} 1 & 0\\ 0& 1\\ 1& 0\\ 0 &1  \end{matrix}\right)
\end{equation}
The indexing of the input of $\delta$, and the output of $\delta^t$, is important; the right hand side of (v) is Not equal to
$$(\delta\circ \delta^t)\oplus (\delta\circ \delta^t), \text{whose } \Phi \text{-image is } \left( \begin{matrix} 2 & 0\\ 0& 2 \end{matrix}\right).$$

\noindent We have the relations,

vi. $\begin{cases}
\underline{(r,\delta)}:&
(r)\circ \delta \equiv \delta \circ ((r)\oplus (r)). \\
\underline{(\delta^t,r)}:&
\delta^t\circ (r)\equiv ((r)\oplus (r))\circ \delta^t \\
\underline{(r_1\cdot r_2)}:&
(r_1)\circ (r_2)\equiv (r_1\cdot r_2) \\
\underline{(r_1+r_2)}:&
\delta\circ ((r_1)\oplus (r_2))\circ \delta^t\equiv (r_1+r_2)
\end{cases}$

\rem{ }When working in the context of $\FR^t$, (with involution!), every relation is equivalent to its transpose, and we should add the relation,
$$(r)^t\equiv (r^t).$$
\thm{2.10.1} \textit{The equivalence ideal $\cK (\Phi)$ is generated by these relations.}
\begin{proof} Let $G\in B$, we shall show that modulo these relations we can bring $G$ to a canonical form which depends only on $\Phi(G)$. By using $(1)$, we can add to $G$ identities $\{   \begin{tikzpicture}[x=1.3cm, y=1cm,
    every edge/.style={
        draw,
        postaction={decorate,
                    decoration={markings,mark=at position 0.5 with {\arrow{>}}}
                   }
        }
]
     \vertex(l) at (0,0){};
     \vertex (r) at (1,0){};
\path

    (r) edge  node[above]{$1$} (l)
;
\end{tikzpicture}
\}$, and assume without loss of generality that $G$ has the form

\begin{equation}
\begin{tikzpicture}

\node (A) at (0.5,0) {$Z_l$};
\node (B) at (1,0) {$\coprod$};
\node (C) at (2.1,0) {$\dots$};
\node (D) at (3.5,0) {$Z_2$};
\node (E) at (4,0) {$\coprod$};
\node (F) at (4.5,0) {$Z_1$};

\node (a1) at (-1,-1) {$Y$};
\node (a2) at (-0.5,-1) {$\supseteq$};
\node (a3) at (-0.1,-1) {$X_l$};
\node (a4) at (0.5,-1) {$\coprod$};
\node (a5) at (1.1,-1) {$X_{l-1}$};
\node (a6) at (1.7,-1) {$\coprod$};
\node (a7) at (2.1,-1) {$\dots$};
\node (a8) at (2.5,-1) {$\coprod$};
\node (a9) at (3,-1) {$X_2$};
\node (a10) at (3.5,-1) {$\coprod$};
\node (a11) at (4,-1) {$X_1$};
\node (a12) at (4.5,-1) {$\coprod$};
\node (a13) at (5,-1) {$X_0$};
\node (a14) at (5.5,-1) {$\subseteq$};
\node (a15) at (6,-1) {$X$};

\path[->,font=\scriptsize,>=angle 90]
(A) edge   (a3)
(A) edge   (a5)
(D) edge   (a9)
(D) edge   (a11)
(F) edge   (a11)
(F) edge   (a13)

;

\end{tikzpicture}
\end{equation}
where each basic graph $\{Z_j\rightrightarrows X_j\coprod X_{j-1}\}$ has the form

\begin{equation}
\begin{tikzpicture}[x=1.3cm, y=1cm,
    every edge/.style={
        draw,
        postaction={decorate,
                    decoration={markings,mark=at position 0.5 with {\arrow{>}}}
                   }
        }
]
     \vertex(l1) at (-0.5,0)[label=above:$X_j$]{};
     \vertex (r1) at (1.5,0)[label=above:$X_{j-1}$]{};

    \vertex(l2) at (-0.5,-0.5){};
     \vertex (r2) at (1.5,-0.5){};

    \node (a) at (0.5,-0.85) {$\vdots$};

    \vertex(l3) at (-0.5,-1.4){};
     \vertex (r3) at (1.5,-1.2){};
    \vertex (r4) at (1.5,-1.6){};

    \node (a) at (0.5,-1.95) {$\vdots$};

    \vertex(l4) at (-0.5,-2.6){};
    \vertex (r5) at (1.5,-2.4){};
    \vertex (r6) at (1.5,-2.8){};

    \vertex (r7) at (1.5,-3.2){};
    \vertex(l5) at (-0.5,-3){};
    \vertex(l6) at (-0.5,-3.4){};

    \vertex(l7) at (-0.5,-3.8){};
    \vertex(l8) at (-0.5,-4.2){};
     \vertex (r8) at (1.5,-4){};

\path

    (r1) edge  node[above]{$r_1$} (l1)
    (r2) edge  node[above]{$r_2$} (l2)
    (r3) edge (l3)
    (r4) edge (l3)
    (r5) edge (l4)
    (r6) edge (l4)
    (r7) edge (l5)
    (r7) edge (l6)
    (r8) edge (l7)
    (r8) edge (l8)
;
\end{tikzpicture}
\end{equation}
\vspace{1.5mm}

i.e. is a direct sum of $\delta,\delta^t,(r)$, with no composition. By further adding identities, we can assume each basic graph is either of the "left"form \\
-a direct sum of $\delta, (r)$'s, (no $\delta^t$), or of the "right" form \\
-a direct sum of $\delta^t$, $(r)$'s, (no $\delta$). By using \underline{$(r,\delta)$}, \underline{$(\delta^t,r)$} and total -comm. we can replace $Z_{j+1}$ of the "right" form, $Z_j$ of the "left" form, by a $Z_{j+1}'$ of the "left" form, and $Z_j'$ of the "right" form:

\begin{equation}
\begin{tikzpicture}[x=1.3cm, y=1cm,
    every edge/.style={
        draw,
        postaction={decorate,
                    decoration={markings,mark=at position 0.5 with {\arrow{>}}}
                   }
        }
]

     \vertex(l1a) at (-1,0){};
     \vertex (l1b) at (0,0){};
     \vertex(l1c) at (1,0){};

    \vertex (l2a) at (-1,-1.25){};
    \vertex (l2b) at (0,-1.25){};
    \vertex (l2c) at (1,-1){};
    \vertex (l2d) at (1,-1.5){};

    \vertex (l3a) at (-1,-2.4){};
    \vertex (l3b) at (-1,-2.9){};
    \vertex (l3c) at (0,-2.65){};
    \vertex (l3d) at (1,-2.65){};

    \vertex (l4a) at (-1,-3.8){};
    \vertex (l4b) at (-1,-4.5){};
    \vertex (l4c) at (0,-4.15){};
    \vertex (l4d) at (1,-3.8){};
    \vertex (l4e) at (1,-4.5){};

    \node (m) at (2,-2.05) {$\implies$};

     \vertex(r1a) at (3,0){};
     \vertex (r1b) at (4,0){};
     \vertex(r1c) at (5,0){};

    \vertex (r2a) at (3,-1.25){};
    \vertex (r2b) at (4,-1){};
    \vertex (r2c) at (4,-1.5){};
    \vertex (r2d) at (5,-1){};
    \vertex (r2e) at (5,-1.5){};

    \vertex (r3a) at (3,-2.4){};
    \vertex (r3b) at (3,-2.9){};
    \vertex (r3c) at (4,-2.4){};
    \vertex (r3d) at (4,-2.9){};
    \vertex (r3e) at (5,-2.65){};

    \vertex (r4a) at (3,-3.8){};
    \vertex (r4b) at (3,-4.5){};
    \vertex (r4c) at (4,-3.6){};
    \vertex (r4d) at (4,-4){};
    \vertex (r4e) at (4,-4.3){};
    \vertex (r4f) at (4,-4.7){};
    \vertex (r4g) at (5,-3.8){};
    \vertex (r4h) at (5,-4.5){};

    \node(n1) at (-1,1){$\underline{X_{j+1}}$};
    \node(n2) at (-0.5,1.5){$\underline{Z_{j+1}}$};
    \node(n3) at (0,1){$\underline{X_j}$};
    \node(n4) at (0.5,1.5){$\underline{Z_j}$};
    \node(n5) at (1,1){$\underline{X_{j-1}}$};
    \node(n6) at (3,1){$\underline{X_{j+1}}$};
    \node(n7) at (3.5,1.5){$\underline{Z_{j+1}'}$};
    \node(n8) at (4.5,1.5){$\underline{Z_j'}$};
    \node(n9) at (5,1){$\underline{X_{j-1}}$};

    \node (d1) at (-0.5,-0.35) {$\vdots$};
    \node (d2) at (-0.5,-1.75) {$\vdots$};
    \node (d3) at (-0.5,-3.1) {$\vdots$};
    \node (d4) at (0.5,-0.35) {$\vdots$};
    \node (d5) at (0.5,-1.75) {$\vdots$};
    \node (d6) at (0.5,-3.1) {$\vdots$};
    \node (d7) at (3.5,-0.35) {$\vdots$};
    \node (d8) at (4.5,-0.35) {$\vdots$};
    \node (d9) at (3.5,-1.75) {$\vdots$};
    \node (d10) at (4.5,-1.75) {$\vdots$};
    \node (d11) at (3.5,-3.1) {$\vdots$};
    \node (d12) at (4.5,-3.1) {$\vdots$};

\path

    (l1b) edge node[above]{$r_1$} (l1a)
    (l1c) edge node[above]{$r_2$} (l1b)

    (l2b) edge node[above]{$r$} (l2a)
    (l2c) edge (l2b)
    (l2d) edge (l2b)

    (l3c) edge (l3a)
    (l3c) edge (l3b)
    (l3d) edgenode[above]{$r$} (l3c)

    (l4c) edge (l4a)
    (l4c) edge (l4b)
    (l4d) edge (l4c)
    (l4e) edge (l4c)

    (r1b) edge node[above]{$r_1$} (r1a)
    (r1c) edge node[above]{$r_2$} (r1b)

    (r2b) edge (r2a)
    (r2c) edge (r2a)
    (r2d) edge (r2b)
    (r2d) edge node[above]{$r$}(r2b)
    (r2e) edge node[above]{$r$}(r2c)

    (r3c) edge  node[above]{$r$}(r3a)
    (r3d) edge  node[above]{$r$}(r3b)
    (r3e) edge (r3c)
    (r3e) edge (r3d)

    (r4c) edge (r4a)
    (r4d) edge (r4a)
    (r4e) edge (r4b)
    (r4f) edge (r4b)
    (r4g) edge (r4c)
    (r4g) edge (r4e)
    (r4h) edge (r4d)
    (r4h) edge (r4f)
;

\end{tikzpicture}
\end{equation}
\end{proof}

Thus we can assume all the basic graphs of the "left" form appear to the left of all the basic graphs of the "right" form. Moreover, we can assume the basics graphs $Z_l,\dots, Z_{j+1}$, are all direct sum of $\delta$ ,$(1)$'s, $Z_j,\dots, Z_{i+1}$ are all direct sum of $(r)$'s, $r\in R$, and $Z_i,\dots, Z_1$ are all direct sum of $\delta^t$, $(1)$'s. We can represent the graph $L=Z_l\circ \dots \circ Z_{j+1}$ (resp $R=Z_1\circ \dots \circ Z_l$) using \underline{Ass}, \underline{Comm} (resp, \underline{$Ass^t$}, \underline{$Comm^t$}) as a direct sum of the graphs $\delta^{(n)}$ (resp. $(\delta^t)^{(n)}$) , $\delta^{(1)}\equiv (1), \delta^{(2)}\equiv \delta$, and $\delta^{(n)}\equiv \delta\circ (\delta^{(n-1)}\oplus (1))= \{ \begin{tikzpicture}[x=1.3cm, y=1cm,
    every edge/.style={
        draw,
        postaction={decorate,
                    decoration={markings,mark=at position 0.5 with {\arrow{>}}}
                   }
        }
]
     \vertex(l) at (0,-0.4){};
     \vertex (r1) at (1,0.7) [label=right:$1$]{};
     \vertex (r2) at (1,0.4)  [label=right:$2$]{};
     \vertex (r3) at (1,-0.4) [label=right:$n$]{};
     \node (r4) at (1,0.1) {$\vdots$};

\path

    (r1) edge (l)
    (r2) edge (l)
    (r3) edge (l)
;
\end{tikzpicture} \}$
, $(\delta^t)^{(n)}\equiv \{ \begin{tikzpicture}[x=1.3cm, y=1cm,
    every edge/.style={
        draw,
        postaction={decorate,
                    decoration={markings,mark=at position 0.5 with {\arrow{>}}}
                   }
        }
]
     \vertex(r) at (1,-0.4){};
     \vertex (l1) at (0,0.7) [label=left:$1$]{};
     \vertex (l2) at (0,0.4)  [label=left:$2$]{};
     \vertex (l3) at (0,-0.4) [label=left:$n$]{};
     \node (l4) at (0,0.1) {$\vdots$};

\path

    (r) edge (l1)
    (r) edge (l2)
    (r) edge (l3)
;
\end{tikzpicture} \}.$ \\

Using the $\underline{(r_1\cdot r_2)}$ relation, the graph $D=Z_j\circ \dots \circ Z_{i+1}$ is equivalent to a "diagonal graph": a direct sum of $(r)$'s, with identifications $Out\; D\overset{\sim}{\leftarrow} D^1\iso In\; D$. Thus, all in all, we can represent $G$ by the data of the set $D$, together with the maps $\pi_1:D\rightarrow Y, \pi_0:D\rightarrow X,\mu:D\rightarrow R:$

\begin{equation}
\begin{tikzpicture}[x=1.3cm, y=1cm,
    every edge/.style={
        draw,
        postaction={decorate,
                    decoration={markings,mark=at position 0.5 with {\arrow{>}}}
                   }
        }
]
    \node (l1) at (0,0) {\underline{Y}};
    \node (m1) at (2,0) {\underline{D}};
    \node (r1) at (4,0) {\underline{X}};
    \node (l0) at (1.05,0.3) {$\pi_0$};
    \node (r0) at (3,0.3) {$\pi_1$};

    \draw [->] (1.7,0) -- (0.2,0);
    \draw [->] (2.2,0) -- (3.7,0);

     \vertex(m2a) at (1.75,-0.5){};
     \vertex (m2b) at (2.25,-0.5){};
     \vertex(m3a) at (1.75,-1){};
     \vertex (m3b) at (2.25,-1){};
      \vertex(m4a) at (1.75,-1.5){};
     \vertex (m4b) at (2.25,-1.5){};
     \vertex(m5a) at (1.75,-2){};
     \vertex (m5b) at (2.25,-2){};
     \vertex(m6a) at (1.75,-3.5){};
     \vertex (m6b) at (2.25,-3.5){};

    \draw  (1.82,-0.5) --  (2.18,-0.5);
    \draw  (1.82,-1) -- (2.18,-1);
    \draw  (1.82,-1.5) --  (2.18,-1.5);
    \draw (1.82,-2) -- (2.18,-2);
    \draw (1.82,-3.5) -- (2.18,-3.5);

     \vertex(l2) at (0,-0.75) [label=left:$y_1$]{};
     \vertex (l3) at (0,-1.75) [label=left:$y_2$]{};
     \vertex(r2) at (4,-0.75) [label=right:$x_1$]{};
     \vertex (r3) at (4,-1.75) [label=right:$x_2$]{};
      \vertex(l4) at (0,-3.5){};
     \vertex (r4) at (4,-3.5){};

    \node (mm1) at (2,-0.35) {$r_1$};
    \node (mm2) at (2,-0.85) {$r_2$};
    \node (mm3) at (2,-1.35){$r_3$};
    \node (mm4) at (2,-1.85){$r_4$};
    \node (mm5) at (2,-3.35) {$r_m$};

    \draw  (0.07,-0.75) --  (1.68,-0.5);
    \draw  (0.07,-0.75) --(1.68,-1);
    \draw  (0.07,-1.75) --  (1.68,-1.5);
    \draw (0.07,-1.75) -- (1.68,-2);
    \draw (0.07,-3.5) -- (1.68,-3.5);

    \draw   (3.92,-0.75) -- (2.32,-0.5);
    \draw   (3.92,-0.75) -- (2.32,-1.5);
    \draw  (3.92,-0.75) --  (2.32,-2);
    \draw (3.92,-1.75) -- (2.32,-1);
    \draw (3.92,-3.5) -- (2.32,-3.5);

    \node (mm5) at (2,-2.55) {$\vdots$};

\end{tikzpicture}
\end{equation}

If $d_1,d_2\in D$, are such that $\pi_1(d_1)=\pi_1(d_2), \pi_0(d_1)=\pi_0(d_2)$, we can use the $\underline{(r_1+r_2)}$ -relation to identify $d_1$ and $d_2$ to a point $d$ with $\pi_1(d)=\pi_1(d_i), \pi_0(d)=\pi_0(d_i), \mu(d)=\mu(d_1)+\mu(d_2)$:

\begin{equation}
\begin{tikzpicture}[scale=0.4,auto=left]

    \node (a) at (-0.5,14) {$\underline{D}$};
    \node (b) at (11.5,14) {$\underline{D}$};
    \node (c) at (-0.5,12.5) {$\vdots$};
    \node (d) at (11.5,12.5) {$\vdots$};
    \vertex(a1) at (-5,9){};
    \vertex (b1) at (-2,10){};
    \vertex(c1) at (-2,8){};
    \vertex (d1) at (1,10){};
    \vertex(e1) at (1,8){};
    \vertex (f1) at (4,9){};
    \node (h1) at (-0.5,10.5) {$r_1$};
    \node (i1) at (-0.5,8.5) {$r_2$};

     \foreach \from/\to in {a1/b1,a1/c1,b1/d1,c1/e1,d1/f1,e1/f1}
    \draw (\from) -- (\to);

    \node (m) at (5.5,9) {$\implies$};

    \vertex(r1) at (7,9){};
    \vertex(r2) at (10,9){};
    \vertex(r3) at (13,9){};
    \vertex(r4) at (16,9){};

     \foreach \from/\to in {r1/r2,r2/r3,r3/r4}
    \draw (\from) -- (\to);

    \node (m) at (5.5,9) {$\implies$};
    \node (i1) at (11.5,9.5) {$r_1+r_2$};
    \node (rr) at (-0.5,6.5) {$\vdots$};
    \node (ll) at (11.5,6.5) {$\vdots$};

\end{tikzpicture}
\end{equation}

Thus we can assume $D$ is a subset of $Y\times X$, $\mu:D\rightarrow R$. Extending $\mu$ by zeros we get $\mu: Y\times X\rightarrow R$, i.e. a $Y\times X$ matrix with values in $R$, which is just $\Phi(G)$. \\

The same proof shows that,

\thm{2.10.2} \textit{For the surjective homomorphism $\Phi:\FF[\delta_{1,2},\delta^t_{1,2}]\sur\FF(\NN)$, the equivalence ideal $\cK \; (\Phi)$ is generated by relations (i)-(v).}

\thm{2.10.3} \textit{For the surjective homomorphism
$\Phi:\FF\{\pm \;1\}[\delta_{1,2},\delta^t_{1,2}]\sur\FF(\ZZ)$,
the relations are (i)-(v), $(-1,\delta),(\delta^t,-1),((-1)\cdot
(-1)),$ and cancellation
$$\delta\circ ((1)\oplus (-1))\circ \delta^t\equiv 0.$$}

\begin{appendices}
\chapter{Proof of Ostrowski's theorem}
\label{appB}
\smallbreak
\begin{proof}[\textit{ Proof of Ostrowski I}]

We can describe the elements of $Val(\FF(\QQ)/\FF\{\pm 1\})$ as collection of mappings
$$|\;|_{Y,X}:\FF(\QQ)_{Y,X}\rightarrow [0,\infty)$$
satisfying (I),(II),(III) and $|\pm 1|_{1,1}=1$, where we identify
the collection $\{|\;|_{Y,X}\}$ with the collection
$\{|\;|_{Y,X}^{\lambda}\}$, for any $\lambda >0$. The "generic
point" $\FF(\QQ)$ corresponds to the trivial valuation
$|y|_{Y,X}=\begin{cases}
1 &  y \neq 0 \\
0 &  y = 0
\end{cases}$.\\
Let $\{|\;|_{Y,X}\}$ be a non - trivial valuation on $\FF(\QQ)$, and let $B\subseteq \FF(\QQ)$ be the associated valuation - $\FF$ - subring, $B_{Y,X}=\{b\in \QQ^{Y\times X}\;,\; |b|_{Y,X}\leq 1\}$.\\
For $q_1,q_2\in \QQ$ we have
$$|q_1+q_2|_{1,1}=\left|(1,1)\circ \left(\begin{matrix}
  q_1 & 0  \\
  0 & q_2
 \end{matrix}\right)\circ \left(\begin{matrix} 1\\ 1 \end{matrix}\right)\right|_{1,1}\leq |(1,1)|_{1,2}\cdot \left|\left( \begin{matrix}
 1 \\ 1
 \end{matrix}\right)\right|_{2,1}\cdot \max\{|q_1|_{1,1},|q_2|_{1,1}\}. $$
Note that $|(1,1)|_{1,2}=\left|\left(\begin{matrix}
1 \\ 1
\end{matrix}\right)\right|_{2,1}$ by III(iii).\\
If we have $|(1,1)|_{1,2}\leq 1$, than $|q_1+q_2|_{1,1}\leq
\max\{|q_1|_{1,1},|q_2|_{1,1}\}$, and it follows that:
$|n|_{1,1}\leq 1$ for all $n\in \ZZ$; $\{n\in \ZZ, |n|_{1,1}<1
\}=p\cdot \ZZ$ is a prime ideal of $\ZZ$; $\ZZ_{(p)}\subseteq
B_{1,1}\underset{\neq}{\subseteq}\QQ$ ; and as $B_{1,1}$ is an
(ordinary) subring of $\QQ$: $\ZZ_{(p)}=B_{1,1}$.\\
For a matrix $b\in B_{Y,X}\subseteq \QQ^{Y\times X}$, its coefficients
$b_{y,x}=j_y^t\circ b\circ j_x\in B_{1,1}=\ZZ_{(p)}$, so
$B\subseteq \FF(\ZZ_{(p)})$. We have by II(i)
$$|(1,1,\dots,1)|_{1,n}=\Sup\bigg\{\left|(1,1,\dots,1)\circ \left(\begin{matrix}
b_1 \\ \vdots \\ b_n
\end{matrix}\right)\right|_{1,1},\left(\begin{matrix}
b_1 \\ \vdots \\ b_n
\end{matrix}\right)\in B_{n,1}\bigg\}$$
$$\leq \sup\{|b_1+\dots+b_n|_{1,1}, b_j\in \ZZ_{(p)}\}\leq 1,$$
and for a vector $b=(b_1,\dots, b_n) \in \FF(\ZZ_{(p)})_{1,n}$, we have
$$|b|_{1,n}=\left|(1,1,\dots,1)\circ \left(\begin{matrix}
b_1&&&0 \\
&b_2&& \\
&&\ddots& \\
0&&&b_n
\end{matrix}\right)\right|_{1,n}\leq |(1,1,...,1)|_{1,n}\cdot \max\{|b_j|_{1,1}\}\leq 1, $$
so $\FF(\ZZ_{(p)})_{1,n}=B_{1,n}$, and $\FF(\ZZ_{(p)})_{n,1}=B_{n,1}$  .\\
For a matrix $b=(b_{y,x})\in \ZZ_{(p)}^{Y\times X}$, we have by
II(i)
$$|b|_{Y,X}=\sup\bigg\{|d\circ b\circ d'|_{1,1},d\in \FF(\ZZ_{(p)})_{1,X}, d'\in \FF(\ZZ_{(p)})_{Y,1}\bigg\}\leq 1 $$
and so $\FF(\ZZ_{(p)})_{Y,X}=B_{Y,X}$. \\
Assume that we have $|(1,1)|_{1,2}>1$. Passing to an equivalent
norm $\{|\;|_{Y,X}^{\lambda}\}$, (with $\lambda\leq \frac{\log
\sqrt{2}}{\log |(1,1)|_{1,1}}$), we can assume that
$|(1,1)|_{1,2}=\left|\left(\begin{matrix} 1 \\ 1
\end{matrix}\right)\right|_{2,1}\leq \sqrt{2}$, and
$$ |q_1+q_2|_{1,1}\leq 2\cdot \max\{|q_1|_{1,1},|q_2|_{1,1}\}\hsm ,q_i\in \QQ. $$
By induction we get
$$|\sum\limits_{i=1}^{2^r} q_i|_{1,1}\leq 2^r\cdot \max\{|q_i|_{1,1}\}\hsm ,q_i\in \QQ$$
hence,
$$\left|\sum\limits_{i=1}^{n} q_i\right|_{1,1}\leq 2\cdot n\cdot \max\{|q_i|_{1,1}\} $$
hence $|n|_{1,1}\leq 2\cdot |n|_{\eta}$, $|n|_{\eta}=\pm n$. the usual absolute value for $n\in \ZZ$. \\ We have for $q_1,q_2\in \QQ$,
$$|q_1+q_2|_{1,1}=|(q_1+q_2)^n|_{1,1}^{1/n}=|\sum\limits_{k=0}^{n} \binom{n}{k} q_1^k\cdot q_2^{n-k}|_{1,1}^{1/n}\leq $$
$$(2(n+1))^{1/n}\cdot (\underset{k}{\max}\{|\binom{n}{k}|_{1,1}\cdot |q_1|_{1,1}^k\cdot |q_2|_{1,1}^{n-k}\})^{1/n}\leq $$
$$ (4(n+1))^{1/n}\cdot (\underset{k}{\max}\{\binom{n}{k}\cdot |q_1|_{1,1}^k\cdot |q_2|_{1,1}^{n-k}\})^{1/n}\leq $$
$$(4(n+1))^{1/n}\cdot ((|q_1|_{1,1}+|q_2|_{1,1})^n)^{1/n}=(4(n+1))^{1/n}\cdot (|q_1|_{1,1}+|q_2|_{1,1}).$$
and taking the limit $n\rightarrow \infty$ we get the triangle inequality
$$|q_1+q_2|_{1,1}\leq |q_1|_{1,1}+|q_2|_{1,1}.$$
We have by II(i),
$$|(1,1,...,1)|_{1,n}=\Sup\{|(1,1,...,1)\circ q|_{1,1},\; |q|_{n,1}\leq 1\}$$
$$=\Sup\{|q_1+\dots+q_n|_{1,1}\;, \; |q|_{n,1}\leq 1\} $$
$$\leq \Sup\{|q_1|_{1,1}+...+|q_n|_{1,1},\; |q|_{n,1}\leq 1\} \leq n $$
Let $a,b\in \NN$, with $a>1$, so we can expand $b$ in the base $a$:
$$b=d_m\cdot a^m+\dots+d_j\cdot a^j+\dots+d_1\cdot a+ d_0 $$
$$0\leq d_j <a \;,\; m<\frac{\log b}{\log a} $$
We get
$$|b|_{1,1}=|(1,\dots,1)\circ \overset{m}{\underset{i=0}{\bigoplus}} a^i\circ \overset{m}{\underset{j=0}{\bigoplus}} d_j\circ (1,\dots, 1)^t|_{1,1}\leq $$
$$|(1,\dots, 1)|_{1,m+1}^2\cdot \underset{0<d<a}{\max} \{|d|_{1,1}\}\cdot \underset{j<\frac{\log b}{\log a}}{\max}\{|a|_{1,1}^j \}\leq $$
$$(1+m)^2\cdot M_a \cdot \max\{1,|a|_{1,1}^{\frac{\log b}{\log a}}\}$$
with $M_a=\underset{0\leq d < a}{\max} \{|d|_{1,1}\}$ a constant independent of $b$. From this follows
$$|b|_{1,1}=|b^n|_{1,1}^{1/n}\leq (1+n\cdot \frac{\log b}{\log a})^{2/n}\cdot M_a^{1/n}\cdot \max\{1,|a|_{1,1}^{\frac{\log b}{\log a}}\}$$
and letting $n\rightarrow \infty$ we obtain
$$|b|_{1,1}\leq \max\{1,|a|_{1,1}^{\frac{\log b}{\log a}}\}$$
It follows that if $|b|_{1,1}> 1 $ for some $b\in \ZZ$, than
$|a|_{1,1}>1$ for all $a\in \ZZ \backslash\{\pm 1,0\}$, in which
case $|b|_{1,1}^{\frac{1}{\log b}}=|a|_{1,1}^{\frac{1}{\log
a}}=e^{\delta}$ is a constant for $a,b\in \ZZ\backslash\{\pm
1,0\}$, or
$$|a|_{1,1}=|a|_{\eta}^{\delta}\;\;\; \text{for }\; a\in \ZZ, \text{hence for}\;\;\; a\in \QQ,\; (|a|_{\eta}=\pm a\; \text{the real absolute value.}) $$
Note that
$$2^{\delta}=|2|_{1,1}\leq |(1,1)|_{1,2}\cdot \left|\left(\begin{matrix}
1 \\ 1
\end{matrix}\right)\right|_{2,1}\leq \sqrt{2}\sqrt{2}=2,\;\;, \delta\leq 1,$$
and passing to an equivalent norms we may assume $\delta=1$, and $|q|_{1,1}=|q|_{\eta}$ is the usual real absolute value. \vspace{1.5mm} \\
For a vector $q=(q_1,\dots, q_n)\in \QQ^n$ we get
$$\sum\limits_{i=1}^{n} q_i^2=|q\circ q^t|_{\eta}=|q\circ q^t|_{1,1}\leq |q|_{1,n}\cdot |q^t|_{n,1}=|q|_{1,n}^2 $$

or $|q|_{\eta}=\left(\sum\limits_{i=1}^{n}q_i^2\right)^{1/2}\leq |q|_{1,n}$.\\
On the other hand, from II(i) we get
$$|q|_{1,n}=\Sup\bigg\{ |q\circ b|_{1,1}\;,\; |b|_{n,1}\leq 1\bigg\}$$
$$\leq \Sup\bigg\{\left|\sum\limits_{i=1}^{n}q_i\circ b_i\right|\;,\; \sum\limits_{i=1}^{n}b_i^2\leq 1\bigg\}=\left(\sum\limits_{i=1}^{n}q_i^2\right)^{1/2}=|q|_{\eta}.$$
and so
$$|(q_1,\dots,q_n)|_{1,n}=\left|\left(\begin{matrix}
q_1 \\ \vdots \\ q_n
\end{matrix}\right)\right|_{n,1}=(\sum\limits_{i=1}^{n}q_i^2)^{1/2}$$
and $B_{1,n}=(\mathcal{O}_{\QQ,\eta})_{1,n}\;,\; B_{n,1}=(\mathcal{O}_{\QQ,\eta})_{n,1}$. \\
Finally, for a matrix $a\in\QQ^{Y\times X}$, from II(i) we get
$$|a|_{Y,X}=\Sup\{|b\circ a\circ b'|\;,\; b=(b_y)\;,\;\sum\limits_{y}^{}|b_y|^2\leq 1\;,\; b'=(b_x')\;,\; \sum\limits_{x}^{}|b_x'|^2\leq 1\}$$
is the usual $l_2$- operator norm, and $B_{Y,X}=(\mathcal{O}_{\QQ,\eta})_{Y,X}$\vspace{3mm}\\
\end{proof}
\begin{proof}[\textit {Proof of Ostrowski II}]
We can describe the elements of $Val(\slfrac{\FF(K)}{\FF\{\mu_K\}})$ as collection of mappings
$$|\;|_{Y,X}=\FF(K)_{Y,X}=K^{Y\times X}\rightarrow [0,\infty).$$
satisfying (I),(II),(III), and $|\mu_K|_{1,1}=1$, identifying $\{|\;|_{Y,X}\}$ with $\{|\;|_{Y,X}^{\lambda}\}$, $\lambda>0$. \\
Let $B_{Y,X}=\{b\in \FF(K)_{Y,X}\;,\; |b|_{Y,X}\leq 1\}$ be the valuation -$\FF$- subring of $\FF(K)$, corresponding to a non- trivial valuation $\{|\;|_{Y,X}\}$. For $q_1,q_2\in K $,
$$|q_1+q_2|_{1,1}=\left|(1,1)\circ \left(\begin{matrix}
q_1&0 \\ 0&q_2
\end{matrix}\right)\circ \left(\begin{matrix}
1 \\ 1
\end{matrix}\right)\right|_{1,1}\leq |(1,1)|_{1,2}^2\cdot \max\{|q_1|_{1,1},|q_2|_{1,1}\}.$$
Thus $|\;|_{1,1}$ is a valuation of $K$ (cf. \cite{CF}),
and passing to equivalent valuation we may assume $|q|_{1,1}=|q|_{\mathfrak{p}},\mathfrak{p}\subseteq \mathcal{O}_K$ a finite prime, or $|q|_{1,1}=|\eta q|$ with $\eta:K\hookrightarrow \CC$ (modulo conjugation) a "real prime". \vspace{3mm}

In the non-archimedean case, $|q|_{1,1}=|q|_{\mathfrak{p}}\;,\;
B_{1,1}=\mathcal{O}_{K,\mathfrak{p}}$, and for any matrix
$b=(b_{y,x})\in B_{Y,X}$,
$$|b_{y,x}|_{\mathfrak{p}}=|b_{y,x}|_{1,1}=|j_y^t\circ b\circ j_x|_{1,1}\leq |j_y^t|_{1,Y}\cdot |b|_{Y,X}\cdot|j_x|_{X,1}\leq |b|_{Y,X}\leq 1, $$
so $B_{Y,X}\subseteq \FF(\mathcal{O}_{K,\mathfrak{p}})_{Y,X}$. Note
that by II(i),
$$|(1,1,\dots,1)|_{1,n}=\Sup\bigg\{|b_1+\dots+b_n|_{\mathfrak{p}}\;,\; \left|\left(\begin{matrix}
b_1 \\ \vdots \\ b_n
\end{matrix}\right)\right|_{n,1}\leq 1\bigg\}$$
$$\leq \Sup\{|b_1+\dots+b_n|_{\mathfrak{p}}\;,\; |b_i|_{\mathfrak{p}}\leq 1\}\leq 1.$$
For a vector $b=(b_1,\dots,b_n)\in \FF(\mathcal{O}_{K,\mathfrak{p}})_{1,n},\\
$ $|b|_{1,n}=\left|(1,1,\dots,1)\circ\left(\begin{matrix}
b_1&&&0\\
&b_2&&&\\
&&\ddots&\\
0&&&b_n
\end{matrix}\right)\right|_{1,n}\leq |(1,1,\dots,1)|_{1,n}\cdot\max\{|b_i|_{\mathfrak{p}}\}\leq 1$,\\
so $B_{1,n}=\FF(\mathcal{O}_{K,\mathfrak{p}})_{1,n}$, and $B_{n,1}=\FF(\mathcal{O}_{K,\mathfrak{p}})_{n,1}$.\\
Finally, for a matrix
$b=(b_{y,x})\in\FF(\mathcal{O}_{K,\mathfrak{p}})_{Y,X}$ we have by
II(i),
$$|b|_{Y,X}=\Sup\{|d\circ b\circ d'|_{\mathfrak{p}}\;,\; d\in \FF(\mathcal{O}_{K,\mathfrak{p}})_{1,Y}\;,\; d'\in \FF(\mathcal{O}_{K,\mathfrak{p}})_{X,1}\}\leq 1,$$ so $B_{Y,X}=\FF(\mathcal{O}_{K,\mathfrak{p}})_{Y,X}$.\\
In the archimedean case, $|q|_{1,1}=|\eta q|\;,\; \eta:K\hookrightarrow \CC$.\\
For a vector $q=(q_1,\dots,q_n)\in K^n$, we get
$$\sum\limits_{i=1}^{n}|\eta q_i|^2=|q\circ \bar{q}^t|_{1,1}\leq |q|_{1,n}^2\;,\; \text{so }\; (\sum\limits_{i=1}^{n}|\eta q_i|^2)^{1/2}\leq |q|_{1,n}, $$
and $B_{1,n}\subseteq (\mathcal{O}_{K,\eta})_{1,n}\;,\;
B_{n,1}\subseteq (\mathcal{O}_{K,\eta})_{n,1}$.\\
Conversely, from II(i)
we get,
$$|q|_{1,n}=\Sup\{|q\circ b|_{1,1}\;,\; |b|_{n,1}\leq 1\}$$
$$\leq \Sup \{|\eta(\sum\limits_{i=1}^{n}q_i\cdot b_i)|\;,\;\sum\limits_{i=1}^{n}|\eta b_i|^2\leq 1 \}=(\sum\limits_{i=1}^{n}|\eta q_i|^2)^{1/2}.$$
and $B_{1,n}=(\mathcal{O}_{k,\eta})_{1,n},$ and similarly
$B_{n,1}=(\mathcal{O}_{K,\eta})_{n,1}$.\\
Finally, for a matrix $a\in
K^{Y\times X}$, we get from II(i),
$$|a|_{Y,X}=\Sup\{|b\circ a\circ b'|_{1,1}\; ;\; b=(b_y),\; \sum\limits_{y}^{}|\eta b_y|^2\leq 1\;, b'=(b_x'), \sum\limits_{x}^{}|\eta b_x'|^2\leq 1\}$$
is the usual $l_2$- operator norm, and $B_{Y,X}=(\mathcal{O}_{K,\eta})_{Y,X}$.
\end{proof}
\end{appendices}

\chapter{Geometry}
\bigskip

In this section $A\in \CFR$ is commutative.

\section{Ideals, maximal ideals and primes}
\smallbreak

According to definition \ref{defA.2.3}.3, a subset $\mathfrak{a}\subseteq A_{1,1}$ is called an \underline{ideal} if for
\begin{equation}
a_1,\dots,a_n\in \mathfrak{a},\;\; b\in A_{1,n},\;\; b'\in A_{n,1}:\;\; b\circ(a_1\oplus \dots \oplus a_n)\circ b'\in \mathfrak{a}.
\end{equation}
Denote the set of ideals of $A$ by $\mathcal{I}(A)$.\\
Given an indexed
set of ideals $\mathfrak{a}_i\subseteq A_{1,1},\; i\in I$, their
intersection $\cap_I \mathfrak{a}_i$ is again an ideal. Their sum
$\Sigma \mathfrak{a}_i$ is an ideal generated by $\cup_I
\mathfrak{a}_i$,
\begin{equation}
\sum_I \mathfrak{a}_i=\bigg\{b\circ (\underset{j}{\oplus} a_j)\circ b'\bigg|\;\; a_j\in \cup \mathfrak{a}_j \bigg\}
\end{equation}
The product $\mathfrak{a}\cdot \mathfrak{a}'$ of two ideals is an ideal generated by the product of elements of these ideals,
\begin{equation}
\mathfrak{a}\cdot \mathfrak{a}'=\bigg\{b\circ(\underset{j}{\oplus}a_j\cdot a_j')\circ b' \bigg| \; a_j\in \mathfrak{a},a_j'\in \mathfrak{a}' \bigg\}
\end{equation}
Let $\varphi:A\rightarrow B$ be a homomorphism of $\frs$. \\
If $\mathfrak{b}\in \mathcal{I}(B)$ then
$\varphi^*(\mathfrak{b})=\varphi^{-1}(\mathfrak{b})\in \mathcal{I}(A)$
,  and we have a map
\begin{equation}
\varphi^*:\mathcal{I}(B)\rightarrow \mathcal{I}(A),\;\;\; \mathfrak{b}\mapsto \varphi^{-1}(\mathfrak{b}).
\end{equation}
If $\mathfrak{a}\in \mathcal{I}(A)$, $\varphi(\mathfrak{a})$ generates
the ideal $\varphi_*(\mathfrak{a})$,
\begin{equation}
\varphi_*:\mathcal{I}(A)\rightarrow \mathcal{I}(B),\;\;\; \mathfrak{a}\mapsto \varphi_*(\mathfrak{a})=\{b\circ(\oplus \varphi(a_i))\circ b'\}.
\end{equation}

\prop{3.1.1 }
\textit{ We have the following:  \vspace{2mm}\\
(1)  $\mathfrak{a}\subseteq \varphi^*\varphi_* \mathfrak{a},\;\;\;\mathfrak{a}\in \mathcal{I}(A)$.\\
(2)  $\mathfrak{b} \supseteq \varphi_*\varphi^* \mathfrak{b},\;\;\; \mathfrak{b}\in \mathcal{I}(B)$. \\
(3)  $\varphi^*\mathfrak{b}=\varphi^*\varphi_*\varphi^*\mathfrak{b},\;\;\;\varphi_*\mathfrak{a}=\varphi_*\varphi^*\varphi_*\mathfrak{a}$. \\
(4)  there is a bijection, via $\mathfrak{a}\mapsto \varphi_*\mathfrak{a}$ (with inverse map $\mathfrak{b}\mapsto \varphi^*\mathfrak{b}$), from the set
\begin{equation}
\{ \mathfrak{a}\in \mathcal{I}(A) |\;\; \varphi^*\varphi_*\mathfrak{a}=\mathfrak{a}\}=\{\varphi^*\mathfrak{b}|\;\; \mathfrak{b}\in \mathcal{I}(B) \}
\end{equation}
to the set
\begin{equation}
\{ \mathfrak{b}\in \mathcal{I}(B) |\;\; \varphi_*\varphi^*\mathfrak{b}=\mathfrak{b} \}=\{\varphi_*\mathfrak{a}|\;\; \mathfrak{a}\in \mathcal{I}(A) \}.
\end{equation}}
\begin{proof}
The proofs of these are straightforward.
\end{proof}

Given an ideal $\mathfrak{a}\in I(A)$, we write $A/\mathfrak{a}$ for the quotient $\fr$ $A/E(\mathfrak{a})$, where $E(\mathfrak{a})$ is the equivalence ideal generated by $\mathfrak{a}$.

\prop{3.1.2}
\textit {We have a one-to-one order-preserving correspondence
\begin{equation}
\pi^*:\mathcal{I}(A/\mathfrak{a})\iso \{\mathfrak{b}\in \mathcal{I}(A)|\;\; \mathfrak{b}\;\text{satisfies  } (*) \}
\end{equation}
where $(*)$ means
\begin{equation}
\text{for any  } a\in \mathfrak{a}:\;\;b\circ(id_Z\oplus a)\circ b'\in \mathfrak{b}\iff b\circ (id_Z\oplus 0)\circ b'\in \mathfrak{b}.
\end{equation}}

\begin{proof}
The proof is clear. (cf. ($\ref{A.2.9}$) for the equivalence ideal $E(\mathfrak{a})$ generated by $(a,0),\; a\in \mathfrak{a}$).\vspace{3mm}\\
\end{proof}

Since the union of a chain of proper ideals is again a proper ideal, an application of Zorn's lemma gives the following result.

\thm{3.1.1 (Zorn)}
\textit{There exists a maximal ideal $\mathfrak{m}\varsubsetneq A_{1,1}$.}

\defin{3.1.1}
\textit{An ideal $\mathfrak{p}\subseteq A_{1,1}$ is
called \underline{prime}:
\begin{equation}
S_{\mathfrak{p}}=A_{1,1}\setminus \mathfrak{p}\; \text{is multiplicatively closed } S_{\mathfrak{p}}\cdot S_{\mathfrak{p}}=S_{\mathfrak{p}} .
\end{equation}
}
We denote by $\Spec A$ the set of prime ideals. \\
For a homomorphism of $\frs$ $\varphi:A\rightarrow B$, the pullback
$\varphi^*=\varphi^{-1}$ induces a map
\begin{equation}
\varphi^*=\Spec(\varphi):\Spec B\rightarrow \Spec A.
\end{equation}

\prop{3.1.3}
\label{3.1.3} \textit{(1)  If $\mathfrak{m}$ is maximal then it is prime.\\
(2) More generally, if $\mathfrak{a}\in \mathcal{I}(A)$, and given $f\in A_{1,1}$ such that,
\begin{equation}
\forall n\in \NN: \;\;\; f^n\notin \mathfrak{a}.
\end{equation}
let $\mathfrak{m}$ be a maximal element of the set
\begin{equation}\label{idealset}
\{\mathfrak{b}\in \mathcal{I}(A) | \mathfrak{b} \supseteq \mathfrak{a}, \mathfrak{b}\not \ni f^n\:
\forall n\in \NN\}
\end{equation}
Then $\mathfrak{m}$ is prime.}
\begin{proof}
(1) If $x,y\in A_{1,1}\setminus \mathfrak{m}$, the ideals
$(x)+\mathfrak{m},(y)+\mathfrak{m}$ are the unit ideals. So we can
write
\begin{equation}
 1=b\circ(\underset{J}{\bigoplus} m_j) \circ d,\hsm \text{with}\hsm m_j=x\hsm \text{or} \hsm m_j\in \mathfrak{m}
\end{equation}
\begin{equation}
 1=b'\circ(\underset{I}{\bigoplus} m'_i)\circ d',\hsm \text{with}\hsm m'_i=y\hsm \text{or} \hsm m'_i\in \mathfrak{m}
\end{equation}
It then follows that,
\begin{equation*}
1=1\cdot 1=(b\circ \underset{J}{\bigoplus} m_j\circ d)\circ (b'\circ \underset{I}{\bigoplus} m_i' \circ d')
\end{equation*}
\begin{equation*}
=(b\circ \underset{J}{\oplus} m_j) \circ \underset{J}{\bigoplus}(b'\circ \underset{I}{\oplus} m'_i\circ d')\circ d
\end{equation*}
\begin{equation*}
=b\circ \underset{J}{\bigoplus} (b'\circ \underset{I}{\oplus} (m_j\circ m'_i)\circ d')\circ d
\end{equation*}
\begin{equation}
=(b\circ \underset{J}{\oplus} b')\circ \underset{J\otimes I}{\oplus}(m_j\circ m_i')\circ (\underset{J}{\oplus}d'\circ d).
\end{equation}
but $m_j\circ m_i'=x\circ y$ or $m_j\circ m_i'\in \mathfrak{m}$, so $1$ is in the ideal generated by $\mathfrak{m}$
and $x\circ y$, and since $1\notin \mathfrak{m}$ then  $x\circ y\notin \mathfrak{m} $. \vspace{2mm} \\
(2) Similarly, if $x\notin \mathfrak{m}$ then $f^n$ is in the ideal generated by $x$ and $\mathfrak{m}$ so
$f^n=b\circ \underset{J}{\bigoplus} m_j \circ d$ , with $m_j=x$ or $m_j\in \mathfrak{m}$. If $y\notin \mathfrak{m}$, $f^{n'}=b'\circ \underset{I}{\bigoplus} m_i'\circ d'$,
with $m_i'=y$ or $m_i'\in \mathfrak{m}$. It then follows that
$$f^{n+n'}=b\circ\underset{J}{\bigoplus} m_j\circ d\circ b'\circ \underset{I}{\bigoplus} m_i'\circ d'$$
\begin{equation}
=b\circ (\underset{J}{\bigoplus} b')\circ \bigg(\underset{J\otimes I}{\bigoplus} m_j\circ m_i'\bigg)\circ (\underset{J}{\bigoplus} d')\circ d
\end{equation}
but $m_j\circ m_i'=x\circ y$ or $m_j\circ m_i'\in \mathfrak{m}$, so $f^{n+n'}$ is in the ideal generated by $\mathfrak{m}$ and $x\circ y$, and since $f^{n+n'}\notin \mathfrak{m}$
then $x\circ y\notin \mathfrak{m}$ .
\end{proof}

\defin{3.1.2}
 \textit{For $\mathfrak{a}\in \mathcal{I}(A)$, the \underline{radical} is
\begin{equation}
\sqrt{\mathfrak{a}}=\{f\in A_{1,1} |\; f^n\in \mathfrak{a}\;\;\text{for some}\;\; n\geq 1 \}
\end{equation}
}
It is easy to see that $\sqrt{\mathfrak{a}}$ is an ideal. This also follows from the following proposition.
\prop{3.1.4}\textit{We have
\begin{equation}
\sqrt{\mathfrak{a}}= \underset{\mathfrak{a}\subseteq \mathfrak{p}}{\bigcap}\mathfrak{p}
\end{equation}
the intersection of prime ideals containing $\mathfrak{a}$.}
\begin{proof}
If $f\in \sqrt{\mathfrak{a}}$, say $f^n\in \mathfrak{a}$, then for all primes $\mathfrak{a}\subseteq \mathfrak{p}$, $f^n\in \mathfrak{p}$ and so $f\in \mathfrak{p}$. If $f\notin \sqrt{\mathfrak{a}}$,
let $\mathfrak{m}$ be a maximal element of the set ($\ref{idealset}$), it exists by Zorn's lemma, and it is prime by proposition \ref{3.1.3}.3(2), $\mathfrak{a}\subseteq \mathfrak{m}$ and $f\notin \mathfrak{m}$.
\end{proof}

\section{The spectrum: $\Spec \; A$}
\smallbreak
\defin{3.2.1} \textit{For a set $\mathfrak{U}\subseteq A_{1,1}$, we let
\begin{equation}
V_A(\mathfrak{U})=\{\mathfrak{p}\in \Spec A \; | \; \mathfrak{U}\subseteq \mathfrak{p}\}.
\end{equation}}
If $\mathfrak{a}$ is the ideal generated by $\mathfrak{U}$, $V_A(\mathfrak{U})=V_A(\mathfrak{a})$; we have
\begin{equation}
V_A(1)=\emptyset,\hsm V_A(0)=\Spec A,
\end{equation}
\begin{equation}
V_A(\Sigma \mathfrak{a})=\underset{i}{\cap} V_A(\mathfrak{a}_i)\;,\; \mathfrak{a}_i\in \mathcal{I}(A),
\end{equation}
\begin{equation}
V_A(\mathfrak{a}\cdot \mathfrak{a}') = V_A(\mathfrak{a})\cup V_A(\mathfrak{a}').
\end{equation}
Hence the sets $\{V_A(\mathfrak{a})\; |\; \mathfrak{a}\in
\mathcal{I}(A)\}$ are the closed sets for the topology on $\Spec
A$, the \textit{Zariski topology}.

\defin{3.2.2} \textit{For $f\in A_{1,1}$ we let
\begin{equation}
D_A(f)=\Spec (A)\setminus V_A(f)=\{\mathfrak{p}\in\Spec A\;|\; f\notin \mathfrak{p}\}.
\end{equation}}
We have
\begin{equation}
D_A(f_1)\cap D_A(f_2)=D_A(f_1\cdot f_2),
\end{equation}
\begin{equation}
\Spec A\setminus V_A(\mathfrak{a})=\underset{f\in \mathfrak{a}}{\bigcup} D_A(f).
\end{equation}
Hence the sets $\{D_A(f)\;|\; f\in A_{1,1}\}$ are the basis for
the open sets in the Zariski topology. We have
\begin{equation}
D_A(f)=\emptyset\;\;\; \iff f\in \underset{\mathfrak{p}\in \Spec A}{\bigcap}\mathfrak{p}=\sqrt{0}\;\;\; \iff f^n=0 \;\text{for some } n
\end{equation}
and we say $f$ is a \textit{nilpotent}. We have
\begin{equation}
D_A(f)=\Spec A\;\;\;\iff (f)=(1)\iff \exists f^{-1}\in A_{1,1}: f\cdot f^{-1}=1
\end{equation}
and we say $f$ is \textit{invertible}. We denote by $GL_1(A)$ the
(commutative) group of invertible elements.

\defin{3.2.3} \textit{For a subset $X\subseteq \Spec A$, we have the associated ideal
\begin{equation}
\mathcal{I}(X)=\underset{\mathfrak{p}\in X}{\bigcap} \mathfrak{p}.
\end{equation}}

\prop{3.2.1} \textit{We have
\begin{equation}\label{3.2.11}
\mathcal{I} V_A{\mathfrak{a}}=\sqrt{\mathfrak{a}},
\end{equation}
\begin{equation}\label{3.2.12}
V_A\mathcal{I}(X)=\bar{X}, \text{the closure of } X \text{in the Zariski topology.}
\end{equation}
}
\begin{proof}
Equation ($\ref{3.2.11}$) is just a restatement of proposition $3.1.4$. For ($\ref{3.2.12}$), $V_A\mathcal{I}(X)$ is clearly a closed set
containing $X$, and if $C=V_A(\mathfrak{a})$ is a closed set containing $X$, then $\sqrt{\mathfrak{a}}=IV_A(\mathfrak{a})\subseteq \mathcal{I}(X),$
hence $C=V_A(\sqrt(\mathfrak{a}))\supseteq V_A\mathcal{I}(X).$
\end{proof}
\cor{3.2.1} \textit{We have a one-to-one order-reversing correspondence between closed sets $X\subseteq \Spec A$, and radical ideals $\mathfrak{a}$,
via $X\mapsto \mathcal{I}(X)$, $V_A(\mathfrak{a})\mapsfrom \mathfrak{a}$
\begin{equation}
\{X\subseteq\Spec A\;|\; \bar{X}=X\}\overset{1:1}{\longleftrightarrow}\{\mathfrak{a}\in \mathcal{I}(A)\;|\; \sqrt{\mathfrak{a}}=\mathfrak{a} \}.
\end{equation}
}
Under this correspondence the closed irreducible subsets correspond to the prime ideals. For $\mathfrak{p}_0,\mathfrak{p}_1\in \Spec A,\mathfrak{p}_0\in \bar{\{\mathfrak{p}_1\}}
\iff \mathfrak{p}_0\supseteq \mathfrak{p}_1$, we say that $\mathfrak{p}_0$ is a \textit{Zariski specialization} of $\mathfrak{p}_1$, or that $\mathfrak{p}_1$
is a \textit{Zariski generalization} of $\mathfrak{p}_0$. The space $\Spec A$ is sober: every closed irreducible subset $C$ has the form $C=V_A(\mathfrak{p})=\bar{\{\mathfrak{p}\}}$,
and we call the (unique) prime $\mathfrak{p}$ the \textit{generic point} of $C$.

\prop{3.2.2} \textit{The sets $D_A(f)$, and in particular $D_A(1)=\Spec A$, are compact (or 'quasi-compact': we do not include Hausdorff in compactness).}
\begin{proof}
Note that $D_A(f)$ is contained in the union $\underset{i}{\bigcup} D_A(g_i)$ if and only if $V_A(f)\supseteq \underset{i}{\bigcap} V_A(g_i)=V_A(\mathfrak{a})$, where
$\mathfrak{a}$ is the ideal generated by $\{g_i\}$, if and only if $\sqrt{f}=\mathcal{I}V_A(f)\subseteq \mathcal{I}V_A(\mathfrak{a})=\sqrt{\mathfrak{a}}$, if and only if
$f^n\in\mathfrak{a}$ for some $n$, if and only if $f^n=b\circ (\underset{i}{\bigoplus} g_i)\circ b'$, and in such expression only a finite number of the $g_i$ are involved.
\end{proof}
Let $\varphi:A\rightarrow B$ be a homomorphism of $\fr$, $\varphi^*:\Spec B\rightarrow \Spec A$ the associated pullback map.

\prop{3.2.3} \textit{We have
\begin{equation}\label{3.2.14}
\varphi^{*-1}(D_A(f))=D_B(\varphi(f)), \;\;\; f\in A_{1,1},
\end{equation}
\begin{equation}\label{3.2.15}
\varphi^{*-1}(V_A(\mathfrak{a}))=V_B(\varphi_*(\mathfrak{a})),\;\;\; \mathfrak{a}\in \mathcal{I}(A),
\end{equation}
\begin{equation}\label{3.2.16}
V_A(\varphi^{-1}\mathfrak{b})=\bar{\varphi^*(V_B(\mathfrak{b}))},\;\;\; \mathfrak{b}\in \mathcal{I}(B).
\end{equation}}
\begin{proof}
The proofs of ($\ref{3.2.14}$) and ($\ref{3.2.15}$) are straightforward:
$$\mathfrak{q}\in \varphi^{*-1}(D_A(f))\iff \varphi^*(\mathfrak{q})\in D_A(f) \iff f\notin \varphi^{-1}(\mathfrak{q})\iff \varphi (f)\notin \mathfrak{q}\iff \mathfrak{q}\in D_B(\varphi(f)), $$
$$\mathfrak{q}\in \varphi^{*-1}(V_A(\mathfrak{a}))\iff \varphi^*(\mathfrak{q})\in V_A(\mathfrak{a}) \iff \mathfrak{a}\subseteq \varphi^{-1}(\mathfrak{q}) \iff \varphi_*(\mathfrak{a}) \subseteq \mathfrak{q}
\iff \mathfrak{q}\in V_B(\varphi_*(\mathfrak{a})). $$
For ($\ref{3.2.16}$) we may assume $\mathfrak{b}=\sqrt{\mathfrak{b}}$ is a radical since $V_B(\mathfrak{b})=V_B(\sqrt{\mathfrak{b}}), \varphi^{-1}(\sqrt{\mathfrak{b}})=\sqrt{\varphi^{-1}(\mathfrak{b})}$.
Let $\mathfrak{a}=\mathcal{I}(\varphi^*(V_B(\mathfrak{b})))$, so that $V_A(\mathfrak{a})=\bar{\varphi^*(V_B(\mathfrak{b}))}$ by ($\ref{3.2.12}$). We have
$$f\in \mathfrak{a} \iff f\in \mathfrak{p}, \; \forall \mathfrak{p} \in \varphi^*(V_B(\mathfrak{b})) \iff f\in \varphi^{-1}(\mathfrak{q}), \;\;\; \forall \mathfrak{q}\supseteq \mathfrak{b}$$
\begin{equation}
\varphi (f)\in \underset{\mathfrak{q}\supseteq \mathfrak{b}}{\bigcap} \mathfrak{q}=\sqrt{\mathfrak{b}}=\mathfrak{b}\iff f\in \varphi^{-1}(\mathfrak{b}).
\end{equation}

\end{proof}
It follows from ($\ref{3.2.14}$), or from ($\ref{3.2.15}$), that
$\varphi^*=\Spec (\varphi)$ is continuous, hence $A\mapsto \Spec
A$ is a contravariant functor from commutative $\frs$ to compact,
sober, topological spaces.

\exmpl{3.2.1} Let $A$ be a commutative ring, $\FF(A)$ the associated $\fr$. An ideal $\mathfrak{a}\subseteq A=\FF(a)_{1,1}$ is also an ideal in our sense, and conversely.
Under this correspondence the primes of $A$ correspond to the primes of $\FF(A)$, and we have a homeomorphism with respect to the Zariski topologies:
\begin{equation}
\Spec A=\Spec \FF (A).
\end{equation}

\exmpl{3.2.2} Let $\eta:\kk\rightarrow \CC$ be a real or complex prime of a number field, and let $\mathcal{O}_{\kk,\eta}$ denote the $\fr$ of real or complex 'integers'.
Then
\begin{equation}
\mathfrak{m}_{\eta}=\{x\in \kk\;|\; |x|_{\eta} <1  \}
\end{equation}
is the (unique) maximal ideal of $\mathcal{O}_{\kk,\eta}$, the closed point of $\Spec \mathcal{O}_{\kk,\eta}$.

\section{Localization $S^{-1}A$}\label{3.3}
\smallbreak

 The theory of localization of an $\fr\;A$, with respect to a
  multiplicative subset $S\subseteq A_{1,1}$, goes exactly as
  in localization of commutative rings - since it is a
  multiplicative theory. We recall this theory next.\par We assume
  $S\subseteq A_{1,1}$ satisfies
\begin{equation}
1\in S
\end{equation}
\begin{equation}
s_1,\;s_2\in S\THEN s_1\cdot s_2\in S
\end{equation}
  On the set
  $$A\times S=\coprod_{Y,X}A_{Y,X}\times S$$

  we define for $a_i\in A_{Y,X},\ s_i\in S$
\begin{equation}
(a_1,s_1)\sim(a_2,s_2)\IFF s\cdot s_2\cdot a_1=s\cdot s_1\cdot a_2\;\text{for some}\;s\in
  S.
\end{equation}
  It follows that $\sim$ is an equivalence relation, and we denote
  by $a/s$ the equivalence class containing $(a,s)$,and by
  $S^{-1}A$ the collection of equivalence classes. On $S^{-1}A$ we
  define the operations:
\begin{equation}
a_1/s_1\circ a_2/s_2=(a_1\circ a_2)/s_1s_2\;,\;a_1\in A_{Z,Y}\;,\;a_2\in A_{Y,X}
\end{equation}
\begin{equation}\label{3.3.5}
a_1/s_1\oplus a_2/s_2=(s_2\cdot a_1\oplus s_1\cdot a_2)/s_1s_2
\end{equation}

  \prop{3.3.1}\textit{ The above operations are well defined, independent of the chosen representatives, and they satisfy the axioms of an $\fr$.}

  \begin{proof}The usual proof works. For example, replacing
  $a_1/s_1$ in (\ref{3.3.5}) by\\ $a_1'/s_1'\sim a_1/s_1,\;\text{say}\;s\cdot s_1'\cdot a_1=s\cdot s_1\cdot
  a_1'$, then $$s\cdot s_1's_2\cdot(s_2a_1\oplus s_1a_2 )=s\cdot s_1s_2\cdot(s_2a_1'\oplus
  s_1'a_2),$$ hence $$(s_2a_1\oplus s_1a_2)/s_1s_2=(s_2a_1'\oplus
  s_1'a_2)/s_1's_2.$$
  \end{proof}

  The $\fr\;S^{-1}A$ comes with a canonical homomorphism
\begin{equation}
\phi=\phi_S:A\to  S^{-1}A\;,\;\;\phi(a)=a/1.
\end{equation}

  \prop{3.3.2}\textit{We have the universal property of $\phi_S$}:
  $$\FR(S^{-1}A,B)=\{\fe\in\FR (A,B)|\;\fe(S)\subseteq GL_{[1]}(B)\}$$
  $$\widetilde{\fe}\longmapsto\widetilde{\fe}\circ\phi_S$$
  $$\widetilde{\fe}(a/s)=\fe(a)\cdot\fe(s)^{-1}\longmapsfrom\fe$$
  \begin{proof}
  Clear.
  \end{proof}
  Note that $S^{-1}A$ is the zero $\fr$ if and only if $0\in S$.\\

  The main examples of localizations are :
\begin{equation}
S_f=\{f^n\}_{n\geq 0},\;f\in A_{1,1}.\;\text{We write $A_f$ for $S^{-1}_fA$}.
\end{equation}
\begin{equation}
S_{\g{p}}=A_{1,1}\minus\g{p},\;\g{p}\in\Spec(A).\;\text{We write } A_{\g{p}}\text{ for }
S^{-1}_{\g{p}}A.
\end{equation}

Consider the canonical homomorphism $\phi=\phi_S:A\to \loc{A},\;\phi(a)=a/1$. \\
If $\g{b}\in \mathcal{I}(\loc{A})$, then $\fe^{-1}(\g{b})\in \mathcal{I}(A)$.\\
If $\g{a}\in \mathcal{I}(A)$ is an ideal of $A$ then
\begin{equation}
\loc{\g{a}}:=\phi_*(\g{a})=\{a/s\in(\loc{A})_{1,1}|\;a\in\g{a},\;s\in S\}
\end{equation}
is an ideal of $\loc{A}$.

\prop{3.3.3} \textit{If $\g{b}\in\mathcal{I}(\loc{A})$, then $\loc{(\fe^{-1}\{\g{b}\})}=\g{b}.$}

\begin{proof}If $a/s\in\g{b},\;a\in \fe^{-1}\{\g{b}\}$, and
$a/s\in\loc{(\fe^{-1}\{\g{b}\})}$; so $\g{b}\subseteq\loc{(\fe^{-1}\{\g{b}\})}$. The
reverse inclusion is clear.
\end{proof}

\prop{3.3.4}\textit{For $\g{a}\in \mathcal{I}(A)$ ,}
\begin{equation}
\fe^{-1}\{\loc{\g{a}}\}=\{a\in A|\;\exists\;s\in S\;:\;s\cdot a\in\g{a}\}.
\end{equation}
\textit{In particular,}
\begin{equation}
\loc{\g{a}}=(1)\IFF \g{a}\cap S\ne\varnothing
\end{equation}
\begin{proof}
$$a\in \fe^{-1}\{\loc{\g{a}}\}\IFF a/1=x/s,\;x\in\g{a},\;s\in S\IFF s\cdot a\in\g{a},\;\text{for some}\;s\in S.$$
\end{proof}

\prop{3.3.5}\textit{The map $\phi^*_S$ induces a bijection}
\begin{equation}
\phi^*_S:\Spec(\loc{A})\iso\{\g{p}\in\Spec A|\;\g{p}\cap
S=\varnothing\},
\end{equation}
\textit{which is a homeomorphism for the Zariski topology.}

\begin{proof}If $\g{q}\in\Spec(\loc{A}),\;\phi^*_S(\g{q})$
belongs to the right-hand-side. Conversely, if $\g{p}$ belongs to
the right-hand-side, $\loc{\g{p}}$ is a (proper) prime of
$\loc{A}$. By propositions $3.3.3, 3.3.4$, these operations are inverses of
each other.
\end{proof}

\cor{3.3.6}\textit{We have homeomorphism for $f\in A_{1,1}$,}
\begin{equation}
\phi^*_f:\Spec(A_f)\iso D_A(f).
\end{equation}

\cor{3.3.7}\textit{We have a homeomorphism for $\g{p}\in\Spec(A),$}
\begin{equation}
\phi^*_{\g{p}}:\Spec(A_{\g{p}})\iso\{\g{q}\in\Spec A|\;\g{q}\subseteq\g{p}\}.
\end{equation}
\textit{In particular, $A_{\g{p}}$ contains a unique maximal
ideal $\g{m}_{\g{p}}=\loc_{\g{p}}{\g{p}}$; we say it is a \underline{local}\;
$\fr$.}

\rem{3.3.8}For $\g{p}\in\Spec(A)$ we let
$\FF_{\g{p}}=A_{\g{p}}/\g{m}_{\g{p}}$ denote the \underline{residue
field} at $\g{p}$. Let $\pi:A\to  A/\g{p}$ be the canonical
homomorphism, and $\bar{S}_{\g{p}}=\pi(S_{\g{p}})$, we have also
$\FF_{\g{p}}=\bar{S}^{-1}_{\g{p}}(A/\g{p})$. The commutative
diagram
\begin{equation}
\xymatrix{
        A \ar^{\phi_{\g{p}}}[r]\ar_{\pi}[d] &  A_{\g{p}}\ar[d]\\
        A/\g{p}\ar[r] & \FF_{\g{p}}}
\end{equation}
is cartesian: $\FF_{\g{p}}=(A/\g{p})\otimes_A A_{\g{p}}$

It is also functorial: given a homomorphism of
$\frs$ $\fe:A\to  B,\;\g{q}\in\Spec B,\;\g{p}=\fe^*(\g{q})$, we
have a commutative cube
\begin{equation}
\xymatrix{
& A\ar[rr]\ar[dl]\ar[dd] & & A_{\g{p}}\ar[dd]\ar[dl]
\\
B \ar[rr]\ar[dd] & & B_{\g{q}}\ar[dd]
\\
& A/\g{p}\ar[dl]\ar[rr] & & \FF_{\g{p}}\ar[dl]
\\
B/\g{q}\ar[rr]& & \FF_{\g{q}}
}
\end{equation}
\section{Structure sheaf $\mathcal{O}_A$}
\smallbreak

Next we define a sheaf $\mathcal{O}_A$ of \;$\frs$ over $\Spec A$.
\defin{3.4.1}\emph{ For an open set $U\subseteq\Spec(A)$, and for $Y,X\in \FF$,
we let $\mathcal{O}_A(U)_{Y,X}$ denote the set of functions
\begin{equation}
s:U\to \bigcup_{\g{p}\in U}(A_{\g{p}})_{Y,X},
\end{equation}
such that $s(\g{p})\in(A_{\g{p}})_{Y,X}$, and $s$ is "locally a fraction":
\begin{equation}\label{spectracondition}
\forall\ \g{p}\in U,\ \exists\ \text{a neighborhood}\ U_{\g{p}}\ \text{of}\
\g{p};\exists\ a\in A_{Y,X};\ \exists f\in A_{1,1}\minus\bigcup_{\g{q}\in U_{\g{p}}}\g{q}
\end{equation}
$$\text{ such that }\; s(\g{q})=a/f\in A_{\g{q}},\;\forall\;\g{q}\in U_{\g{p}}\en{\bigstar}$$}

It is clear that
\begin{equation}
\mathcal{O}_A(U)=\bigcup_{Y,X}\mathcal{O}_A(U)_{Y,X}
\end{equation}
is an $\fr$. If $U'\subseteq U$, the natural restriction map $s\mapsto s|_{U'}$,
is a homomorphism of $\frs\;\;\mathcal{O}_A(U)\to \mathcal{O}_A(U')$,
thus  $\mathcal{O}_A$ is a presheaf of $\frs$. From the local nature
of $(\bigstar)$ we see that $\mathcal{O}_A$ is in fact a sheaf of $\frs$ over $\Spec A$, in the sense that for any $X,Y\in \FF$, $U\mapsto \mathcal{O}_A(U)_{Y,X}$ is a sheaf of (pointed) sets.

\prop{3.4.2}\textit{For $\g{p}\in\Spec(A)$, the stalk}
\begin{equation}
\mathcal{O}_{A,\g{p}}=\varinjlim\limits_{\g{p}\in U}\mathcal{O}_A(U)
\end{equation}
\textit{ of the sheaf $\mathcal{O}_A$ is isomorphic to $A_{\g{p}}$.}

\begin{proof}The map taking a local section $s$ in a
neighborhood of $\g{p}$ to $s(\g{p})\in A_{\g{p}}$, induces a
homomorphism $\mathcal{O}_{A,\g{p}}\to  A_{\g{p}}$, which is clearly
surjective. It is also injective:\\ Let
$s_1,s_2\in\mathcal{O}_A(U)_{Y,X}$ have the same value at
$\g{p},\;s_1(\g{p})=s_2(\g{p})$. Shrinking $U$ we may assume
$s_i=a_i/f_i$ on $U,\;a_i\in A_{Y,X},\;f_i\in
A_{1,1}.\;\;a_1/f_1=a_2/f_2\;\text{in}\;A_{\g{p}}$ means $h\cdot f_2\cdot a_1=h\cdot f_1\cdot a_2,\;h\in
A_{1,1}\minus\g{p}$, but then $a_1/f_1=a_2/f_2$ in $A_{\g{q}}\;\forall\g{q}\in U\cap
D_A(h)$.
\end{proof}

\prop{3.4.3}\textit{For $f\in A_{1,1}$, the $\fr\;\mathcal{O}_A(D_A(f))$ is isomorphic to $A_f$.\\
 In particular, the global sections $\Gamma(\Spec(A),\mathcal{O}_A)\Def\mathcal{O}_A(D_A(1))\cong A.$}

\begin{proof}Define the homomorphism
$\psi:A_f\to \mathcal{O}_A(D_A(f))$ by sending $a/f^n$ to the
section whose value at $\g{p}$ is the image of $a/f^n$ in
$A_{\g{p}}$.\par We shall show that $\psi$ is injective:\\If
$\psi(a_1/f^{n_1})=\psi(a_2/f^{n_2})$ then $\forall\;\g{p}\in
D_A(f)$ there is $h_{\g{p}}\in A_{1,1}\minus\g{p}$ with
\begin{equation}
h_{\g{p}}f^{n_2}a_1=h_{\g{p}}f^{n_1}a_2.
\end{equation}
Let
$\g{a}=ann_A(f^{n_2}a_1,f^{n_1}a_2)$, it is an ideal of
$A$, and $\forall\;\g{p}\in D_A(f),\;\g{p}\notin V_A(\g{a})$,
so $D_A(f)\cap V_A(\g{a})=\varnothing$, hence $V_A(\g{a})\subseteq
V_A(f)$, hence $f\in IV_A(\g{a})=\sqrt{\g{a}}$, hence
$f^n\in\g{a}$ for some $n\geq 1$, showing that $a_1/f^{n_1}=a_2/f^{n_2}$ in
$A_f$.\par
We show next that $\psi$ is surjective:\\
Let  $s\in\mathcal{O}_A(D_A(f))_{Y,X}$. By proposition $3.2.2,\ D_A(f)$ is compact,
so there exists a finite open covering
\begin{equation}
D_A(f)=\bigcup_{1\leq i\leq N}D_A(h_i),
\end{equation}
such that for all $\g{p}\in D_A(h_i):\ s(\g{p})=a_i/g_i \in A_{\g{p}}$,
where $a_i\in A_{Y,X}$ and $g_i\in A_{[1],[1]}$ is such that
$ D_A(g_i)\supseteq D_A(h_i)$ for $1\leq i\leq N$.\\
We have $V_A(g_i)\subseteq
V_A(h_i)$, hence
\begin{equation}
\sqrt{(g_i)}=IV_A(g_i)\supseteq IV_A(h_i)=\sqrt{(h_i)},
\end{equation}
hence $h_i\in \sqrt{(g_i)}$ so that for some $n_i\geq 1$ we have
$h_i^{n_i}=c_i\cdot g_i$, hence  $s(\g{p})=c_ia_i/h_i^{n_i}$.
So we can replace $h_i$ by $g_i$. On the set
\begin{equation}
D_A(g_i)\cap D_A(g_j)=D_A(g_ig_j)
\end{equation}
we have $a_i/g_i=s(\g{p})=a_j/g_j$, hence by the injectivity of
$\psi$
\begin{equation}
a_i/g_i=a_j/g_j\;\text{in}\; A_{g_ig_j}.
\end{equation}
This means $(g_ig_j)^n\cdot g_ja_i=(g_ig_j)^n\cdot
g_ia_j$, and we can choose $n$ big enough to work for all
$i,j$. We can replace $g_i$ by $g_i^{n+1}$\;(since
$D_A(g_i)=D_A(g_i^{n+1})$), and replace $a_i$ by $g_i^n\cdot
a_i$\;(since $s(\g{p})\equiv g_i^na_i/g_i^{n+1}$), and then
have the simpler equation
\begin{equation*}
g_j\cdot a_i=g_i\cdot
a_j\;\forall\;i,j.\end{equation*}
Since the sets $D_A(g_i)$ cover $D_A(f)$ we
have, ($cf$., \emph{Proposition} $3.2.2$),
\begin{equation}
f^m=b\circ(\oplus_ig_i)\circ b'.
\end{equation}
Set
\begin{equation}
a=(\underset{Y}{\oplus} b)\circ(\underset{i}{\oplus} a_i)\circ(\underset{X}{\oplus} b').
\end{equation}
Then
$$g_j\cdot a = (\underset{Y}{\oplus} b)\circ(\underset{i}{\oplus} g_ja_i)\circ(\underset{X}{\oplus} b')$$
$$=(\underset{Y}{\oplus} b)\circ(\oplus_i g_i a_j)\circ(\underset{X}{\oplus}b')$$
\begin{equation}
=\underset{Y}{\oplus}(b\circ (\underset{i}{\oplus}g_i)\circ b')\circ a_j \text{ (by commutativity!)}
\end{equation}
$$=f^m\cdot a_j .$$
Hence $a_j/g_j=s(\g{p})=a/f^m$ and $s=\psi(a/f^m)$.
\end{proof}

\section{Grothendieck $\FF$ - Schemes and locally-$\FF$-ringed spaces}
\smallbreak

We define the categories of $\FF$-(locally)-ringed-spaces, and its full subcategory of
(Grothendieck) $\FF$-schemes.

\defin{3.5.1}\textit{An \underline{\frsp} $(X,\mathcal{O}_X)$ is a topological space with a sheaf $\mathcal{O}_X$
of $\frs$: $U\mapsto \mathcal{O}_X(U)$ is a pre-sheaf of $\frs$
such that for any $W,Z\in \FF, U\mapsto \mathcal{O}_X(U)_{W,Z}$ is
a sheaf of (pointed) sets. That is, a collection
$\mathcal{O}_X(U)$ together with restriction maps $\rho^U_V\in
\FR(\mathcal{O}_X(U),\mathcal{O}_X(V))$ for every inclusion of
open sets $V\subseteq U$, such that:}

$$\begin{tikzcd}
&\mathcal{O}_X(U)\arrow{d}{\rho^U_V}\arrow{ddl}[swap]{\rho^U_W}\\
& \mathcal{O}_X(V)\arrow{dl}{\rho^V_W}\\
\mathcal{O}_X(W)
\end{tikzcd} \text{is commutative for }W\subseteq V\subseteq U.  $$
\textit{And for all $Y,X\in \FF$: $U\mapsto \mathcal{O}_X(U)_{Y,X}$ is a sheaf,\\ i.e. for any open covering $U=\underset{i}{\bigcup}U_i$:}


$$ 0\longrightarrow \mathcal{O}_X(U)_{Y,X}\longrightarrow \underset{i}{\prod}\mathcal{O}_X (U_i)_{Y,X} \rightrightarrows \underset{i,j}{\prod}\mathcal{O}_X(U_i\cap U_j)_{Y,X} $$
is exact.

 \textit{A map of \frsps  $f:X\to  Y$ is a continuous map of the
   underlying topological spaces together with a map of sheaves of $\frs$ on $Y,\\
   f^{\#}:\mathcal{O}_Y\to  f_*\mathcal{O}_X$, i.e. for $U\subseteq Y$ open we have
   $f^{\#}_U:\mathcal{O}_Y(U)\to \mathcal{O}_X(f^{-1}U)$ a map of $\frs$, such that for $U'\subseteq U$:}
\begin{equation}
f^{\#}_U(s)|_{f^{-1}U'}=f^{\#}_{U'}(s|_{U'}).
\end{equation}
   \textit{The \frsp $X$ is \underline{\lfrsp}   if for all $\g{p}\in X$ the stalk
   $\mathcal{O}_{X,\g{p}}$ is a local $\fr$, i.e. contains  a unique maximal ideal
   $\g{m}_{X,\g{p}}$.\\
For a map of \frsps $f:X\to  Y$, and for $\g{p}\in X$,
   we get an induced homomorphism of $\frs$ on the stalks}

\begin{equation}
f^{\#}_{\g{p}}:\mathcal{O}_{Y,f(\g{p})}=\varinjlim\limits_{f(\g{p})\in V}\mathcal{O}_Y(V)\to
    \varinjlim\limits_{\g{p}\in f^{-1}V}\mathcal{O}_X(f^{-1}V)\to
    \varinjlim\limits_{\g{p}\in U}\mathcal{O}_X(U)=\mathcal{O}_{X,\g{p}}
\end{equation}
   \textit{A map $f:X\to  Y$ of \lfrsps is a map of \frsps such that
    $f^{\#}_{\g{p}}$ is a \underline{local homomorphism} for all
    $\g{p}\in X$, i.e.}
\begin{equation}
f^{\#}_{\g{p}}(\g{m}_{Y,f(\g{p})})\subseteq\g{m}_{X,\g{p}}\;\text{\textit{or equivalently}}\;
    (f^{\#}_{\g{p}})^{-1}\g{m}_{X,\g{p}}=\g{m}_{Y,f(\g{p})}.
\end{equation}
   \textit{We let $\FRS$ (resp. $\LFRS$) denote the category of
   $\FF$-(resp. locally)-ringed-spaces.} \\

   For a homomorphism of commutative $\frs$ $\fe:A\to  B$, for
   $\g{p}\in\Spec(B)$, we have a unique homomorphism $\fe_{\g{p}}:A_{\fe^{-1}\g{p}}\to
   B_{\g{p}}$, such that we have a commutative diagram,
   \begin{equation}
\xymatrix{
        A \ar^{\fe}[r]\ar[d] &  B\ar[d]\\
        A_{\fe^{-1}\g{p}}\ar^{\fe_{\g{p}}}[r] & B_{\g{p}}}
\end{equation}
   $\fe_{\g{p}}(a/s)=\fe(a)/\fe(s)$, and $\fe_{\g{p}}$ is a local homomorphism. \\
Thus $A\mapsto\Spec(A)$ is a contravariant functor from $\catfr$ to $\LFRS$. \\
It is the adjoint of the functor $\Gamma$ of
   taking global sections
\begin{equation}
\Gamma(X,\mathcal{O}_X)=\mathcal{O}_X(X),\;\Gamma(f)=f^{\#}_Y:\mathcal{O}_Y(Y)\to \mathcal{O}_X(X).
\end{equation}

   \prop{3.5.2}
\begin{equation}
\LFRS(X,\Spec(A))=\FR(A,\mathcal{O}_X(X)).
\end{equation}

   \begin{proof}For an \lfrsp $X$, and for a point $x\in X$, the
   canonical homomorphism $\phi_x:\mathcal{O}_X(X)\to \mathcal{O}_{X,x}$
   gives a prime
   $\mathcal{P}(x)=\phi^{-1}_x(\g{m}_{X,x})\in\Spec\mathcal{O}_X(X)$. The
   map $\mathcal{P}:X\to \Spec\mathcal{O}_X(X)$ is continuous:
\begin{equation}
\mathcal{P}^{-1}(D(f))=\{x\in X|\;\phi_x(f)\notin\g{m}_{X,x}\}
\end{equation}
   is open for $f\in\mathcal{O}_X(X)$. We have an induced homomorphism
\begin{equation}
\mathcal{P}^{\#}_{D(f)}:\mathcal{O}_X(X)_f\to \mathcal{O}_X(\{x\in
   X|\;\phi_x(f)\notin\g{m}_{X,x}\}),
\end{equation}
   making $\mathcal{P}$ a map of \frsps, and taking the direct limit over $f$ with\\
   $\phi_x(f)\notin\g{m}_{X,x}$ we get
\begin{equation}
\mathcal{P}^{\#}_x:\mathcal{O}_X(X)_{\mathcal{P}(x)}\to \mathcal{O}_{X,x},
\end{equation}
   showing $\mathcal{P}$ is a map of \lfrsps.\\
   To a homomorphism of $\frs\ \fe:A\to \mathcal{O}_X(X)$ we associate the map of \lfrsps
\begin{equation}
X\stackrel{\mathcal{P}}{\to }\Spec\mathcal{O}_X(X)\stackrel{\Spec(\fe)}{\longrightarrow}\Spec
   A.
\end{equation}
   Conversely, to a map $f:X\to \Spec A$ of \lfrsps (as in definition $3.5.1$) we associate
   its action on global sections
\begin{equation}
\Gamma(f)=f^{\#}_{\Spec A}:A=\mathcal{O}_A(\Spec
   A)\to \mathcal{O}_X(X).
\end{equation}
Clearly, $\Gamma(\Spec(\fe)\circ\mathcal{P})=\fe$.\\ Conversely, given a map $f:X\to \Spec A$ (as in definition $3.5.1$), for $x\in X$ we have a commutative diagram,
\begin{equation}
\xymatrix{
        A=\mathcal{O}_A(\Spec A)\ar^{\hspace{7mm}\Gamma(f)}[r]\ar^{\phi_{f(x)}}[d] &  \mathcal{O}_X(X)\ar_{\phi_x}[d]\\
        A_{f(x)}\ar^{f^{\#}_x}[r] & \mathcal{O}_{X,x}}.
\end{equation}
   \\
   Since $f^{\#}_x$ is assumed to be local,
   $(f^{\#}_x)^{-1}(\g{m}_{X,x})=\g{m}_{f(x)}$, and by the
   commutativity of the diagram we get
   $\Gamma(f)^{-1}(\mathcal{P}(x))=f(x)$, i.e. $f=(\Spec\Gamma(f))\circ\mathcal{P}$ is
   the continuous map associated to the homomorphism $\Gamma(f)$.
   Similarly, for $g\in A$, the commutativity of the diagram
\begin{equation}
\xymatrix{
        A \ar^{\Gamma(f)}[r]\ar[d] &  \mathcal{O}_X(X)\ar[d]\\
        A_{g}\ar^{\hspace{-8mm}f^{\#}_{D(g)}}[r] & \mathcal{O}_X(D_X(f^{\#}g))}
\end{equation}
   \\
   gives $f^{\#}_{D(g)}(a/g^n)=\Gamma(f)(a)/(\Gamma(f)(g))^n$, hence
   $f=(\Spec\Gamma(f))\circ\mathcal{P}$ as a map of \lfrsps.
   \end{proof}

   \cor{3.5.3}\textit{ For $\frs\;A,B$:}

\begin{equation}
    \LFRS(\Spec B,\Spec A)=\FR(A,B).
\end{equation}

   \defin{3.5.4}\textit{ A (Grothendieck) $\FF$-scheme is an \lfrsp $(X,\mathcal{O}_X)$,
   such that there is a covering by open sets $X=\cup_iU_i$, and the canonical
   maps}
\begin{equation}
\mathcal{P}:(U_i,\mathcal{O}_X|_{U_i})\to \Spec\mathcal{O}_X(U_i)
\end{equation}
  \textit{are isomorphisms
   of \lfrsps. A morphism of $\FF$-schemes is a map of  \lfrsps. \vspace{3mm}\\
 We denote the category of $\FF$-schemes by \underline{$\FS$}}\\

   $\FF$-schemes can be glued:

   \prop{3.5.5}\textit{Given a set of indices $I$, and for $i\in I$ given
   $X_i\in\FS$, and for $i\ne j,\ i,j\in I$, an
   isomorphism $\fe_{ij}:U_{ij}\iso U_{ji}$, with $U_{ij}\subseteq X_i$ open
   (and hence $U_{ij}$ are $\FF$-schemes), such that}
\begin{equation}
\fe_{ji}=\fe^{-1}_{ij}
\end{equation}
\begin{equation}
\fe_{ij}(U_{ij}\cap U_{ik})=U_{ji}\cap U_{jk},\;\text{\textit{and}}\;
   \fe_{jk}\circ\fe_{ij}=\fe_{ik}\;\text{\textit{on}}\;U_{ij}\cap
   U_{ik}.
\end{equation}
   \textit{There exists $X\in
    \FS$, and maps $\psi_i:X_i\to  X$, such that}
\begin{equation}
\psi_i\;\text{\textit{is an isomorphism of $X_i$ onto the open set $\psi_i(X_i)\subseteq X$}},
\end{equation}
\begin{equation}
X=\bigcup_i\psi_i(X_i),
\end{equation}
\begin{equation}
\psi_i(U_{ij})=\psi_i(X_i)\cap\psi_j(X_j),
\end{equation}
\begin{equation}
\psi_i=\psi_j\circ\fe_{ij}\;\text{\textit{on}}\; U_{ij}.
\end{equation}

   \begin{proof} Clear: glue the topological spaces and glue the
   sheaves of $\frs$. For $V\subseteq X$ open
\begin{equation}
\mathcal{O}_X(V)=\ker\left\{\prod_i\mathcal{O}_{X_i}(\psi_i^{-1}V)\rightrightarrows
   \prod_{i,j}\mathcal{O}_{X_i}(\psi_i^{-1}V\cap U_{ij})\right\}.
\end{equation}
   \end{proof}


\chapter{Symmetric Geometry}
\bigskip

In this section $A\in \CFR^t$ is commutative and has involution.

\section{Symmetric ideals and symmetric primes}
\smallbreak

\defin{4.1.1} \textit{Consider the set $A^+:=\{a\in A_{1,1}\;| a=a^t\}$.\\
An ideal $\mathfrak{a}\subseteq A_{1,1}$ is called \underline{symmetric}
if it is generated by $\mathfrak{a}^+:=\mathfrak{a}\cap A^+$. (In particular $\mathfrak{a}=\mathfrak{a}^t$).\\
Denote by $I^+(A)$ the symmetric ideals of $A$.}\vspace{3mm}\\
Given an indexed set of symmetric ideals $\mathfrak{a}_i\subseteq A,i\in I$, their symmetric intersection $\cap_I^+ \mathfrak{a}_i$ is the ideal generated
by the mutual symmetric elements:
$$\underset{I}{\cap}^+ \mathfrak{a}_i=\{b\circ(\underset{J}{\oplus} a_j)\circ b' | a_j\in \underset{I}{\cap} \mathfrak{a_i}^+   \}.\en{4.1.1}$$
Their sum is again a symmetric ideal:
$$\underset{I}{\Sigma} \mathfrak{a}_i=\{b\circ(\underset{J}{\oplus} a_j)\circ b' | a_j\in \underset{I}{\cup} \mathfrak{a_i}^+   \}.\en{4.1.2}$$
Their product $\mathfrak{a}\cdot \mathfrak{a}'$ is again the symmetric ideal:
$$\mathfrak{a}\cdot \mathfrak{a}'=\{b\circ(\underset{J}{\oplus} a_j\cdot a_j')\circ b' | a_j\in \mathfrak{a}^+, a_j'\in \mathfrak{a'}^+  \}.\en{4.1.3}$$

Let $\varphi:A\rightarrow B$ be a homomorphism of $\frs$. If $\mathfrak{b}\in \mathcal{I}^+(B)$, we define:
$$\varphi^*(\mathfrak{b})= \text{ideal generated by } \{a=a^t\in A^+,\; \varphi(a)\in \mathfrak{b}\}\equiv \langle \varphi^{-1}(\mathfrak{b}^+) \rangle \en{4.1.4}$$
and if $ \mathfrak{a}\in \mathcal{I}^+(A)$, we define,
$$\varphi_*(\mathfrak{a})=\text{ideal generated by } \varphi(\mathfrak{a}^+)\subseteq B^+\equiv \langle \varphi(\mathfrak{a}^+) \rangle .\en{4.1.5}$$

\prop{4.1.1} \textit{For any $\mathfrak{a}\in\mathcal{I}^+(A),\mathfrak{b}\in \mathcal{I}^+(B)$, we have the following:\vspace{1.5mm} \\
(1) $\mathfrak{a}\subseteq \varphi^*\varphi_*\mathfrak{a}.$\\
(2) $\mathfrak{b}\supseteq \varphi_*\varphi^*\mathfrak{b}.$\\
(3)\;$\varphi^*\mathfrak{b}=\varphi^*\varphi_*\varphi^*\mathfrak{b},\;\;\;\varphi^*\mathfrak{a}=\varphi^*\varphi_*\varphi^*\mathfrak{a}. $\\
(4) there is a bijection, via $\mathfrak{a}\mapsto \varphi_*\mathfrak{a}$ (with inverse map $\mathfrak{b}\mapsto \varphi^*\mathfrak{b}$), from the set
$$\{\mathfrak{a}\in \mathcal{I}^+(A)|\; \varphi^*\varphi_*\mathfrak{a}=\mathfrak{a}\}=\{\varphi^*\mathfrak{b}|\; \mathfrak{b}\in \mathcal{I}^+(B)\}\en{4.1.6}$$
to the set
$$\{\mathfrak{b}\in \mathcal{I}^+(B)|\; \varphi_*\varphi^*\mathfrak{b}=\mathfrak{b}\}=\{\varphi_*\mathfrak{a}|\; \mathfrak{a}\in \mathcal{I}^+(A)\}.\en{4.1.7}$$ }

Since the union of a chain of (proper) symmetric ideals is a (proper) symmetric ideal, we have

\thm{4.1.1 (Zorn)} \textit{There exists a maximal symmetric ideal $\mathfrak{m}\varsubsetneq A_{1,1}$.}

\defin{4.1.2} \textit{A symmetric ideal $\mathfrak{p}\subseteq A_{1,1}$ is called \underline{symmetric prime}:
$$S_{\mathfrak{p}}^+=A^+ \setminus \mathfrak{p} \text{ is multiplicatively closed}\;\;\; S_{\mathfrak{p}}^+\cdot S_{\mathfrak{p}}^+=S_{\mathfrak{p}}^+. \en{4.1.8}$$}
We denote by $\Spec^+A$ the set of symmetric prime ideals.\\
For a homomorphism of $\frs$ $\varphi:A\rightarrow B$, the pullback $\varphi^*$ induce a map
$$\varphi^* =\Spec (\varphi): \Spec^+\; B\rightarrow \Spec^+ A.\en{4.1.9}$$

\prop{4.1.2} \textit{(1) If $\mathfrak{m}\in \mathcal{I}^+(A)$ is a maximal symmetric ideal then it is symmetric prime. \\
(2) More generally, if $\mathfrak{a}\in \mathcal{I}^+(A)$, and given $f=f^t\in A^+$ such that,
$$\forall n\in \NN: \hsm f^n\notin \mathfrak{a}.$$
let $\mathfrak{m}$ be a maximal element of the set
$$\{\mathfrak{b}\in \mathcal{I}^+(A)|\mathfrak{b}\supseteq \mathfrak{a}, f^n \not \in \mathfrak{b} \; \forall n\in \NN \}\en{4.1.10}$$
Then $\mathfrak{m}$ is symmetric prime.}
\begin{proof}
(1) If $x=x^t,y=y^t\in A_{1,1}\setminus \mathfrak{m}$, the ideals $(x)+\mathfrak{m},(y)+\mathfrak{m}$ are the unit ideals. So we can write
$$1=b\circ (\underset{J}{\bigoplus} m_j)\circ d, \;\;\; \text{with }\; m_j=x\;\; \text{or }\; m_j\in \mathfrak{m}.\en{4.1.11a}$$
$$1=b'\circ (\underset{J}{\bigoplus} m_j')\circ d', \;\;\; \text{with }\; m_j'=y\;\; \text{or }\; m_j'\in \mathfrak{m}.\en{4.1.11b}$$
It then follows that,
$$1=1\cdot 1=b\circ \underset{J}{\bigoplus} m_j\circ d\circ b'\circ \underset{I}{\bigoplus} m_i'\circ d'\en{4.1.12}$$
$$=b\circ(\underset{J}{\bigoplus}b')\circ \bigg(\underset{I\otimes J}{\bigoplus}m_j\circ m_i'\bigg)\circ (\underset{I}{\bigoplus}d)\circ d'.$$
but $m_j\circ m_i'=x\circ y$ or $m_j\circ m_i'\in \mathfrak{m}$, so 1 is in the ideal generated by $\mathfrak{m}$ and $x\circ y$, and since $1\notin \mathfrak{m}$ then
$x\circ y\notin \mathfrak{m}$. \\
(2) Similarly, if $x\notin \mathfrak{m}$ then $f^n$ is in the ideal generated by $x$ and $\mathfrak{m}$ so $f^n=b\circ \underset{J}{\oplus}m_j \circ d$, with $m_j=x$ or $m_j\in \mathfrak{m}$.
If $y\notin \mathfrak{m}, f^{n'}=b'\circ \underset{I}{\oplus} m_i'\circ d'$, with $m_i=y$ or $m_i\in \mathfrak{m}$. It then follows that
$$f^{n+n'}=b\circ \underset{J}{\bigoplus}m_j\circ d\circ b'\circ \underset{I}{\bigoplus} m_i'\circ d'  \en{4.1.13}$$
$$=b\circ (\underset{J}{\bigoplus}b')\circ \bigg(\underset{J\otimes I}{\bigoplus} m_j\circ m_i' \bigg) \circ (\underset{I}{\bigoplus}d)\circ d' $$
but $m_j\circ m_i'=x\circ y$ or $m_j\circ m_i'\in \mathfrak{m}$, so $f^{n+n'}$ is in the ideal generated by $\mathfrak{m}$ and $x\circ y$, and since $f^{n+n'}\notin \mathfrak{m}$
then $x\circ y\notin \mathfrak{m}$ .
\end{proof}

\defin{4.1.3}\textit{For $\mathfrak{a}\in \mathcal{I}^+(A)$, the \underline{symmetric radical} is
$$\sqrt{\mathfrak{a}}^+=\text{ideal generated by } \{ f=f^t\in A^+, f^n\in \mathfrak{a}\; \text{for some } \; n\geq 1\} \en{4.1.14}$$
$$=\{b\circ \underset{J}{\oplus} f_j\circ d \;|\; b\in A_{1,J},d\in A_{J,1}, f_j=f_j^t\in A^+, \text{and } \hsm f_j^n\in \mathfrak{a} \}.$$ }
Note that $\sqrt{\mathfrak{a}}^+\subseteq \sqrt{a}$, and for any $a\in \sqrt{\mathfrak{a}}^+$, we have $a^n\in \mathfrak{a}$ for $n>>1$.

\prop{4.1.4} \textit{We have
$$ \sqrt{\mathfrak{a}}^+=\underset{\mathfrak{a}\subseteq \mathfrak{p}}{\bigcap}^+\mathfrak{p} =\text{i.e. the symmetric ideal generated by } \underset{\mathfrak{a}\subseteq \mathfrak{p}}{\bigcap} \mathfrak{p}^+ \en{4.1.16} $$
$$\text{where $\mathfrak{p}$ runs over symmetric primes containing } \mathfrak{a}.$$}
\begin{proof}
If $f=f^t, \; f^n\in \mathfrak{a}$, then for all symmetric primes $\mathfrak{a}\subseteq \mathfrak{p}$: $f\in \mathfrak{p}$. If $f=f^t$ and
$f^n\notin \mathfrak{a},\forall \; n$, let $\mathfrak{m}$ be a maximal element of the set (4.1.10), it exists by Zorn's lemma, and it is symmetric prime by proposition (4.1.2b), $\mathfrak{a}\subseteq \mathfrak{m}$
and $f\notin \mathfrak{m}$.
\end{proof}

\section{The symmetric spectrum: $\Spec^+ (A)$}
\smallbreak

\defin{4.2.1} \textit{For a set $\mathfrak{U}\subseteq A^+$, we let
$$V_A^+(\mathfrak{U})=\{\mathfrak{p}\in \Spec^+(A)| \mathfrak{U}\subseteq \mathfrak{p} \} \en{4.2.1}$$}
If $\mathfrak{a}$ is the ideal generated by $\mathfrak{U}$, $V_A^+(\mathfrak{U})=V_A^+(\mathfrak{a})$; we have
$$V_A^+(1)=\emptyset,\hsm V_A^+(0)=\Spec^+ (A),\en{4.2.2}$$
$$V_A^+(\Sigma \mathfrak{a})=\underset{i}{\cap} V_A^+(\mathfrak{a}_i)\;,\; \mathfrak{a}_i\in \mathcal{I}^+(A),$$
$$V_A^+(\mathfrak{a}\cdot \mathfrak{a}') = V_A^+(\mathfrak{a})\cup V_A^+(\mathfrak{a}').$$
Hence the sets $\{V_A^+(\mathfrak{a})\; |\; \mathfrak{a}\in
\mathcal{I}^+(A)\}$ are the closed sets for the topology on $\Spec^+
(A)$, the \textit{Zariski topology}.

\defin{4.2.2} \textit{For $f=f^t\in A^+$ we let
$$D_A^+(f)=\Spec^+(A)\setminus V_A^+(f)=\{\mathfrak{p}\in \Spec^+ (A)| f\notin \mathfrak{p} \}.\en{4.2.3}$$
We have
$$D_A^+(f_1)\cap D_A^+(f_2)=D_A^+(f_1\cdot f_2),\en{4.2.4}$$
$$\Spec^+ A\setminus V_A^+(\mathfrak{a})=\underset{f\in \mathfrak{a}^+}{\bigcup} D_A^+(f).$$}
Hence the sets $\{D_A^+(f)|f\in A^+\}$ are the basis for the open sets in the Zariski topology. We have
$$D_A^+(f)=\emptyset \iff f\in \underset{\mathfrak{p}\in \Spec^+ A}{{\bigcap}^+} \mathfrak{p}= \sqrt{0}^+\iff f^n=0 \text{ for some } n\en{4.2.5} $$
and we say $f$ is a \textit{nilpotent}. We have
$$D_A^+(f)=\Spec^+ A \iff (f)=(1) \iff \exists f^{-1}\in A^+ : f\cdot f^{-1}=1 \en{4.2.6}$$
and we say $f$ is \textit{invertible}. We denote by $GL^+_1(A)$ the (commutative) group of symmetric invertible elements.

\defin{4.2.3} \textit{For a subset $X\subseteq \Spec^+(A)$, we have the associated ideal
$$\mathcal{I}^+(X)=\underset{\mathfrak{p}\in X}{\bigcap}^+ \mathfrak{p}=\text{ideal generated by }\underset{\mathfrak{p}\in X}{\bigcap}\mathfrak{p}^+ .\en{4.2.7}$$}

\prop{4.2.1} \textit{We have
$$ \mathcal{I}^+ V_A^+ \mathfrak{a}= \sqrt{\mathfrak{a}}^+, \en{4.2.8a}$$
$$V^+_A\mathcal{I}^+(X)=\bar{X},\text{ the closure of } X \text{ in the Zariski topology.} \en{4.2.8b}$$}

\begin{proof}
The first equation is just a restatement of proposition $4.1.4$. Indeed,
$$\mathcal{I}^+ V_A^+ \mathfrak{a}=\text{ideal generated by } \underset{\mathfrak{a}\subseteq \mathfrak{p}\in \Spec^+(A)}{\bigcap} \mathfrak{p}^+=\sqrt{\mathfrak{a}}^+\en{4.2.9}$$
For the second, $V_A^+\mathcal{I}^+(X)$ is clearly a closed set containing $X$, and if $C=V_A^+(\mathfrak{a})$ is a closed set containing $X$,
then $\sqrt{\mathfrak{a}}^+=\mathcal{I}^+V_a^+(\mathfrak{a})\subseteq \mathcal{I}^+(X) $,
hence $C=V_A^+(\sqrt{(\mathfrak{a})}^+)\supseteq V_A^+\mathcal{I}^+(X)$
\end{proof}

\cor{4.2.1} \textit{We have a one-to-one order-reversing correspondence between closed sets $X\subseteq \Spec^+(A)$, and radical symmetric ideals $\mathfrak{a}$, via $X\mapsto \mathcal{I}^+(X),\; V_A^+(\mathfrak{a})\mapsfrom \mathfrak{a}$
$$\{ X\subseteq \Spec^+(A)|\bar{X}=X\}\overset{1:1}{\longleftrightarrow} \{\mathfrak{a}\in \mathcal{I}^+(A)| \sqrt{\mathfrak{a}}^+=\mathfrak{a}\}.  \en{4.2.10}$$}
Under this correspondence the closed irreducible subsets corresponds to symmetric prime ideals. For $\mathfrak{p}_0,\mathfrak{p}_1\in \Spec^+(A)$, $\mathfrak{p}_0 \in \bar{\{\mathfrak{p}_1\}}\iff \; \mathfrak{p}_0\supseteq\mathfrak{p}_1$,
we say that $\mathfrak{p}_0$ is a \textit{Zariski specialization} of $\mathfrak{p}_1$, or that $\mathfrak{p}_1$ is the \textit{Zariski generalization} of $\mathfrak{p}_0$.
The space $\Spec^+(A)$ is sober: every closed irreducible subset $C$ has the form $C=V^+_A(\mathfrak{p})=\bar{\{\mathfrak{p}\}}$, and we call the (unique) prime $\mathfrak{p}$ the \textit{generic point} of $C$.

\prop{4.2.2} \textit{The sets $D_A^+(f)$, and in particular $D_A^+(1)=\Spec^+(A)$, are compact.}
\begin{proof}
Note that $D_A^+(f)$ is contained in the union $\underset{i}{\bigcup} D_A^+(g_i)$ if and only if $V_A^+(f)\supseteq \underset{i}{\bigcap} V_A^+(g_i)=V_A^+(\mathfrak{a})$, where
$\mathfrak{a}$ is the ideal generated by $\{g_i\}$, if and only if $\sqrt{f}^+=\mathcal{I}^+V_A^+(f)\subseteq \mathcal{I}^+V_A^+(\mathfrak{a})=\sqrt{\mathfrak{a}}^+$, if and only if
$f^n\in \mathfrak{a}^+$ for some $n$, if and only if $f^n=b\circ (\underset{i}{\bigoplus} g_i)\circ b'$, and in such expression only a finite number of the $g_i$ are involved.
\end{proof}

Let $\varphi:A\rightarrow B$ be a homomorphism of $\fr$ with involution, $\varphi^*:\Spec^+(B)\rightarrow \Spec^+(A)$ the associated pullback map.

\prop{4.2.3} \textit{We have
$$\varphi^{*-1}(D_A^+(f))=D_B^+(\varphi(f)), \;\;\; f\in A^+, \en{4.2.11}$$
$$\varphi^{*-1}(V_A^+(\mathfrak{a}))=V_B^+(\varphi_*(\mathfrak{a})),\;\;\; \mathfrak{a}\in \mathcal{I}^+(A), \en{4.2.12}$$
$$V_A^+(\varphi^{-1}\mathfrak{b})=\bar{\varphi^*(V_B^+(\mathfrak{b}))},\;\;\; \mathfrak{b}\in \mathcal{I}^+(B).\en{4.2.13}$$}
\begin{proof}
The proofs of $(4.2.11)$ and $(4.2.12)$ are straightforward:
$$\mathfrak{q}\in \varphi^{*-1}(D_A^+(f))\iff \varphi^*(\mathfrak{q})\in D_A^+(f) \iff f\notin \varphi^{-1}(\mathfrak{q})
\iff \varphi (f)\notin \mathfrak{q}\iff \mathfrak{q}\in D_B^+(\varphi(f)), $$
$$\mathfrak{q}\in \varphi^{*-1}(V_A^+(\mathfrak{a}))\iff \varphi^*(\mathfrak{q})\in V_A^+(\mathfrak{a}) \iff \mathfrak{a}^+\subseteq \varphi^{-1}(\mathfrak{q})
 \iff \varphi_*(\mathfrak{a}) \subseteq \mathfrak{q} \iff \mathfrak{q}\in V_B^+(\varphi_*(\mathfrak{a})). $$
For $(4.2.13)$ we may assume $\mathfrak{b}=\sqrt{\mathfrak{b}}^+$ is a radical since $V_B^+(\mathfrak{b})=V_B^+(\sqrt{\mathfrak{b}}^+), \varphi^*(\sqrt{\mathfrak{b}}^+)=\sqrt{\varphi^* (\mathfrak{b})}^+$.
Let $\mathfrak{a}=\mathcal{I}^+(\varphi^*(V_B^+(\mathfrak{b})))$, so that $V_A^+(\mathfrak{a})=\bar{\varphi^*(V_B^+(\mathfrak{b}))}$ by $(4.2.8b)$. We have
$$f\in \mathfrak{a} \iff f\in \mathfrak{p}, \; \forall \mathfrak{p} \in \varphi^*(V_B^+(\mathfrak{b})) \iff f\in \varphi^{-1}(\mathfrak{q}), \;\;\; \forall \mathfrak{q}\supseteq \mathfrak{b}$$
$$\varphi (f)\in \underset{\mathfrak{q}\supseteq \mathfrak{b}}{\bigcap}^+ \mathfrak{q}=\sqrt{\mathfrak{b}}^+=\mathfrak{b}\iff f\in \varphi^*(\mathfrak{b}). \en{4.2.14}$$
\end{proof}

It follows from $(4.2.11)$, or from $(4.2.12)$, that $\varphi^*=\Spec^+ (\varphi)$ is continuous, hence $A\mapsto \Spec^+(A)$ is a contravariant functor from $\Cfrs^t$
to compact, sober, topological spaces.

\exmpl{4.2.1} Let $A$ be a commutative ring with (nontrivial) involution $a \mapsto a^t$, $\FF(A)$ the associated $\fr$ with involution $(a^t)_{x,y}=(a_{y,x})^t,a=(a_{y,x})\in A_{Y,X}$. An ideal $\mathfrak{a}\subseteq A=\FF(a)_{1,1}$ is also an ideal in our sense, and conversely.
Under this correspondence the primes of $A$ correspond to the primes of $\FF(A)$, and we have a homeomorphism with respect to the Zariski topologies: $\Spec\; \FF(A)\simeq \Spec\;A$.
The symmetric primes correspond to primes of the subring
$$A^+=\{a\in A\;|\; a^t=a\} \en{4.2.15}$$
consisting of symmetric elements. (note that $A/A^+$ is always integral, $\alpha \in A$ is a root of $x^2-(\alpha+\alpha^t)x+\alpha\alpha^t\in A^+[x]$).
We have as well a homeomorphism with respect to the topology defined in $(4.2.1-2)$: $\Spec^+ \;\FF (A)\simeq \Spec \;A^+.$ \\

Given $A\in \CFR^t$, we can forget the involution and consider $\Spec\; A$. There is a canonical map
$$\pi_A^+:\Spec\; A\rightarrow \Spec^+\;A$$
$$\mathfrak{p}\mapsto \pi_A^+(\mathfrak{p})=\text{ ideal generated by } \mathfrak{p}^+.\en{4.2.16}$$
and this map is continuous:
$$(\pi_A^+)^{-1}(V_A^+(\mathfrak{a}))=V_A(\mathfrak{a})\en{4.2.17a}$$
$$(\pi_A^+)^{-1}(D_A^+(f))=D_A(f). \en{4.2.17b}$$

\section{Symmetric localization}
\smallbreak
A set $S\subseteq A_{1,1}^+$ (consisting of symmetric elements!) is called "multiplicative" when
\begin{equation}
1\in S
\end{equation}
\begin{equation}
S\cdot S=S
\end{equation}
For such $S$ the localized $\fr$ $S^{-1}A$ has involution:
\begin{equation}
(a/s)^t=a^t/s.
\end{equation}
The following universal property holds:
\begin{equation}
\FR^t(S^{-1}A,B)=\{\varphi\in \FR^t(A,B), \varphi(S)\subseteq GL_1^+(B)\}.
\end{equation}
The main examples of localizations are:
\begin{equation}
\text{For } \mathfrak{p}\in \Spec^+A,\hsm\hsm S_{\mathfrak{p}}=A^+\setminus \mathfrak{p}\;,\; S^{-1}_{\mathfrak{p}}A:=A_{\mathfrak{p}}
\end{equation}
\begin{equation}
\text{For } f=f^t\in A^+,\hsm\hsm S_f=\{1,f,f^2,\dots, f^n,\dots\}\;,\; S^{-1}_fA:=A_f=A[\frac{1}{f}].
\end{equation}

\section{Structure sheaf $\mathcal{O}_A^t/ \Spec^+\;(A)$}
\smallbreak
\defin{4.4.1}\textit{ A sheaf of $\FR^t$ (with involution!) over a topological space $X$, $\mathcal{O}\in \FR^t/X$ is a pre-sheaf of $\frs$ with involutions $U\mapsto \mathcal{O}(U)$, such that for any $W,Z\in \FF$, $U\mapsto \mathcal{O}(U)_{W,Z}$ is a sheaf.}
\\

Next we define a sheaf $\mathcal{O}_A^t$ of $\FR^t$ over $\Spec^+\;A$.

\defin{4.4.2}\textit{For an open set $U\subseteq \Spec^+\;A,\; X,Y\in \FF:$
$$\mathcal{O}_A (U)_{Y,X}:= \{s:U\rightarrow \underset{\mathfrak{p}\in U}{\coprod} (A_{\mathfrak{p}})_{Y,X},\; s(\mathfrak{p})\in (A_{\mathfrak{p}})_{Y,X}\;,\; \text{and satisfy } (*) \}$$
$$(*)\hspace{5mm} \forall \mathfrak{p}\in U, \exists f=f^t\notin \mathfrak{p},\; \mathfrak{p}\in D^+(f)\subseteq U, \exists a\in A_{Y,X},$$
\begin{equation}
s(\mathfrak{q})=a/f,\hsm \forall \mathfrak{q}\in D^+(f).
\end{equation}}

\thm{4.4.1}\textit{(i) For $\mathfrak{p}\in \Spec^+A$ we have
\begin{equation}
\mathcal{O}_{A,\mathfrak{p}} =\underset{\mathfrak{p}\in U}{\underset{\longrightarrow}{\lim}} \mathcal{O}_A\;(U)\iso A_{\mathfrak{p}}.
\end{equation}
(ii)   For $f=f^t\in A^+$,
\begin{equation}
\mathcal{O}_A(D^+(f))\overset{\sim}{\leftarrow}A_f.
\end{equation}   }

\begin{proof}
The proofs are exactly as in propositions $3.4.2$ and $3.4.3$.
\end{proof}

\section{Schemes with involution $\FSt$ and locally- $\FR^t$- spaces $\LFRtS$}
\smallbreak

 \defin{4.5.1}  \textit{ An \underline{$\FR^t$-space}
   $(X,\mathcal{O}_X)$ is a topological space $X$ with a sheaf $\mathcal{O}_X$
   of $\frs$ with involutions. A map of $\FR^t$-spaces  $f:X\to  Y$ is a continuous map of the
   underlying topological spaces together with a map of sheaves of $\FR^t$ on $Y,\\
   f^{\#}:\mathcal{O}_Y\rightarrow  f_*\mathcal{O}_X$   , i.e. for $U\subseteq Y$ open we have
   $f^{\#}_U:\mathcal{O}_Y(U)\rightarrow \mathcal{O}_X(f^{-1}U)$ a map of $\FR^t$, such that for
   $$U'\subseteq U\;:\;\;\;\;\;f^{\#}_U(s)|_{f^{-1}U'}=f^{\#}_{U'}(s|_{U'}).$$
The $\FR^t$-space $X$ is \underline{a locally-$\FR^t$-space} if for all $\g{p}\in X$ the stalk
   $\mathcal{O}_{X,\g{p}}$ is a local $\FR^t$, i.e. contains  a unique maximal \underline{symmetric} ideal
   $\g{m}_{X,\g{p}}$.\\ For a map of $\FR^t$ -spaces $f:X\to  Y$, and for $\g{p}\in X$,
   we get an induced homomorphism of $\FR^t$ on the stalks
\begin{equation}
f^{\#}_{\g{p}}:\mathcal{O}_{Y,f(\g{p})}=\varinjlim_{f(\g{p})\in V}\mathcal{O}_Y(V)\rightarrow
    \varinjlim_{\g{p}\in f^{-1}(V)}\mathcal{O}_X(f^{-1}V)\rightarrow
    \varinjlim_{\g{p}\in U}\mathcal{O}_X(U)=\mathcal{O}_{X,\g{p}}
\end{equation}
   A map $f:X\to  Y$ of locally- $\FR^t$-spaces is a map of $\FR^t$-spaces such that
    $f^{\#}_{\g{p}}$ is a \underline{local homomorphism} for all
    $\g{p}\in X$, i.e.
\begin{equation}
f^{\#}_{\g{p}}(\g{m}_{Y,f(\g{p})})\subseteq\g{m}_{X,\g{p}}\;\text{\textit{ or equivalently}}\;
    (f^{\#}_{\g{p}})^*(\g{m}_{X,\g{p}})=\g{m}_{Y,f(\g{p})}.
\end{equation}   }
  We let $\LFRtS$ denote the category of $\FF$-locally-ringed-spaces with involution.
  \defin{4.5.2}\textit{A pair $(X,\mathcal{O}_X)\in \LFRtS$ is an  (Grothendieck)\hsm  \underline{$\FF-\mathcal{S}cheme$} with involution},
  if there exists an open covering $\{U_i\}_{i\in I}$ of the topological space $X$ s.t.
\begin{equation}
(U_i,\mathcal{O}_X|_{U_i})\iso \Spec^+(\mathcal{O}_X(U_i)).
\end{equation}
  We denote the category of Grothendieck $\FF$-$\mathcal{S}cheme$ with involution by $\FSt$.
Note that we have the following full embedding of categories:

\begin{equation}
(\FR^t)^{op}\subseteq \FSt\subseteq \LFRtS.
\end{equation}

\chapter{Pro - limits}
\bigskip
We can work with general filtered small category $J$ (cf. \cite{AM}), or restrict attention to the case of $(J,\leq)$ a partial ordered set that is directed,
$$\forall j_1,j_2\in J, \;\; \exists j\in J,\; j\geq j_1,j\geq j_2.$$
and co-finite,
$$\forall j\in J,\hspace{3mm} \bigg| \{i\in J|\;\; i\leq j\}\bigg| <\infty.$$
The following inverse limits over $J$-indexed inverse systems in $(\FR^t)^{op}$ and $\LFRtS$ exist:
$$\begin{tikzcd}
\underset{J}{\varprojlim}:&(\LFRtS)^J\arrow{r}&\LFRtS \\
&\begin{turn}{90}$\subseteq$ \end{turn}&\begin{turn}{90}$\subseteq$ \end{turn} \\
\underset{J}{\varprojlim}:&((\FR^t)^{op})^J\arrow{r}&(\FR^t)^{op}
\end{tikzcd}$$
Moreover in the affine case the following limit can be interchanged:
$$\underset{J}{\varprojlim}\Spec^+(A_j)=\Spec^+(\underset{J}{\varinjlim} A_j).$$
But inverse limits over $(\FSt)^J$ do not always exist!. This leads us to define the "pro" category of $\FSt$ as our category of schemes.

\section{Pro - Schemes }
\smallbreak

\defin{5.1} \textit{Define the $\text{pro-}\FS$ (respectfully $\text{pro-}\FSt$) category to be the category with objects all inverse systems in $\FS$ (respectfully $\FSt$ ) over arbitrary directed co- finite partially ordered set.
The morphisms between objects $X=\{X_j\}_{j\in J},Y=\{Y_i\}_{i\in I}$ are given by,
\begin{equation}
\text{pro-} \FS (X,Y)=\underset{I}{\varprojlim}\;\;\underset{J}{\varinjlim}\FS (X_j,Y_i),
\end{equation}
$$ (respectfully.   \hsm  \text{pro-} \FSt (X,Y)=\underset{I}{\varprojlim}\;\;\underset{J}{\varinjlim}\FSt (X_j,Y_i)).$$}
For every $i'\leq i,j\leq j'$ we have a commutative diagram in
$Set$:
\begin{equation}
\begin{tikzcd}
\FS (X_j,Y_i)\arrow{d}{\pi_{i'}^{i}\circ \underline{\;\;}}\arrow{r}{\underline{\;\;}\circ \pi_j^{j'}}& \FS (X_{j'},Y_i)\arrow{d}{\pi_{i'}^{i}\circ \underline{\;\;}}\\
\FS (X_j,Y_{i'})\arrow{r}{\underline{\;\;}\circ  \pi_j^{j'}}& \FS (X_{j'},Y_{i'})
\end{tikzcd}
\end{equation}
(here $\pi_{i'}^i:Y_i\rightarrow Y_{i'}$ and $\pi_j^{j'}:X_{j'}\rightarrow X_j $). \\
We can describe a morphism between $X$ and $Y$
in the pro- category as collections of maps
$\{\varphi^j_i:X_j\rightarrow Y_i\}$ defined for every $i$ and for
every $j\geq \sigma(i)$ large enough (depending on $i$). The maps
are such that the following conditions must hold for every $i'\leq
i$ and $\sigma(i) \leq j\leq j'$,
$$\varphi^j_i\circ \pi^{j'}_j=\varphi^{j'}_i$$
\begin{equation}
\pi^i_{i'}\circ \varphi^j_i=\varphi^j_{i'}
\end{equation}
The maps $\{\varphi^j_i\},\{\tilde{\varphi}^j_i\}$ are considered equivalent if for all $i\in I$ and all $j\in J$ large enough (depending on $i$): $\varphi_i^j=\tilde{\varphi}_i^j$. \\
We have a full and faithful embedding $\FSt\inj \; \text{pro-}
\FSt$, (taking the indexing set to be a point). We have a functor
$\varprojlim$ (which is generally not full, and not faithful, but
is such on finitely presented objects), making a commutative
diagram
\begin{equation}
\begin{tikzcd}
\text{pro- }\FSt\arrow{rr}{\varprojlim} &&  \LFRtS\\
&\FSt \arrow[hook]{lu}\arrow[hook]{ru}&
\end{tikzcd}
\end{equation}
\section{The compactified $\Spec \ZZ$.}
\smallbreak

Let $N\geq 2$ be a square free integer.
Let $A_N$ be the $\fr$ defined by:
$$A_N=\FF(\ZZ[\frac{1}{N}])\cap\mathcal{O}_{\QQ,\eta},$$
\begin{equation}
(A_N)_{Y,X}=\bigg\{a\in \FF\bigg(\ZZ\bigg[\frac{1}{N}\bigg]\bigg)_{Y,X}\;\;\bigg|\;\; |a|_{\eta}\leq 1\bigg\}.
\end{equation}
The map $j:A_N\rightarrow \FF(\ZZ[\frac{1}{N}])$ is a localization, $(A_N)_{\frac{1}{N}}\cong \FF(\ZZ[\frac{1}{N}])$, and so defines an injection:
\begin{equation}
j^*:\Spec \;\ZZ\bigg[\frac{1}{N}\bigg]\cong \Spec\; \FF\bigg(\ZZ\bigg[\frac{1}{N}\bigg]\bigg)\cong \Spec\; (A_N)_{\frac{1}{N}}\cong D_{A_N}\bigg(\frac{1}{N}\bigg)\inj \Spec\;A_N.
\end{equation}
The space $\Spec A_N$ also contains the closed point,
\begin{equation}
\eta= i^*(\mathfrak{m}_{\QQ,\eta})=\bigg\{a\in (A_N)_{1,1}\bigg| |a|_{\eta} <1\bigg\}=\ZZ\bigg[\frac{1}{N}\bigg]\cap (-1,1).
\end{equation}
which is the real prime given by the injection $i:A_N\hookrightarrow \mathcal{O}_{\QQ,\eta}$. The prime ideal $\eta$ contains any other ideal of $A_N$ thus it is the unique maximal ideal for $A_N$, which is a local $\fr$ (of Krull dimension $2$).

Note that $\Spec A_N=\Spec \ZZ[\frac{1}{N}]\cup \{\eta\}$ as sets. The point $\eta\in \Spec A_N$ is very closed in the sense that the only open set containing it is the whole space. For any non-trivial basic open set $D_{A_N}(f)$, say $f=\frac{a}{N^k},a\in\ZZ, |a|_{\eta}<N^k$, we have $(A_N)_f=\FF(\ZZ[\frac{1}{N\cdot a}])$, and so
\begin{equation}
D_{A_N}(f)=\Spec(A_N)_f=\Spec \FF\bigg(\ZZ\bigg[\frac{1}{N\cdot a}\bigg]\bigg)\simeq \Spec \ZZ\bigg[\frac{1}{N\cdot a}\bigg]
\end{equation}
does not contain $\eta$.

Let $X_N\in \FS$ be the Grothendieck-$\FF$-scheme defined by:
\begin{equation}
X_N=\Spec \FF(\ZZ)\underset{\Spec \FF(\ZZ[\frac{1}{N}])}{\coprod} \Spec  A_N.
\end{equation}
i.e. gluing $\Spec  \FF(\ZZ)$ with $\Spec A_N$ along the common
basic open set $\Spec \FF(\ZZ[\frac{1}{N}])$.
\begin{equation}
\begin{tikzcd}
&X_N& \\
\Spec\FF(\ZZ)\arrow[hook]{ru}{i_1}&&\Spec A_N\arrow[hook]{lu}{i_2} \\
&\Spec\FF(\ZZ[\frac{1}{N}])\arrow[hook]{lu}\arrow[hook]{ru}&
\end{tikzcd}
\end{equation}
The open sets of $X_N$ are open sets of $\Spec \ZZ$ and sets of the form $U\cup\{\eta\}$ with $\Spec \ZZ[\frac{1}{N}]\subseteq U\subseteq \Spec \ZZ$.

The points and specializations of $X_N$, $N=p_1\dots p_k$,

\begin{equation}
\begin{tikzpicture}[scale=0.3,auto=left]
\tikzstyle{no}=[circle,draw,fill=black!100,inner sep=0pt, minimum width=4pt]
\node (n1) at (28,10) {$\Spec A_N$  };

\node[no] (n2) at (17,11)  {};
\node (n12) at (17,10){$\eta$};

\node[no] (n3) at (6,7)  {};
\node (n13) at (6,8){$p_1$};
\node[no] (n4) at (9,7)  {};
\node (n14) at (9,8){$p_2$};
\node (n5) at (12,7)  {$\dots$};
\node[no] (n6) at (15,7)  {};
\node (n15) at (15,8){$p_k$};

\node[no] (n7) at (22,7)  {};
\node (n16) at (23,7){$q_1$};
\node[no] (n8) at (25,7)  {};
\node (n17) at (26,7){$q_2$};
\node[no] (n9) at (28,7)  {};
\node (n18) at (29,7){$q_3$};
\node (n10) at (31,7)  {$\dots$};

\node[no] (n11) at (17,3)  {};
\node (n16) at (17,2){$(0)$};
\node(n17) at (27,2.5){($q_i$ the primes $\nmid N$). };

\draw [-] (n2) --  (n7);
\draw [-] (n2) --  (n8);
\draw [-] (n2) --  (n9);

\draw [-] (n11) --  (n7);
\draw [-] (n11) --  (n8);
\draw [-] (n11) --  (n9);

\draw [-] (n3) --  (n11);
\draw [-] (n4) --  (n11);
\draw [-] (n6) --  (n11);

\end{tikzpicture}
\end{equation}

The structure sheaf $\mathcal{O}_{X_N}$ is defined as follows: for any open set $U\subseteq X_N$,

$$\mathcal{O}_{X_N}(U)=\bigg\{<s_1,s_2>\bigg| s_1\in\mathcal{O}_{\FF(\ZZ)}(i_1^{-1}(U)) \text{ and } s_2\in \mathcal{O}_{A_N}(i_2^{-1}(U)) \text{ and } $$
\begin{equation}
  s_1|_{i_1^{-1}(U)\cap \Spec \FF(\ZZ[\frac{1}{N}])}\simeq s_2|_{i_2^{-1}(U)\cap \Spec \FF(\ZZ[\frac{1}{N}])} \bigg\}
\end{equation}

For an open set $U=\Spec \ZZ[\frac{1}{M}]= D(M)\subseteq \Spec \ZZ$ , we have
\begin{equation}
s_1\in \mathcal{O}_{\FF(\ZZ)}(D(M))\simeq \FF(\ZZ[\frac{1}{M}]),
\end{equation}
\begin{equation}
s_2\in \mathcal{O}_{A_N} (i_2^{-1}(U)), \;\; i_2^{-1}(U)\subseteq \Spec \FF(\ZZ[\frac{1}{N}]).
\end{equation}
the isomorphism condition gives us in that case that
\begin{equation}
s_1|_{i_1^{-1}(U)\cap \FF(\ZZ[\frac{1}{N}])}=s_2,
\end{equation}
and therefore $\mathcal{O}_{X_N}(U)= \FF(\ZZ[\frac{1}{M}])$  .

For a set $U= \Spec \ZZ [\frac{1}{M}]$ with $M|N$ , we have,

\begin{equation}
i_1^{-1}(U\cup \{\eta\})=\Spec \FF\bigg(\ZZ[\frac{1}{M}]\bigg),
\end{equation}
\begin{equation}
i_2^{-1}(U\cup \{\eta\})=\Spec A_N.
\end{equation}
and so,

\begin{equation}
 s_1\in \mathcal{O}_{\FF(\ZZ)}(D(M))\simeq \FF(\ZZ[\frac{1}{M}])
\end{equation}
\begin{equation}
s_2\in \mathcal{O}_{A_N}(D(1)) \simeq A_N
\end{equation}
where $s_1|_{\Spec \FF(\ZZ[\frac{1}{N}])}=s_2|_{\Spec \FF(\ZZ[\frac{1}{N}])}$.\\
Thus $\mathcal{O}_{X_N}(U\cup \{\eta\})$ is the pullback of the diagram

\begin{equation}
\begin{tikzcd}
A_N\arrow[hook]{rd}&&\FF(\ZZ[\frac{1}{M}])\arrow[hook]{ld}\\
&  \FF(\ZZ[\frac{1}{N}])&
\end{tikzcd}
\end{equation}
which is $A_M$, and so, $\mathcal{O}_{X_N}(U\cup \{\eta\}) = A_M$.\vspace{3mm}\\

Alternatively, $X_N$ is "integral", and $\mathcal{O}_{X_N}(U)$,
for $U\subseteq X_N$ open, are all $\FF$- sub-rings of the stalk
at the generic point $\mathcal{O}_{X_N,(0)}= \FF(\QQ)$, given by
\begin{equation}
\mathcal{O}_{X_N}(U)=
\begin{cases}
  \underset{p\in U}{\bigcap} \FF(\ZZ_{(p)})  & \eta \notin U, \\
  \underset{p\in U}{\bigcap} \FF(\ZZ_{(p)})\cap \mathcal{O}_{\QQ,\eta}  & \eta\in U
\end{cases}
\end{equation}
For $N$ dividing $M$ we have commutative diagrams
\begin{equation}
\xymatrix{
        A_N\ar@{^(->}[r]\;\ar@{_(->}[d]&A_M\;\ar@{_(->}[d]&\Spec A_N&\Spec A_M\ar@{_(->}[l]\\
       \FF\left(\ZZ\left[\frac{1}{N}\right]\right )\ar@{^(->}[r]&\FF\left(\ZZ\left[\frac{1}{M}\right]\right )&\Spec\FF\left(\ZZ\left[\frac{1}{N}\right]\right )\ar@{^(->}[u] &\Spec\FF\left(\ZZ\left[\frac{1}{M}\right]\right )\;\;\ar@{_(->}[l]\;\ar@{^(->}[u]}
\end{equation}
and we obtain a map $\pi_N^M:X_M\rightarrow X_N$ which is a
bijection on points and further
$(\pi_M^N)_*\mathcal{O}_{X_M}=\mathcal{O}_{X_N}$, i.e.
$(\pi_N^M)^{\#}$ is the identity. But note that $X_M$ has more
open neighborhoods of $\eta$ than $X_N$ such as $\Spec
\ZZ[\frac{1}{M}]\cup\{\eta\}=\Spec A_M\subseteq \Spec A_N.$
The "compactified $\Spec \ZZ$" is the inverse system $\{X_N\}$ with indices square free integers $N\geq 2$, maps $\pi^M_N$, and partial order given by divisibility. We denote it by $\bar{\Spec \ZZ}$. \\
Note that the $\FF$-locally-ringed-space $\mathcal{L}(\bar{\Spec\ZZ})=\varprojlim\limits_NX_N$
   has for points $\Spec\ZZ\cup\{\eta\}$, with open sets of the form
   $U$ or $U\cup\{\eta\}$ with $U$ an arbitrary open set of
   $\Spec\ZZ$ (hence $\bar{\Spec\ZZ}$ is of "Krull" dimension 1).
   Note that each $X_N$ is compact, and
   hence $\mathcal{L}(\bar{\Spec\ZZ})$ is compact.   Furthermore, the local
   $\fr$  $\mathcal{O}_{\bar{\Spec\zz},\eta}$ is just $\mathcal{O}_{\qq,\eta}$
   (while the local $\fr$  $\mathcal{O}_{X_N,\eta}$ is only $A_N$).\\
   For an open set $U=\Spec\ZZ\left[\frac{1}{N}\right]$ we have
\begin{equation}
\mathcal{O}_{\bar{\Spec\zz}}(U)=\FF\left(\ZZ\left[\frac{1}{N}\right]\right),
\end{equation}
   and
\begin{equation}
\mathcal{O}_{\bar{\Spec\zz}}(U\cup\{\eta\})=A_N.
\end{equation}
   The global sections $\mathcal{O}_{\bar{\Spec\zz}}(\bar{\Spec\ZZ})$ is the $\fr\  \FF\{\pm
   1\}$.

\section {The compactified $\Spec \mathcal{O}_K$.}
\smallbreak
Similarly, for a number field $K$, with ring of integers $\mathcal{O}_K$, with real primes $\{\eta_i\},i=1,\dots,r(=r_{\tiny \RR}+r_{\tiny \CC}),$ let $A_{N,i}=\FF(\mathcal{O}_K[\frac{1}{N}])\cap \mathcal{O}_{K,\eta_i}$, be the $\fr$ with
\begin{equation}
(A_{N,i})_{Y,X}=\bigg\{ a\in \FF(\mathcal{O}_K[\frac{1}{N}]):\;\;\;\; |a|_{\eta_i}\leq 1 \bigg\}
\end{equation}
the $Y\times X$ matrices with values in $\mathcal{O}_K\left[\frac{1}{N}\right]$ and with $\eta_i$-operator
$L_2$- norm bounded by 1.\\
Let $X_N$ be the Grothendieck-$\FF$-scheme obtain
by gluing $\{\Spec A_{N,i}\}_{i=1}^r$, and $\{\Spec\FF(\mathcal{O}_K)\}$, along the common open set
$\Spec\FF\left(\mathcal{O}_K\left[\frac{1}{N}\right]\right)$. For $N$
dividing $M$ we obtain a map $\pi^M_N:X_M\to X_N$, with $\pi^M_N|_{\Spec
A_{M,i}}$ induced by $A_{N,i}\subseteq A_{M,i}$. The inverse
system $\{X_N\}$ is the pro- $\FF$-scheme $\bar{\Spec\mathcal{O}_K}$, the
compactification of $\Spec\mathcal{O}_K$.\\
The space $\mathcal{L}(\bar{\Spec\mathcal{O}_K})=\varprojlim\limits_NX_N$
has for points $\Spec\mathcal{O}_K\cup\{\eta_i\}_{i\leq r}$, and open
sets are of the form $U\cup\{\eta_i\}_{i\in I}$ with $U$ open in
$\Spec\mathcal{O}_K$, and $I\subseteq\{1,...,r\}$ a subset
(and hence it is of "Krull" dimension $1$). The local
$\frs \; \mathcal{O}_{\bar{\Spec\mathcal{O}_K},\eta_i}$ is the ring $\mathcal{O}_{K,\eta_i}$. The global
sections $\mathcal{O}_{\bar{\Spec\mathcal{O}_K}}(\bar{\Spec\mathcal{O}_K})$
is the $\fr \; \FF\{\mu_K\}$, $\mu_K$ the group
of roots of unity in $\mathcal{O}_K^*$.\\

\newpage
\chapter{Vector bundles}
\bigskip

\section{Meromorphic functions $\mathcal{K}_{X_N}$.}
\smallbreak

Let $X=\{X_N,\pi_N^M\}_{M\geq N\in \mathcal{N}}$ be a pro- object in the category $\mathfrak{g}\FF\mathcal{S}c^{(t)}$,\\
let $U\subseteq X_N$ be an open set and define a multiplicative set,
$$S_{N}(U)=\bigg \{s \in \mathcal{O}^{+}_{X_N}(U)_{1,1}\; |\; \forall M\geq N, \forall V\subseteq (\pi_N^M)^{-1}(U)$$
\begin{equation}
\forall a,a'\in \mathcal{O}_{X_M}(V)_{Y,X}, (\pi_N^M)^{\#}(s)\cdot a=(\pi_N^M)^{\#}(s)\cdot a'\implies a=a' \bigg\}.
\end{equation}
Define the sheaf of $\frs$ with involution, the "meromorphic functions" $\mathcal{K}_{N}$, to be the sheaf associated to the pre- sheaf $U\mapsto S_N(U)^{-1}\cdot \mathcal{O}_{X_N}(U)_{Y,X}$
\begin{equation}
\mathcal{K}_N=S_N^{-1}\mathcal{O}_{X_N} \in \FR^t/X_N.
\end{equation}

For $M\geq N$ we have the following commutative diagram:
\begin{equation}
\begin{tikzcd}
\mathcal{O}_{X_N}\arrow{d}\arrow[hook]{r}&\mathcal{K}_N\arrow{d}&S_N\arrow{d}{\pi_N^M} \\
(\pi_N^M)_*\mathcal{O}_{X_M}\arrow[hook]{r}&(\pi_N^M)_*\mathcal{K}_M &S_M
\end{tikzcd}
\end{equation}

\section{Rank-d Vector Bundles at finite layer.}
\smallbreak

For any $d\in \FF$, we have an injection (of sheaves) of groups
$$\bigg( GL_d(\mathcal{O}_{X_N})\inj GL_d(\mathcal{K}_N)\bigg)\in \;\slfrac{Grps}{X_N}\en{6.2.1}$$
define,
$$\mathcal{D}_d(X_N)=\Gamma(X_N,\slfrac{GL_d(\mathcal{K}_N)}{GL_d(\mathcal{O}_{X_N})})\in Sets_*.\en{6.2.2}$$
Elements of $\mathcal{D}_d(X_N)$ are represented as $\mathcal{D}:=(u_{\alpha},f_{\alpha})$, where
$$X_N=\underset{\alpha}{\cup}u_{\alpha},\;\;\; f_{\alpha} \in GL_d(\mathcal{K}_N)(u_{\alpha}),\;\;\; f_{\alpha}^{-1}|_{u_{\alpha}\cap u_{\beta}}\circ f_{\beta}|_{u_{\alpha}\cap u_{\beta}}\in GL_d(\mathcal{O}_{X_N}),\en{6.2.3}$$
and two such elements $D:=(u_{\alpha},f_{\alpha}) ,D':= (v_{\beta},g_{\beta})$ represent the same element of $\mathcal{D}_d(X_N)$:
$$D=D'$$
$$\iff $$
$$ \exists X_N=\underset{\gamma}{\bigcup}w_{\gamma}\text{ a common refinement of } \{u_{\alpha}\} \text{ and } \{v_{\beta}\}, $$
$$\text{ and }\exists \; u_{\gamma}\in GL_d(\mathcal{O}_{X_N})(w_{\gamma}) \text{ such that } f_{\alpha}|_{w_{\gamma}}=g_{\beta}|_{w_{\gamma}} \circ u_{\gamma},  \text{ for } w_{\gamma}\subseteq u_{\alpha}\cap v_{\beta}. \en{6.2.4}$$


\noindent We obtain a category $\Vec(X_N)$ of "vector- bundles over $X_N$", with objects $\underset{\small d\in \FF}{\bigcup} \mathcal{D}_d(X_N)$, and with arrows from $\mathcal{D}=(u_{\alpha},f_{\alpha})\in \mathcal{D}_d(X_N)$ to $\mathcal{D}'=(v_{\beta},g_{\beta})\in \mathcal{D}_{d'}(X_N)$ given by

\begin{equation*}
\begin{split}
\Vec(X_N)(\mathcal{D},\mathcal{D}')= \{h\in \Gamma (X_N,(\mathcal{K})_{d',d}),\; g_{\beta}^{-1}\circ h|_{v_{\beta}\cap u_{\alpha}}\circ f_{\alpha} \\ \in \mathcal{O}_{X_N} (v_{\beta}\cap u_{\alpha})_{d',d}\subseteq \mathcal{K}_N(v_{\beta}\cap u_{\alpha})_{d',d} \}
\end{split}
\end{equation*}
$$\en{6.2.5}$$
\noindent For such $\mathcal{D}=(u_{\alpha},f_{\alpha}), \mathcal{D}'=(v_{\beta},g_{\beta})$, we have the well defined element $D\oplus D'\in \mathcal{D}_{d\oplus d'}(X_N)$,
$$D\oplus D':= (u_{\alpha}\cap v_{\beta}, \; f_{\alpha}|_{u_{\alpha}\cap v_{\beta}}\oplus g_{\beta}|_{u_{\alpha}\cap v_{\beta}}) \en{6.2.6}$$
We have associativity, commutativity, and unit isomorphisms, the unit object $(0)=(X_N,id_{[0]})$ (where: $\mathcal{D}_{[0]}(X_N)=\{(0)\}$) is the initial and final object of $\Vec(X_N)$. For $h_i\in \Vec(\mathcal{D}_i ,\mathcal{D}'_i)$, we have $h_0\oplus h_1\in \Vec (\mathcal{D}_0\oplus \mathcal{D}_1, \mathcal{D}_0'\oplus \mathcal{D}_1')$, thus $\Vec(X_N)$ is a symmetric monoidal category. \vspace{3mm} \\
For $d\geq 0$, we have the isomorphism classes of rank- $d$ vector bundles

$$Pic_d(X_N):= \Gamma(X_N,GL_d(\mathcal{K}_N)) \bigg \backslash \mathcal{D}_d(X_N) \en{6.2.7}$$
and we get from $\oplus$ an induced commutative monoid structure on $Pic_{\infty}(X_N)=\underset{d\geq 0}{\bigcup} Pic_d (X_N)$ . \\
For $\mathcal{D}_i,\mathcal{D}'_i\in Pic_{\infty}(X_N)$, we let $(\mathcal{D}_0, \mathcal{D}_0')\sim (\mathcal{D}_1, \mathcal{D}_1')$ iff $\mathcal{D}_0\oplus \mathcal{D}_1'\oplus \mathcal{D}\cong \mathcal{D}_1\oplus \mathcal{D}_0'\oplus \mathcal{D}$ for some $\mathcal{D}\in Pic_{\infty}(X_N)$. This defines an equivalence relation on pairs $(\mathcal{D},\mathcal{D}')$. We let $\mathcal{D}-\mathcal{D}'$ denote the equivalence class of $(\mathcal{D},\mathcal{D}')$, and we let
$K(X_N)= \slfrac{Pic_{\infty}(X_N)\times Pic_{\infty}(X_N)}{\sim}$, the Grothendieck group of stable isomorphism classes of (virtual) vector bundles. \\
We have a homomorphism of monoids $Pic_{\infty}(X_N)\to K(X_N)$, which is universal into an abelian group. \\
For $M\geq N$ we have pull-back of vector bundles
\begin{equation*}
\begin{split}
(\pi_N^M)^{\#}: \mathcal{D}_d(X_N)\longrightarrow\mathcal{D}_d(X_M) \hspace{41.5mm} \\
\mathcal{D}=(u_{\alpha}, f_{\alpha})\longmapsto (\pi_N^M)^{\#}\mathcal{D}=((\pi_N^M)^{-1}u_{\alpha}, (\pi_N^M)^{\#}(f_{\alpha}))
\end{split}
\end{equation*}
$$\en{6.2.8}$$

It is (the object part of) a strict- symmetric- monoidal functor $\Vec(X_N)\to \Vec(X_M)$,
and it induces a homomorphism of commutative monoids$\slash$abelian groups

$$\begin{tikzcd}
Pic_{\infty}(X_N)\arrow{r}\arrow{d} & Pic_{\infty}(X_M)\arrow{d} \\
K(X_N)\arrow{r} & K(X_M)
\end{tikzcd}\en{6.2.9}$$

\rem{6.2.1} We have a partial order on $\mathcal{D}_d(X_N)$
$$\mathcal{D}=(u_{\alpha},f_{\alpha})\leq \mathcal{D}'=(v_{\beta},g_{\beta})\iff id_{[d]}\in \Vec(X_N)(\mathcal{D},\mathcal{D}')$$
$$ \iff g_{\beta}^{-1}\circ f_{\alpha}\in \mathcal{O}_{X_N}(v_{\beta}\cap u_{\alpha})_{d',d}\;\; \text{for all } \beta, \alpha . \en{6.2.10}$$
The action of the group $\Gamma(X_N,GL_d(\mathcal{K}_N))$ on $\mathcal{D}_d(X_N)$ preserves this partial order.
The maps $(\pi_N^M)^{\#}$ of $(6.2.8)$ is order preserving, and is covariant with respect to $\Gamma(X_N,GL_d(\mathcal{K}_N))\to \Gamma(X_M,GL_d(\mathcal{K}_M))$ .

\section{$\mathcal{D}_d(X)$, Rank-d Vector Bundles in the limit.}
\smallbreak

Passing to the limit $\varprojlim_{N} X_N$, we have a symmetric monoidal category
$\Vec(X)= \varinjlim_{N} \Vec(X_N)$. \\
It has objects $\underset{d\in \FF}{\bigcup} \varinjlim_N \mathcal{D}_d(X_N)$, and for $\mathcal{D}\in \mathcal{D}_d(X_{N_0}), \mathcal{D}'\in \mathcal{D}_{d'}(X_{N_0'})$,
$$\Vec(X)(\mathcal{D},\mathcal{D}')=\varinjlim_{N\geq N_0,N_0'} \Vec(X_N)((\pi_{N_0}^{N})^{\#}\mathcal{D},(\pi_{N_0'}^N)^{\#}\mathcal{D}') \en{6.3.1} $$
The isomorphism classes of rank-$d$ vector-bundles are given by
$$Pic_d(X)=\varinjlim_N Pic_d(X_N)=\varinjlim_N \Gamma(X_N,GL_d(\mathcal{K}_N))\bigg \backslash \varinjlim_N \mathcal{D}_d(X_N)$$
The direct sum $\oplus$ induces a commutative monoid structure on $Pic_{\infty}(X)=\underset{d\geq 0}{\bigcup}Pic_d(X)$ , and passing to stable isomorphism classes of virtual vector- bundles we obtain the Grothendieck group
$$K(X)= \varinjlim_N K(X_N)\in Ab. \en{6.3.2}$$
We describe next another passage to the limit (which is kind of dual to the one above, as we consider "$\varprojlim_N\mathcal{D}_d(X_N)$") . Set


$$  \mathcal{B}_d(X)=
\begin{Bmatrix}
  D=\{D_N\},\;\mathcal{D}_N \in \mathcal{D}_d(X_N),\hspace{3mm}     (*) \hsm \forall M\geq N, \hsm (\pi_N^M)^{\#}\mathcal{D}_N\geq \mathcal{D}_M, \\
    \hspace{10mm} (**) \hsm \exists \delta_0\in \mathcal{D}_d(X_{N_0}), \hsm \forall N\geq N_0, \hsm \mathcal{D}_N\geq (\pi_{N_0}^N)^{\#}\delta_0
\end{Bmatrix}\en{6.3.3}$$
i.e. $B_d(X)$ is the collection of bounded below $(**)$, monotone decreasing $(*)$, filters of vectors bundles. \\
Define a transitive reflexive relation on $\mathcal{B}_d(X)$,
$$\{\mathcal{D}_N\}\geq \{\mathcal{D'}_N\} $$
$$\iff$$
$$\bigg\{ \forall\delta_0\in \mathcal{D}_d(X_{N_0})\text{ such that } \hsm \forall N\geq N_0\;\;\; \mathcal{D}'_N\geq (\pi_{N_0}^N)^{\#}\delta_0 \en{6.3.4}$$
$$\implies \forall N\geq N_0' (\geq N_0) \hsm D_N\geq (\pi_{N_0}^N)^{\#}\delta_0\bigg \} \en{6.3.5}$$
and an equivalence relation,
$$\{D_N\} \sim \{D_N'\} \en{6.3.6}$$
$$\text{if }\{\mathcal{D}_N\}\leq \{\mathcal{D}_N'\} \text{ and } \{\mathcal{D}_N\}\geq \{\mathcal{D}_N'\} \en{6.3.7}$$
Finally, we define $\widehat{\mathcal{D}}_d(X)$ to be the $\sim$ equivalence classes:
$$\widehat{\mathcal{D}}_d(X):= \slfrac{\mathcal{B}_d(X)}{\sim}. \en{6.3.8}$$
It is a partially ordered set, with an (order preserving) action of the group $GL_d(\mathcal{K}(X)):= \varinjlim_{N} GL_d(\mathcal{K}_N)(X_N).$

\exmpl{6.3.1}
We have
$$\widehat{\mathcal{D}}_d(\bar{\Spec \ZZ})\cong \slfrac{GL_d(\mathbb{A}_{\QQ})}{O_d\times \underset{\sigma}{\Pi} GL_d(\ZZ_p)},\en{6.3.9}$$
where $\mathbb{A}_{\QQ}$ is the ring of Adeles of $\QQ$, and $O_d=GL_d(\mathcal{O}_{\RR,\eta})$ the orthogonal group.\vspace{1.5mm} \\

For a finite $p$, the symmetric space $\Xi_p=\slfrac{GL_d(\QQ_p)}{GL_d(\ZZ_p)}$ can be identified with the $GL_d(\QQ_p)$- set
$$\underline{\underline{\mathcal{L}}}_p=\{\mathcal{L}\subseteq \QQ_p^{\oplus d},\; \mathcal{L} \text{ is a }\ZZ_p\text{ -lattice} \}\en{6.3.10}$$
via,
$$g\cdot GL_d(\ZZ_d)\mapsto g(\ZZ_p^{\oplus d}),\en{6.3.11}$$
since $GL_d(\QQ_p)$ acts transitively on $\underline{\underline{\mathcal{L}}}_p$, and the stabilizer of $\ZZ^{\oplus d}_p$ is precisely $GL_d(\ZZ_p)$.\\
Similarly, for the real prime $\eta$, the symmetric space $\;\Xi_{\eta}= \slfrac{GL_d(\RR)}{O_d}$, can be identified with the $GL_d(\RR)$ -set $\underline{\underline{\mathcal{L}}}_{\eta}=\{Q\in Mat_{d,d}(\RR), \text{ symmetric and positive definite}\}$, of positive definite quadratic forms $Q$, or alternatively by the associated ellipsoids $\mathcal{L}_{Q}:= \{x\in \RR^{\oplus d},\hsm x\circ Q\circ x^t\leq 1\}$, so that the correspondence is given (more like the case of finite p's) as
$$\Xi_{\eta}=\slfrac{GL_d(\RR)}{O_d}\iso \underline{\underline{\mathcal{L}}}_{\eta}=\{\mathcal{L}\subseteq \RR^{\oplus d} \text{ ellipsoid}\}\en{6.3.12} $$
$$g\cdot O_d\mapsto g(\ZZ_{\eta}^{\oplus d}).\en{6.3.13}$$
with $\ZZ^{\oplus d}_{\eta}=\{x\in \RR^{\oplus d},\;\;\;\; x\circ x^t\leq 1\}$ the $d$- dimensional unit ball. \vspace{3mm}

The Adelic symmetric space $\Xi_{\mathbb{A}}={GL_d(\mathbb{A})}\bigg /{O_d\times \Pi_p GL_d(\ZZ_p)}$ can be identified with the restricted- product $\underline{\underline{\mathcal{L}}}_{\mathbb{A}}=\underline{\underline{\mathcal{L}}}_{\eta} \times \prod_p \underline{\underline{\mathcal{L}}}_p$, it is the subset of the product consisting of $l=(l_1,l_2,l_3,l_5,\dots )$ such that $l_p=\ZZ_p^{\oplus d}$ for all but finitely many p's. Each such $l$ is a compact subset of $\mathbb{A}^{\oplus d}$.\vspace{3mm}\\

In the space $X_N=\Spec \; \FF(\ZZ)\underset{\Spec \FF(\ZZ[\frac{1}{N}])}{\coprod} A_N,\; N=p_1\dots p_l $, we have the system
$\{V_{\eta}^M,V_{p_1}^M,\dots,V_{p_l}^M\}$ of cofinal coverings for $M=q_1\dots q_k$ prime to $N$, with
$$V_{p_j}^M=\Spec \FF(\ZZ[\frac{p_j}{N\cdot M}])\;,\; j=1,\dots, l \en{6.3.14}$$
$$V_{\eta}^M=\Spec A_N\hspace{5mm} \text{(independent of  } M\text{!)}.\en{6.3.15}$$
(for every open cover $\{u_{\alpha}\}$ of $X_N$, there exists $M$ such that $\{V_{\eta}^M,V_{p_i}^M\}$ is a refinement of $\{u_{\alpha}\}$). \\
The sheaf $\mathcal{K}_N$ of meromorphic functions is the constant sheaf $\FF(\QQ)$, for all $N$. Thus we can represent an element of $\mathcal{D}_d(X_N)$ by a sequence $(g_{\eta},g_{p_1},\dots, g_{p_l})\in GL_d(\QQ)^{(l+1)}$, with $g_{p_i}^{-1}\circ g_{\eta}$, and $g_{p_i}^{-1}\circ g_{p_j}$ in $\underset{p \nmid N\cdot M}{\prod} GL_d(\ZZ_p)$ , (a vacuous condition since $M$ is arbitrary).

Two such sequences $(g_{\eta},g_{p_1},\dots, g_{p_l})$, and $(h_{\eta},h_{p_1},\dots, h_{p_l})$, represent the same element of $\mathcal{D}_d(X_N)$, if and only if $h_{p_i}^{-1}\circ g_{p_i}\in GL_d(\ZZ_{p_i})$, and $h_{\eta}^{-1}\circ g_{\eta}\in O_d\cap GL_d(\ZZ[\frac{1}{N}])$.

Thus we have the well defined map

$$div_N:\mathcal{D}_d(X_N)\rightarrow \underline{\underline{\mathcal{L}}}_{\mathbb{A}}\en{6.3.16}$$
$$div_N(g_{\eta},g_{p_1},\dots,g_{p_l})=\{\mathcal{L}_p\}\;,\; \mathcal{L}_p= \begin{cases}
    g_{p_i}(\ZZ_{p_i}^{\oplus d}) & p=p_i\mid N \\
   g_{\eta}(\ZZ_p^{\oplus d}) & p\nmid N \text{ or } p=\eta.
  \end{cases}\en{6.3.17}$$

For $N=p_1\cdots p_l$, and for $N'=p_1'\cdots p_{l'}'$ prime to $N$, the map $\pi_N^{N\cdot N'}: X_{N\cdot N'}\rightarrow X_N$, induce pullback $(\pi_N^{N\cdot N'})^{\#}: \mathcal{D}_d(X_N)\rightarrow \mathcal{D}_d(X_{N\cdot N'})$ , and
$$(\pi_N^{N\cdot N'})^{\#}(g_{\eta},g_{p_1},\dots, g_{p_l})\slash\sim=\slfrac{(g_{\eta},g_{p_1},\dots, g_{p_l},g_{\eta},\dots ,g_{\eta}) }{\sim},\en{6.3.18}$$

i.e. $g_{\eta}$ placed in the $p_i'$- spots. Thus we have a commutative diagram:

$$\xymatrix{
\mathcal{D}_d(X_{N\cdot N'})\ar@{^{(}->}[rrd]^{div_{N\cdot N'}}&& \\
&&\underline{\underline{\mathcal{L}}}_{\mathbb{A}} \\
\mathcal{D}_d(X_N)\ar@{^{(}->}[rru]_{div_N}\ar@{^{(}->}[uu]&&}
\en{6.3.19}$$

Note that $\underset{N}{\bigcup} div_N (\mathcal{D}_d(X_N))$ is the dense subset of $\underline{\underline{\mathcal{L}}}_{\mathbb{A}}$ consisting of $\{\mathcal{L}_p\}$ with arbitrary $\mathcal{L}_p$ at finite primes $p$, and with $\mathcal{L}_{\eta}$ defined over $\QQ$ (i.e. $\mathcal{L}_{\eta}\in GL_d(\QQ)\cdot (\ZZ_{\eta}^{\oplus d})\subseteq \underline{\underline{\mathcal{L}}}_{\eta}$, while $GL_d(\QQ)\cdot (\ZZ_p^{\oplus\;d})\equiv \underline{\underline{\mathcal{L}}}_p$). \\
Given $\mathcal{D}=\{\mathcal{D}_N\}\in \mathcal{B}_d(X)$, $X=\bar{\Spec\;\ZZ}=\{X_N\}$, monotone decreasing bounded below sequence, we let $\hat{div}{\mathcal{D}}:=\underset{N}{\bigcap} div_N(\mathcal{D}_N)$, the intersection taken in $\mathbb{A}^{\oplus d}$.
Thus $\widehat{div}(\mathcal{D})_p= div_N(\mathcal{D}_N)_p$, for all finite $p$, and all $N$ divisible by (some fixed) $N_0$, and $\widehat{div}(D)_{\eta}=\underset{N}{\bigcap} g_{\eta}^N(\ZZ_{\eta}^{\oplus d})$ (where $g_{\eta}^N$ is the $\eta$- component of $\mathcal{D}_N$) is an (arbitrary real) ellipsoid. We have $\mathcal{D}=\{D_N\}\geq D'=\{D_N'\}$, if and only if, $\widehat{div}(\mathcal{D})\supseteq \widehat{div}(D')$, and so $\widehat{div}(\;)$ induces a bijection
$$\widehat{\mathcal{D}}_d(X)=\slfrac{\mathcal{B}_d(X)}{\sim}\underset{\widehat{div}}{\iso} \mathcal{L}_{\mathbb{A}} \underset{\sim}{\leftarrow} X_{\mathbb{A}}. \en{6.3.20}$$

\vspace{3mm}

In exactly the same way, we obtain for any number field $K$, a $GL_d(K)$- covarient identification
$$\widehat{\mathcal{D}}_d(\bar{\Spec \mathcal{O}_K})\cong \slfrac{ GL_d(\mathbb{A}_K)}{\prod_{\nu} GL_d(\mathcal{O}_{K,\nu})},\en{6.3.21}$$
where $\mathbb{A}_{K}$ is the ring of Adeles of $K$, $\mathcal{O}_{K,\nu}$ the local ring at $\nu$ for finite $\nu$'s, and $GL_d(\mathcal{O}_{K,\nu})\cong O_d$ (resp. $U_d$) the orthogonal (resp. unitary) group for $\nu$ real (resp. complex). \vspace{3mm} \\
For $d=1$, $\slfrac{GL_1(\mathcal{K}_N)}{GL_1(\mathcal{O}_{X_N})}\cong \slfrac{\mathcal{K}_N^*}{\mathcal{O}_{X_N}^*}$ is a sheaf of abelian groups, and we have the ordered abelian group $\widehat{\mathcal{D}}_1(X)$, and its quotient by $GL_1(\mathcal{K}(X))=\underset{N}{\varinjlim} \; GL_1(\mathcal{K}_N (X_N)),$ the completed- Picard group
$$\widehat{Pic}_1(X):= GL_1(\mathcal{K}(X))\bigg\backslash \widehat{\mathcal{D}}_1(X). \en{6.3.22}$$
For a number field $K$,
$$\widehat{Pic}_1(\bar{\Spec \mathcal{O}_K})= K^*\Bigg \backslash \slfrac{\mathbb{A}_K^*}{\prod_{\nu}\mathcal{O}_{K,\nu}^*}.\en{6.3.23}$$


\numberwithin{equation}{section}
\chapter{Modules}
\bigskip
\section{Definitions}
\smallbreak

\defin{7.1.1}\textit{Let $A\in \FR$. We denote by $\Amod$ the full subcategory of the functor category $(Ab)^{A\times A^{op}}$ given by
\begin{equation}
\Amod\equiv \bigg\{ M=\{ M_{Y,X} \}\in (Ab)^{A\times A^{op}}\;,\; M_{0,X}=\{0\}=M_{Y,0}\bigg\}.
\end{equation}   }
Thus an $A$- module $M$ is a collection of abelian groups $M_{Y,X}$, for $X,Y\in \FF$, together with maps:
\begin{equation}
A_{Y',Y}\times M_{Y,X}\times A_{X,X'}\rightarrow M_{Y',X'}
\end{equation}
\begin{equation}
a,m,b\mapsto a\circ m\circ b
\end{equation}
such that,
\begin{equation}
a\circ (m+m')\circ b=a\circ m\circ b+a\circ m'\circ b,\hsm \text{(homomorphism)}
\end{equation}
\begin{equation}
(a_1\circ a_2)\circ m \circ (b_2\circ b_1)=a_1\circ (a_2\circ m\circ b_2)\circ b_1, \hsm \text{(associativity)}
\end{equation}
\begin{equation}
id_Y\circ m\circ id_X=m.\hsm \text{(identity)}
\end{equation}
\underline{Notation:} for $m\in M_{Y,X}$ we will write
\begin{equation}
a\circ m:=a\circ m\circ id_X,
\end{equation}
\begin{equation}
m\circ b:=id_Y\circ m\circ b.
\end{equation}

\prop{7.1.1}\label{7.1.1}\textit{$\Amod$ is complete and co-complete abelian category. It has enough projectives and injectives.}
\begin{proof} (Well known). All (co)limits can be taken pointwise $(\lim M)_{Y,X}=\lim (M_{Y,X})$. The evaluation functor $i^{Y,X}:\Amod\rightarrow Ab$, $i^{Y,X}M=M_{Y,X}$, has a left (resp. right) adjoint $i_!^{Y,X}$ (resp. $i_*^{Y,X}$). Taking a surjection $p^{Y,X}\sur M_{Y,X}$ (resp. injection $M_{Y,X}\inj p^{Y,X}$) with $p^{Y,X}$ projective (resp. injective) in $Ab$, we obtain a surjection (resp. injection) $\underset{Y,X}{\bigoplus} i_!^{Y,X}p^{Y,X}\sur M$ (resp. $M\inj \underset{Y,X}{\prod} i_*^{Y,X}p^{Y,X}$) from a projective (resp. into an injective) $A$- module.
\end{proof}

We have the following injective and surjective maps:
\begin{equation}
\begin{tikzcd}
M_{Y_0,X_0}\oplus M_{Y_0,X_1}\oplus M_{Y_1,X_0}\oplus M_{Y_1,X_1}\arrow[bend right,hook]{r}{f}&M_{Y_0\oplus Y_1,X_0\oplus X_1}\arrow[bend right,two heads]{l}[swap]{g}
\end{tikzcd}
\end{equation}
given by,
$$\begin{pmatrix}
  m_{00} & m_{01} \\
  m_{10} & m_{11}
\end{pmatrix}\longmapsto \sum_{i,j=0,1}1_{Y_i}\circ m_{i,j}\circ 1_{X_j}^t,$$
\begin{equation}
(1_{Y_i}^t\circ m\circ 1_{X_j})\longmapsfrom m.
\end{equation}
where,
\begin{equation}
1_{Y_i}\in \FF_{Y_0\oplus Y_1,Y_i},1_{X_j}^t\in \FF_{X_j,X_0\oplus X_1}.
\end{equation}
$$ (1_{X_j}^t\circ 1_{X_j}=id_{X_j}, 1_{X_j}\circ 1_{X_j}^t=id_{X_j}\oplus 0_{X_{1-j}}\hsm,j=0,1).$$
Note that the composition $g\circ f=id$, therefore it is a direct summand, but in general $f\circ g\neq id$, in fact,
\begin{equation}
f\circ g=id \iff m=\sum_{i,j=0,1}1_{Y_i}\circ 1_{Y_i}^t\circ m\circ 1_{X_j}\circ 1_{X_j}^t,
\end{equation}
$$\forall m\in M_{Y_0\oplus Y_1,X_0\oplus X_1}.$$
There is a similar correspondence :
\begin{equation}
\begin{tikzcd}
(M_{1,1})^{X\otimes Y}\arrow[bend right,hook]{r}{f}&M_{Y,X}\arrow[bend right,two heads]{l}[swap]{g}
\end{tikzcd}
\end{equation}
where $(M_{1,1})^{X\otimes Y}$ is a direct summand,
$$\bigg(1_y^t\circ m\circ 1_x\bigg)_{(x,y)\in Y\otimes X}\longmapsfrom m$$
\begin{equation}
m_{y,x}\longmapsto \sum_{x\in X,y\in Y} 1_y\circ m_{y,x}\circ 1_x^t.
\end{equation}
If these maps are isomorphisms for all $X,Y\in \FF$, we say $M$ is a "matrix $A$- module". \\
In particular, the map $f$ of $(7.1.9)$, made up from the $\FF$- action of $1_{Y_i}$ and $1_{X_j}^t$, gives the "direct- sum" for $M$:
\begin{equation}
M_{Y_0,X_0}\times M_{Y_1,X_1}\rightarrow M_{Y_0\oplus Y_1,X_0\oplus X_1}
\end{equation}
\begin{equation}
(m_0,m_1)\mapsto m_0\oplus m_1:=1_{Y_0}\circ m_0\circ 1^t_{X_0}+1_{Y_i}\circ m_1\circ 1_{X_1}^t.
\end{equation}
Note: associativity, commutativity, unit -isomorphisms are the
canonical $\FF$- isomorphisms. The map $(m_0,m_1)\mapsto m_0\oplus
m_1$ is strongly natural in the sense that
\begin{equation}
a\circ (m_0\oplus m_1)\circ a'=(a\circ 1_{Y_0})\circ m_0\circ (1_{X_0}^t\circ a')+(a\circ 1_{Y_1})\circ m_1\circ (1_{X_1}^t\circ a'),
\end{equation}
In particular it is natural,
\begin{equation}
(a_0\oplus a_1)\circ (m_0\oplus m_1)\circ (a_0'\oplus a_1')=(a_0\circ m_0\circ a_0')\oplus (a_1\circ m_1\circ a_1')
\end{equation}

\exmpl{7.1.1} Define $M^{k,l}\in \FF\text{-}mod$ by
\begin{equation}
\bigg(M^{k,l}\bigg)_{Y,X}:=\text{free abelian group on } (Y_0,X_0),
\end{equation}
\begin{equation}
\begin{array}{lcl} Y_0\subseteq Y&,& \#Y_0:=k \\ X_0\subseteq X&,& \# X_0=l \end{array},
\end{equation}
and the $\FF$- action of $a\in \FF_{Y',Y}, b\in \FF_{X,X'}$, is defined on the generators $(Y_0,X_0)$ of $(M^{k,l})_{Y,X}$ by:
\begin{equation}
a\circ (Y_0,X_0)\circ b=
\begin{cases}
 (a(Y_0),b^t(X_0)) & Y_0\subseteq D(a)\;,\;X_0\subseteq I(b)\\
 0 & \text{otherwise}
\end{cases}.
\end{equation}
The rank at every degree $(Y,X)$ is easy to calculate and is given by,
\begin{equation}
Rank_{\ZZ}\bigg(M_{Y,X}^{k,l}\bigg)=\dbinom{\# Y}{k}\cdot \dbinom{\# X}{l}.
\end{equation}
where in case $k>1$ or $l>1$ and $X=Y=1$ we have: $\bigg( M^{k,l}\bigg)_{1,1}=\{0\}.$ \\
It is an example of a non- matrix $\FF$- module.

\exmpl{7.1.2}  Define $M^{k,l}\in \FF(\ZZ)\text{-}mod$ by,
\begin{equation}
\bigg( M^{k,l} \bigg)_{Y,X}=\Hom_{\ZZ} \bigg(\bigwedge_{\ZZ}^l(\ZZ \cdot X),\bigwedge^k_{\ZZ}(\ZZ\cdot Y)\bigg)
\end{equation}
and the $\FF(\ZZ)$ -action of
\begin{equation}
a\in \FF(\ZZ)_{Y',Y}=\Hom_{\ZZ}(\ZZ\cdot Y,\ZZ\cdot Y'), \;\;\; b\in \FF(\ZZ)_{X,X'}=\Hom_{\ZZ}(\ZZ\cdot X',\ZZ\cdot X ),
\end{equation}
is given by,
\begin{equation}
a\circ m\circ b:= \bigwedge^k(a)\circ m\circ \bigwedge^l(b).
\end{equation}
Note that it is non- matrix. \vspace{3mm}\\
{\bf \large Tensor product:} \vspace{3mm} \\
For $M,N\in \Amod$, we have their tensor product $M\underset{A}{\otimes} N\in \Amod$,
$$(M\underset{A}{\otimes} N)_{Y,X} := \slfrac{\underset{Z}{\bigoplus} M_{Y,Z}\underset{\ZZ}{\otimes} N_{Z,X}}{\{(m\circ a)\otimes n\sim m\otimes (a\circ n), m\in M_{Y,Z'},a\in A_{Z',Z}, n\in N_{Z,X}\}}, $$
and we have their left (resp. right) inner hom $\Hom_A^{l/r}(M,N)\in \Amod$
$$\Hom_A^l(M,N)_{Y,X}:= \{\varphi =\{\varphi_Z\},\; \varphi_Z\in Ab(M_{Z,X},N_{Z,Y}),\; \varphi(a\circ m)=a\circ \varphi (m)\} $$
$$\text{with $A$- action: \hspace{20mm} }(a_1\circ \varphi\circ a_2)(m):= \varphi (m\circ a_1)\circ a_2. $$
resp.,
$$\Hom_A^r(M,N)_{Y,X}:= \{\varphi =\{\varphi_Z\}, \; \varphi_Z\in Ab(M_{X,Z},N_{Y,Z}), \; \varphi(m\circ a)=\varphi (m)\circ a \}, $$
$$\text{with $A$- action:  \hspace{20mm} }(a_1\circ \varphi\circ a_2)(m):= a_1\circ \varphi(a_2\circ m).$$
We have the adjunction:
$$\Amod(M\underset{A}{\otimes} N, K) =\Amod(M,\Hom^r_A(N,K))=\Amod(N,\Hom^l_A(M,K)).$$
Note that $\otimes_A$ is associative, but Not commutative, and has unit $A$ only if $A$ is an $A$- module, i.e. $A\in (\FR)^{Ab}$.
\defin{7.1.2. $A$- modules with involution: $\Amod^t$} \textit{For any $A\in \FRt$ define
$\Amod^t$ to be the category of $M=\{M_{Y,X}\}\in \Amod$ with an involution,
\begin{equation}
(\;)^t:M_{Y,X}\rightarrow M_{X,Y}
\end{equation}
\begin{equation}
m\mapsto m^t
\end{equation}
satisfying:
\begin{equation}
(m_0+m_1)^t=m_0^t+m_1^t
\end{equation}
$$(m^t)^t=m$$
\begin{equation}
(a\circ m\circ b)^t=b^t\circ m^t\circ a^t.
\end{equation}
The morphisms in $\Amod^t$ between two such modules $M,N$, are given by the set:
\begin{equation}
\Amod^t(M,N)\equiv \{\varphi\in \Amod (M,N)\;\;\; \varphi(m^t)=\varphi(m)^t\}.
\end{equation} }
\prop{7.1.2} \textit{The category $\Amod^t$ is complete and co-complete abelian category with enough injectives and projectives objects.}
\begin{proof}
Same as proposition 7.1.1.\vspace{3mm}\\
\end{proof}
\noindent \textbf{\large Free $A$- modules}\vspace{3mm}\\
For $X,Y\in \FF$ and an abelian group $N\in Ab$, we have the adjunction:
\begin{equation}
\Amod\bigg(i_!^{Y,X}N,M\bigg)\equiv Ab\bigg(N,M_{Y,X}\bigg)
\end{equation}
where the $i_!$ functor defined by,
\begin{equation}
(i_!^{Y,X}N)_{W,Z}=\underset{A_{W,Y}\times A_{X,Z}}{\bigoplus} N.
\end{equation}
An element of  $(i_!^{Y,X}N)_{W,Z}$ has the form $\sum_{i=1}^{k} n_i\cdot (a_i,b_i)$ with  $n_i\in N,a_i\in A_{W,Y}, b_i\in A_{X,Z}$.
\begin{equation}
(\text{where } \;0\cdot (a,b)\equiv n\cdot (0,b)\equiv n\cdot (a,0)\equiv 0,)
\end{equation}
The $A$- action on such an element of $a\in A_{W',W},b\in A_{Z,Z'}$ is given by:
\begin{equation}
a\circ (\Sigma n_i(a_i,b_i))\circ b=\Sigma n_i (a\circ a_i,b_i\circ b).
\end{equation}
In particular,  we have the free $A$- module of degree $(Y,X)$:
\begin{equation}
A^{Y,X}:=i_!^{Y,X}\ZZ,
\end{equation}
\begin{equation}
\Amod (A^{Y,X},M)\equiv M_{Y,X}.
\end{equation}
For $A\in \FRt \text{ (with involution!) } , X,Y\in \FF, N\in Ab,M\in \Amod^t$, we have the adjunction:
\begin{equation}
\Amod^t(i_!^{Y,X},M)\equiv Ab(N,M_{Y,X})
\end{equation}
\begin{equation}
\text{where now:}\;(i_!^{Y,X}N)_{W,Z}=\bigg( \underset{A_{W,Y}\times A_{X,Z}}{\bigoplus} N \bigg) \oplus \bigg( \underset{A_{W,X}\times A_{Y,Z}}{\bigoplus} N\bigg)
\end{equation}
has $A$- action as above and has involution, interchanging the summands above:
\begin{equation}
(\Sigma n_i(a_i,b_i) )^t=\Sigma n_i(b_i^t,a_i^t).
\end{equation}
The free $\A$- module with involution of degree $(Y,X)$, $A^{Y,X}_t$, is obtained by taking $N=\ZZ$:
\begin{equation}
A^{Y,X}_t:=i_!^{Y,X}\ZZ \;\; \text{with involution},
\end{equation}
\begin{equation}
\Amod^t(A_t^{Y,X},M)\equiv M_{Y,X}\end{equation}   \vspace{1.5mm} \\
Similarly, for $I\in \slfrac{Set}{\FF\times \FF}$, i.e. $I$ is a set together with a map $I\rightarrow \FF\times \FF,i\mapsto (Y_i,X_i)$, we have the free $A$- module, $A^I=\underset{i\in I}{A^{Y_i,X_i}}$, giving the adjunction:
\begin{equation}
\xymatrix{\Amod \ar@/^/@{->}[d]^{U} \\ \slfrac{Set}{\FF\times \FF}\ar@/^/@{->}[u]^{A^I}}\hspace{5mm} \Amod(A^I,M)\equiv \slfrac{Set}{\FF\times \FF}(I,UM).
\end{equation}
If we let $(\slfrac{Set}{\FF\times \FF})^t$ denote the category of sets over $\FF\times \FF,I\rightarrow \FF\times \FF$, together with involution i.e. a bijection $I\rightarrow I, i\mapsto i^t,i^{tt}=i, (I^{-1}(Y,X))^t=I^{-1}(X,Y)$, than we have the adjunction:
\begin{equation}
\xymatrix{\Amod^t \ar@/^/@{->}[d]^{U} \\ (\slfrac{Set}{\FF\times \FF})^t\ar@/^/@{->}[u]^{A^I_t}}\hspace{5mm} \Amod^t(A^I_t,M)\equiv (\slfrac{Set}{\FF\times \FF})^t(I,UM).
\end{equation}
\section{Commutativity for Modules}\label{7.2}
\smallbreak
Let $M\in \Amod^{(t)}$.
\defin{7.2.1}\textit{We say $M$ is \underline{Total commutative}:
$$\forall m\in M_{Y,X},\;\;\; b\in A_{I,J},$$
\begin{equation}
(\underset{Y}{\bigoplus} b)\circ (\underset{J}{\bigoplus}m)=(\underset{I}{\bigoplus}m)\circ (\underset{X}{\bigoplus}b),
\end{equation}
 \underline{Commutative}:
$$ \forall m\in M_{Y,X},\;\;\; b\in A_{1,J},\;\;\; d\in A_{J,1},$$
\begin{equation}
(\underset{Y}{\bigoplus}b)\circ (\underset{J}{\bigoplus}m)\circ (\underset{X}{\bigoplus}d)=(\underset{Y}{\oplus}b\circ d)\circ m= m\circ  (\underset{X}{\oplus} b\circ d),
\end{equation}
\underline{Central}:
$$ \forall m\in M_{Y,X},\;\;\; b\in A_{1,1}$$
\begin{equation}
(\underset{Y}{\bigoplus} b)\circ m=m\circ (\underset{X}{\bigoplus}b):=b\cdot m.
\end{equation}
i.e. the monoid $A_{1,1}$ acts centrally on $M_{Y,X}$, and we denote this actions by $b\cdot m$.} \vspace{0mm}\\
Note that,
$$ \text{total commutativity}\implies \text{commutativity} \implies \text{centrality}. $$
We let $CA$-$mod$ (resp. $A$-$mod_{tot-com}$, $A$-$mod_{cent'l}$) denote the full subcategory of $A$-$mod$ consisting of commutative (resp. totally commutative, central) $A$- modules. \\
We have the following adjunctions between the categories,
$$\xymatrix{\Amod_{tot-com}\ar@{^{(}->}@<-2pt>[r] &C\Amod\ar@{->>}@<-2pt>[l] \ar@{^{(}->}@<-2pt>[r] &\Amod_{centr'l}\ar@{->>}@<-2pt>[l] \ar@{^{(}->}@<-2pt>[r] &\Amod\ar@{->>}@<-2pt>[l] }$$
the left adjoint given by the quotient maps: \vspace{3mm} \\
\begin{equation}
\xymatrix{M_{tot-com}&CM\ar@{->>}[l]&M_{cent'l}\ar@{->>}[l]&M\ar@{->>}[l] }
\end{equation}
\vspace{3mm} \\
For $M\in \Amod_{cent'l}^{(t)}$, and $S\subseteq A_{1,1}^{(+)}$ multiplicative, we have the localization
\begin{equation}
S^{-1}M\in S^{-1}A\text{-mod}^{(t)}.
\end{equation}
We have
\begin{equation}\label{localization26}
(S^{-1} M)_{Y,X}:= (M_{Y,X}\times S)/\sim ,
\end{equation}
$$(m_0,s_0)\sim (m_1,s_1) \iff \exists s\in S \text{ such that } (s\cdot s_1)\cdot m_0=(s\cdot s_0)\cdot m_1,$$
and denoting by $m/s$ the equivalence class of $(m,s)$, we have
\begin{equation}\label{localization27}
(a'/s')\circ (m/s)\circ (a''/s''):=(a'\circ m\circ a'')/(s'\cdot s\cdot s'')
\end{equation}
In particular we have localizations
$$M_f:=S^{-1}_f M, \; S_f=\{1,f,f^2,\dots,f^n,\dots \}, \; f\in A^{(+)}_{1,1}$$
$$M_{\frak{p}}:=S_{\frak{p}}^{-1} M, \; S_{\frak{p}}=A_{1,1}^{(+)}\setminus \frak{p}, \; \frak{p}\in \Spec^{(+)} A. $$
\rem{0}{If $A\in \FR$ is such that there exists $b\in A_{1,J},d\in A_{J,1}, b\circ d=1, b_j=b\circ 1_j=0,$ (or $d_j=1_j^t\circ d=0$) then $CA$-mod$\equiv \{0\}$: If $M\in \Amod$ is commutative, $m\in M_{Y,X},$
\begin{equation}
m=(b\circ d)\cdot m=(\underset{Y}{\bigoplus}b)\circ (\sum_{j\in J} 1_{Y_j}\circ m\circ 1_{X_j}^t)\circ (\underset{X}{\bigoplus} d)=
\sum_{j\in J}(\underset{Y}{\bigoplus} b_j)\circ m \circ (\underset{X}{\bigoplus}d_j)=\sum_j (b_j\circ d_j)\cdot m=0.
\end{equation}
E.g. when $A=\FF_{\eta}$.}
\rem{1}{If $M\in \FF$ -mod, then $M$ is automatically commutative.}
\rem{2}{If $M\in \FF\{S\}$ -mod, $M$ commutative $\iff$ $M$ central: $(\underset{Y}{\oplus}s)\circ m=m\circ (\underset{X}{\oplus}s)$, $\forall\; s\in S,m\in M_{Y,X}$. For $S=\{\pm 1\}$, $(\underset{Y}{\oplus} (-1))\circ m=m\circ (\underset{X}{\oplus}(-1))$; if this is $-m$ (=the inverse of $m$ in the abelian group $M_{Y,X}$), we shall say $M$ is "$(-1)$- true". For $S=\NN$, $\underset{Y}{\oplus}(p)\circ m=m\circ \underset{X}{\oplus}(p)$ for all (prime) $p\in \NN$; if this is equal to $p\cdot m$ ($=(m+m+\dots+m)$, $p$ times), we shall say $M$ is " $\NN$- true". If it is both $(-1)$- true and $\NN$- true we shall say it is "$\ZZ$- true".  }
\rem{3}{For $R\in\;CRig$, $M\in \FF(R)$ -mod, $m_0,m_1\in M_{Y,X}$, the element $m_0\star m_1\in M_{Y,X}$
\begin{equation}
m_0\star m_1= (\underset{Y}{\oplus}(1,1))\circ (m_0\oplus m_1)\circ (\underset{X}{\oplus}  \left(\begin{matrix}
1 \\ 1
\end{matrix}\right)) =(\underset{Y}{\oplus}(1,1))\circ (1_{Y_0}\circ m_0\circ 1_{X_0}^t+1_{Y_1}\circ m_1\circ 1_{X_1}^t)\circ (\underset{X}{\oplus}  \left(\begin{matrix}
1 \\ 1
\end{matrix}\right)).
\end{equation}
is equal to $m_0+m_1$ (indeed, $(m_0+m_0')\star(m_1+m_1')=(m_0\star m_1)+(m_0'\star m_1')$ ).} \\

Assume that $M$ is commutative as $\FF(R)$ -module. For $r\in R,m\in M_{Y,X}$, put
\begin{equation}
r\cdot m:= (\underset{Y}{\oplus}(r))\circ m=m\circ (\underset{X}{\oplus}(r))\hspace{10mm} (M \;\;\; \text{central}).
\end{equation}
We have $(r_1\cdot r_2)\cdot m=r_1\cdot (r_2\cdot m), 1\cdot m=m,0\cdot m=0,r\cdot (m+m')\equiv r\cdot m+r\cdot m',$ (by definition!),
and moreover,

\begin{equation}
\begin{array} {lcl}  r\cdot m+r'\cdot m & = &(\underset{Y}{\oplus}(1,1))\circ (1_{Y_0}\circ r\cdot m\circ 1_{X_0}^t+1_{Y_1}\circ r'\cdot m\circ 1_{X_1}^t)\circ (\underset{X}{\oplus} \left(\begin{matrix}
1 \\ 1
\end{matrix}\right)) \\ & = &(\underset{Y}{\oplus} (r,r'))\circ (m\oplus m)\circ (\underset{X}{\oplus} \left(\begin{matrix}
1 \\ 1
\end{matrix}\right)) \\  & = &\bigg(\underset{Y}{\oplus} (r,r')\circ \left(\begin{matrix}
1 \\ 1
\end{matrix}\right) \bigg)\circ m \hspace{10mm} (M \text{ commutative}) \\ & = & (r+r')\cdot m.  \end{array}
\end{equation}
\\
Thus $M_{Y,X}$ is an $R$ -module. For $a\in \FF(R)_{Y',Y}, b\in \FF(R)_{X,X'}$, the map $m\mapsto a\circ m\circ b$ is an $R$ -module homomorphism, so $M\in (R$-mod$)^{\FF(R)\times \FF(R)^{op}}$. \\
Note that $M$ is $\NN$- true, and if $R\in CRing$, has negatives, $(-1)\cdot m=(\underset{Y}{\oplus}(-1))\circ m=m\circ (\underset{X}{\oplus}(-1))$ is equal to $-m$, the inverse of $m$, and $M$ is $\ZZ$- true.  \\
Conversely, if $M\in \; (R$-mod$)^{\FF(R)\times \FF(R)^{op}}$, is
an $\FF(R)$-module with values in $R$ -modules such that for $r\in
R, m\in M_{Y,X}, (\underset{Y}{\oplus}(r))\circ m=m\circ
(\underset{X}{\oplus}(r))=r\cdot m$, (i.e. $M$ is central, and the
action of the monoid $\FF(R)_{1,1}=R$ is the given $R$ -module
structure on $M_{Y,X})$, then $M$ is commutative as $\FF(R)$
-module:
$$\left(\underset{Y}{\oplus}(r_1,\dots, r_n) \right)\circ \bigg( \overset{n}{\underset{i=1}{\oplus}} m  \bigg) \circ \bigg( \underset{X}{\oplus}  \left(\begin{matrix}
r_1' \\ \vdots \\ r_n'
\end{matrix}\right)  \bigg)=$$
$$\left(\underset{Y}{\oplus}(1,\dots, 1) \right)\circ \bigg( \overset{n}{\underset{i=1}{\oplus}} r_i\cdot r_i' \cdot m  \bigg) \circ \bigg( \underset{X}{\oplus}  \left(\begin{matrix}
1 \\ \vdots \\ 1
\end{matrix}\right)  \bigg)=$$
$$r_1\cdot r_1'\cdot m+\dots +r_n r_n'\cdot m=$$
$$[r_1 r_1'+\dots +r_n r_n']\cdot m= \hspace{3mm} (M_{Y,X}\in R\text{-mod})$$
\begin{equation}
\bigg[(r_1,\dots,r_n)\circ \left(\begin{matrix}
r_1' \\ \vdots \\ r_n'
\end{matrix}\right)\bigg]\cdot m
\end{equation}

\section{Sheaves of $\mathcal{O}_X$- modules}
\smallbreak

\defin{7.3.1} \textit{Given $(X,\mathcal{O}_X)\in \FR^{(t)}\mathcal{S}p$, define an $\mathcal{O}_X\text{-mod}^{(t)}$ $M$ to be a functor
\begin{equation}
U\mapsto M(U), \;\text{for } U\subseteq X\;\;\;\text{open},
\end{equation}
such that,
\begin{equation}
M(U)\in C\mathcal{O}_X(U)\text{-mod}^{(t)},
\end{equation}
for $U\subseteq U'\subset X$ open:
\begin{equation}
b|_U\circ m|_U\circ d|_U=(b\circ m\circ d)|_U,
\end{equation}
\begin{equation}
(\text{resp.} \; (m|_U)^t=m^t|_U)  ,
\end{equation}
and $\forall X,Y\in \FF$,
\begin{equation}
U\mapsto M(U)_{Y,X} \;\;\; \text{is a sheaf}.
\end{equation}}
\defin{7.3.2}\textit{Let $A\in \CFR^{(t)}$, $X=\Spec^{(+)}A$, $M\in C\Amod^{(t)}$. Define $\tilde{M}\in \mathcal{O}_X\text{-mod}^{(t)}$,
\begin{equation}
\tilde{M}(U)_{Y,X}\equiv \{s:U\rightarrow \underset{\mathfrak{p}\in U}{\coprod} (M_{\mathfrak{p}})_{Y,X}\;, s(\mathfrak{p})\in (M_{\mathfrak{p}})_{Y,X}\;\;\; \text{such that locally   } s(\mathfrak{p})=m/f.\}
\end{equation}}
We have,
\prop{7.3.1} \textit{For $\mathfrak{p}\in \Spec^{(+)} A$,
\begin{equation}
(\tilde{M})_{\mathfrak{p}}=M_{\mathfrak{p}}.
\end{equation}
}
\begin{proof} see Proposition 3.4.2. \end{proof}
\thm{7.3.2} For $f\in A_{1,1}^{(+)}$,
\begin{equation}
\tilde{M}(D^{(+)}(f))\equiv M_f.
\end{equation}
\begin{proof}
Replace $a\in A_{Y,X}$ by $m\in M_{Y,X}$ in Proposition 3.4.3.
\end{proof}
\thm{7.3.3}\textit{Let $X\in \gFRst$, $M\in \mathcal{O}_X\text{-mod}^{(t)}$, the following conditions are equivalent,
\begin{equation}
(1)\;\;\; X=\bigcup \Spec^{(+)} A_i, \;\;\; \exists M_i\in A_i\text{-mod}^{(t)}, M|_{\Spec^{(+)} A_i}=\tilde{M_i},
\end{equation}
\begin{equation}
(2)\;\;\; \forall \Spec^{(+)} A\subseteq X, \;\;\; M(\Spec^{(+)} A)\iso M|_{\Spec^{(+)} A},
\end{equation}
$$(3)\;\;\;\forall \; U\subseteq X\; \text{open, } \forall g\in \mathcal{O}_X(U)_{1,1}^{(+)},\; \text{letting } D(g)=\{\mathfrak{p}\in U, g|_{\mathfrak{p} } \in GL_1(\mathcal{O}_{X,\mathfrak{p}})\}\subseteq U,$$
\begin{equation}
\text{restriction induces isomorphism: } M(U)_g\iso M(D(g)).
\end{equation}}
We say $M$ is "quasi- coherent" $\mathcal{O}_X$-module, and we denote by $q.c.\; \mathcal{O}_X$-$\text{mod}^{(t)}$ the full subcategory of $\mathcal{O}_X$-$\text{mod}^{(t)}$ consisting of the quasi coherent $\mathcal{O}_X$-modules. \vspace{1mm}\\
For an affine $X=\Spec^{(+)} A$, we have an equivalence
\begin{equation}
C\Amod^{(t)}\iso q.c.\;\; \mathcal{O}_X\text{-mod}^{(t)}\subseteq \mathcal{O}_X\text{-mod}^{(t)},
\end{equation}
$$M\longmapsto \tilde{M}, $$
\begin{equation}
M(X)\longmapsfrom M.
\end{equation}

\section{Extension of scalars}
\smallbreak
For $\varphi \in \FR^{(t)}(B,A)$, we have an induced pair of adjoint functors: We use geometric notations
$$\xymatrix{\Amod^{(t)} \ar@<1ex>@{->}[r]^{\varphi_*}& B\text{-}mod^{(t)}\ar@<1ex>@{->}[l]^{\varphi^*} }$$
The right adjoint is:
\begin{equation}
\xymatrix{ N\ar@{|->}[r] &\varphi_*N\equiv N_B:=N}
\end{equation}
with $B$- action given via $\varphi$:
\begin{equation}
b_1\circ n\circ b_2:= \varphi(b_1)\circ n\circ \varphi(b_2).
\end{equation}
The left adjoint will be denoted by $M\mapsto \varphi^*M=M^A$. The abelian group $(M^A)_{Y,X}$ is obtained from the free sum:
\begin{equation}
\underset{W,Z\in \FF}{\bigoplus} \underset{A_{Y,W}\times A_{Z,X}}{\bigoplus} M_{W,Z}
\end{equation}
whose elements can be written as sums
\begin{equation}
\sum_{i=1}^k a_i\circ [m_i]\circ a_i',\; a_i\in A_{Y,W_i},\; a_i'\in A_{Z_i,X},\; m_i\in M_{W_i,Z_i}.
\end{equation}

\begin{equation}
\text{with} \hsm\;\; a\circ [m+m']\circ a'=a\circ [m]\circ a'+a\circ [m']\circ a',
\end{equation}
by dividing by the subgroup generated by all elements of the form
\begin{equation}
a\circ [b\circ m\circ b']\circ a'-(a\circ \varphi(b))\circ [m] \circ (\varphi(b')\circ a').
\end{equation}
If $\varphi\in \FRt(B,A)$, $M\in B\text{-}mod^t\text{ then } M^A\in \Amod^t$ has automatically an involution,
\begin{equation}
(\sum a_i\circ [m_i]\circ c_i)^t=\sum c_i^t\circ [m_i^t]\circ a_i^t.
\end{equation}

\section{Infinitesimal extensions}
\smallbreak Let $A\in \FR,M\in \Amod$, We define the infinitesimal
extension $A\prod M\in (\slfrac{\FR}{A})^{ab}$, an abelian group
object of the category $\slfrac{\FR}{A}$ of $\frs$ over $A$:
\begin{equation}
(A\Pi M)_{Y,X}:=A_{Y,X}\Pi M_{Y,X}
\end{equation}
\begin{equation}
(a,m)\circ (b,n):=(a\circ b,a\circ n+m\circ b)
\end{equation}
\begin{equation}
(a_0,m_0)\oplus (a_1,m_1):=(a_0\oplus a_1,m_0\oplus m_1)
\end{equation}
We have homomorphism
\begin{equation}
\pi \in \FR(A\Pi M,A)\;,\; \pi(a,m)=a,
\end{equation}
Furthermore the map,
\begin{equation}
\mu\in \slfrac{\FR}{A}\bigg( (A\Pi M) \prod_A (A\Pi M), A\Pi M\bigg)\;,\; \mu((a,m),(a,m'))=(a,m+m').
\end{equation}
satisfy associativity, commutativity, unit: $\epsilon\in \slfrac{\FR}{A}\bigg(A,A\Pi M\bigg),\; \epsilon (a)=(a,0)$, antipode:  $S(a,m)=(a,-m)$ and so makes $A\prod M$ into an abelian group object in $\slfrac{\FR}{A}$. \vspace{3mm}

In the case where $A\in \FRt,M\in \Amod^t, A\prod M$, has a natural involution
$$(a,m)^t:=(a^t,m^t),$$
 and so $A\prod M\in  (\slfrac{\FRt}{A})^{ab}$.\vspace{3mm}\\

Note that we have the (strict) implications
$$\begin{matrix}  A \text{ totally- commutative}\\ M \text{ commutative}\hspace{12mm}\end{matrix} \implies A\Pi M \text{ commutative}\implies \begin{matrix} A \text{  commutative } \;(\text{as } \fr)\hspace{2mm} \\ M \text{ commutative} \;(\text{as } A\text{-module}) \end{matrix}.$$

\section{Derivations and differentials}
\smallbreak
\defin{7.6.1} \textit{Let $\varphi\in \FR^{(t)}(C,A)$ , $M\in \Amod^{(t)}$.\\
Define the $C$ -derivations from $A$ to $M$ to be the set:
$$\Der^{(t)}_C(A,M):=\bigg\{\delta= \{\delta_{Y,X}:A_{Y,X}\rightarrow M_{Y,X}\}\text{ such that }$$
$$\hspace{37mm}(*) \;\; Leibnitz:\; \delta(a\circ a')=\delta(a)\circ a'+a\circ \delta(a')$$
$$\hspace{37mm}(**)\;\; C\text{-}\; linear \; \delta(\varphi(c))\equiv 0. \hspace{30mm} $$
$$\hspace{17mm}(***) \;\delta (a_0\oplus a_1)= \delta (a_0)\oplus \delta (a_1)\bigg\}$$
\begin{equation}
(\text{resp.} \;\;\; \delta(a)^t=\delta(a^t))
\end{equation}}
It is a functor $\Der^{(t)}_C(A,\underline{ \;\;\; }): \Amod^{(t)}\rightarrow Ab$, representable by the $A$- module of K\"{a}hler differentials $\Omega(\slfrac{A}{C})\in \Amod^{(t)}$.
The abelian group $\Omega(\slfrac{A}{C})_{Y,X}$ has elements that are sums of the form
\begin{equation}
\sum_{i=1}^{k}m_i\cdot a_i'\circ d(a_i)\circ a_i''
\end{equation}
\begin{equation}
m_i\in \ZZ,a_i'\in A_{Y,W_i}, a_i\in A_{W_i,Z_i}, a_i''\in A_{Z_i,X}.
\end{equation}
with the relations,
$$(*) \;\; Leibnitz:\; a'\circ d(b\circ b')\circ a''=(a'\circ b)\circ d(b')\circ a''+a'\circ d(b)\circ (b'\circ a'')$$
\begin{equation}
(**)\;\; C\text{-}\; linear \; a'\circ d(\varphi(c))\circ a''\equiv 0.\hspace{60mm}
\end{equation}
(or equivalently,
$$ a'\circ \varphi(c')\circ d(a)\circ \varphi(c'')\circ a''=a'\circ d(\varphi(c')\circ a\circ \varphi (c''))\circ a''.)$$
$$(***)\;\;a'\circ d(a_0\oplus a_1)\circ a''=(a'\circ 1_{Y_0})\circ d(a_0)\circ (1_{X_0}^t\circ a'')+(a'\circ 1_{Y_1})\circ d(a_1)\circ (1_{X_1}^t\circ a''),\; a_i\in A_{Y_i,X_i}.$$
If $\varphi\in \FRt(C,A)$, i.e. has involution, then $\Omega(\slfrac{A}{C})\in \Amod^t$:
\begin{equation}
\bigg( \Sigma m_i\cdot a_i'\circ d(a_i)\circ a_i''\bigg)^t= \Sigma m_i\cdot (a_i'')^t\circ d(a_i^t)\circ (a_i')^t.
\end{equation}
\vspace{3mm} \\
E.g. for $A=C[\delta_{Y,X}]=C\otimes_{\FF} \FF[\delta_{Y,X}]$ we have,
\begin{equation}
\Omega(\slfrac{C[\delta_{Y,X}]}{C})\equiv \text{free } C[\delta_{Y,X}]\text{- module of degree }(Y,X) \text{ with generator }  d(\delta_{Y,X})  .
\end{equation}
and similarly, for $A=C[\delta_{Y,X},\delta_{Y,X}^t]$ the free - $C$ - $\fr$ with involution,
$$\Omega(\slfrac{C[\delta_{Y,X},\delta^t_{Y,X}]}{C})\equiv \text{free } C[\delta_{Y,X},\delta_{Y,X}^t]\text{- mod with involution of degree }(Y,X) $$
\begin{equation}
\text{ generated by } d(\delta_{Y,X}).
\end{equation}
\vspace{3mm} \\
For $\varphi \in \FR^{(t)}(C,A), \; B\in C \diagdown \FR^{(t)}\diagup A$, i.e. we have $\FR^{(t)}$- maps
$$\epsilon:C\rightarrow B,$$
$$\pi:B\rightarrow A,$$
\begin{equation}
\pi\circ \epsilon =\varphi,
\end{equation}
we have the following identifications, for $M\in \Amod^{(t)}$:
$$ \bigg(C \diagdown \FR^{(t)}\diagup A\bigg)(B,A\Pi M)\equiv
\begin{Bmatrix}
\psi:B\rightarrow A\prod M, \hsm \psi(b)=(\pi(b),\delta(b)) \\
 \delta_{Y,X}:B_{Y,X}\rightarrow M_{Y,X}, \delta(b\circ b')=\delta (b)\circ \pi(b')+\pi(b)\circ \delta(b')\\
\delta(\epsilon(C))\equiv 0 \\
\delta(b_0\oplus b_1)=\delta(b_0)\oplus \delta (b_1) \\
(\text{respectfully},\; \delta(b^t)=\delta(b)^t)
\end{Bmatrix}$$
\begin{equation}
\equiv \Der^{(t)}_C(B,M_B)\equiv B\text{-mod}^{(t)}(\Omega(B/C),M_B)\equiv \Amod^{(t)}(\Omega(B/C)^A,M).
\end{equation}
Thus we have the adjunction:
\begin{equation}\label{7.6.10}
\xymatrix{\Omega(\slfrac{B}{C})^A&\Amod^{(t)} \ar@/^/@{|->}[d]& M \ar@{|->}[d] \\
 B\ar@{|->}[u]&C \diagdown \FR^{(t)}\diagup A\ar@/^/@{|->}[u]& A\Pi M}
\end{equation}
Restricting to \emph{commutative} $A$-modules $CA$-$mod^{(t)}$, we have similar adjunction, with $\Omega(B/C)^A$ replaced by its commutative quotient $C\Omega(B/C)^A$.

\section{Properties of differentials}
\smallbreak

Given a homomorphism $k\rightarrow A$ of $\FR^{(t)}$, and $M\in \Amod^{(t)}$, we have the bijection:
$$\Der^{(t)}_k(A,M)\equiv \Amod^{(t)}\bigg(\Omega(A/k),M\bigg).$$
\begin{equation}
\varphi\circ d_{A/k}\mapsfrom \varphi,
\end{equation}
where $d_{A/k}:A\rightarrow \Omega(A/k)$ is the universal derivation.


\property{0} Given any commutative diagram in $C\FR^{(t)}$, and $M'\in A'\text{- mod}^{(t)}$,
\begin{equation}
\xymatrix{A\ar[r]& A'\\ k\ar[r]\ar[u]&k'\ar[u]}
\end{equation}
we have a sequence of homomorphisms:
\begin{equation}
\xymatrix{ \mathcal{D}er_{k'}^{(t)}(A',M')\ar@{^{(}->}[r] & \mathcal{D}er^{(t)}_k(A',M')\ar[r] &  \mathcal{D}er^{(t)}_k(A,M') }
\end{equation}
represented by the $A'$- mod homomorphisms (with commutative
diagrams of derivations)

\begin{equation}
\xymatrix{\Omega(A/k)^{A'}\ar[r] &\Omega(A'/k)\ar@{->>}[r]&\Omega(A'/k') \\ A\ar[u]_{(d_{A/k})^{A'}}\ar[r]&A'\ar[u]_{d_{A'/k}}\ar[ru]_{d_{A'/k'}}}
\end{equation}

\property{1 (First exact sequence)}
Given a commutative diagram in $\FR^{(t)}$
\begin{equation}
\xymatrix{A\ar[rr]^{\varphi}&&A' \\ &k\ar[lu]\ar[ru]&}
\end{equation}
we have an exact sequence:
\begin{equation}\label{firstexactsequence}
\xymatrix{\Omega(A/k)^{A'}\ar[r]&\Omega(A'/k)\ar[r]&\Omega(A'/A)\ar[r]&0}
\end{equation}
\begin{proof} {Applying $A'\text{-} mod(\underline{\;\;},M')$ this is equivalent to
\begin{equation}
\xymatrix{0\ar[r]& \mathcal{D}er_A(A',M')\ar[r]& \mathcal{D}er_k(A',M')\ar[r]& \mathcal{D}er_k(A,M')}
\end{equation}
is exact for every $M'\in A'$-mod, which is clear. }
\end{proof}
The first exact sequence (\ref{firstexactsequence}) will be exact on the left, if and only if,  any derivation $D:A\rightarrow M'$ into an $A'$- mod has an extension to a derivation $D'$ of $A'$
\begin{equation}
\xymatrix{A\ar^D[rd]\ar[dd]& \\ &M' \\A'\ar@{-->}_{D'}[ru]&}
\end{equation}
e.g. this holds if $\varphi:A\rightarrow A'$ is a retract: have $\psi:A'\rightarrow A$, $\psi\circ \varphi =id_A$, and can take $D'=D\circ \psi$.

\property{2 (Second exact sequence)} For a surjective $\varphi:A\sur A'$,
\begin{equation}
\xymatrix{\mathcal{KER}(\varphi)= A\prod_{A'} A\ar@<2pt>[r]\ar@<-2pt>[r]&A\ar@{->>}^{\varphi}[rr]&&A' \\
&&k\ar[ru]\ar[lu]&}
\end{equation}
we have an exact sequence
\begin{equation*}
\Omega(A\prod_{A'} A/k)^{A'}\xrightarrow{\delta} \Omega(A/k)^{A'}\rightarrow \Omega(A'/k)\rightarrow 0
\end{equation*}
\begin{equation}
d(a_0,a_1)\overset{\delta}{\longmapsto} d(a_0)-d(a_1).\hspace{20mm}
\end{equation}
\begin{proof} {Applying the functor $A'$- mod$(\underline{\;\;},M')$ this is equivalent to the left exact sequence,
\begin{equation*}
0\rightarrow  \mathcal{D}er_k(A',M')\rightarrow  \mathcal{D}er_k(A,M')\rightarrow  \mathcal{D}er_k(A\prod_{A'} A, M')\hspace{10mm}
\end{equation*}
\begin{equation}
\hspace{52mm} D\longmapsto D(a_0,a_1)=D(a_0) -D(a_1).
\end{equation}
which is clear.}
\end{proof}
Moreover, we can replace $\Omega^{(2)}= \Omega(A\underset{A'}{\prod} A/k)^{A'}$ on the left of $(7.7.11)$ by its odd quotient (w.r.t. the involution permuting the factors), $\Omega^{(-)}=\slfrac{\Omega^{(2)}}{d(a_0,a_1)+d(a_1,a_0)}$.
Furthermore, we have a map $\delta: \Omega^{(3)}= \Omega (A\underset{A'}{\prod} A\underset{A'}{\prod} A/k)^{A'}\rightarrow \Omega^{(2)},$
$$ d(a_0,a_1,a_2)\mapsto d(a_1,a_2)-d(a_0,a_2)+d(a_0,a_1),$$
and since $\delta\circ \delta =0$, we can replace $\Omega^{(2)}$ on the left of $(7.7.11)$ by $\slfrac{\Omega^{(-)}}{\delta (\Omega^{(3)})}$.

\property{3}
\begin{equation}
\Omega(A_1\otimes_k A_2/k)\cong \Omega(A_1/k)^{A_1\otimes_k A_2}\oplus  \Omega(A_2/k)^{A_1\otimes_k A_2}.
\end{equation}
\property{4}
\begin{equation}
\Omega(A\otimes_k k'/k')\cong \Omega(A/k)^{A\otimes_k k'}.
\end{equation}

\property{5} Differentials commute with direct limits: we have isomorphisms of $A=\varinjlim A_i$- modules,
\begin{equation}
\Omega(\varinjlim A_i/ \varinjlim k_i)\cong \varinjlim(\Omega(A_i/k_i)),
\end{equation}
\begin{proof}
Taking $\Amod(\underline{\;\;},M)$ this is equivalent to
\begin{equation}\Amod(\Omega(A/\varinjlim k_i),M)\equiv\mathcal{D}er_{\varinjlim k_i}(A,M)\equiv \varprojlim \mathcal{D}er_{k_i}(A_i,M) \equiv \varprojlim A_i\text{-mod}(\Omega(A_i/k_i),M)$$
$$\equiv \Amod (\varinjlim \Omega(A_i/k_i),M).
\end{equation}
\end{proof}
Properties $(0)-(5)$ all hold with $\Omega$ replaced by its commutative quotient $C\Omega$. \vspace{3mm} \\
Given multiplicative sets $\sigma\subseteq k_{1,1}, S\subseteq A_{1,1}$, such that $\varphi:k\rightarrow A$, takes
$\varphi_{1,1}(\sigma)\subseteq S$, we have
\begin{equation}
C\Omega (S^{-1}A/\sigma^{-1} k )\cong S^{-1}C\Omega (A/k).
\end{equation}
\vspace{1.5mm}\\
\exmpl{7.7.1} For a map of monoids $\varphi\in Mon(M_0,M_1)$, and the associated map of $\frs$ $\Phi= \FF\{\varphi\}\in \FR(\FF\{M_0\},\FF\{M_1\})$, the $\FF\{M_1\}$- module $\Omega =\Omega(\slfrac{\FF\{M_1\}}{\FF\{M_0\}})$ is generated by
$d(m)\in \Omega_{1,1}, m\in M_1$, and so is a matrix- module $\Omega_{Y,X}=(\Omega_{1,1})^{Y\times X}$, and $\Omega_{1,1}=\Omega(\slfrac{\ZZ[M_1]}{\ZZ[M_0]})$ is the usual bi- module of differentials of the associated rings. For $M_i$ commutative, $C\Omega_{Y,X}=(C\Omega_{1,1})^{Y\times X}$, and $C\Omega_{1,1}$ the usual module of K\"{a}hler differentials of $\slfrac{\ZZ[M_1]}{\ZZ[M_0]}$.
\exmpl{7.7.2} For a map of rings $\varphi\in Mon(R_0,R_1)$, and the associated map of $\frs$ $\Phi= \FF(\varphi)\in \FR(\FF(R_0),\FF(R_1))$, the $\FF(R_1)$- module $\Omega =\Omega(\slfrac{\FF(R_1)}{\FF(R_0)})$ is generated by
$d(r)\in \Omega_{1,1}, r\in R_1$, and so is a matrix- module $\Omega_{Y,X}=(\Omega_{1,1})^{Y \times X}$, and $\Omega_{1,1}=\Omega(\slfrac{\ZZ[R_1]}{\ZZ[R_0]})$ is the usual bi- module of differentials of the associated rings. For $R_i$ commutative, $C\Omega_{Y,X}=(C\Omega_{1,1})^{Y\times X}$, and $C\Omega_{1,1}$ the usual module of K\"{a}hler differentials of $\slfrac{\ZZ[R_1]}{\ZZ[R_0]}$.

\section{Differentials of $\FF(\ZZ)$ and $\FF(\NN)$. }
\smallbreak


\thm{7.8.1}\label{thm7.8.1}\textit{The module $\Omega=\Omega(\slfrac{\FF(\NN)}{\FF})$, (respectively, $\Omega(\slfrac{\FF(\ZZ)}{\FF\{\pm 1\}})$ ) is defined as:
\begin{equation}
\begin{split}
&\Omega_{n,m}= \text{free abelian group on generators:}\\
&[a | \begin{matrix}
b \\ b'
\end{matrix}]=a\circ d(1,1)\circ \left(\begin{matrix}
b \\ b'
\end{matrix}\right) \text{ and } [a,a'|b]=(a,a')\circ d\left(\begin{matrix}
1 \\ 1
\end{matrix}\right)\circ b, \hsm \forall b,b'\in \NN^m,\forall a,a'\in \NN^n,
\end{split}
\end{equation}
(respectfully, $b,b'\in \ZZ^m, a,a'\in \ZZ^n$), \\
modulo relations:
\begin{equation*}
\begin{matrix}
(0): & \left[0\bigg|\begin{matrix}
b \\ b'
\end{matrix}\right]=\left[a\bigg |\begin{matrix}
0 \\ b'
\end{matrix}\right]=\left[a\bigg|\begin{matrix}
b \\ 0
\end{matrix}\right]=[a,a'|0]=[0,a'|b]=[a,0|b]=0, \vspace{3mm} \\
(Comm): & \left[a\bigg|\begin{matrix}
b \\ b'
\end{matrix}\right]=\left[a\bigg|\begin{matrix}
b' \\ b
\end{matrix}\right]\;,\; [a,a'|b]=[a',a|b],   \vspace{3mm}  \\
(\text{Ass}): & \left[a\bigg|\begin{matrix}
b_1+b_2 \\ b_3
\end{matrix}\right]+ \left[a\bigg |\begin{matrix}
b_1 \\ b_2
\end{matrix}\right]=\left[a\bigg|\begin{matrix}
b_1 \\ b_2+b_3
\end{matrix}\right]+\left[a\bigg|\begin{matrix}
b_2 \\ b_3
\end{matrix}\right]   \vspace{3mm}    \\
 & [a_1+a_2,a_3|b]+[a_1,a_2|b]=[a_1,a_2+a_3|b]+[a_2,a_3|b], \vspace{3mm}  \\
(\text{Almost linear}): & [a_1,a_2|b_1+b_2] +\left[a_1+a_2\bigg |\begin{matrix}
b_1 \\ b_2
\end{matrix}\right] = \left[a_1\bigg |\begin{matrix}
b_1 \\ b_2
\end{matrix}\right]+\left[a_2\bigg|\begin{matrix}
b_1 \\ b_2
\end{matrix}\right]+[a_1,a_2|b_1]+[a_1,a_2|b_2].
\end{matrix}
\end{equation*}
respectively and the relations,
\begin{equation}
\begin{matrix}
(\text{Cancellation}): &\left[a\bigg |\begin{matrix} b \\ -b
\end{matrix}\right]+[a,-a|b]=0    \vspace{3mm}        \\
(-1): & \left[-a\bigg|\begin{matrix}
b \\ b'
\end{matrix}\right]=\left[a\bigg|\begin{matrix}
-b \\ -b'
\end{matrix}\right]\;,\; [a,a'|-b]=[-a,-a'|b],
\end{matrix}
\end{equation}}
\begin{proof}
By the description of \S \ref{genandrel}, $\FF(\NN)=\slfrac{\FF[\delta_{1,2},\delta_{1,2}^t]}{\text{relations}}$, (respectively  $\FF(\ZZ)= \FF\{\pm 1\}[\delta_{1,2},\delta_{2,1}^t]/\text{relations}$), the second exact sequence (7.7.11) gives that $\Omega$ is the free $\FF(\NN)$ (respectively $\FF(\ZZ)$)- module on $d\delta_{1,2}=d(1,1)$ and $d\delta_{1,2}^t=d\left( \begin{matrix} 1 \\ 1 \end{matrix}\right)$, modulo the derived - relations.\\
The $(0)$ relation follows, from the implication,
\begin{equation}
(1,1)\circ \left( \begin{matrix}
1 \\ 0
\end{matrix}\right)=(1) \implies d(1,1)\circ \left(\begin{matrix}
1 \\ 0
\end{matrix}\right)=0.
\end{equation}
The $(Comm)$ relation follows since,
\begin{equation}
(1,1)\circ \left( \begin{matrix}
0 & 1 \\ 1 & 0
\end{matrix}  \right)  = (1,1) \implies d(1,1)\circ \left(\begin{matrix}
0 & 1 \\ 1 & 0
\end{matrix} \right)=d(1,1).
\end{equation}
The $(Ass)$ relation follows since,
$$(1,1)\circ ((1,1)\oplus id_1)=(1,1)\circ (id_1\oplus (1,1))$$
$$\Downarrow $$
$$d(1,1)\circ  ((1,1)\oplus (1))+ (1,1)\circ  (d(1,1)\oplus 0)=d(1,1)\circ ((1)\oplus (1,1))+(1,1)\circ (0\oplus d(1,1))$$
$$\Downarrow $$
\begin{equation}
[a|\begin{matrix} b_1+b_2\\ b_3\end{matrix} ]+[a|\begin{matrix} b_1 \\ b_2\end{matrix}]=[a|\begin{matrix} b_1 \\ b_2+b_3\end{matrix}]+[a|\begin{matrix} b_2\\ b_3\end{matrix}].
\end{equation}
Thus $\left[ a\bigg | \begin{matrix} - \\ - \end{matrix}\right], [-,-|b]$ are symmetric normalized 2-cycles  $\forall b, \forall a$. \\
The $(Almost \; linear)$ relation follows since by total commutativity,
$$\bigg (\begin{matrix} 1 \\ 1 \end{matrix} \bigg) \circ (1,1)=\bigg (\begin{matrix} 1 & 1 \\ 1 & 1 \end{matrix} \bigg )= \bigg ( \begin{matrix} 1 & 1 & 0 & 0 \\ 0 & 0& 1& 1 \end{matrix}\bigg )\circ \left( \begin{matrix} 1 & 0 \\ 0 & 1 \\ 1 & 0 \\ 0 & 1 \end{matrix} \right)= ((1,1)\oplus (1,1))\circ (\bigg ( \begin{matrix} 1 \\ 1 \end{matrix} \bigg ) \oplus \bigg ( \begin{matrix} 1 \\ 1 \end{matrix}\bigg ) ) $$
$$\Downarrow$$
$$ d\left(\begin{matrix}1 \\ 1 \end{matrix}\right)\circ (1,1)+\left( \begin{matrix} 1 \\ 1 \end{matrix}\right)\circ d(1,1)=(d(1,1)_{1/1,2} +d(1,1)_{2/3,4})\circ  \left( \begin{matrix} 1 & 0 \\ 0 & 1 \\ 1 & 0 \\ 0 & 1 \end{matrix} \right)+$$
\begin{equation}\label{7.8.6}
+ \bigg ( \begin{matrix} 1 & 1 & 0 & 0 \\ 0 & 0& 1& 1 \end{matrix}\bigg ) \circ (d\bigg ( \begin{matrix} 1 \\ 1 \end{matrix} \bigg )_{1,3/1}+d \bigg ( \begin{matrix} 1 \\ 1 \end{matrix} \bigg )_{2,4/2} ) .
\end{equation}
Here the subindices indicate how the matrices act on the differentials (left/right), so that multiplying (\ref{7.8.6}) by $(a_1,a_2)$ on the left, and by $\left(\begin{matrix} b_1 \\ b_2 \end{matrix}\right)$ on the right, we obtain the \emph{almost linear} form of the theorem. \vspace{3mm}\\
Respectively for $\ZZ$, the $(Cancelation)$ relation follows since,
\begin{equation}
(1,1) \circ \bigg ( \begin{matrix} 1 & 0 \\ 0& -1 \end{matrix} \bigg ) \circ \bigg( \begin{matrix} 1 \\ 1 \end{matrix}\bigg)=(0) \implies d(1,1)\circ \bigg ( \begin{matrix} 1 \\ -1 \end{matrix}\bigg )+ (1,-1)\circ d\bigg ( \begin{matrix} 1 \\ 1 \end{matrix}\bigg ) = 0 .
\end{equation}
\end{proof}
The commutative quotient $C\Omega=C\Omega(\FF(\NN)/\FF)$ (respectively, $C\Omega(\FF(\ZZ)/\FF\{\pm 1\})$) is obtained by adding the relations for all $\lambda\in \NN (\text{respectively, }\ZZ)$, (cf. Remark 3 of \S \ref{7.2}):
\begin{equation}
\begin{matrix}  \lambda\cdot \left[ a\bigg |\begin{matrix} b\\ b'\end{matrix}\right]=\left[ \lambda\cdot a\bigg |\begin{matrix} b\\ b'\end{matrix}\right]=\left[ a\bigg |\begin{matrix} \lambda\cdot b\\ \lambda\cdot b'\end{matrix}\right]\\
\lambda\cdot [a,a'|b]=[\lambda\cdot a,\lambda\cdot a'|b]=[a,a'|\lambda \cdot b].  \end{matrix}
\end{equation}
\vspace{3mm} \\
Define $\tilde{N}_{Y,X}:=$ free abelian group on generators $[a|b], a\in \ZZ^Y=\FF(\ZZ)_{Y,1}, b\in \ZZ^X=\FF(\ZZ)_{1,X}$, (i.e. we think of the $a$'s as column vectors, the $b$'s as row vectors), modulo the relations $[a|0]=[0|b]=0, [-a|b]\equiv [a|-b]=-[a|b]$. \\
$\tilde{N}=\{\tilde{N}_{Y,X}\}$ is an $\fr$ with involution,
\begin{equation}
\tilde{N}_{Z,Y}\times \tilde{N}_{Y,X}\rightarrow \tilde{N}_{Z,X},
\end{equation}
\begin{equation}
(\sum_j n_j [c_j | a_j'])\circ (\sum_i m_i [a_i|b_i])=\sum_{i,j} n_j\cdot m_i \cdot a_j'\circ a_i \cdot [c_j|b_i],
\end{equation}
\begin{equation}
(\sum_i m_i [a_i|b_i])^t=\sum_i m_i [b_i^t|a_i^t].
\end{equation}
and $\tilde{N}$ is also an $\FF(\ZZ)$ -module, via
\begin{equation}
\FF(\ZZ)_{Y',Y}\times \tilde{N}_{Y,X}\times \FF(\ZZ)_{X,X'}\rightarrow \tilde{N}_{Y',X'}
\end{equation}
\begin{equation}
A\circ (\sum_i m_i [a_i|b_i])\circ B:=\sum_i m_i[Aa_i|b_i B]
\end{equation}
(and the monoidal structure, $\oplus: \tilde{N}_{Y_0,X_0}\times \tilde{N}_{Y_1,X_1}\rightarrow \tilde{N}_{Y_0\oplus Y_1,X_0\oplus X_1}$, which is part of the $\fr$ structure, comes from the $\FF(\ZZ)$- module structure). \\
We have a surjective map
\begin{equation}
\tilde{\pi}: \tilde{N}\sur \FF(\ZZ)
\end{equation}
\begin{equation}
\tilde{\pi}_{Y,X}: \tilde{N}_{Y,X}\rightarrow \FF(\ZZ)_{Y,X}
\end{equation}
\begin{equation}
\tilde{\pi} (\sum m_i [a_i|b_i]) =\sum m_i a_i\otimes b_i\in \ZZ^Y\otimes \ZZ^X\equiv \ZZ^{Y\otimes X}\equiv \FF(\ZZ)_{Y,X}
\end{equation}
and $\tilde{\pi}$ is both a homomorphism of $\frs$ with involution, and a homomorphism of $\FF(\ZZ)$ -modules. \\
We have a surjective map $\partial: \Omega(\FF(\ZZ)/\FF\{\pm 1\})\sur \Ker (\tilde{\pi})$, define on the generators by
\begin{equation}
(\star) \hspace{5mm} \begin{matrix} \partial [a_1,a_2|b]=[a_1|b]+[a_2|b]-[a_1+a_2|b] \\ \partial [a|\begin{matrix} b_1 \\b_2\end{matrix}]=-[a|b_1]-[a|b_2]+[a|b_1+b_2]\end{matrix}
\end{equation}
(indeed it is easy to verify that $\partial (relations)=0$).\\
Similarly, let $N_{Y,X}:=$ free abelian group on generators $[a|b], a\in \ZZ^Y,b\in \ZZ^X$ modulo
\begin{equation}
\lambda \cdot [a|b]=[\lambda\cdot a|b]= [a|\lambda \cdot b], \text{for all } \lambda\in \ZZ. \
\end{equation}
$N=\{N_{Y,X}\}$ is an $\fr$ with involution, and a commutative $\FF(\ZZ)$ -module, and again we have a surjection
\begin{equation}
\pi:N\sur \FF(\ZZ)
\end{equation}
which is both a homomorphism of $\frs$ with involution, and a homomorphism of $\FF(\ZZ)$ -modules. We have a surjective map
$\bar{\partial}: C\Omega (\FF(\ZZ)/\FF\{\pm 1\})\sur \Ker (\pi)$, defined as in $(\star)$.

Thus we have exact sequence (of $\FF(\ZZ)$- modules)
\begin{equation}
\xymatrix{\Omega (\FF(\ZZ)/\FF\{\pm 1\})\ar@{->>}[d]\ar[r]^-{\partial}&\tilde{N}\ar@{->>}[d]\ar[r]^-{\tilde{\pi}}&\FF(\ZZ)\ar[r]\ar@{=}[d]&0 \\ C\Omega(\FF(\ZZ)/\FF\{\pm 1\})\ar[r]^-{\bar{\partial}}&N\ar[r]^-{\pi}&\FF(\ZZ)\ar[r]&0 }
\end{equation}
The derivation $d=\bar{\partial} \circ d_{\FF\{\ZZ\}/\FF\{\pm 1\}}:\FF(\ZZ)\rightarrow N$, is given by :

\begin{equation}
d_{1,1}\equiv 0,
\end{equation}
\begin{equation}
d_{1,2}(a,b)=[a|1,0]+[b|0,1]-[1|a,b]
\end{equation}
\begin{equation}
d_{2,1} \left(\begin{matrix} a \\ b \end{matrix}\right)= \bigg [\begin{matrix} a \\ b \end{matrix}\bigg | 1 \bigg]- \bigg [\begin{matrix} 1 \\ 0 \end{matrix}\bigg | a \bigg]-  \bigg [\begin{matrix} 0 \\ 1 \end{matrix}\bigg | b \bigg]
\end{equation}
\begin{equation}
d_{2,2} \left( \begin{matrix}a & b\\ c & d \end{matrix}\right)= \bigg [\begin{matrix} a \\ c \end{matrix}\bigg | 1,0 \bigg]
+\bigg [\begin{matrix} b \\ d \end{matrix}\bigg | 0,1 \bigg]- \bigg [\begin{matrix} 1 \\ 0 \end{matrix}\bigg | a,b \bigg] -\bigg [\begin{matrix} 0 \\ 1 \end{matrix}\bigg | c,d \bigg]
\end{equation}
For $A\in \FF(\ZZ)_{n,m}$,with columns $A^{(1)},\dots ,A^{(m)}\in \ZZ^{n}$, rows $A_1,\dots, A_n\in \ZZ^m$, letting $f_1^n,\dots, f_n^n$ be the standard column basis, $e_1^m,\dots, e_m^m$ the standard row basis,
\begin{equation}
f_i^n=\left(\begin{matrix} 0 \\ 0 \\ 1^{(i)} \\ 0 \\ \vdots \\ 0 \end{matrix}\right), e_j^m=(0,\dots,0,1^{(j)}, 0,\dots,0)
\end{equation}
we have:
\begin{equation}
d_{n,m}(A)=\sum_{j=1}^{m}[A^{(j)}|e_j^m]-\sum_{i=1}^{n}[f_i^n|A_i].
\end{equation}

\section{Differentials for commutative rings}
\smallbreak

Let $R\in C\catring$, $S\subseteq R$, a multiplicative set. \\
We have the exact sequences of $\FF(R)$- modules
\begin{equation}
\xymatrix{ \Omega(\slfrac{\FF(\ZZ)}{\FF\{\pm 1\}})^{\FF(R)}\ar[r]&\Omega(\slfrac{\FF(R)}{\FF\{\pm 1\}})\ar[r]&\Omega(\slfrac{\FF(R)}{\FF(\ZZ)})\ar[r]&0\\
\Omega(\slfrac{\FF\{S\}}{\FF\{\pm 1\}})^{\FF(R)}\ar[r]&\Omega(\slfrac{\FF(R)}{\FF\{\pm 1\}})\ar[r]&\Omega(\slfrac{\FF(R)}{\FF\{S\}})\ar[r]&0}
\end{equation}
\thm{7.9.1} For the $\FF(R)$- module $\Omega = \Omega (\slfrac{\FF(R)}{\FF})$,
the abelian group $\Omega_{n,m}$ can be described as the free $\ZZ$ -module with generators
\begin{equation}
[a|\begin{matrix} b \\ b' \end{matrix}]=\left(\begin{matrix} a_1 \\ \vdots \\ a_n \end{matrix}\right) \circ d(1,1) \circ \left(\begin{matrix}b_1 &\dots& b_m\\ b_1'&\dots& b_m'   \end{matrix}\right)
\end{equation}
\begin{equation}
[a,a'|b]=\left(\begin{matrix} a_1 & a_1' \\ \vdots & \vdots \\ a_n& a_n' \end{matrix}\right) \circ d\left(\begin{matrix} 1\\1 \end{matrix}\right) \circ (\begin{matrix}b_1 &\dots& b_m \end{matrix})
\end{equation}
\begin{equation}
[a|b]^{(r)}=\left(\begin{matrix} a_1 \\ \vdots \\ a_n \end{matrix}\right) \circ d(r) \circ (b_1,\dots,b_m),
\end{equation}
modulo the relations:\hspace{25mm} $a_i\in R^n, b_i\in R^m,r\in R$
\begin{equation}
\begin{matrix}
 \text{Zero}: &  [a|\begin{matrix} 0 \\ b\end{matrix}]=[a|\begin{matrix} b \\ 0 \end{matrix}]=0 \\
 \text{Zero}^t: & [0,a|b]= [a,0|b] =0 \\
(0): & [a,b]^{(0)}=0\\
(1): & [a,b]^{(1)}=0\\
\text{Comm:} &[a|\begin{matrix} b_1 \\ b_2 \end{matrix}]=[a|\begin{matrix}b_2 \\ b_1\end{matrix}] \\
 \text{Comm}^t: & [a_1,a_2|b] = [a_2,a_1|b]\\
 \text{Ass:} & [a|\begin{matrix} b_1+b_2\\ b_3\end{matrix}]+[a|\begin{matrix} b_1\\ b_2\end{matrix}] = [a| \begin{matrix} b_1\\ b_2+b_3\end{matrix}]+[a|\begin{matrix} b_2\\ b_3\end{matrix}]\\
 \text{Ass}^t: & [a_1+a_2,a_3|b]+[a_1,a_2|b]=[a_1,a_2+a_3|b]+[a_2,a_3|b]\\
 \text{tot-com:} & [a_1,a_2|b_1+b_2]-[a_1,a_2|b_1]  -[a_1,a_2|b_2] \\
& +[a_1+a_2|b_1,b_2]-[a_1|b_1,b_2]-[a_2|b_1,b_2]=0 \\
(r,\delta): & [a|b_1+b_2]^{(r)}=[a|b_1]^{(r)}+[a|b_2]^{(r)}\\
(\delta^t,r): & [a_1+a_2|b]^{(r)}=[a_1|b]^{(r)}+[a_2|b]^{(r)}\\
(r_1\cdot r_2): & [a|b]^{(r_1\cdot r_2)}=[a\cdot r_1 |b]^{(r_2)}+[a|r_2\cdot b]^{(r_1)}\\
(r_1+r_2): & [a|b]^{(r_1+r_2)}=[a|b]^{(r_1)}+[a|b]^{(r_2)}+[a|\begin{matrix}r_1\cdot b\\ r_2\cdot b\end{matrix}]+[r_1\cdot a,r_2\cdot a|b].
\end{matrix}
\end{equation}
The $\FF(R)$- module $\Omega(\slfrac{\FF(R)}{\FF\{S\}})$ is the quotient of $\Omega(\slfrac{\FF(R)}{\FF})$ obtained by adding the relations $[a|b]^{(s)}=0\; \forall s\in S$.
\begin{proof} Same as the proof of \ref{thm7.8.1}.1, by derivation of the relations \S \ref{genandrel} of $\FF(R)$.
\end{proof}
\rem{(cf. Remark 3\; \S \ref{7.2})}
For the commutative quotient $C\Omega=C\Omega(\FF(R)/\FF\{S\})$, $C\Omega_{n,m}$ is obtained from the $R$ -module $\Omega_{n,m}\underset{\ZZ}{\otimes} R$ by adding the relations
\begin{equation}
\begin{matrix}
  r\cdot [a| \begin{matrix} b_1 \\ b_2\end{matrix} ]= [r\cdot a| \begin{matrix} b_1\\b_2\end{matrix}]= [a| \begin{matrix} r\cdot b_1\\  r\cdot b_2\end{matrix}]\\
  r\cdot [a_1,a_2|b]= [r\cdot a_1, r\cdot a_2| b]=[a_1,a_2| r\cdot b]\\
 r\cdot [a|b]^{(r')}=[r\cdot a| b]^{(r')}=[a|r\cdot b]^{(r')}
\end{matrix}
\end{equation}\vspace{3mm} \\

Let $N_{Y,X}$ denote the $R$ -module obtained from the free $R$ -module with generators $[a|b], a\in R^Y,b\in R^X,$ modulo the relations $r\cdot [a|b]=[r\cdot a|b]=[a| r\cdot b]$. \\
We think of the $b$'s (resp. $a$'s) as row (resp. column) vector, i.e. $b\in \FF(R)_{1,X}$ (resp. $a\in \FF(R)_{Y,1}$). \\
The collection $N=\{N_{Y,X}\}$ forms an $\fr$ with involution with respect to the operations
\begin{equation}
\circ : N_{Z,Y}\times N_{Y,X}\rightarrow N_{Z,X}
\end{equation}
\begin{equation}
(\sum_j [c_j| \bar{a}_j])\circ (\sum_i [a_i|b_i])=\sum_{i,j}(\bar{a}_j\circ a_i)\cdot [c_j|b_i]
\end{equation}
\begin{equation}
(\sum_i [a_i|b_i])^t =\sum_i [b_i^t |a_i^t].
\end{equation}
It is also a commutative $\FF(R)$ -module with respects to the operation
\begin{equation}
\FF(R)_{Y',Y}\times N_{Y,X}\times \FF(R)_{X,X'}\rightarrow N_{Y',X'}
\end{equation}
\begin{equation}
A\circ (\sum_i [a_i|b_i])\circ B= \sum_i [A\circ a_i| b_i\circ B].
\end{equation}
We have a surjection of $\frs$ with involution
\begin{equation}
\pi:N\sur \FF(R),
\end{equation}
\begin{equation}
N_{n,m}\overset{\pi}{\sur} R^n\underset{R}{\otimes} R^m \equiv R^{n\times m}\equiv \FF(R)_{n,m},
\end{equation}
\begin{equation}
\pi(\Sigma r_i [a_i; b_i])= \Sigma r_i\cdot a_i\otimes b_i.
\end{equation}
We have a surjection $\partial: C\Omega(\FF(R)/\FF)\sur \Ker (\pi)$, defined on the generators by
\begin{equation}
\partial [a|\begin{matrix} b_1\\ b_2\end{matrix} ]=[a|b_1]+[a|b_2]-[a|b_1+b_2]
\end{equation}
\begin{equation}
\partial [a_1,a_2|b]= -[a_1|b]-[a_2|b]+[a_1+a_2|b]
\end{equation}
$$\partial [a|b]^{(r)}=0$$
(Indeed, it is easy to check that $\partial$ takes our relations to zero. The minus signs comes from the need to have $\partial(\text{tot-com})=0$ and $\partial ((r_1+r_2))=0$). Thus we have,

\thm{7.9.2}\textit{For $R\in C\catring$, we have an exact sequence of commutative $\FF(R)$-modules
\begin{equation}
C\Omega(\FF(R)/\FF)\xrightarrow{\partial} N(R)\xrightarrow{\pi} \FF(R)\rightarrow 0.
\end{equation}
The map $\pi$ is a homomorphism of $\frs$ with involution.}

\section{Quillen model structures}
\smallbreak

Define a Quillen model structure on the category of simplicial $\frs$ (with involution) under $C\in \FR^{(t)}$, $\Delta (C \diagdown \FR^{(t)}):= (C \diagdown \FR^{(t)})^{\Delta^{op}}$:
\begin{equation}
\begin{matrix} \text{Fibrations:} & \mathcal{F}=\bigg \{\varphi:B\rightarrow B' \; |\; \forall X,Y\in\FF, \; \{\varphi_{Y,X}:B_{Y,X}\rightarrow B_{Y,X}'\} \in \mathcal{F}_{Set_0} \bigg \}, \\ \text{Weak\; equivalences:} &  \mathcal{W}= \bigg \{\varphi:B\rightarrow B' \; |\; \forall X,Y\in\FF \; \{\varphi_{Y,X}:B_{Y,X}\rightarrow B_{Y,X}'\}\in \mathcal{W}_{Set_0} \bigg \}, \\  \text{Cofibrations:} & \mathcal{C}=L(\mathcal{W}\cap\mathcal{F}). \hspace{75mm} \end{matrix}
\end{equation}
where $\mathcal{F}_{Set_0}$ (resp. $\mathcal{W}_{Set_0}$) denote the fibrations (resp. weak - equivalences) of simplicial pointed sets, and
where $L(\mathcal{W}\cap \mathcal{F})$ denotes maps satisfying the left lifting property with respects to all the acyclic fibrations $\mathcal{W}\cap \mathcal{F}$.\\
\thm{7.10.1}  \textit{This is a closed model structure.}
\begin{proof}
 (\cite{Q67}, theorem 4,II,$\S$4 )
\end{proof}
For $C\in \FR^{(t)}, \mathcal{I}\in (\slfrac{Set}{\FF\times \FF})$ define,
\begin{equation}
C[\mathcal{I}]=C[\delta_{Y_i,X_i}]= C\; \underset{\FF}{\bigotimes}\underset{i\in \mathcal{I}}{\otimes} \FF[\delta_{Y_i,X_i}].
\end{equation}
We have the adjunction,
\begin{equation}
\Delta (\slfrac{Set}{\FF\times \FF})^{(t)} \underset{U}{\overset{F=C[\;]}{\rightleftarrows}} \Delta (C \diagdown \FR^{(t)})
\end{equation}
The model structure on $\Delta(C\setminus \FR^{(t)})$ is cofibrantly generated. The set of generating cofibrations is:
\begin{equation}
\mathcal{I}=\bigg\{C[\{\partial \Delta (n)\rightarrow \Delta(n)\}^{X,Y} ]\bigg\}\equiv \bigg\{C[(\partial \Delta (n))_{Y,X}]\inj C[(\Delta (n))_{Y,X}], n\geq 1, Y,X\in \FF.\bigg\}
\end{equation}
The set of generating acyclic fibrations is:
\begin{equation}
\mathcal{J}=\bigg\{C[\{ \Lambda(n,k)\rightarrow \Delta(n)\}^{X,Y} ]\bigg\}\equiv \bigg\{ C[\Lambda(n,k)]_{Y,X}\rightarrow C[\Delta(n)]_{Y,X}, 0\leq k\leq n, Y,X\in \FF.\bigg \}
\end{equation} \vspace{3mm} \\
A model structure on simplicial $A$ -modules (with involution),
\begin{equation}
\Delta (\Amod^{(t)}):= (\Amod^{(t)})^{\Delta^{op}}
\end{equation}
is given similarly by,
\begin{equation}
\begin{matrix} \text{Fibrations:} & \mathcal{F}=\bigg \{\varphi:M \rightarrow M' \; |\; \forall X,Y\in\FF,\{\varphi_{Y,X}:M_{Y,X}\rightarrow M_{Y,X}'\}\in \mathcal{F}_{Set_0} \bigg \}, \\ \text{Weak\; equivalences:} &  \mathcal{W}=\bigg \{\varphi:M \rightarrow M' \;|\; \forall X,Y\in\FF,\{\varphi_{Y,X}:M_{Y,X}\rightarrow M_{Y,X}'\}\in \mathcal{W}_{Set_0}\bigg \}, \\  \text{Cofibrations:} & \mathcal{C}=L(\mathcal{W}\cap\mathcal{F}). \hspace{80mm} \end{matrix}
\end{equation}
For $\mathcal{I}\in \slfrac{Set}{\FF\times \FF}$, define the functor $A^-_{(t)}:\Delta  (\slfrac{Set}{\FF\times \FF})\rightarrow \Delta (\Amod^{(t)})$,
\begin{equation}
A^{\mathcal{I}}_{(t)}\equiv \underset{i\in \mathcal{I}}{\bigoplus} A^{Y_i,X_i}_{(t)}.
\end{equation}
$A^{Y_i,X_i}_{(t)}$ the free $A$- module (with involution) on a generator of degree $(Y_i,X_i)$.\vspace{1.5mm} \\
This model structure on $\Delta (\Amod^{(t)})$ is cofibrantly generated with,
\begin{equation}
\mathcal{I}=\bigg \{A^{(\partial \Delta (n))_{Y,X}}_{(t)} \rightarrow A^{\Delta(n)_{Y,X}}_{(t)} \;,\; n\geq 1\;,\; Y,X\in \FF \bigg\},
\end{equation}
\begin{equation}
\mathcal{J}=\bigg\{ A^{\Lambda(n,k)_{Y,X}}_{(t)}\rightarrow A^{\Delta(n)_{Y,X}}_{(t)}\;,\; 0\leq k\leq n\geq 1 , Y,X\in\FF \bigg\}.
\end{equation}
Under the Dold-Puppe equivalence $\Delta (\Amod^{(t)})\iso ch(\Amod^{(t)})$, these model structures corresponds to the projective model structure on chain complexes,
\begin{equation}
\begin{split}
&\text{Fibrations:} \; \mathcal{F}=\bigg \{ \varphi:M_{\bullet}\rightarrow N_{\bullet}\;| \; \varphi_n:M_n\sur N_n \; \text{surjective } \forall n\geq 1 \bigg \}, \\
&\text{Weak\; equivalences:} \; \mathcal{W}=\bigg \{\varphi:M_{\bullet} \rightarrow N_{\bullet} \;|\; H_n(\varphi): H_n(M_{\bullet})\iso H_n(N_{\bullet})\; \text{isomorphism } \forall n\geq 0 \bigg \} \\
&\text{Cofibrations:} \; \mathcal{C}=\bigg \{ \varphi:M_{\bullet}\rightarrow N_{\bullet} \;,\;\text{ injective, with }\Coker \{\varphi_n:M_n\rightarrow N_n\} \text{ projective } \forall n\geq 0 \bigg \}.
\end{split}
\end{equation} \\
For a free $C\text{-}\FR^{(t)},B=C[\delta_{Y_i,X_i}]_{i\in \mathcal{I}}\; (\text{resp. } B=C[\delta_{Y_i,X_i},\delta_{Y_i,X_i}^t]_{i\in \mathcal{I}/\sim}$), there is an identification,
\begin{equation}
\Omega (\slfrac{B}{C}) \equiv B^{(Y_i,X_i)_{i\in \mathcal{I}}}_{(t)}.
\end{equation}
Thus the left adjoint functor of (\ref{7.6.10}), $B\mapsto \Omega(B/ C)^A$, takes (acyclic) cofibrations to (acyclic) cofibrations, and is a left Quillen functor. \vspace{3mm} \\

For $\varphi\in \FR^{(t)}(C,A)$, Quillen's cotangent bundle is the element of the derived category of $\Amod^{(t)}$,$\mathbb{D}(A\text{-mod}^{(t)})=HO(\Delta (\Amod^{(t)}))$, given by,
\begin{equation}
\mathbb{L}\Omega(\slfrac{A}{C}):=\Omega(\slfrac{P_C(A)}{C})^A\in \mathbb{D}(A\text{-mod}^{(t)}).
\end{equation}
where $P_C(A)\rightarrow A$ is a cofibrant replacement of $A$.\vspace{1.5mm}\\
$\exists$ standard resolution associative to the pair of adjoint
functors,
\begin{equation}
\xymatrix{C\diagdown \FR^{(t)} \ar@/^/[d]^{U} \\
 (\slfrac{Set}{\FF\times \FF})^{(t)}\ar@/^/[u]^{F}}
\end{equation}

\begin{equation}
\xymatrix@R=0.5em{P_C(A)\equiv \dots FUFU(A) \ar@<1ex> [r] \ar@<-1ex>[r]&FU(A)=C[\underset{=}{a};a\in \underset{Y,X}{\coprod} A_{Y,X} ]\ar@{->>}[r]\ar[l]& A \\
&\underset{=}{a}\ar@{|->}[r]&a}
\end{equation}
given by the unit and counit of the adjunction.

We have,
\begin{equation}
\mathbb{L}\Omega(\slfrac{\varinjlim B_i}{\varinjlim C_i})\equiv \varinjlim \mathbb{L} \Omega (\slfrac{B_i}{C_i}).
\end{equation}
Similarly, for a map of simplicial objects $C_{\bullet}\to A_{\bullet}$ in $\Delta (\FR^{(t)})$, applying $P_{C_n}(A_n)$ at each dimension $n$, we obtain a bi- simplicial object $P_{C_n}(A_n)_m=(F_{C_n}U)^m(A_n)$, and taking the diagonal object $(n=m)$ we obtain the resolution $P_C^{\nabla}(A_{\bullet})\to A_{\bullet}$, cf. \cite{I}.\\
Given $\FR^{(t)}$- homomorphisms
\begin{equation}
\xymatrix{C\ar[r]&B\ar[r]&A \\ C\ar[r]\ar@{=}[u]&Q=P_C(B)\ar[r]\ar[u]&P_{Q}^{\nabla}(A)\ar[u]},
\end{equation}
we have compatible resolutions, hence an exact (also on the left!) sequence of $A$- modules,
\begin{equation}
0\to \Omega(Q/C)^A\to \Omega(P_Q^{\nabla}(A)/C)^A\to \Omega(P_Q^{\nabla}(A)/Q)^A\to 0.
\end{equation}
This can be interpreted (see \cite{I}) as the exact triangle in $\mathbb{D}(\Amod)$,
\begin{equation}
\mathbb{L}\Omega(B/C)^A\to \mathbb{L}\Omega(A/C)\to \mathbb{L}\Omega(A/B),
\end{equation}
or as the long exact sequence of $A$- modules

\begin{equation}
\begin{tikzpicture}[>=angle 90,scale=0.3,text height=1.5ex, text depth=0.25ex]

\node (a4) at (45,13) {$0$};
\node (a3) at (35,13) {$\Omega(A/B)$};
\node (a2) at (25,13) {$\Omega(A/C)$};
\node (a1) at (15,13) {$\Omega(B/C)^A$};

\node (b3) at (33,9) {$L_1\Omega(A/B)$};
\node (b2) at (23,9) {$L_1\Omega(A/C)$};
\node (b1) at (13,9) {$L_1\Omega(B/C)^A$};

\node (c3) at (30,5) {$L_2\Omega(A/B)$};
\node (c2) at (20,5) {$L_2\Omega(A/C)$};
\node (c1) at (10,5) {$L_2\Omega(B/C)^A$};
\node (c0) at (1,5) {$\dots$};

\draw[->,font=\scriptsize]
(a1) edge (a2)
(a2) edge (a3)
(a3) edge (a4)

(b1) edge (b2)
(b2) edge (b3)

(c0) edge (c1)
(c1) edge (c2)
(c2) edge (c3) ;

\draw[->]

(c3) edge[out=0,in=180] (b1)
(b3) edge[out=0,in=180] (a1) ;

\end{tikzpicture}
\end{equation} \vspace{3mm}

Let $X$ be a finite (Krull) dimensional, noetherian, topological space. \\
Define a model structure on $\Delta (\mathcal{O}_X\text{-mod}^{(t)})$, $(X,\mathcal{O}_X)\in \FR^{(t)}\mathcal{S}p$,
\begin{equation}
\begin{matrix} \text{Fibrations:} &\mathcal{F}=& \begin{Bmatrix}
   \varphi:m\rightarrow m', \;\;\; \forall X,Y\in \FF,\forall V\subseteq U\subseteq X,\\
m(U)_{Y,X}\rightarrow m'(U)_{Y,X}\underset{m'(V)_{Y,X}}{\prod} m(V)_{Y,X}\in \mathcal{F}_{Set_0} \end{Bmatrix}, \\ \text{Weak\; equivalences:} &\mathcal{W}=& \{\varphi: m\rightarrow m',\;\;\; \forall X,Y\in \FF, \forall \mathfrak{p}\in X, (m_{\mathfrak{p}})_{Y,X}\rightarrow (m_{\mathfrak{p}}')_{Y,X}\in \mathcal{W}_{Set_0}\} \\
\text{Cofibrations:}& \mathcal{C}=& L(\mathcal{W}\cap \mathcal{F}).\hspace{75mm}  \end{matrix}
\end{equation} \vspace{3mm}

Define a model structure on $\Delta (\mathcal{O}_X\setminus (\slfrac{\FR^{(t)}}{X}))$, the category of simplicial sheaves of $\FR^{(t)}$ over $X$, together with a map $\mathcal{O}_X\rightarrow B$, with arrows $B\rightarrow B'$ are homomorphisms of simplicial $\FR^{(t)}$ -sheaves over $\mathcal{O}_X$.
\begin{equation}
\begin{matrix} \text{Fibrations:} & \\ &\mathcal{F}= &\begin{Bmatrix}
   \varphi: B\rightarrow B', \text{such that } \forall X,Y\in \FF,\forall V\subseteq U\subseteq X,\\
\{B(U)_{Y,X}\rightarrow B'(U)_{Y,X}\underset{B'(V)_{Y,X}}{\prod} B(V)_{Y,X}\}  \in \mathcal{F}_{Set_0} \end{Bmatrix}, \\ \text{Weak\; equivalences:} & \\ &\mathcal{W}=&\{\varphi: B\rightarrow B',\;\;\; \forall X,Y\in \FF, \forall \mathfrak{p}\in X, \{ (B_{\mathfrak{p}})_{Y,X}\rightarrow (B_{\mathfrak{p}}')_{Y,X}\}  \in \mathcal{W}_{Set_0}\} \\
\text{Cofibrations:}& \\ & \mathcal{C}=&L(\mathcal{W}\cap \mathcal{F}).\hspace{75mm} \end{matrix}
\end{equation}
That these constitute a Quillen model structure on $\Delta (\mathcal{O}_X\setminus (\slfrac{\FR^{(t)}}{X}))$, and on $\Delta (\mathcal{O}_X\text{-mod}^{(t)})$, follows as in \cite{Q67} (theorem 4.II,$\S$ 4), with the aid of the Brown- Gerstein lemma \cite{BG}: For $X$ finite dimensional, noetherian,
\begin{equation}
\mathcal{F}\cap \mathcal{W}\subseteq \left \{ \begin{matrix} \varphi:B\to B', \text{ such that }\forall X,Y\in \FF, \forall V\subseteq U\subseteq X,\\ \{B(U)_{Y,X}\to B'(U)_{Y,X}\underset{B'(V)_{Y,X}}{\prod} B(V)_{Y,X}\}\in \mathcal{W}_{set_0} \end{matrix}\right \} \text{ the \underline{global} weak equivalences.}
\end{equation}
For a map $f\in \FS^{(t)}(X,Y)$, the cotangent bundle $\mathbb{L} \Omega (\slfrac{X}{Y})$ is the element of the derived category of $\mathcal{O}_X$- modules, $\mathbb{D}(\mathcal{O}_X\text{- mod})=HO(\Delta (\mathcal{O}_X\text{- mod}))$,
\begin{equation}
\mathbb{L}\Omega (\slfrac{X}{Y}):= \Omega (\slfrac{P_{f^*\mathcal{O}_Y}(\mathcal{O}_X)}{f^*\mathcal{O}_Y}).
\end{equation}
Given another map $g\in \FS^{(t)}(Y,Z)$, we have the exact triangle,
\begin{equation}
\mathbb{L}\Omega(\slfrac{Y}{Z})^X\rightarrow \mathbb{L}\Omega(\slfrac{X}{Z})\rightarrow \mathbb{L}\Omega(\slfrac{X}{Y}).
\end{equation}
\numberwithin{equation}{section}
\let\cleardoublepage\clearpage
\part{Generalized Rings}
\let\cleardoublepage\clearpage
\chapter{Generalized Rings}
\let\cleardoublepage\clearpage
\bigskip
\section{Definitions}
\smallbreak
For $Y\in\FF$ and $X\subseteq Y$, define the operation of contracting $X$ to a point $*_X$ :
\begin{equation}
Y/X:=(Y\setminus X)\amalg \{*_X\}
\end{equation}
equipped with a map,
\begin{equation}
\pi:Y\sur Y/X, \hspace{3mm} \pi(X)=\{*_X\}.
\end{equation}
We have an ''inverse'' operation: for $X,Z\in \FF$ and  $z_0\in Z$, , set
\begin{equation}
Z\underset{z_0}{\lhd} X:= (Z\setminus \{z_0\})\amalg
X,\end{equation} equipped with a map,
\begin{equation}
\pi: Z\underset{z_0}{\lhd} X\sur Z,\hspace{3mm} \pi(X)=z_0.
\end{equation}
We have the following two trivial identifications:
\begin{equation}
(Y/X)\underset{*_X}{\lhd} X\equiv Y\hspace{3mm} \text{and} \hspace{3mm} (Z\underset{z_0}{\lhd}X)/ X\equiv Z.
\end{equation}
\defin{8.1.1}\textit{A generalized ring is a functor $A\in (Set_0)^{\FF}$, such that $A_{[0]}=\{0\}$ with the two operations:
}
$$\begin{tikzcd}
\textbf{multiplication:} \text{ for } z_0\in Z\in \FF, &&\textbf{contraction:} \text{ for } X\subseteq Y\in \FF,&\\
A_Z\times A_X\arrow{r} & A_{Z\underset{z_0}{\lhd} X}&A_Y\times A_X\arrow{r}&A_{Y/X} \\
c,a\arrow[mapsto]{r} & (c\underset{z_0}{\lhd} a)&b,a\arrow[mapsto]{r}&(b\sslash a)
\end{tikzcd}$$
\textit{such that the following axioms holds:} \\
\begin{enumerate}[O.]
\item \textbf{Functoraility of operations.}\\
For $\varphi\in \FF(Z,Z'), \varphi(z_0)=z_0'$ and $\psi\in \FF(X,X')$, we have a map defined in the obvious way $\varphi\lhd \psi \in \FF(Z\underset{z_0}{\lhd} X, Z'\underset{z_0'}{\lhd} X')$ and we require the following diagram to commute:
\begin{equation}
\begin{tikzcd}
A_Z\times A_X\ar{r}\ar{d}{\varphi\times \psi}& \;\;A_{Z\underset{z_0}{\lhd} X}\ar{d}{\varphi\lhd\psi}\\
A_{Z'}\times A_{X'}\ar{r}& A_{Z'\underset{z_0'}{\lhd} X'}
\end{tikzcd}
\end{equation}
that is,
\begin{equation}
\varphi(a_Z)\underset{z_0'}{\lhd}\psi(a_X)=\varphi\lhd \psi (a_Z\underset{z_0}{\lhd} a_X).
\end{equation}
Secondly, for $\varphi\in \FF(Y/X,Y'/X'), \varphi(*_X)=*_{X'}$ and $\psi\in \FF(X,X')$, we have a diagram:\\
\begin{equation}
\xymatrix{(A_Y\times A_X)\ar@/_/[d]_{\varphi\lhd\psi}\ar[rrr]&&&A_{Y/X}\ar[d]^{\varphi}\\
(A_{Y'}\times A_{X'})\ar@/_/[u]_{\psi^t}\ar[rrr]&&&A_{Y'/X'}}
\end{equation}
and we require,
\begin{equation}
\varphi\lhd\psi(a_Y)\sslash a_{X'}=\varphi(a_Y\sslash \psi^t(a_{X'}))
\end{equation}
e.g. for the contraction,
\begin{equation}
\label{eq8.1.9}
\xymatrix@R=1mm{ A_X\times A_X\ar[r]& A_{[1]} \\
(a_1,a_2)\ar@{|->}[r]& a_1\sslash a_2
}
\end{equation}
and for $\varphi\in \FF(X,X')$ we have,
\begin{equation}
\varphi(a_X)\sslash a_{X'}=a_X\sslash \varphi^t(a_{X'})
\end{equation}
(up-to the identifications, $A_{*_X}=A_{[1]}=A_{*_{X'}}$).
\end{enumerate}
And we have the four \emph{zero- axioms}:
\begin{enumerate}[(i)]
\item $\;0_Z\underset{z_0}{\lhd} a_X= 0_{Z\underset{z_0}{\lhd} X},$
\item $\;0_Y\sslash a_X=0_{Y/X},$
\item $\;a_Z\underset{z_0}{\lhd} 0_X=\varphi(a_Z) \;\; \text{with} \;\; (\varphi: \overbrace{Z\setminus \{z_0\}}^{\subseteq Z}\xrightarrow{=} \overbrace{Z\setminus \{z_0\}}^{\subseteq Z \underset{z_0}{\lhd} X})\in \FF_{Z\underset{z_0}{\lhd} X,Z},$
Also for $a_Z\in A_{Z\setminus z_0}\subseteq A_Z$, viewed as an element of $A_Z$, we have for any $a_X\in A_X$,
$$a_Z\underset{z_0}{\lhd} a_X=\varphi(a_Z),$$
for $\varphi\in \FF(Z\setminus \{z_0\},Z\underset{z_0}{\lhd}X)$ defined as above.
\item $\;a_Y\sslash 0_X=\varphi(a_Y) \;\; \text{with} \;\; (\varphi: \overbrace{Y\setminus X}^{\subseteq Y}\xrightarrow{=} \overbrace{Y\setminus X}^{\subseteq Y/X})\in \FF_{Y/X,Y}.$
\end{enumerate} \vspace{3mm}
\begin{enumerate}[I.]
\item \textbf{Disjointness Axiom.}\\
Suppose $X_0\amalg X_1\subseteq Y\in \FF$ is a \emph{disjoint} union of subsets. Then for any $b\in A_Y, a_i\in A_{X_i}$,
\begin{equation}
(b\sslash a_1)\sslash a_0\equiv (b\sslash a_0)\sslash a_1\;.
\end{equation} \vspace{0.3cm}\\
Schematically we have, \\
\begin{equation}
\begin{tikzpicture}[scale=0.3]

\tikzstyle{no}=[circle,draw,fill=black!100,inner sep=0pt, minimum width=4pt]
\tikzstyle{no2}=[circle,draw,fill=black!10,inner sep=0pt, minimum width=4pt]


\node[no] (n1) at (15,12) {};
\node[no] (n2) at (12,9)  {};
\node[no] (n3) at (14,9)  {};
\node[no] (n4) at (16,9) {};
\node[no] (n5) at (18,9)  {};
\node[no2] (n6) at (20,9){};
\node[no] (n7) at (12,5) {};
\node[no] (n8) at (14,5){};
\node[no] (n9) at (16,5) {};
\node[no] (n10) at (18,5) {};
\node[no] (n11) at (17,3) {};
\node[no] (n12) at (13,3){};
\node at (13,2) {$a_0$};
\node at (17,2) {$a_1$};

\draw  (12,8.4)-- (12,8) -- (14,8) --(14,8.4);
\draw  (16,8.4)-- (16,8) -- (18,8) --(18,8.4);

\node at (15,13) {$b$};
\node at (13,7) {$X_0$};
\node at (17,7) {$X_1$};

\foreach \from/\to in {n2/n1,n3/n1,n4/n1,n5/n1,n6/n1,n11/n9,n11/n10,n12/n7,n12/n8}
\draw[thick] (\from) -- (\to);
\foreach \from/\to in {n2/n7,n3/n8,n4/n9,n5/n10}
\draw[dashed] (\from) -- (\to);
\end{tikzpicture}
\end{equation}
\vspace{0.1mm} \\
and we can contract from $b$, $a_0$ first and then $a_1$ or the other way round.

\begin{equation}
\begin{tikzcd}
&X/X_0\ar{rd}&\\
Y\ar{rd}\ar{ru}&&\hspace{7mm} \begin{matrix}(Y/X_0)/X_1 \\ {\begin{turn}{90}
$\equiv$
\end{turn}}\\(Y/X_1)/X_0 \end{matrix}&:=\begin{matrix} Y\sslash\{X_0,X_1\}\\ {\begin{turn}{90}
$=$
\end{turn}}\\ (Y\setminus (\{X_0\amalg X_1\}))\amalg \{*_{X_0},*_{X_1}\}\end{matrix}\\
&X/X_1\ar{ru}&
\end{tikzcd}
\end{equation}
Generally, we denote such multiple contractions as,
\begin{equation}
b\sslash (a_i)\in A_{Y\slash  \{X_i\}_{\tiny I}}, \;\; b\in A_{Y},a_i\in A_{X_i}, \underset{I}{\amalg} X_i\subseteq Y.
\end{equation}
\defin{8.1.2} For a map of sets $f\in Set_{\bullet}(Y,Z)$, ($Y,Z\in \FF$), we put
\begin{equation}
A_f:= \underset{z\in Z}{\prod} A_{f^{-1}(z)}.
\end{equation}
We have the (multiple) contraction \\
\begin{equation}
\xymatrix@R=1mm{A_Y\times A_f\ar[r]& A_Z \\
a_Y,a^f\ar@{|->}[r]& a_Y\sslash a^f.}
\end{equation}
for $a_Y\in A_Y, a^f=(a^f_z),a^f_z\in A_{f^{-1}(z)}$, $z\in Z$.
\item \textbf{Disjointness Axiom.} \\
Suppose $z_0\neq z_1\in Z$ and $X_0,X_1\in \FF$. We have for $c\in A_Z$ and $a_i\in A_{X_i}$,
\begin{equation}
(c\underset{z_0}{\lhd} a_0)\underset{z_1}{\lhd} a_1=(c\underset{z_1}{\lhd} a_1)\underset{z_0}{\lhd}a_0,
\end{equation} \vspace{1.5mm}\\
Schematically we have, \\
\begin{equation}
\begin{tikzpicture}[scale=0.3]

\tikzstyle{no}=[circle,draw,fill=black!100,inner sep=0pt, minimum width=4pt]
\tikzstyle{no2}=[circle,draw,fill=black!10,inner sep=0pt, minimum width=4pt]


\node[no] (n1) at (15,12) {};
\node[no] (n2) at (12,9)  {};
\node at (11,8.2) {$z_0$};
\node[no2] (n3) at (14,9)  {};
\node[no] (n4) at (16,9) {};
\node at (17,8.2) {$z_1$};
\node[no2] (n5) at (18,9)  {};
\node[no2] (n6) at (20,9){};
\node[no] (n7) at (11,4) {};
\node[no] (n8) at (13,4){};
\node[no] (n9) at (15,4) {};
\node[no] (n10) at (17,4) {};
\node[no] (n11) at (16,6) {};
\node[no] (n12) at (12,6){};
\node at (12,3) {$a_0$};
\node at (16,3) {$a_1$};
\node at (15,13) {$b$};

\foreach \from/\to in {n2/n1,n3/n1,n4/n1,n5/n1,n6/n1,n11/n9,n11/n10,n12/n7,n12/n8}
\draw[thick] (\from) -- (\to);
\foreach \from/\to in {n11/n4,n12/n2}
\draw[dashed] (\from) -- (\to);
\end{tikzpicture}
\end{equation}

and a corresponding diagram of sets:
\begin{equation}
\begin{tikzcd}
&Z\underset{z_0}{\lhd}X_0\ar[two heads]{rd}& \\
Z\underset{z_i}{\lhd} X_i:=\begin{matrix}(Z\underset{z_0}{\lhd}X_0)\underset{z_1}{\lhd} X_1\\ {\begin{turn}{90}
$\equiv$
\end{turn}}\ar[two heads]{ru}\ar[two heads]{rd}\hspace{2mm}&& \\
(Z\underset{z_1}{\lhd}X_1)\underset{z_0}{\lhd} X_0 \end{matrix}\hspace{10mm}&& \hspace{5mm}Z\ni z_0,z_1\\
&Z\underset{z_1}{\lhd}X_1\ar[two heads]{ru}&
\end{tikzcd}
\end{equation}
Generally we denote such a multiple multiplications, indexed by a subset $Z_0\subseteq Z$ with $\{X_z\}_{z\in Z_0}$ by,
\begin{equation}
c\underset{z\in Z_0}{\lhd} (a_z).
\end{equation}
for $c\in A_Z$ and $\{a_z\in A_{X_z}\}_{z\in Z_0}$. \\
\defin{8.1.3} For a map of sets $f\in Set_{\bullet}(Y,Z)$, ($Y,Z\in \FF$), we get (multiple) multiplication
\begin{equation}
\xymatrix@R=0.5em{A_Z\times A_f\ar[r]& A_Y \\
a_Z,a^f\ar@{|->}[r]& a_Z\lhd a^f}
\end{equation}

\item \textbf{Disjointness Axiom.}\\
For $b\in A_Z$, $a_i\in A_{X_i}$ and $X_0\subseteq Z, z_1\in Z\setminus X_0$, we have,
\begin{equation}
(b\sslash a_0)\underset{z_1}{\lhd} a_1\equiv (b\underset{z_1}{\lhd} a_1)\sslash a_0.
\end{equation}
Schematically,

\begin{equation}
\begin{tikzpicture}[scale=0.3]

\tikzstyle{no}=[circle,draw,fill=black!100,inner sep=0pt, minimum width=4pt]
\tikzstyle{no2}=[circle,draw,fill=black!10,inner sep=0pt, minimum width=4pt]


\node[no] (n1) at (15,12) {};
\node[no2] (n2) at (12,9)  {};
\node at (11,8.2) {$z_1$};
\node[no] (n3) at (14,9)  {};
\node[no2] (n4) at (16,9) {};
\node at (17,8.2) {$X_0$};
\node[no2] (n5) at (18,9)  {};
\node[no] (n6) at (20,9){};
\node[no] (n7) at (11,4) {};
\node[no] (n8) at (13,4){};
\node[no] (n9) at (16,6) {};
\node[no] (n10) at (18,6) {};
\node[no] (n11) at (17,4) {};
\node[no] (n12) at (12,6){};
\node at (12,3) {$a_1$};
\node at (16,3) {$a_0$};
\node at (15,13) {$b$};

\foreach \from/\to in {n2/n1,n3/n1,n4/n1,n5/n1,n6/n1,n11/n9,n11/n10,n12/n7,n12/n8}
\draw[thick] (\from) -- (\to);
\foreach \from/\to in {n2/n12,n9/n4,n10/n5}
\draw[dashed] (\from) -- (\to);
\end{tikzpicture}
\end{equation}

\item \textbf{Associativity Axiom.}\\
For $b\in A_Z, \; a_i\in A_{X_i},\; z_0\in Z,\; x_0\in X_0, $
\begin{equation}
(b\underset{z_0}{\lhd}a_0)\underset{x_0}{\lhd} a_1=b\underset{z_0}{\lhd}(a_0\underset{x_0}{\lhd}a_1).
\end{equation}

Schematically,
\begin{equation}
\begin{tikzpicture}[scale=0.3]

\tikzstyle{no}=[circle,draw,fill=black!100,inner sep=0pt, minimum width=4pt]
\tikzstyle{no2}=[circle,draw,fill=black!10,inner sep=0pt, minimum width=4pt]


\node[no] (n1) at (14,12) {};
\node[no2] (n2) at (12,9)  {};
\node at (11,9) {$z_1$};
\node[no] (n3) at (14,9)  {};
\node[no] (n4) at (16,9){};
\node[no] (n7) at (11,4) {};
\node[no2] (n8) at (13,4){};
\node at (14.1,4) {$x_0$};
\node[no] (n9) at (12,0) {};
\node[no] (n10) at (14,0) {};
\node[no] (n11) at (13,2) {};
\node[no] (n12) at (12,6){};
\node at (10.5,5.5) {$a_1$};
\node at (11.3,1.5) {$a_0$};
\node at (14,13) {$b$};

\foreach \from/\to in {n2/n1,n3/n1,n4/n1,n11/n9,n11/n10,n12/n7,n12/n8}
\draw[thick] (\from) -- (\to);
\foreach \from/\to in {n2/n12,n11/n8}
\draw[dashed] (\from) -- (\to);
\end{tikzpicture}
\end{equation}

\item \textbf{Left adjunction Axiom.}\\
For $X_0\subseteq X_1\subseteq Y$, $b\in A_Y, a_0\in A_{X_0}, a_1\in A_{X_1/X_0}$  we have,
\begin{equation}
(b\sslash a_0)\sslash a_1\equiv b\sslash (a_1\underset{*_{X_0}}{\lhd}a_0).
\end{equation}

Schematically seen,

\begin{equation}
\begin{tikzpicture}[scale=0.3]

\tikzstyle{no}=[circle,draw,fill=black!100,inner sep=0pt, minimum width=4pt]
\tikzstyle{no2}=[circle,draw,fill=black!10,inner sep=0pt, minimum width=4pt]


\node[no] (n1) at (8,12) {};
\node[no2] (n2) at (5,9)  {};
\node[no2] (n3) at (7,9)  {};
\node[no] (n4) at (9,9) {};
\node[no] (n6) at (11,9){};
\node[no] (n7) at (6,7) {};
\node[no] (n9) at (6,5) {};
\node[no] (n10) at (9,5) {};
\node[no] (n11) at (7.5,4) {};
\node[no2] (n12) at (5,9){};
\node at (3,8.1) {$a_0$};
\node at (5,4.7) {$a_1$};
\node at (6.5,12.8) {$b$};
\node at (14,8.8) {$\equiv$};
\node[no] (m1) at (20,12) {};
\node[no] (m2) at (17,9)  {};
\node[no] (m3) at (19,9)  {};
\node[no] (m4) at (21,9) {};
\node[no] (m5) at (23,9){};
\node[no] (m6) at (17,7)  {};
\node[no] (m7) at (19,7)  {};
\node[no2] (m8) at (18,5)  {};
\node[no] (m9) at (21,5)  {};
\node[no] (m10) at (19.5,4)  {};

\foreach \from/\to in {n2/n1,n3/n1,n4/n1,n6/n1,n11/n9,n11/n10,n12/n7,n3/n7}
\draw[thick] (\from) -- (\to);
\foreach \from/\to in {n7/n9,n10/n4}
\draw[dashed] (\from) -- (\to);

\foreach \from/\to in {m1/m2,m1/m3,m1/m4,m1/m5,m6/m8,m7/m8,m8/m10,m9/m10}
\draw[thick] (\from) -- (\to);
\foreach \from/\to in {m2/m6,m3/m7,m9/m4}
\draw[dashed] (\from) -- (\to);

\end{tikzpicture}
\end{equation}

with a corresponding identity of sets:
\begin{equation}
(Y/X_0)/(X_1/X_0)\equiv Y/X_1.
\end{equation}

\item \textbf{Right adjunction Axiom.}\\
For $X\subseteq Y_0\subseteq Y_1, b\in A_{Y_0}, a\in A_X,c\in A_{Y_1/X}$,  we have,
\begin{equation}
c\sslash (b\sslash a)\equiv (c\underset{*_X}{\lhd} a)\sslash b
\end{equation}

Schematically seen,

\begin{equation}
\begin{tikzpicture}[scale=0.3]

\tikzstyle{no}=[circle,draw,fill=black!100,inner sep=0pt, minimum width=4pt]
\tikzstyle{no2}=[circle,draw,fill=black!10,inner sep=0pt, minimum width=4pt]


\node[no] (n1) at (22,12) {};
\node[no2] (n2) at (19,9)  {};
\node at (18,9.2) {$*_X$};
\node[no] (n3) at (21,9)  {};
\node[no] (n4) at (23,9) {};
\node[no] (n5) at (25,9)  {};
\node[no] (n6) at (27,9){};
\node[no] (n7) at (16,7) {};
\node[no] (n8) at (18,7){};
\node[no] (n13) at (23,6){};
\node[no] (n9) at (16,6) {};
\node[no] (n10) at (18,6) {};
\node[no] (n11) at (17,4) {};
\node[no2] (n12) at (19,9){};
\node at (16,8.8) {$a$};
\node at (15.5,5) {$b$};
\node at (22,13) {$c$};

\node at (13,8.8) {$\equiv$};

\node[no] (m1) at (5,12) {};
\node[no] (m2) at (4,9)  {};
\node[no] (m3) at (6,9)  {};
\node[no] (m4) at (8,9) {};
\node[no] (m5) at (10,9)  {};
\node[no] (m6) at (1,9){};
\node[no] (m7) at (1,7) {};
\node[no2] (m8) at (0,5.5){};
\node[no2] (m9) at (2,5.5) {};
\node[no] (m10) at (1,4) {};
\node[no] (m11) at (6,5) {};

\foreach \from/\to in {n2/n1,n3/n1,n4/n1,n5/n1,n6/n1,n11/n9,n11/n10,n12/n7,n12/n8,n13/n11}
\draw[thick] (\from) -- (\to);
\foreach \from/\to in {n7/n9,n8/n10,n13/n4}
\draw[dashed] (\from) -- (\to);

\foreach \from/\to in {m2/m1,m3/m1,m4/m1,m5/m1,m6/m1,m7/m8,m7/m9,m8/m10,m9/m10,m10/m11}
\draw[thick] (\from) -- (\to);
\foreach \from/\to in {m6/m7,m11/m3}
\draw[dashed] (\from) -- (\to);

\end{tikzpicture}
\end{equation}

with a corresponding identity of sets:
\begin{equation}
(Y_1/X)/(Y_0/X)\equiv Y_1/Y_0.
\end{equation}

\item \textbf{Left linear Axiom.}\\
For $a\in A_X,b\in A_Y, c\in A_Z$ where $X\subset Y, z_0\in Z$, we have,
\begin{equation}
c\underset{z_0}{\lhd}(b\sslash a)\equiv (c\underset{z_0}{\lhd}b)\sslash a
\end{equation}
Schematically seen,

\begin{equation}
\begin{tikzpicture}[scale=0.3]

\tikzstyle{no}=[circle,draw,fill=black!100,inner sep=0pt, minimum width=4pt]
\tikzstyle{no2}=[circle,draw,fill=black!10,inner sep=0pt, minimum width=4pt]


\node[no] (n1) at (8,12) {};
\node[no2] (n2) at (5,9)  {};
\node at (4,9) {$z_0$};
\node[no] (n3) at (7,9)  {};
\node[no] (n4) at (9,9) {};
\node[no] (n5) at (11,9)  {};
\node[no] (n6) at (13,9){};
\node[no] (n7) at (2,7) {};
\node[no] (n8) at (4,7){};
\node[no] (n13) at (6,7){};
\node[no] (n9) at (2,4) {};
\node[no] (n10) at (4,4) {};
\node[no] (n11) at (3,2) {};
\node[no2] (n12) at (5,9){};
\node at (2,8.8) {$b$};
\node at (1.5,3) {$a$};
\node at (8,13) {$c$};

\node at (16,9) {$\equiv$};

\node[no] (m1) at (22,12) {};
\node[no] (m2) at (19,9)  {};
\node[no] (m3) at (21,9)  {};
\node[no] (m4) at (23,9) {};
\node[no] (m5) at (25,9)  {};
\node[no] (m6) at (27,9){};
\node[no] (m7) at (19,7){};
\node[no2] (m8) at (18,5){};
\node[no2] (m9) at (20,5){};
\node[no] (m10) at (22,5){};
\node[no] (m11) at (19,3){};

\foreach \from/\to in {n2/n1,n3/n1,n4/n1,n5/n1,n6/n1,n11/n9,n11/n10,n12/n7,n12/n8,n13/n2}
\draw[thick] (\from) -- (\to);

\foreach \from/\to in {n9/n7,n10/n8}
\draw[dashed] (\from) -- (\to);

\foreach \from/\to in {m1/m2,m1/m3,m1/m4,m1/m5,m1/m6,m7/m8,m7/m9,m7/m10,m8/m11,m9/m11}
\draw[thick] (\from) -- (\to);

\foreach \from/\to in {m7/m2}
\draw[dashed] (\from) -- (\to);

\end{tikzpicture}
\end{equation}

\item \textbf{Unit Axiom.}\\
We have an element $1\in A_1$. We obtain for any $x\in X\in \FF$,
\begin{equation}
1_x\equiv1_x^X\in A_X, 1_x^X=j_x^X(1)
\end{equation}
for $j_x^X\in \FF_{X,1}, j_x^X(1)=x$. \\
We require that the following identities holds: for any $a\in A_X$, and any $x\in X$,
\begin{equation}
\begin{split}
1\underset{1}{\lhd}a\equiv a,  \\
a\underset{x}{\lhd}1\equiv a,  \\
a\sslash 1_x \equiv a.
\end{split}
\end{equation}
with a corresponding identities of sets:
\begin{equation}
\begin{split}
[1]\underset{1}{\lhd}X\equiv X,  \\
X\underset{x}{\lhd}[1]\equiv X,  \\
X/\{x\} \equiv X.
\end{split}
\end{equation}

Schematically it can be seen as follows:

\begin{equation}
\begin{tikzpicture}[scale=0.3]

\tikzstyle{no}=[circle,draw,fill=black!100,inner sep=0pt, minimum width=4pt]
\tikzstyle{no2}=[circle,draw,fill=black!10,inner sep=0pt, minimum width=4pt]


\node[no] (b1a) at (-1,11) {};
\node at (-1.5,12) {$1$};
\node[no] (b2a) at (-3,9)  {};
\node[no] (b3a) at (-1,9)  {};
\node[no] (b4a) at (1,9) {};
\node[no] (n0) at (-1,13) {};

\node at (2,10) {$\equiv$};

\node[no] (b1b) at (6,11) {};
\node[no] (b2b) at (4,9)  {};
\node[no2] (b3b) at (6,9)  {};
\node at (5.3,9)  {$\tiny x$};
\node[no] (b4b) at (8,9) {};
\node[no] (n2) at (6,7) {};
\node at (5.5,8) {$1$};

\node at (9.5,10) {$\equiv$};

\node[no] (b1c) at (13,11) {};
\node[no] (b2c) at (11,9)  {};
\node[no2] (b3c) at (13,9)  {};
\node at (12.3,8.8)  {$\tiny x$};
\node[no] (b4c) at (15,9) {};
\node at (12.5,8) {$1$};
\node[no] (n4) at (13,7) {};

\node at (16,10) {$\equiv$};

\node[no] (b1d) at (20,11) {};
\node[no] (b2d) at (18,9)  {};
\node[no] (b3d) at (20,9)  {};
\node[no] (b4d) at (22,9) {};

\draw[dashed]  (-3,6)-- (-3,5) -- (1,5) --(1,6);
\node at (-1,4) {$1\lhd a$};
\draw[dashed]  (4,6)-- (4,5) -- (8,5) --(8,6);
\node at (6,4) {$a\underset{x}{\lhd} 1$};
\draw[dashed]  (11,6)-- (11,5) -- (15,5) --(15,6);
\node at (13,4) {$a/1_x$};
\draw[dashed]  (18,6)-- (18,5) -- (22,5) --(22,6);
\node at (20,4) {$a$};

\foreach \from/\to in
 {b1a/b2a,b1a/b3a,b1a/b4a,b1b/b2b,b1b/b3b,b1b/b4b,b1c/b2c,b1c/b3c,b1c/b4c,b1d/b2d,b1d/b3d,b1d/b4d,n0/b1a,b3b/n2,
b3c/n4}
\draw[thick] (\from) -- (\to);
\end{tikzpicture}
\end{equation}



\item \textbf{Commutativity.}
\end{enumerate}
\defin{8.1.4}
\textit{We say a generalized ring $A$ is \emph{commutative} if it satisfies for any $b\in A_Y, X_0\subseteq Y, a_i\in A_{X_i}$,
\begin{equation}
(b\sslash a_0)\underset{*_{X_0}}{\lhd}a_1\equiv \underbrace{(b\underset{x_0\in X_0}{\lhd}(a_1))}_{\in A_{(X_0\times X_1)\amalg (Y\setminus X_0)}}\sslash (a_0)_{x_1\in X_1}\label{comm_identity}
\end{equation}
(where the $x_1$' \ the copy of $a_0$ is attached through the indices $X_0\times \{x_1\}$). \vspace{1.5mm} \\
We also call this identity ''\emph{Right- linear}''.}\\
Schematically seen,

\begin{equation}
\begin{tikzpicture}[scale=0.3]

\tikzstyle{no}=[circle,draw,fill=black!100,inner sep=0pt, minimum width=4pt]
\tikzstyle{no2}=[circle,draw,fill=black!10,inner sep=0pt, minimum width=4pt]


\node[no] (b1a) at (-1,11){};
\node at (-1.5,12) {$b$};
\node[no] (b2a) at (-3,9){};
\node[no] (b3a) at (-1,9){};
\node[no] (b4a) at (1,9){};
\node[no] (n1) at (-3,5){};
\node[no] (n2) at (-1,5){};
\node[no] (n3) at (-2,7){};

\node at (2.3,10) {$\equiv$};

\node[no] (b1b) at (6,11) {};
\node at (6,12) {$b$};
\node[no] (b2b) at (4,9)  {};
\node[no] (b3b) at (6,9)  {};
\node[no] (b4b) at (8,9) {};
\node[no] (m1) at (3.5,7){};
\node[no] (m2) at (4.5,7){};
\node[no] (m3) at (5.5,7){};
\node[no] (m4) at (6.5,7){};
\node[no] (m5) at (4,5){};
\node[no] (m6) at (6,5){};

\foreach \from/\to in {b1a/b2a,b1a/b3a,b1a/b4a,b1b/b2b,b1b/b3b,b1b/b4b,n3/b2a,n3/b3a,n1/n3,n2/n3,b2b/m1,b2b/m3, b3b/m2,b3b/m4,m1/m5,m2/m5,m3/m6,m4/m6}
\draw[thick](\from) -- (\to);
\end{tikzpicture}
\end{equation}

\noindent We are mainly interested in the commutative generalized
rings. We are \emph{not} interested in the vertical-commutative
generalized rings of the following:

\defin{8.1.5}
\textit{We say a generalized ring $A$ is \emph{vertical- commutative} if it satisfies, for any $a\in A_X, b\in A_Y$,}
\begin{equation}
a\underset{x\in X}{\lhd}(b)\equiv b\underset{y\in Y}{\lhd}(a)\hspace{5mm} \in A_{X\underset{x\in X}{\lhd}Y}\equiv A_{X\times Y}\equiv
 A_{Y\underset{y\in Y}{\lhd} X}
\end{equation}
Schematically seen,
\begin{equation}
\begin{tikzpicture}[scale=0.3]

\tikzstyle{no}=[circle,draw,fill=black!100,inner sep=0pt, minimum width=4pt]
\tikzstyle{no2}=[circle,draw,fill=black!10,inner sep=0pt, minimum width=4pt]

\node[no] (b1a) at (-1,11){};
\node at (-1.5,12) {$b$};
\node[no] (b2a) at (-3,9){};
\node[no] (b3a) at (-1,9){};
\node[no] (b4a) at (1,9){};
\node[no] (n1) at (-3.5,7){};
\node[no] (n2) at (-2.5,7){};
\node[no] (n3) at (-1.5,7){};
\node[no] (n4) at (-0.5,7){};
\node[no] (n5) at (0.5,7){};
\node[no] (n6) at (1.5,7){};

\node at (2.3,10) {$\equiv$};

\node[no] (b1b) at (6,11) {};
\node at (6,12) {$b$};
\node[no] (b2b) at (4,9)  {};
\node[no] (b4b) at (8,9) {};
\node[no] (m1) at (3.5,7){};
\node[no] (m2) at (4.5,7){};
\node[no] (m3) at (4,7){};
\node[no] (m4) at (7.5,7){};
\node[no] (m5) at (8.5,7){};
\node[no] (m6) at (8,7){};

\foreach \from/\to in {b1a/b2a,b1a/b3a,b1a/b4a,b1b/b2b,b1b/b4b,n1/b2a,n2/b2a,n3/b3a,n4/b3a,n5/b4a,n6/b4a,b2b/m1,b2b/m2,b2b/m3,b4b/m4,b4b/m5,b4b/m6}
\draw[thick](\from) -- (\to);
\end{tikzpicture}
\end{equation}

\noindent We say that $A$ is "totally- commutative" if it is both commutative and vertical- commutative.

\defin{8.1.6} \textit{A homomorphism
 of generalized rings
$\varphi : A \rightarrow A'$
is a natural transformation of functors
(so for $X \in \FF$,
we have $\varphi_X \in Set_0 (A_X , A'_X)$,
and for $f \in \FF (X,Y)$,
we have $\varphi_Y \circ f_A = f_{A'} \circ \varphi_X$),
such that $\varphi$ preserves multiplication, contraction, and the unit:
\begin{equation}\label{8.1.43a}
\varphi(a \lhd b) = \varphi (a) \lhd \varphi (b),
\end{equation}
\begin{equation}\label{8.1.43b}
\varphi((a\sslash b)) = (\varphi (a)\sslash  \varphi (b)) ,
\end{equation}
\begin{equation}\label{8.1.43c}
\varphi (1_A) = 1_{A'}
\end{equation}
We remark that for generalized rings $A, A'$,
a collection of maps
$\varphi = \{\varphi_X \in Set_0 (A_X , A'_X )\}$
satisfying
(\ref{8.1.43a}-\ref{8.1.43c})
is a homomorphism;
i.e. it is automatically a natural transformation of functors by functoriality (O) and unit axiom (VIII).}

Thus we have a category of generalized rings and homomorphisms
which we denote by $\cG\cR$. \\
We denote by $\cG\cR_C$ the full subcategory of $\cG\cR$ consisting of the \emph{commutative} generalized rings. \vspace{3mm}\\

\noindent There are three equivalent ways to describe the operations of a generalized ring $A$:
\begin{enumerate}
\item The elementary operations of multiplications and contraction
as above:

\begin{equation}
\begin{cases}
a\lhd b: A_Z\times A_X\rightarrow A_{Z\underset{z_0}{\lhd} X}, \;\; z_0\in Z \\
a\sslash b: A_Y\times A_X\rightarrow A_{Y/X}, \;\; X\subseteq Y \\
\end{cases}
\end{equation}

\item Using the disjointness axioms we can view these operations as maps: \\
For $f\in Set_{\bullet}(Y,Z)$ a map of sets, ($Y,Z\in \FF$), and $A_f=\prod_{z\in Z} A_{f^{-1}(z)}$,
\begin{equation}
\begin{cases}
a\lhd b: A_Z\times A_f\rightarrow A_Y \\
(a_Z, (b^f))\longmapsto a_Z\underset{z}{\lhd} (b_z^f)\\
\\
a\sslash b: A_Y\times A_f\rightarrow A_Z \\
(a_Y, (b^f))\longmapsto a_Y\sslash (b^f_z) \\
\end{cases}
\end{equation}

\item We can further extend the operations "fiber by fiber", so that for map of sets $Y\xrightarrow{f} Z\xrightarrow{g} W$ we obtain:

\begin{equation}
\begin{cases}
a\lhd b: A_g\times A_f\rightarrow A_{g\circ f} \\
[(a_w^g)\lhd (a_z^f)]_{w_0}:=a_{w_0}^g\underset{z\in g^{-1}(w_0)}{\lhd} (a_z^f) \\
\\
a\sslash b: A_{g\circ f}\times A_f\rightarrow A_g \\
[(a_w^{g\circ f})\sslash (a_z^f)]_{w_0}:=a_{w_0}^{g\circ f}\sslash (a_z^f)_{ z\in g^{-1}(w_0)}
\end{cases}
\end{equation}
\end{enumerate}

\noindent \textbf{\large The axioms in the fiber-extended form.}\vspace{3mm}\\
We can now write the axioms again in the fiberwise extended form:

\noindent{\bf Associativity:} For $  W \xleftarrow{h} Z \xleftarrow{g} Y \xleftarrow{f} X  $ in
$ Set_{\bullet}$,
and for $d \in A_h$, $c \in A_g$, $b \in A_f$, we have in
$ A_{h \circ g\circ f}$:
\begin{equation}
    d \lhd (c \lhd b) = (d \lhd c ) \lhd b
\end{equation}

\begin{center}
\begin{diagram}
 & & Z &&\lTo^{c \lhd b}&& X\\
 &\ldTo^{d}& & \luTo^{c}& & \ldTo^{b}\\
W& &\lTo_{d \lhd c}&& Y
\end{diagram}
\end{center}

\vspace{10pt}
\noindent {\bf Left-Adjunction:}
  \label{sec1.4}
For  $  W \xleftarrow{h} Z \xleftarrow{g} Y \xleftarrow{f} X  $ in
$ Set_{\bullet}$,  and for $d \in A_{h \circ g\circ f}$, $ a \in A_g$,
$c \in A_f$, we have in $ A_h$:
\begin{equation}
 ( d \sslash (a \lhd c)) = \left( (d \sslash c ) \sslash a \right)
\end{equation}

\begin{center}
\begin{diagram}
& &X& & \\
&\ldTo(2,4)^d&\dTo_c&\rdTo(2,4)^{a \lhd c}& \\
& & Y & & \\
&\ldTo_{(d\slash c)}& & \rdTo_{a}& \\
W  & & \lTo& &  Z
\end{diagram}
\end{center}

\vspace{10pt}
\noindent{\bf Right-Adjunction:}
\label{sec1.5}
 For $  W \xleftarrow{h} Z
\xleftarrow{g} Y \xleftarrow{f} X$ in $ Set_{\bullet}$,
and for $d \in A_{h \circ g}$, $ a \in A_{g\circ f}$, $c \in A_f$, we have
in $ A_h$:
\begin{equation}
  (d \lhd c) \sslash a = d \sslash (a\sslash c)
\end{equation}

\begin{center}
\begin{diagram}
  & &X& & \\
  & \ldTo^{c}& &\rdTo^{a}(2,4)& \\
Y  & & & & \\
\dTo^{d}&\rdDashto^{a\sslash c}(4,2) & & & \\
W  &\lTo  & & & Z
\end{diagram}
\end{center}

\vspace{10pt}
\noindent{\bf Left-Linear:}
 For $  W \xleftarrow{h} Z \xleftarrow{g} Y \xleftarrow{f} X$
     in $ Set_{\bullet}$,
and for $d \in A_h$, $ a \in A_{ g\circ f}$, $c \in A_f$,
we have in $ A_{h \circ g}$:
\begin{equation}
   (d \lhd a) \sslash c= d \lhd (a\sslash c)
\end{equation}

\begin{center}
\begin{diagram}
  & &X& & \\
  & \ldTo^{a}& &\rdTo^{c}(2,4)& \\
Z  & & & & \\
\dTo^{d}&\luDashto^{a\sslash c}(4,2) & & & \\
W  &\lTo  & & & Y
\end{diagram}
\end{center}

The axiom of right-linearity is the delicate one, as it captures the commutativity we need, and it will require the following notations for generalizing the indices adjustment we had in (\ref{comm_identity}).

\vspace{10pt}
 For $  Z \xrightarrow{g} Y \xleftarrow{f} X  $
 in $ Set_{\bullet}$, we form the cartesian square,
\begin{equation}
%
%
  Z  \prod\limits_{Y}  X = \left\{ (z,x) \in D(g) \times D(f) ,
           g(z) = f(x)  \right\}
\end{equation}
\begin{center}
\begin{diagram}
 & &Z \prod\limits_Y X & &\\
 &\ldTo^{\tilde{f}}& &\rdTo^{\tilde{g}}& \\
Z& & & &X\\
 &\rdTo_g&&\ldTo_f&\\
 & & Y& & &
\end{diagram}
\end{center}

with $  \tilde{f}$, $\tilde{g}$ the natural projections. Note that
we have an identification of fibers
\[ \tilde{f}^{-1}(z) \cong f^{-1} ( g (z) )\,  \mbox{for} \, z \in Z
\]
 and we get a map
\begin{equation}
 A_f \rightarrow A_{\tilde{f}} ,\quad c \mapsto
  \ g^*c \in A_{\tilde{f}}
\quad \mbox{with} \quad g^*c^{ (z) } = c^{ ( g (z))}.
\end{equation}
Similarly, $   \tilde{g}^{-1} (x) \cong g^{-1} ( f (x) )$  for
$x\in X$, and we get a map
\[
  A_g \rightarrow A_{\tilde{g}} \quad,\quad
  a \mapsto f^*a \in A_{\tilde{g}} \quad \mbox{with}
  \quad f^*a^{ (x) } = a^{ ( f (x) ) }
\]
\noindent{\bf  Right-Linear:}
\label{sec1.7}
For
$f, g, \tilde{f}, \tilde{g}$ as above, and $  W \xleftarrow{h} Y$ in
$ Set_{\bullet}$, and for $   d \in A_{h \circ f}$, $a \in A_g$,
$ c \in A_f$,
we have in $  A_{h \circ g}$:
\begin{equation}
  ( d \sslash c ) \lhd a = (d \lhd f^*a) \sslash g^*c
\end{equation}

\begin{center}
\begin{diagram}
& & & & X \prod\limits_{Y} Z& & \\
& & & \ldTo^{f^*a}& & \rdTo^{g^*c}& \\
& & X& & & & Z \\
 & \ldTo^{d}& & \rdTo_{c}& &\ldTo_{a}& \\
W  && \lTo_{(d\slash c)} & &Y& &
\end{diagram}
\end{center}

\vspace{10pt}
\noindent{\bf Unit axioms:}
  We have a distinguished element
$   1 = 1_A \in A_{[1]}$.

Hence for each singleton $  \{ x \} \in \FF$,
using the unique isomorphism $   [1] \xrightarrow{ \sim} \{ x \}$
we get $    1_x \in A_{ \{ x \}}$;
and for $ f \in \FF (X,Y)    $
we have $   1_f = ( 1_{ f^t(y) } )_{ y \in f(X) } \in A_f$.

We have for all $       h \in Set_{\bullet}(X,W)$,
$d \in A_h$:
\begin{equation}
  \label{eq1.8.1}
\begin{array}{c}
d \lhd 1_{id_X}  = d , \\
 1_{id_W} \lhd d \; = d ,  \\
 (d \,\sslash \, 1_{id_X}) = d
\end{array}
\end{equation}

\numberwithin{equation}{section}
\setcounter{equation}{0}
\section{Remarks}

\rem{8.2.1. \ Functoriality}

In $   A_{[1]}  $ we have $    0 \lhd 1 = 1 \lhd 0 = 0
\lhd 0 = 0$, and $    1 \lhd 1 = 1$, hence for $   f \in
\FF (X,Y)$,  $g \in \FF (Y,Z)$, we have
\begin{equation}
   1_g \lhd 1_f = 1_{g \circ f}
\end{equation}
 Also we have
$0\sslash 1 = 1\sslash 0 = 0\sslash 0 =0$,
and $1\sslash 1=1$,
hence
\begin{equation}
\label{8.2.2}
 1_{id_Y}\sslash  1_{f^t} = 1_f
\end{equation}
We obtain for all $a \in A_X$, $f \in  \FF (X,Y)$,
\begin{equation}
\label{8.2.3}
  a \lhd 1_{f^t} =  (a \lhd 1_{f^t}) \sslash 1_{id_Y} = ( a\sslash
(1_{id_Y}\sslash 1_{f^t} ) ) = a\sslash 1_f.
\end{equation}
This gives a structure of a functor $\FF \rightarrow Set_0$ on $X
\mapsto A_X$, which is the given structure by the functoriality
and the unit axiom: for all $a \in A_X$, $f \in \FF (X,Y)$,
\begin{equation}
\label{8.2.4}
   a \lhd 1_{f^t} = (a\sslash 1_f ) = f_A(a).
\end{equation}

\rem{8.2.2. \ The involution}

Note that $     (1\sslash a)   $ makes sense only for $     a \in A_{[1]}$,
we define
\begin{equation}
      a^t = 1\sslash a \quad, \; a \in A_{[1]}.
\end{equation}
It is an involution of $    A_{[1]}$:
\begin{equation}
   (a^t)^t = 1\sslash (1\sslash a) =  (1 \lhd a) \sslash 1  = a\sslash 1=a
\end{equation}
 It preserves the operation of multiplication, and the
special elements $0$, $1$:
\begin{equation}\label{eq8.2.7}
   (a \lhd b)^t = 1\sslash (a \lhd b) = (1\sslash b)\sslash a =
  \Big( 1 \lhd (1\sslash b)\Big) \sslash a =
\end{equation}
\[= (1\sslash a) \lhd (1\sslash b) = a^t \lhd b^t
\]
\[
   0^t = (1\sslash 0)=0 \qquad , \quad 1^t = (1\sslash 1)=1  \]

\defin{8.2.3}  We shall say that $A$
 is \emph{self-adjoint} if
 \[    a^t=a
\mbox{ for all} \, a \in A_{[1]}.\]

In general, we let
\begin{equation}\label{eq8.2.8}
A_{[1]}^+=\{a\in A_{[1]},a^t=a\}
\end{equation}
denote the subset of symmetric elements.

\noindent \rem{8.2.4. \ The structure of $A_{[1]}$}

It follows from the associativity and unit axioms
  that $A_{[1]}$  is an associative monoid
  with unit $1$. It also follows from the commutativity axiom IX, that it is commutative:
\begin{equation}\label{eq8.2.9}
     a \lhd b = (a \lhd b)\sslash1 = (1 \lhd a)\sslash (1\sslash b) \overset{\text{comm.}}{=}
(1\sslash (1\sslash b))\lhd a = b \lhd a
\end{equation}
For $  X \in \FF$, the commutative monoid $A_{id_X} =
(A_{[1]})^X $ acts on the right on $A_X$.
\begin{equation}
\xymatrix@R=2mm{A_X\times (A_{(1)})^X\ar[r]&A_X \\
 b \ \ ,(a^{(x)})\ar@{|->}[r]& b\underset{x\in X}{\lhd} (a^{(x)})}
\end{equation}
Note that for $  a=(a^{(x)}) \in A_{id_X}$, and $  b \in A_X$,
\begin{equation}
  b \lhd a = (b \lhd a) \sslash 1_{id_X}  = (b\sslash (1_{id_X}\sslash a)) = (b\sslash a^t)
\end{equation}
The monoid $  A_{[1]} $ also acts on the left on $ A_X$, via
\begin{equation}
\xymatrix@R=2mm{A_{[1]}\times A_X\ar[r]&A_X \\
 a \ \ ,b \ar@{|->}[r]& a\lhd b}
\end{equation}
  and this is the diagonal
right action: For $  a \in A_{[1]}$, $b \in A_X$, putting
$\tilde{a} \in A_{id_X}$, $\tilde{a}^{(x)}\equiv a$ for all $  x \in X$,
\begin{equation}
\label{eq_up2}
  a \lhd b = (1 \sslash a^t) \lhd b \overset{\text{comm.}}{=} (1 \lhd b) \sslash \tilde{a}^t =
b\sslash \tilde{a}^t = b \lhd \tilde{a}.
\end{equation}
More generally, for $  f \in Set_{\bullet}(X,Y)$, $b \in A_f$,
$a= (a^{(y)})
\in A_{id_Y} = (A_{[1]})^Y$, we get
\begin{equation*}
  f^*a \in A_{id_X} = (A_{[1]})^X \;,\;
f^*a^{(x)} =  a^{( f(x) )} \, \mbox{for} \, x \in X,
\end{equation*}
and we have in $  A_f$,
\begin{equation}
\label{eq_up1}
  a \lhd b = b \lhd f^*a.
\end{equation}
Indeed checking (\ref{eq_up1}) at a given fiber
$   A_{f^{ -1 } (y) }$, $y \in Y$,
reduces to (\ref{eq_up2}).

The action of $(A_{[1]})^X$ on $  A_X  $ is self-adjoint
with respect to the pairing (\ref{eq8.1.9}):
for $a \in (A_{[1]})^X$,  $b,d \in A_X$, we have
\begin{equation}\label{eq_up3}
   ( b \lhd a)\sslash d  =  b\sslash (d\sslash a)  =  b\sslash (d \lhd a^t )
\end{equation}
We denote the group of invertible elements of $A_{[1]}$ by $A^*$:
\begin{equation}
A^*=\{a\in A_{[1]}, \exists a^{-1}\in A_{[1]}, a\lhd a^{-1}=a^{-1}\lhd a=1\}
\end{equation}
\noindent \rem{8.2.5. "one contraction suffices"}

The axioms of a \emph{commutative} generalized ring allow to transform any formula in the operations of multiplication and contraction, to an equivalent formula with only
\underline{one} contraction. Thus expressions of the form
\begin{equation}
 (a\sslash b) = (a_1 \lhd a_2 \lhd \cdots \lhd a_n )\sslash ( b_1 \lhd \cdots \lhd b_m)
\end{equation}
are closed under multiplication and contraction.
We have the following lemma:
\lem{8.2.6} Suppose we have a commutative diagram in $Set_{\bullet}$ and elements of $A$ over it ($b_i\in A_{f_i},b\in A_h$):
\begin{equation}
\begin{tikzcd}
&\tilde{X}_2\ar{rd}{\tilde{f_2}}\ar{ld}{\tilde{\tilde{h}}}&&\\
X_2\ar{rd}{f_2}&&\tilde{X}_1\ar{rd}{\tilde{f_1}}\ar{ld}{\tilde{h}}& \\
&X_1\ar{rd}{f_1}&&\tilde{X}_0\ar{ld}{h}\\
&&X_0&
\end{tikzcd}
\begin{tikzcd}
&\tilde{X}_2\ar{rd}{\tilde{h}^*(b_2)}\ar{ld}{f_2^*(f_1^*(b))}&&\\
X_2\ar{rd}{b_2}&&\tilde{X}_1\ar{rd}{h^*(b_1)}\ar{ld}{f_1^*(b)}& \\
&X_1\ar{rd}{b_1}&&\tilde{X}_0\ar{ld}{b}\\
&&X_0&
\end{tikzcd}
\end{equation}
where the inner squares are fibre products. Then we have the following identities:
\begin{equation}
h^*(b_1\lhd b_2)=h^*(b_1)\lhd \tilde{h}^*(b_2),
\end{equation}
\begin{equation}
f_1^*(f_2^*(b))=(f_2\circ f_1)^*(b).
\end{equation}
\begin{proof}
For any $z\in \tilde{X}_0$ we get:
$$h^*(b_1\lhd b_2)^{(z)}=(b_1\lhd b_2)^{(h(z))}=b_1^{(h(z))}\underset{f_1|_{h(z)}}{\lhd}{b_2|_{h(z)}}=$$
$$ =b_1^{(h(z))}\underset{f_1|_{h(z)}}{\lhd}(b_2)_{x\in f^{-1}(h(z))}, $$
and
$$(h^*(b_1)\lhd \tilde{h}^*(b_2))^{(z)}= h^*(b_1)^{(z)}\underset{\tilde{f}_1|_z}{\lhd}\tilde{h}^*(b_2)|_z=$$
$$=b_1^{(h(z))}\underset{\tilde{f}_1|_z}{\lhd} \left( \tilde{h}^*(b_2)\right)_{\tilde{x}\in \tilde{f}_1^{-1}(z)}=b_1^{(h(z))}\underset{\tilde{f}_1|_z}{\lhd}(b_2^{(\tilde{h}(\tilde{x}))}).$$
and we get the equality by the identification of fibers for any $z\in \tilde{X}_0$ of $f_1^{-1}(h(z))$ and $\tilde{f_1}^{-1}(z)$ given by
$$(x,z)\in \tilde{f_1}^{-1}(z) \iff x\in f_1^{-1}(h(z)).$$
\end{proof}

In the general case where there are $n$ consecutive fiber products and an element $b_1\lhd b_2\lhd \dots b_n$, we have formulas:
\begin{equation}
\begin{split}
h^*(b_1\lhd b_2 \lhd \dots \lhd b_n)=h^*(b_1)\lhd \tilde{h}^*(b_2)\lhd \tilde{\tilde{h^*}}(b_3)\lhd \dots \\
(f_1\circ \dots f_n)^*(b_1\lhd\dots \lhd b_n)= f_n^*(f_{n-1}^*(\dots\circ f_1^*(b)\dots).
\end{split}
\end{equation}

We have the formulas:
\begin{equation}
\label{multiplication}
\textbf{\large Multiplication 8.2.7:}\;\;(a\sslash b) \lhd (c\sslash d) = (a \lhd g^*c )\sslash ( d \lhd f^*b ),
\mbox{for} \ b \in A_g , \ c \in A_f
\end{equation}
Indeed, for $c=c_1\lhd \dots \lhd c_l$ an element over $g:Y\rightarrow Z$, and $b=b_1\lhd\dots \lhd b_m$ an element over
$f:X\rightarrow Z$, we have,
\begin{equation}
(a\sslash b)\lhd (c\sslash d)= ((a\sslash b)\lhd c)\sslash d\overset{\text{comm.}}{\large =}((a\lhd f^*c)\sslash g^* b)\sslash d=(a\lhd f^*c)\sslash (d\lhd g^* b).
\end{equation}

\begin{center}
\begin{equation}
\begin{diagram}
&& & &\bullet & & & &\\
&& &\ldTo^{g^*c}& &\rdTo^{f^*b}& & &\\
&& \bullet & &\large{\prod}& &\bullet&& \\
&\ldTo^{a}& &\rdTo^{b}& &\ldTo^{c}& &\rdTo^{d}&\\
\bullet&&\lTo& &\bullet&\lTo& &&\bullet
\end{diagram}
\end{equation}
\end{center}

\noindent Similarly we have:

\begin{equation}\label{contraction}
\textbf{\large Contraction 8.2.8:}\;\;
(a\sslash b) \sslash (c\sslash d)  = (a \lhd g^*d)\sslash (
c \lhd f^*b), \ \mbox{for} \ b \in A_g , \ d \in A_f
\end{equation}

\begin{equation}
\begin{diagram}
& & \bullet & & \\
&\ldTo^{g^*d}& & \rdTo^{f^*b}& \\
\bullet& & \large{\prod}& & \bullet \\
 &\rdTo^{b}& &\ldTo^{d}& \\
\dTo^{a}& & \bullet& &\dTo_{c}\\
 & \ldTo & &\rdTo & \\
\bullet&\lTo& & &\bullet
\end{diagram}
\end{equation}

We can also see this in a more explicit way: for elements $d_1\lhd d_2\lhd \dots \lhd d_m$ over $Y_m\rightarrow \dots \rightarrow Y_0$ and $b_1\lhd b_2\lhd \dots \lhd b_n$ over $X_m\rightarrow \dots X_0$, with an isomorphism $Y_m\cong X_n$, and an element $b$ over $h:\tilde{X}_0\rightarrow X_0$ we form the diagram:

\begin{equation}
\begin{tikzpicture}[scale=0.3]

\tikzstyle{no}=[circle,draw,fill=black!100,inner sep=0pt, minimum width=4pt]
\tikzstyle{no2}=[circle,draw,fill=black!10,inner sep=0pt, minimum width=4pt]


\node at (9.4,17) {$Y_m$};
\node at (11,17) {$\cong$};
\node at (13,17) {$X_n$};
\node at (5,9.1) {$Y_2$};
\node at (3,5.5) {$Y_1$};
\node at (1,1.9) {$Y_0$};
\node at (17.5,9.1) {$X_2$};
\node at (19.5,5.5) {$X_1$};
\node at (21.5,1.9) {$X_0$};
\node at (18,20) {$\tilde{X}_n$};
\node at (22.5,12.3) {$\tilde{X}_2$};
\node at (24.5,8.6) {$\tilde{X}_1$};
\node at (26.5,4.9) {$\tilde{X}_0$};

\draw
(13,16.5) edge[->,>=angle 90] node[left] {\tiny $b_n$} (14,14.5) ;
\draw
(17.8,8.5) edge[->,>=angle 90] node[left] {\tiny $b_2$} (18.8,6.5) ;
\draw
(19.9,4.7) edge[->,>=angle 90] node[left] {\tiny $b_1$} (21,2.6) ;
\draw
(22.7,11.5) edge[->,>=angle 90] node[right] {\tiny $\tilde{h}^*(b_2)$} (24,9.2) ;
\draw
(25,7.7) edge[->,>=angle 90] node[right] {\tiny $h^*(b_1)$} (26,5.5) ;
\draw
(8.8,16.3) edge[->,>=angle 90] node[left] {\tiny $d_m$} (7.3,13.5) ;
\node at (7.2,13.1) {$\cdot$};
\node at (6.5,12.1) {$\cdot$};
\node at (5.8,11.1) {$\cdot$};
\draw
(4.4,8.5) edge[->,>=angle 90] node[left] {\tiny $d_2$} (2.8,6.2) ;
\draw
(2.4,5) edge[->,>=angle 90] node[left] {\tiny $d_1$} (.7,2.7) ;
\draw
(20.4,1.9) edge[->,>=angle 90] node[left] {} (1.9,1.9) ;
\draw
(17,20) edge[->,>=angle 90] node[left] {} (13,17.7) ;
\draw
(21.5,11.7) edge[->,>=angle 90] node[above] {\tiny $f_2^*(f_1^*(b))$} (17.7,9.5) ;
\draw
(23.5,7.8) edge[->,>=angle 90] node[right] {\tiny $f_1^*(b)$} (20,5.7) ;
\draw
(25.5,4.5) edge[->,>=angle 90] node[right] {\tiny $b$} (22.5,2.2) ;
\draw
(18.2,19.3)  edge[->,>=angle 90] node[right] {\tiny $\tilde{h^*}^{(n)}(b_n)$} (19.5,17) ;
\draw
(18.5,16.5)  edge[->,>=angle 90] node[left] {} (14.5,14.4) ;

\node at (15,13) {$\cdot$};
\node at (15.7,12) {$\cdot$};
\node at (16.3,11) {$\cdot$};

\node at (21,14) {$\cdot$};
\node at (20,16) {$\cdot$};
\node at (20.5,15) {$\cdot$};





\end{tikzpicture}
\end{equation}

By commutativity we have now,
\begin{equation}
\bigg ((d_1\lhd d_2\lhd\dots\lhd d_m)\sslash (b_1\lhd\dots \lhd b_n)\bigg )\lhd b=(d_1\lhd\dots\lhd (f_n^*\circ\dots\circ f_1^*(b))\sslash (h^*(b_1)\lhd\dots \lhd \tilde{h^*}^{(n)}(b_n)).
\end{equation}

\noindent \rem{8.2.9. Matrixness and Tameness}
We say that $A$ is a "matrix" generalized ring if for $X\in \FF$ we have an \emph{injection:}
\begin{equation}
A_X\inj (A_{[1]})^X, \;\; a\mapsto (a\sslash 1_x)_{x\in X}.
\end{equation}
We say $A$ is "tame"if for $X\in \FF$, and any $a,a'\in A_X$, we have
\begin{equation}
a\sslash d= a'\sslash d\; \text{ for all }\; d\in A_X\implies a=a'.
\end{equation}
Note that we have the implication
$$A\; \text{matrix} \implies \; A\text{ tame }.$$

\section{Examples of generalized Rings}
\smallbreak
\setcounter{subsection}{-1}
\subsection{Generalized Rings arising from $\frs$.}
If $A=\{A_{Y,X}\}$ is an $\fr$ with involution, we get a generalized ring $\mathcal{G}(A)$, with
$$\mathcal{G}(A)_X:=A_{1,X}.$$
and operations, \\
\textbf{multiplication:}
\begin{equation}
\begin{split}
z_0\in Z:\hspace{3mm} \mathcal{G}(A)_Z\times \mathcal{G}(A)_X= & A_{1,Z}\times A_{1,X}\rightarrow \mathcal{G}(A)_{Z\underset{z_0}{\lhd}X}=A_{1,Z\underset{z_0}{\lhd}X} \\
&a_Z\underset{z_0}{\lhd} a_X:= a_Z\circ (a_X\oplus \underset{z\in Z\setminus \{z_0\}}{\bigoplus}1).
\end{split}
\end{equation}
\textbf{contraction:}
\begin{equation}
\begin{split}
\mathcal{G}(A)_Y\times \mathcal{G}(A)_X= &A_{1,Y}\times A_{1,X}\rightarrow \mathcal{G}(A)_{Y/X}=A_{1,Y/X} \\
&a_Y\sslash a_X:= a_Y\circ (a_X^t\oplus \underset{y\in Y\setminus X}{\bigoplus} 1).
\end{split}
\end{equation}
Equivalently, for a map of sets $f:Y\rightarrow Z$, we have \\
\textbf{multiplication:}
\begin{equation}
\begin{split}
A_{1,Z}\times \prod_{z\in Z} A_{1,f^{-1}(z)}\rightarrow A_{1,Y}\\
a_Z\lhd (a_z^f):= a_Z\circ (\underset{z\in Z}{\bigoplus} a_z^f)
\end{split}
\end{equation}
\textbf{contraction:}
\begin{equation}
\begin{split}
A_{1,Y}\times  \prod_{z\in Z} A_{1,f^{-1}(z)}\rightarrow A_{1,Z} \\
a_Y\sslash (a_z^f):= a_Z\circ (\underset{z\in Z}{\bigoplus} a_z^f)^t
\end{split}
\end{equation}
When $A$ is $\times$- commutative (resp. totally- commutative) $\fr$ with involution, the generalized ring $\mathcal{G}(A)$ is commutative (resp. totally- commutative). \vspace{1.5mm}\\
We have the following diagram:

\begin{equation}
\xymatrix{\fr^t\ar[r]&\cG\cR \\
\fr^t_{\times\text{-com}}\ar[r]\ar@{^{(}->}[u]&\cG\cR_C\ar[u] \\
\fr^t_{\text{tot-com}}\ar[r]\ar@{^{(}->}[u]&\mathcal{GR}_{\text{tot-com}}\ar@{_{(}->}[u]}
\end{equation}

In particular we have the totally- commutative matrix generalized rings:
$$\begin{matrix}
\text{I.   }\;\;\;\;\;\;\FF &\text{the initial object of } \mathcal{GR}: \; \FF_X=X\cup\{0\}\\
\text{II. }\mathcal{G}(R)& R\text{\;\;a commutative rig}\\
\text{III. }\mathcal{O}_{K,\eta}& \eta:K\rightarrow \CC, \;\;\;\;\; (\ell_2)\text{!}\\
\text{IV. }\FF\{M\}& M\text{ commutative monoid.}
\end{matrix}$$
We specify now each of these examples explicitly. \\
\subsection{The field with one element $\FF$}
We write $\FF$ for the functor
$\FF \stackrel{\thicksim}{\rightarrow} \FF_0 \subseteq Set_0$.
Thus
\begin{equation}
\FF_X=X_0=X \coprod  \{0_X \} , \, \mbox{and} \,
\FF_f = \prod_y\left(f^{-1}(y) \right)_0 \, \mbox{for} \,
f \in Set_{\bullet}(X,Y).
\end{equation}

The multiplication is given by
\begin{equation}
\begin{array}{ll}
\FF_Y \times \FF_f & \rightarrow  \FF_X\\
y_0 , (x_f^{(y)}) & \mapsto y_0 \lhd x_f = x_f^{(y_0)} \in
          \left(f^{-1}(y_0)\right)_0\subseteq X_0
\end{array}
\end{equation}
The contraction is given by
\begin{equation}
\begin{array}{ll}
\FF_X \times \FF_f & \rightarrow  \FF_Y\\
x_0 , (x_f^{(y)}) & \mapsto (x_0 \sslash x_f) =
   \left\{\begin{array}{ll}
y & \, \, x_0 = x_f^{(y)}\\
0 & \mbox{otherwise}.
\end{array}\right.
\end{array}
\end{equation}
It is easy to check that $\FF$ is a totally-commutative,
self-adjoint matrix generalized ring.\\
\noindent For $A \in {\mathcal GR}$, and for $x \in X \in \FF$,
put
\begin{equation}
\varphi_X(x) = (1_x  \sslash 1_{j_x}) = 1_x \lhd 1_{j_x^t} \in A_X.
\end{equation}
\noindent with $j_x\in \FF(\{x\},X)$ the inclusion.\\
We have $\varphi_X \in Set_0 (\FF_X , A_X)$,
and the collection of $\varphi_X$ define a homomorphism
$\varphi \in {\mathcal GR}(\FF, A)$.
It is easy to check this is the only possible homomorphism,
and $\FF$ is the initial object of ${\mathcal GR}$.

\subsection{Commutative Rigs}
For a rig $A$, let $\mathcal{G}(A)_X = A \cdot X = A^X$
be the free $A$-module with basis $X$.
It forms a functor $\mathcal{G}(A): \FF \rightarrow A \cdot mod \subseteq Set_0$.
We define the multiplication for $f \in Set_{\bullet}(X, Y)$
by
\begin{equation}
{\cal G}(A)_Y \times {\cal G}(A)_f = A^Y \times \prod_{y \in Y} A^{f^{-1}(y)}
        \rightarrow A^X = {\cal G}(A)_X
\end{equation}
\[ a = (a_y) , b = (b_x^{(y)})_{x \in f^{-1}(y)} \mapsto (a \lhd b )_x
       = a_{f(x)} \cdot b_x^{f(x)}
\]
We define the contraction by
\begin{equation}
{\cal G}(A)_X \times {\cal G}(A)_f = A^X \times \prod_{y \in Y} A^{f^{-1}(y)}
        \rightarrow A^Y = {\cal G}(A)_Y
\end{equation}
\[
a = (a_x) , b = (b_x^{(y)})_{x \in f^{-1}(y)} \mapsto (a \sslash b )_y =
        \sum_{x \in f^{-1}(y)} a_x \cdot b_x^{(y)}
\]
It is straightforward to check that $\mathcal{G}(A)$ is a
self-adjoint, Matrix generalized ring. When $A$ is a commutative
rig, $\mathcal{G}(A)$ is a totally- commutative generalized ring.
A homomorphism of (commutative) rigs $\varphi \in Rig (A,B)$ gives
a homomorphism $\mathcal{G}(\varphi) \in
\mathcal{GR}(\mathcal{G}(A), \mathcal{G}(B))$, thus we have a
functor
\begin{equation}\label{embedding}
\mathcal{G} : Rig \rightarrow \mathcal{GR}
\end{equation}
It is fully-faithful: if $\varphi \in {\cal GR}({\cal G}(A), {\cal G}(B))$,
and $a = (a_x)\in {\cal G}(A)_X$,
then $\varphi_X(a)_x =\left(\varphi_{[1]}(a_x)\right)$ by functoriality
over $\FF$, so $\varphi$ is determined by \\
$\varphi_{[1]}: A \rightarrow B$;
the map $\varphi_{[1]}$ is multiplicative (and preserves $1$),
but it is also additive
\begin{equation}
\begin{split}
\varphi_{[1]} (a_1 + a_2) &= \varphi_{[1]} \left( (a_1, a_2) \sslash (1,1)\right) =
(\varphi_{[1]}(a_1) , \varphi_{[1]}(a_2)) \sslash (1,1) = \\
&=\varphi_{[1]}(a_1) + \varphi_{[1]}(a_2)
\end{split}
\end{equation}
Thus $\varphi_{[1]} \in Rig(A, B)$ , and
$\varphi = \mathcal{G}(\varphi_{[1]})$; and we have
\begin{equation}
Rig(A,B) = \mathcal{GR}(\mathcal{G}(A) ,\mathcal{G}(B) )
\end{equation}

Note that for every $X \in \FF$,
we have a distinguished element
\begin{equation}
\mathbbm{1}_X  \in {\cal G}(A)_X , \quad
\left(\mathbbm{1}_X\right)_x=1 \, \mbox{for all} \, x \in X
\end{equation}
hence for $f \in Set_{\bullet}(X,Y)$, the element
$\mathbbm{1}_f = (\mathbbm{1}_{f^{-1}(y)}) \in {\cal G}(A)_f$.
These elements satisfy
\begin{equation}
\label{eq2.2.6}
\mathbbm{1}_g \lhd \mathbbm{1}_f = \mathbbm{1}_{g \circ f} ; \quad
\mathbbm{1}_{[1]} =1 ;
\end{equation}
which imply naturality
\begin{equation}
\mathbbm{1}_X \lhd 1_{f^t} = \mathbbm{1}_{I(f)} \in {\cal G}(A)_{I(f)}
\subseteq  {\cal G}(A)_Y \, \mbox{for} \, f \in \FF(X,Y)
\end{equation}
Note that any element
$a  = (a_x) \in  {\cal G}(A)_X$, gives an element of the monoid
\begin{equation}
\langle a \rangle = (a_x) \in  {\cal G}(A)_{id_X}= A^X
\end{equation}
and the vector $\mathbbm{1}_X$ is "cyclic" in the sense that
\begin{equation}
\label{8.3.19}
a = \mathbbm{1}_X \lhd \langle a \rangle
\end{equation}

We have the generalized rings ${\cal G} (\mathbb{N})$ and ${\cal
G}([0, \infty ))$, as well as the "tropical'' examples ${\cal G}
(\mathbb{N}_0)$ and ${\cal G}([0, \infty )_0)$ where we replace
addition by the operation of taking the maximum $\max \{x,y\}$.
For $\sigma = 1/p \in (0,1]$ we have the generalized ring $\cG
([0, \infty )_{\sigma})$ where we replace addition by the
operation $x+_{\sigma} y := (x^p + y^p)^{\sigma}$.

\subsection{Real primes}\label{generealprimes}

Let $|\, | :K \rightarrow [0, \infty)$ be a non-archimedean
absolute value on a field $K$. Let
\begin{equation}\label{realprimes} {\cal O}_X = \left\{ a = (a_x)
\in \mathcal{G}(K)_X \, ,\, \sum\limits_{x \in X}|a_x|^2 \leq 1
\right\}
\end{equation}
Note that $X \mapsto {\cal O}_X$ is a subfunctor of ${\cal G}(K):
\FF \rightarrow Set_0$. But it is also a
\emph{sub-generalized-ring}, in the sense that it is closed under
the operations of multiplication and contraction, and $1 \in
\mathcal{O}_{[1]}$. The proof for multiplication is
straightforward:

\noindent For $a =(a_y) \in {\cal O}_Y$,
    $b =(b_x^{(y)}) \in {\cal O}_f$, we have
\begin{equation}
\sum\limits_{x \in X} |(a \lhd b)_x|^2 =
\sum\limits_{x \in X} |a_{f(x)}|^2 \cdot |b_x^{(f(x))}|^2 =
\end{equation}
\[
=\sum\limits_{y \in Y} |a_y|^2 \sum\limits_{x \in f^{-1}(y)}|b_x^{(y)}|^2 \leq
\sum\limits_{y \in Y}  |a_y|^2 \leq 1
\]
and so $a \lhd b \in {\cal O}_X$.

The proof for contraction is just the Cauchy-Schwartz inequality: \\
For $a =(a_x) \in {\cal O}_X$,
    $b =(b_x^{(y)}) \in {\cal O}_f$ , we have
\begin{equation}\begin{array}{ll}
\sum\limits_{y \in Y} |(a \sslash b)_y|^2 &=
\sum\limits_{y \in Y}\left|\sum\limits_{x \in f^{-1}(y)} a_x \cdot b_x^{(y)}
                     \right|^2 \leq\\

 &\leq \sum\limits_{y \in Y} \left(\sum\limits_{x \in f^{-1}(y)}|a_x|^2\right) \cdot
                     \left(\sum\limits_{x \in f^{-1}(y)}|b_x^{(y)}|^2\right)\\
&\leq \sum\limits_{y \in Y}\sum\limits_{x \in f^{-1}(y)}|a_x|^2 \leq
 \sum\limits_{x \in X}|a_x|^2 \leq 1
\end{array}\end{equation}
and so $a\sslash b \in {\cal O}_Y$.

\noindent Thus $\mathcal{O}$ is a generalized ring (totally- commutative, matrix, self- adjoint). \\
 Note that
\begin{equation}
m_X = \left\{ a = (a_x) \in \mathcal{O}_X , \sum\limits_{x \in X}|a_x|^2  <1 \right\}
\end{equation}
forms a subfunctor of $\mathcal{O}$, and it is (the unique maximal) functorial-
\emph{ideal} of $\mathcal{O}$,
in the sense that we have
\begin{equation}\label{eq8.3.24}
{\cal O} \lhd m \,  , \,  m \lhd {\cal O} \, , \,   (m\sslash {\cal O})\, ,\, ({\cal O}\sslash m) \subseteq m
\end{equation}
By collapsing $m_X$ to zero we obtain the quotient (in $Set_0$):
\begin{equation}
k_X = {\cal O}_X / m_X =
\left\{  a = (a_x) \in {\cal O}_X , \sum\limits_{x \in X}|a_x|^2  =1\right\} \amalg
\{ 0_X \}
\end{equation}
There is a canonical projection map $\pi_X:{\cal O}_X -\hspace{-5pt}-\hspace{-10pt}\gg k_X$,
with $\pi_X(m_X) = 0_X$.
Now (\ref{eq8.3.24}) imply that there is a (unique) structure of
 a generalized ring on $k$,
such that $\pi$ is a homomorphism,
$\pi \in {\cal GR}({\cal O}, k)$. It is given by
\begin{equation*}
a \lhd b = \left\{\begin{array}{cl}a \lhd b & \, \mbox{if} \
                                             ||a \lhd b||=1\\
0 & \, \mbox{if}\ ||a \lhd b||< 1
 \end{array}\right.
\quad,
(a\sslash b)=\left\{\begin{array}{cl} a\sslash  b & \, \mbox{if} \
                                             ||a\sslash b||=1\\
0 & \, \mbox{if}\ ||a\sslash b||< 1
 \end{array}\right.
\end{equation*}
with the $l_2$ -norm
\begin{equation*}
||(a_x) || = \left( \sum\limits_{x \in X} |a_x|^2 \right)^{\frac{1}{2}}
\end{equation*}
Note that the (totally- commutative, self- adjoint) generalized ring $k$ is Not matrix, but is tame.

\subsection{Ostrowski's theorem}

Let $K\in \mathcal{GR}_C$. We say $K$ is a \emph{field} when any element in $K_1\setminus \{0\}$ is invertible:
\begin{equation}
K_1^*=K_1\setminus \{0\}.
\end{equation}
For $A\subseteq K$ a sub-generalized ring of $K$ we say $A$ is \emph{full} when:
\begin{equation}
\forall y\in K_X, \exists d\in A_1\setminus \{0\}, d\lhd y\in A_X.
\end{equation}
(i.e. $K=A_{(0)}$ the localization of $A$ with respects to the prime $\{0\}$). \\
We say $A$ is \emph{tame} when:  for any $X\in \FF$
\begin{equation}
A_X=\{y\in K_X, y\sslash a\in A_1, \forall a\in A_X\}.
\end{equation}
We call $A$ a valuation of $K$ when it is full and tame and
\begin{equation}
\forall y\in K_1^*, y\in A_1 \text{ or } y^{-1}\in A_1.
\end{equation}
For $k\subseteq K$ we denote the set of $K/k$- valuations by
\begin{equation}
\Val(K/k):=\{A\subseteq K \text{ a valuation }, k\subseteq A\}.
\end{equation}
For $B\in \Val(K)$, the ordered-abelian-group $K_1^*/B^*_1$ can be embedded in a complete ordered-abelian-group $\Gamma$, and the quotient map,
\begin{equation}
|\;|:K_1\rightarrow \Gamma\cup\{0\},\hspace{5mm} |x|:=x\lhd B_1^*,
\end{equation}
satisfies:
\begin{enumerate}[(I)]
\item\label{i1}
\begin{equation}
\begin{split}
&|x|=0 \iff x=0 \\
&|x_1\lhd x_2|=|x_1|\cdot |x_2| \\
&|1|=1\;\; =(\text{unit of }\Gamma).
\end{split}
\end{equation}
For  $y\in K_X$ we have the equality (cf. Claim(II) of \S\ref{2.5})
\item \label{i2}
\begin{equation}
|y|_X:= \sup\{|y\sslash b|,\;  b\in B_X\}= \inf\{|d^{-1}|, d\in K_1^*,\; d\lhd y\in B_X\}
\end{equation}
For $f\in \ \Set(Y,Z)$, $a=(a_z)\in K_f=\prod_{z\in Z} K_{f^{-1}(z)}$, we put:
\begin{equation}
|a|_f=\Max \{|a_z|_{f^{-1}(z)}\}
\end{equation}
and we have,
\item\label{i3}
\begin{equation}
\begin{split}
&|z\lhd a|_Y\leq |z|_Z\cdot |a|_f,\hspace{5mm} ,z\in K_Z. \\
&|y\sslash a|_Z\leq |y|_Y\cdot |a|_f,\hspace{5mm} ,y\in K_Y.
\end{split}
\end{equation}
\end{enumerate}
Conversely, given the maps for $X\in \FF$,
\begin{equation}
|\;|_X:K_X\rightarrow \Gamma\cup\{0\}
\end{equation}
 satisfying (\ref{i3}), with $|\;|=|\;|_1$ satisfying (\ref{i1}), $B=\{B_X\}$, with $B_X=\{b\in K_X,\; |b|_X\leq 1\}$ is a sub-generalized-ring of $K$. If further (\ref{i2}) is satisfied, than $B$ is tame and full in $K$, and is a valuation of $K$.\\
When $\Gamma$ is rank-1, we can replace it by $\RR^+=(0,\infty)$, and two valuations
\begin{equation}
|\;|_X,|\;|'_X:K_X\rightarrow [0,\infty)
\end{equation}
are equivalent ($|x|_X\leq 1\iff |x|'_X\leq 1$) if and only if
$|x|_X'=|x|_X^{\sigma}$ for some $\sigma>0$. We have exactly as in
\S \ref{2.5}, the Ostrowski theorems. The proofs are the same as
in Appendix \ref{appB}, only written in the language of
generalized rings. Thus e.g. $(1,1)\circ \left( \begin{matrix}q_1
&0 \\0 &q_2 \end{matrix}\right)\circ \left( \begin{matrix} 1 \\1
\end{matrix}\right)$ should be replaced by $((1,1)\lhd
(q_i))\sslash (1,1)$. The proof for generalized ring is even
shorter: we only need to show $B_X= (\ZZ_{(\mathfrak{p})})^X$ or
$(\mathcal{O}_{\QQ,\eta})_X$, the rest of the proof (that
$B_{Y,X}=(\ZZ_{(\mathfrak{p})})^{X\times Y}$ or
$(\mathcal{O}_{\QQ,\eta})_{Y,X})$ is irrelevant.

\thm{Ostrowski  I}
\begin{equation}
\Val(\mathcal{G}(\QQ))=\Val(\slfrac{\mathcal{G}(\QQ)}{\FF\{\pm 1\}})=
\begin{cases}
  \mathcal{G}(\QQ) &  \\
  \mathcal{G}(\ZZ_{(\mathfrak{p})})  & \mathfrak{p} \mbox{ is prime} \\
  \mathcal{O}_{\QQ,\eta}& \mbox{"the real prime"}
\end{cases}
\end{equation}
\thm{Ostrowski  II} Let $K$ be a number field, $\mathcal{O}_K$ its ring of integers,
\begin{equation}
\Val(\mathcal{G}(K))=\Val(\slfrac{\mathcal{G}(K)}{\FF\{\mu_K\}})=
\begin{cases}
  \mathcal{G}(K) &  \\
  \mathcal{G}(\mathcal{O}_{K,\mathfrak{p}})  & \mathfrak{p}\subseteq \mathcal{O}_K \mbox{ a prime ideal, } \mathcal{O}_{K,\mathfrak{p}}=S_{\mathfrak{p}}^{-1}\mathcal{O}_K  \\
  \mathcal{O}_{K,\eta}, &  \eta:K\inj \CC \text{ mod conjugation}
\end{cases}
\end{equation}
\subsection{Commutative Monoids}
\label{s8.3.5}
For $A \in {\cal GR}$, we have
$A_{[1]} \in Mon^t$, cf. Remarks $8.2.2$-$8.2.4$, giving rise to a functor
\begin{equation}
{\cal GR} \rightarrow Mon^t \, , \hspace{1cm}  A \mapsto A_{[1]}
\end{equation}
This functor has a left adjoint. For $M \in Mon^t$, let
\begin{equation}
\mathbb{F} \{M\}_X = \coprod\limits_{x \in X} M =
\left( X \times (M \setminus \{ 0 \})\right) \amalg \{ 0 \} =
\left[X \times M \right] / [x,0] \sim 0
\end{equation}
It forms a functor $X \mapsto \mathbb{F} \{M\}_X : \FF\rightarrow Set_0$.
We define the operation of multiplication by
\begin{equation}
\lhd :\mathbb{F}\{M\}_Y \times \mathbb{F}\{M\}_f \rightarrow \mathbb{F}\{M\}_X
\end{equation}
\[[y_0 , m_0 ] \, ,\, [x^{(y)} , m^{(y)} ] \mapsto
[x^{(y_0)} , m_0 \cdot m^{(y_0)}]
\]
We define the operation of contraction by
\begin{equation}
( \,\, \sslash \,\, ) :\mathbb{F}\{M\}_X \times \mathbb{F}\{M\}_f
\rightarrow \mathbb{F}\{M\}_Y
\end{equation}
\[ [x_0 , m_0 ] \, ,\, [x^{(y)} , m^{(y)} ] \mapsto
[ y , m_0 ( m^{(y)})^t] \, \mbox{if } \, x_0 =x^{(y)}, \,
\mbox{otherwise} \, 0.
\]
It is straightforward to check that with these operations
$\mathbb{F}\{M\}$ forms a generalized ring, and for $A \in {\cal
GR}$, we have adjunction

\begin{equation}\label{mon_adj}
\begin{array}{cll}
{\cal GR} ( \mathbb{F}\{M\} , A )& = & Mon^t(M, A_{[1]})\\
\varphi &\mapsto &\varphi_{[1]}\\
\widetilde{\psi}_X ([x, m] ) = \psi(m)\sslash 1_{j_x}&
     \mapsfrom & \psi \\
\mbox{with} \, j_x \in \mathbb{F}([1], X ) , j_x (1) =x& &
\end{array}
\end{equation}

We let $CMon \subseteq Mon^t$ denote the full-subcategory of
commutative monoids (with trivial involution). For $M \in CMon$,
the generalized ring $\mathbb F\{M\}$ is
totally-commutative,matrix and self-adjoint.

\subsection{The free commutative generalized ring $\Delta^W$}

In this section we give a description of the commutative
generalized ring $\Delta^W$, characterized by the universal
property that homomorphisms of commutative generalized rings
$\varphi: \Delta^W \rightarrow A$
correspond bijectively with elements of $A_W$. \\
The elements of $(\Delta^W)_X$, of "degree" $X\in \FF$, are the formulas of degree $X$ that can be written using multiplication, contraction, the elements of $\FF$, and one variable $\delta_W$ of degree $W$, modulo the identifications of formulas that are consequences of the axioms of a commutative generalized ring. By Remark 8.2.5 - "one contraction suffice", these formulas have a concrete combinatorial shape which we describe next. \\

By a \emph{tree} $F$
we shall mean a finite set with a distinguished element $0_F \in F$,
the \emph{root},
and a map ${\cal S} = {\cal S}_F :F \setminus  \{ 0_F \} \rightarrow F $,
such that for all $a \in F$ there exists $n$ with
$\mathcal{S}^n (a) = 0_F$;
we write $n = ht(a)$, and put $ht(0_F) =0$.
For $a \in F$, we put $\nu (a) = \sharp \mathcal{S}^{-1}(a)$.
The boundary of $F$ is the set
\begin{equation}
\partial F = \{ a \in F \, ,\, \nu (a) =0 \}
\end{equation}
The unit tree is the tree with just a root,
$\{ 0 \}$, and
$\partial \{ 0 \} = \{ 0 \} $.
The zero tree is the empty set, $\emptyset$, and $\partial \emptyset = \emptyset$.
Given a subset $B \subseteq \partial F$, we have the reduced tree
$F|_B$,
obtained from $F$ by omitting all the elements of $\partial F \setminus B$,
and then omitting all elements $a \in F $ such that all the elements of
$\mathcal{S}^{-1}(a)$ have been omitted and so on;
we have $\partial (F|_B)=B$.

\begin{center}
\includegraphics[scale=1]{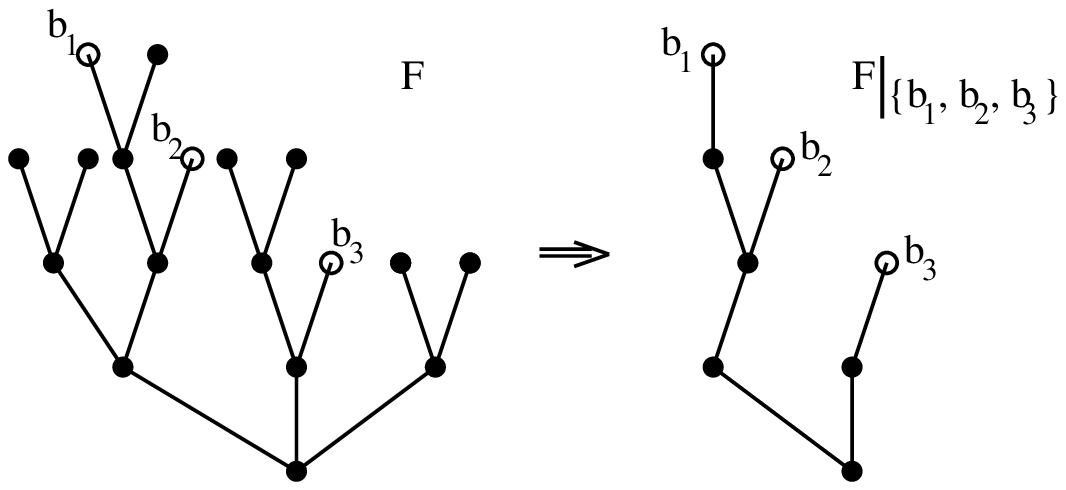}
\end{center}
\begin{equation}
\label{eq8.3.39}
\end{equation}
Given for each $b \in \partial F $, a tree $G_b$,
let $B =\{b \in \partial F \, , \, G_b \neq \emptyset \}$, and form the tree
\begin{equation}
\label{eq2.5.3}
F \lhd G := F|_B \amalg \coprod\limits_{b \in B} \left( G_b \setminus \{ 0_{G_b}\} \right)
\end{equation}
with $\mathcal{S} (a) = b $ if $a \in G_b$ and $\mathcal{S}_{G_b}(a) = 0_{G_b}$,
and otherwise $\mathcal{S}$ being the restriction of the given $\mathcal{S}_F$ and $\mathcal{S}_{G_b}$.

\vspace{10pt} Fix $W\in \FF$. We say that a tree $F$ is $W$-
\emph{labelled}, if we are given for all $b\in F$ an injection
\begin{equation}
\mu_b:S^{-1}(b)\inj W
\end{equation}
We view $\mu$ as a map $\mu:F\setminus \{0_F\}\rightarrow W$, $\mu(b)=\mu_{S(b)}(b)$, injective on fibers of $S_F$.

Let for $X \in \FF$,
\begin{equation}
\tilde{\Delta}_X = \{ F = (F_1 ; \{ \bar{F}_x \} ; \sigma )\}
\end{equation}
consist of the data of a $W$- labelled tree $F_1$, a $W$- labelled
tree $\bar{F}_x$ for each $x \in X$, and a bijection $\sigma :
\partial F_1 \stackrel{\thicksim}{\rightarrow} \coprod\limits_{x
\in X} \partial \bar{F}_x$. We view such data $F$ as being the
same as $F' = (F_1' ; \{\bar{F}_x'\} ; \sigma'   )$ if there are
isomorphisms of $W$- labelled trees $\tau_1 : F_1
\stackrel{\thicksim}{\rightarrow} F'_1 $, $\tau_x : \bar{F}_x
\stackrel{\thicksim}{\rightarrow} \bar{F}'_x $, with $\sigma'
\circ \tau_1 (b) = \tau_x \circ \sigma (b)$ for $b \in \partial
F_1$, $\sigma (b) \in \partial \bar{F}_x $.

Note that for such data $F=(F_1 ; \{\bar{F}_x\} ; \sigma)$,
we have an associated map
\begin{equation}
\underline{\sigma} : \partial F_1 \rightarrow X
\end{equation}
\[\underline{\sigma}(b) = x \, \mbox{iff} \, \sigma(b) \in \partial \bar{F}_x
\]
For $f \in Set_{\bullet} (X,Y)$ , we have
$\tilde{\Delta}_f = \prod\limits_{y \in Y} \tilde{\Delta}_{f^{-1}(y)}$,
its elements are isomorphism classes of data
$F=\left( \{ F_y \}_{y \in f(X)} ; \{ \bar{F}_x \}_{x \in D(f)} ; \sigma_y \right)$,
with bijections \\
$\sigma_y : \partial F_y \stackrel{\sim}{\rightarrow} \coprod\limits_{x \in f^{-1}(y)}
  \partial \bar{F}_{x}$,
for all $y \in Y$.

We define the operation of \emph{multiplication} by cf. (\ref{multiplication}):
\begin{equation}
\label{eq8.3.44}
\lhd : \tilde{\Delta}_Y \times \tilde{\Delta}_f \rightarrow \tilde{\Delta}_X
\end{equation}
\[G \lhd F = (G_1 ; \bar{G}_y ; \tau) \lhd (F_y ; \bar{F}_x ; \sigma_y) =
   (G_1 \lhd F_{\underline{\tau}(\, \, )} ;
\bar{F}_x \lhd \bar{G}_{f(x)} ;
 \tau \circ \sigma)
\]
with
\begin{equation}\label{eq8.3.45}
 \tau \circ \sigma : \partial (G_1 \lhd F_{\underline{\tau}(\, \,)}) =
\hspace{2.5in}
\end{equation}
\[=
\partial G_1 \times \partial F_{\underline{\tau}(\, \,)}
\stackrel{\stackrel{\textstyle \tau}{\textstyle \sim}}{\rightarrow}
\coprod\limits_{y \in Y} \partial \bar{G}_y \times \partial F_y
\stackrel{\textstyle \amalg \sigma_y}{\stackrel{\textstyle \sim}{\textstyle \rightarrow}}
\coprod\limits_{y \in Y} \coprod\limits_{x \in f^{-1}(y)} \partial \bar{G}_y  \times
\partial \bar{F}_x =\]
\[\hspace{1.8in}=
\coprod\limits_{x \in X} \partial \bar{F}_x  \times
\partial \bar{G}_{f(x)} =
 \coprod\limits_{x \in X} \partial ( \bar{F}_x \lhd \bar{G}_{f(x)} )
\]
We define the operation of \emph{contraction} by, cf. (\ref{contraction}):
\begin{equation}
\label{eq8.3.46}
( \, \, \sslash \,\, ) : \tilde{\Delta}_X \times \tilde{\Delta}_f \rightarrow
\tilde{\Delta}_Y
\end{equation}
\[(G \sslash F ) = \left( (G_1 ; \bar{G}_x ; \tau)\sslash (F_y ; \bar{F}_x ; \sigma_y)\right) =
   \left( G_1 \lhd \bar{F}_{\underline{\tau}(\, \, )}; F_y \lhd \bar{G}_{\underline{\sigma_y}(\, )} ; (\tau ,\sigma) \right)
\]
with
\begin{equation}\label{eq8.3.47}
(\tau , \sigma): \partial (G_1 \lhd \bar{F}_{\underline{\tau}(\,)})=
\partial G_1 \times \partial \bar{F}_{\underline{\tau}(\,)}
\stackrel{\tau}{\stackrel{\sim}{\rightarrow}}
\coprod\limits_{x \in X} \partial \bar{G}_x \times \partial \bar{F}_x =
\end{equation}
\[=
\coprod\limits_{y \in Y} \coprod\limits_{x \in f^{-1}(y)} \partial \bar{F}_x  \times
\partial \bar{G}_x
\stackrel{\amalg \sigma_y^{-1}}{\stackrel{\sim}{\rightarrow}}
\coprod\limits_{y \in Y}  \partial F_y \times \partial \bar{G}_{\underline{\sigma_y}(\,)}=
\coprod\limits_{y \in Y}  \partial (F_y \lhd \bar{G}_{\underline{\sigma_y}(\,)})
\]

It is straightforward to check that $\tilde{\Delta}$ with these
operations is a (non-commutative!) \emph{generalized-ring}, in the
sense that all the axioms of a generalized ring \emph{except} the
axioms of right-linearity
 (\ref{comm_identity}) are satisfied.

\vspace{10pt}

We let $\Delta_X = \tilde{\Delta}_X / \varepsilon$ denote the
quotient commutative generalized ring, where $\varepsilon$ is the
equivalence ideal generated by right-linearity: $(F,F') \in
\varepsilon_X$ if and only if there exists a path $F = F^0, \ldots
, F^l=F'$, with $\{F^{j-1},F^j\}$ of the form
\begin{equation}\label{eq8.3.48}
\big\{\left((A_j\sslash B_j) \lhd C_j)\sslash D_j,
\big((A_j \lhd f_j^*C_j \big) \sslash g_j^*B_j \big)\sslash D_j \right\}
\end{equation}
where the diagram is
\begin{center}
\begin{equation}
\begin{diagram}
& & & &W_j \prod\limits_{Z_j} Y_j & & &  &\\
 & &&\ldTo^{f^*_jC_j}& &\rdTo^{g^*_jB_j} & & &\\
& & W_j& & & & Y_j & &\\
& \ldTo^{A_j}& &\rdTo^{B_j}_{f_j}& & \ldTo^{C_j}_{g_j}& &\rdTo^{D_j} &\\
[1]&\lTo &  & &Z_j & & & &X \\
\end{diagram}
\end{equation}
\end{center}

\noindent or of the form
\begin{equation}\label{eq8.3.50}
\left\{D_j\sslash ( (A_j\sslash B_j)\lhd C_j ) \ , \
D_j\sslash ((A_j\lhd f^*_j C_j)\sslash g_j^*B_j)
\right\}
\end{equation}
where the diagram is
\begin{center}
\begin{equation}
\begin{diagram}
& & & &W_j \prod\limits_{Z_j} Y_j & & &  &\\
 & &&\ldTo^{f^*_jC_j}& &\rdTo^{g^*_jB_j} & & &\\
& & W_j& & & & Y_j & &\\
& \ldTo^{A_j}& &\rdTo^{B_j}_{f_j}& & \ldTo^{C_j}_{g_j}& &\rdTo^{D_j} &\\
X&\lTo &  & &Z_j & & & &[1] \\
\end{diagram}
\end{equation}
\end{center}

Note that if $\{F,F'\} \in \varepsilon$,
then $\{E \lhd F , E \lhd F' \}$, \ \
$\{F \lhd E , F' \lhd E\}$, \ \
$\{E\sslash F) , E\sslash  F'\}$, \ \
$\{F\sslash E,F'\sslash E\}$
are in $\varepsilon$,
so that $\varepsilon$ is indeed an \emph{equivalence ideal},
and the operations of multiplication and contraction
descend to well defined operations on the $\varepsilon$-equivalence
classes $\Delta = \tilde{\Delta}/\varepsilon$.\vspace{1.5mm}\\

We get in this way a commutative generalized ring $\Delta^W$,
\begin{equation}
\label{eq2.6.4}
\Delta_X^W = \left\{\begin{array}{r} ( F_1 ; \{\bar{F}_x \}_{x \in X} ; \sigma )
             /\hspace{-3pt}\approx \, , \,  F_1 , \bar{F}_x \, \,  W\mbox{-labeled trees} \, , \\
              \sigma : \partial F_1 \stackrel{\sim}{\rightarrow}
                \coprod \partial \bar{F}_x \, \mbox{bijection}
             \end{array}\right\}
\end{equation}
The element $\delta^W = \left( [\{  0 \} \amalg W] ; \{ 0_w \}_{w \in W} ;
            \sigma  \right) \in \Delta_W^W$
generates $\Delta^W$.\\
For any commutative generalized ring $A$, we have a bijection
\begin{equation}
\label{eq8.3.53}
\cG\cR_C (\Delta^W , A) = A_W
\end{equation}
\[\begin{array}{c}\varphi \mapsto \varphi (\delta^W)\\
             \varphi^{(a)} \mapsfrom a
\end{array}
\]

\vspace{10pt}

Given an injection $j: W \hookrightarrow W'$, every $W$-labelled
tree $F$ is naturally $W'$-labelled, and we have an injective
homomorphism
\begin{equation}
\Delta^j \in \cG \cR (\Delta^W , \Delta^{W'} ) ,
\Delta^j_X (\{F_1\}_W ; \{\bar{F}_x\}_W ; \sigma ) =
(\{F_1\}_{W'} ; \{\bar{F}_x\}_{W'} ; \sigma )
\end{equation}
It is dual via (\ref{eq8.3.53}) to the map
$A_{W'} \twoheadrightarrow A_W$, $a \mapsto a \lhd 1_j$.

Conversely, if $F$ is a $W'$-labelled tree, and
\begin{equation}
B_F = \{ b \in \partial F , \mu(\mathcal{S}^n (b) ) \in W \mbox{for} \,
         n=0, \ldots , ht(b)-1 \}
\end{equation}
than the reduced tree $F|_{B_F}$, c.f. (\ref{eq8.3.39}), is
$W$-labelled. We have a surjective homomorphism
\begin{equation}
\Delta^{j^t} \in \cG \cR (\Delta^{W'} , \Delta^W ) ,
\end{equation}
\[
\Delta^{j^t}_X (\{F_1\}_{W'} ; \{\bar{F}_x\}_{X'} ; \sigma ) =
(\{F_1|_B\}_{W} ; \{\bar{F}_x|_{\sigma(B) \cap \partial \bar{F}_x}\} ;
           \sigma|_B )
\]
\[ \mbox{with} \,  B=B_{F_1} \cap \sigma^{-1} (\coprod\limits_{x \in X} B_{\bar{F}_x})
\]
It is dual to the map $A_W \hookrightarrow A_{W'}$,
$a \mapsto a \lhd 1_{j^t}$, and we have \\
$\Delta^{j^t} \circ \Delta^j = id_{\Delta^W}$. Thus the map $W\mapsto \Delta^W$ is a functor
\begin{equation}
\Delta:\FF\rightarrow \cG\cR_C
\end{equation}

Fixing a family $\{W_i\}_{i\in I} W_i\in \FF$, we say that a tree
$F$ is \emph{$\{W_i\}$- labelled}, if for all $b \in F \setminus
\partial F$ we are given $j(b) \in I$, and an injection $\mu_b :
\mathcal{S}^{-1}_F(b) \hookrightarrow W_{j(b)}$.

Repeating the above construction, with $\{W_i \}$-labelled trees,
we get a generalized ring $\Delta^{\{W_i\}}$. The elements
$\delta^{W_i} = \left( [\{0\} \amalg W_i] ; \{ 0_w \}_{w \in W_i}
; \sigma \right) \in \Delta_{W_i}^{\{W_i\}}$ generate
$\Delta^{\{W_i\}}$, and we have for any commutative generalized
ring $A$,
\begin{equation}
\label{8.3.54}
\cG \cR_C (\Delta^{\{W_i\}} , A) = \prod_{i\in I} A_{W_i}
\end{equation}
\[\varphi \mapsto \left(\varphi ( \delta^{\{W_i\}})\right)
\]
(i.e. $\Delta^{W_1 \cdots W_N} = \Delta^{W_1} \bigotimes\limits_{\mathbb{F}}
                             \cdots
                        \bigotimes\limits_{\mathbb{F}} \Delta^{W_N}$
cf. \S \ref{13.1} for the tensor product $\otimes$ of commutative
generalized rings.).

Thus we have the adjunction:
\begin{equation}
\xymatrix{\Delta^{\{W_i\}}&\cG\cR_C \ar@/^/@{->}[d]^U& A \ar@{|->}[d] \\
 {\begin{matrix} I\rightarrow \FF \\i\mapsto W_i \end{matrix}}\ar@{|->}[u]& Set/\FF \ar@/^/@{->}[u]^{\Delta}& UA=\coprod_{X\in \FF}A_X}
\end{equation}
For an element $F=(F_1; \{\bar{F}_x\}; \sigma)\in \Delta_X^{\{W_i\}}$, we have its (well-defined!) degree
\begin{equation}
\deg F= \underset{b\in \partial F_1}{\Max} \{\Ht(b)+\Ht(\sigma(b))\}\in \NN
\end{equation}
\begin{equation}
\deg F=0\iff F\in \FF_X\subseteq \Delta_X^{\{W_i\}}
\end{equation}
For $f\in \Set(Y,Z)$, $G=\{G^{(z)}\} \in \Delta_f^{\{W_i\}}$, putting $\deg G=\underset{z\in Z}{\Max}\{\deg G^{(z)}\}$, we have: \begin{equation}
\begin{split}
\deg(F_Z\lhd G)\leq \deg F_Z+\deg G,\; F_Z\in \Delta_Z^{\{W_i\}} \\
\deg(F_Y\sslash G)\leq \deg F_Y+\deg G,\;  F_Y\in \Delta_Y^{\{W_i\}}
\end{split}
\end{equation}

\vspace{20pt}

Given a map $f \in Set_{\bullet}(Z,W)$, we have the element
\begin{equation}
\begin{array}{r}
\delta^f = ( [ \{0\} \amalg f(Z) \amalg D(f) ] ;
\{ [0_z]\}_{z \in D(f)} ; \sigma = id_{D(f)} ) \in
\Delta_Z^{W, f^{-1}(w)} =\\
 = (\Delta^W \otimes \bigotimes_{w \in W}
\Delta^{f^{-1}(w)})_Z
\end{array}
\end{equation}
where the tree $F = \{0\} \amalg f(Z) \amalg D(f)$, has
$\mathcal{S}_F|_{D(f)}=f$, $\mathcal{S}_F(f(Z)) \equiv 0$, and is
labelled by $\mu_0 : f(Z) \hookrightarrow W$, and $\mu_w =
id_{f^{-1}(w)}$ for $w \in f(Z)$. This element gives a
homomorphism of generalized rings, co-multiplication,
\begin{equation}
\Delta_{\lhd}^f \in \cG\cR_C ( \Delta^Z , \Delta^{W,f^{-1}(w)}),\;\;\;\;\delta^Z\mapsto \delta^f,
\end{equation}
which is dual to multiplication
$A_W \times \prod\limits_{w \in W} A_{f^{-1}(w)} \rightarrow A_Z$.

On the other hand we have the element
\begin{equation}
\begin{array}{r}
\varepsilon^f = ( [ \{0\} \amalg D(f) ] ;
\{ [ \{0\} \amalg f^{-1}(w) ] \}_{w \in f(Z)};
\sigma = id_{D(f)}) \in
\Delta_W^{Z, f^{-1}(w)}=\\
=
( \Delta^Z \otimes \bigotimes_{w \in W} \Delta^{f^{-1}(w)} )_W
\end{array}
\end{equation}
giving rise to a homomorphism of generalized rings,
co-contraction,
\begin{equation}
\Delta^f_{\sslash} \in \cG \cR_C (\Delta^W , \Delta^{Z,f^{-1}(w)}),\;\;\;\;\delta^W\mapsto  \varepsilon^f
\end{equation}
which is dual to contraction
$A_Z \times \prod\limits_{w \in W} A_{f^{-1}(w)} \rightarrow A_W$.

\vspace{10pt}

The functor $\Delta : \FF\rightarrow  \cG \cR_C  $, with its
structure (of co-multiplication, co-contraction, and co-unit) is
thus a co-generalized-ring-object in the tensor category $( \cG
\cR_C , \bigotimes\limits_{\FF})$, i.e. the dual of our axioms are
satisfied. (Just as the polynomial ring $\mathbb{Z} [X]$, with
co-multiplication $\Delta_{\bullet} (X) = X \otimes X$, and
co-addition $\Delta_+ (X) = X \otimes 1 + 1 \otimes X$, is a
co-ring object in the tensor category of commutative rings and
$\bigotimes\limits_{\mathbb{Z}}$).

\exmpl{8.3.1} Taking for $W =[1]$, the unit set, a
$[1]$-labelled-tree is just a ''ladder'' $\{x_0, x_1, \ldots , x_n
\}$, $S(x_j) = x_{j-1}$, and is determined by its length $n$. Thus
the element $F = ([F_1]; [\bar{F}_x] ; \sigma ) \in
\Delta_X^{[1]}$, is determined by the length  $n$ of $F_1$, the
point $x \in X$ such that $\bar{F}_x \neq 0$, and the length $m$
of $\bar{F}_x$. We write $F = (x, z^n \cdot (z^t)^m)$, and we have
an isomorphism
\begin{equation}
\Delta^{[1]} \xrightarrow{\sim}{\FF}\{z^{\mathbb{N}} \cdot (z^t)^{\mathbb{N}}\}
\end{equation}
with the generalized ring associated with the free
monoid-with-involution on one element, $M = z^{\mathbb{N}} \cdot
(z^t)^{\mathbb{N}} \cup \{ 0 \}$.

The self-adjoint quotient $\Delta_+^{[1]}$ of $\Delta^{[1]}$ is
isomorphic to the generalized ring associated with
the free monoid on one element \\
$M = z^{\mathbb N} \cup \{0 \}$,
\begin{equation}
\label{eq8.3.63}
\Delta_+^{[1]} \xrightarrow{\sim}
{\mathbb F} \{z^{\mathbb N}\}
\end{equation}

\subsection{Limits}
\label{Limits}
Given a partially ordered set $I$, a functor $A \in (\cG \cR)^I$ is given by objects
$A^{(i)}\in \cG \cR $ for $i \in I$,
and homomorphism $\varphi_{i,j}: A^{(j)} \rightarrow A^{(i)}$ for $i \leq j , i,j \in I$,
such that $\varphi_{i,j} \circ \varphi_{j,k} = \varphi_{i,k}$ for $i \leq j \leq k $, and
$\varphi_{ii} = id_{A^{(i)}}$.
The inverse limit of $A$ exists,
and can be computed in $Set_0$. We put
\begin{equation}
(\lim\limits_{\stackrel{\longleftarrow}{I}} A^{(i)})_X =
\{ a= (a^{(i)} ) \in \prod\limits_{i \in I } A_X^{(i)} ,
     \varphi_{i,j}(a_j ) = a_i \mbox{for all} \, i \leq j  \}
\end{equation}
With the operations of componentwise multiplication and
contraction $\lim\limits_{\longleftarrow} A^{(i)}$ is a
generalized-ring (sub-ring of $\prod\limits_{i \in I} A^{(i)}$).
We have the universal property
\begin{equation}
\cG \cR (B, \lim\limits_{\stackrel{\longleftarrow}{I}} A^{(i)} ) =
\lim\limits_{\stackrel{\longleftarrow}{I}} \cG \cR (B, A^{(i)}) =
\end{equation}
\[=
\{ (\psi_i ) \in \prod\limits_{I} \cG \cR  (B, A^{(i)}) ,
\varphi_{i,j} \circ \psi_j = \psi_i \, \mbox{for} \, i \leq j \}
\]

\vspace{10pt}

If the set $I$ is \emph{directed} (for $j_1 , j_2 \in I$ there is
$i \in I$ with $i \leq j_1, j_2$) the direct limit of $A$ exists,
and can be computed in $Set_0$. We have
\begin{equation}
\label{eq2.7.3}
(\lim\limits_{\stackrel{\longrightarrow}{I}} A^{(i)})_X =
\lim\limits_{\stackrel{\longrightarrow}{I}} A^{(i)}_X =
(\coprod\limits_{i \in I} A_X^{(i)}) / \{a \sim \varphi_{i,j} (a)\}
\end{equation}
There are well defined operations of multiplication and
contraction making $\lim\limits_{\stackrel{\longrightarrow}{I}}
A^{(i)}$ into a generalized ring. We have the universal property
\begin{equation}
\cG \cR (\lim\limits_{\stackrel{\longrightarrow}{I}} A^{(i)} , B) =
\lim\limits_{\stackrel{\longleftarrow}{I}} \cG \cR (A^{(i)}, B) =
\end{equation}
\[=
\{ (\psi_i ) \in \prod\limits_{I} \cG \cR  (A^{(i)}, B) ,
\psi_i \circ \varphi_{i,j} = \psi_j  \, \mbox{for} \, i \leq j \}
\]
We shall see below that $\mathcal{GR}_C$ has sums and push-outs \S 13.1, and quotients by equivalence-ideals \S 9.1, hence arbitrary co- limits, and we have:
\thm{8.3.8} $\mathcal{GR}_C$ is complete and co-complete.

\chapter{Ideals}
\smallbreak
\setcounter{equation}{0}
\section{Equivalence ideals}

\defin{9.1.1}
For  $A \in \cG \cR$ an \emph{equivalence ideal} $\varepsilon$  is
a sub-generalized-ring $\varepsilon \subseteq A \times A$, such
that $\varepsilon_X$ is an equivalence relation on $A_X$, (we
write $a \sim a'$ for $(a,a') \in \varepsilon_X$; and for $a =
(a_y)$, $a' = (a'_y) \in A_f$, we write $a \sim a'$ for $(a_y ,
a'_y ) \in \varepsilon_{f^{-1}(y)}$
 for all $y \in Y$).
Thus, we have an equivalence relation $\sim$ that respects the
operations: if $a \sim a'$ and $b \sim b'$ then
\begin{equation}
\label{eq3.1.2}
\begin{array}{l}
a \lhd b \sim a' \lhd b' , \\
a\sslash  b \sim a'\sslash  b' .
\end{array}
\end{equation}
(whenever these operations are defined.)

We let $eq(A)$ denote the set of equivalence ideals of $A$.

For $ \varepsilon \in eq(A)$ we can form the quotient $A/\varepsilon$, with
\begin{equation}
(A/\varepsilon)_X = A_X/\varepsilon_X = \varepsilon_X\mbox{-equivalence classes in}
\, A_X
\end{equation}
There is a natural surjection $\pi_X : A_X \rightarrow A_X/\varepsilon_X$,
and there is a unique structure of a generalized ring on
$A/\varepsilon$ such that $\pi \in \cG\cR (A,A/\varepsilon)$.

\vspace{10pt}

\defin{9.1.2} For $\varphi \in \cG\cR (A,A')$, we let
$KER (\varphi) = A \prod\limits_{A'} A$,
\begin{equation}
KER(\varphi)_X = \{ (a_1, a_2) \in A_X \times A_X , \varphi (a_1 ) = \varphi_X(a_2) \}
\end{equation}
$KER (\varphi )$ is an equivalence ideal.

\vspace{10pt}

\noindent We have the universal property of the quotient
\begin{equation}
\cG\cR (A/\varepsilon, A' ) = \{ \varphi \in \cG\cR (A, A') ,
KER ( \varphi)  \supseteq \varepsilon \}
\end{equation}

\vspace{10pt}

\noindent Every homomorphism
$ \varphi \in \cG\cR (A, A')$ has a

 canonical factorization
\begin{equation}
\varphi = j \circ \tilde{\varphi} \circ \pi
\end{equation}
\begin{center}
\begin{diagram}
A & \rTo^{\varphi}& A'\\
\dOnto^{\pi}& &\uInto_{j}\\
A/KER(\varphi) & \rTo_{\sim}^{\widetilde{\varphi}}& \varphi(A)
\end{diagram}
\end{center}

with $\pi$ surjection, $j$ injection, and $\tilde{\varphi}$ an isomorphism.

\setcounter{equation}{0}
\section{Functorial Ideals}
\defin{9.2.1} For $A \in \cG\cR_C$,
a \emph{functorial ideal} $\sqa$ is a collection of subsets $\sqa_X \subseteq A_X$,
with $0_X \in \sqa_X$, and with $A \lhd \sqa$, $\sqa \lhd A$, $A\sslash \sqa$, $\sqa\sslash A \subseteq \sqa$. \\
We let $fun\cdot il(A)$ denote the set of
functorial  ideals of $A$.

\vspace{10pt}

For $\varepsilon \in eq(A)$, we have the associated
functorial ideal $Z(\varepsilon)$:
\begin{equation}
Z(\varepsilon)_X = \left\{a \in A_X , (a,0_X) \in \varepsilon_X \right\}
\end{equation}
For $\sqa \in fun \cdot il(A)$, we have the equivalence ideal $E(\sqa)$ generated by $\sqa$,
it is the intersection of all equivalence ideals containing $(a,0)$ for all
$a \in \sqa$.

These give a Galois correspondence:
\begin{equation}\label{galois_ideals}
fun \cdot il(A)
\begin{array}{c}
  \stackrel{\scriptstyle E}{\longrightarrow} \\
  \stackrel{\textstyle \longleftarrow}{\scriptstyle Z}
\end{array}
eq(A)
\end{equation}
It is monotone,
\begin{equation}\label{eq9.2.3}
\begin{array}{ll}
\sqa_1 \subseteq \sqa_2 & \Rightarrow E(\sqa_1) \subseteq E(\sqa_2)\\
\varepsilon_1 \subseteq \varepsilon_2 & \Rightarrow Z(\varepsilon_1) \subseteq Z(\varepsilon_2)
\end{array}
\end{equation}
and we have
\begin{equation}
\sqa \subseteq ZE(\sqa) \quad , \quad EZ(\varepsilon) \subseteq \varepsilon
\end{equation}
It follows that we have
\begin{equation}
ZEZ(\varepsilon) = Z(\varepsilon) \quad , \quad E(\sqa) = EZE(\sqa) ,
\end{equation}
and the maps $E$, $Z$ induce inverse bijections
\begin{equation}
E\mbox{-}il(A) \stackrel{\textstyle \sim}{\longleftrightarrow}
Z\mbox{-}eq(A)
\end{equation}
with
\begin{equation}
E\mbox{-}il(A) = \{ \sqa \in fun \cdot il(A) , \sqa = ZE(\sqa)\} =
          \{ Z(\varepsilon) , \varepsilon \in eq(A)\} ,
\end{equation}
the \emph{stable} functorial ideals, and
\begin{equation}
Z\mbox{-}eq(A) = \{ \varepsilon \in eq(A) , \varepsilon = EZ(\varepsilon)\} =
         \{E(\sqa) , \sqa \in il(A)\}
\end{equation}

Let $\sqa \in fun \cdot il(A)$, and let $\varepsilon= (\varepsilon_X)$ be defined by
\begin{equation}\label{eq9.2.9}
\begin{split}
\varepsilon_X=& \\
=&\left\{\begin{array}{l}
  (a, a') \in A_X \times A_X, \mbox{there exists a \emph{path}}\,
a = a_1 , a_2, \ldots , a_n =a',\\
 \mbox{with} \,
\{a_j , a_{j+1} \} \,  \mbox{of the form
   either} \, \{ (d_j \lhd c_j) \sslash b_j, (d_j \lhd \bar{c}_j)\sslash b_j \},\\
\mbox{or} \ \{ d_j\sslash (b_j \lhd c_j ) , d_j \sslash (b_j \lhd \bar{c}_j ) \}
\ \mbox{with} \
    d_j \in A_{Y_j} , b_j \in A_{f_j} ,\\ f_j \in Set_{\bullet}(Z_j , X) ,
    c_j , \bar{c}_j \in A_{g_j} ,
 g_j \in Set_{\bullet}(Z_j , Y_j ) \ \mbox{or resp.} \ \\
g_j \in Set_{\bullet}(Y_j , Z_j )\
   \mbox{and with} \ c_j^{(z)} = \bar{c}_j^{(z)}, \mbox{ or }\
     c_j^{(z)},  \bar{c}_j^{(z)} \in \sqa , \\
\mbox{for all} \ z \in Y_j
    (\mbox{resp.} \, z \in Z_j)
\end{array}\right\}
\end{split}
\end{equation}
\cl{9.2.2}
$$E(\sqa) = \varepsilon$$
\begin{proof}
It is clear  that $\varepsilon_X$ is an equivalence relation on $A_X$,
and that $\varepsilon_X \subseteq E(\sqa)_X$,
and we need to show that $\varepsilon$ respects the operations (\ref{eq3.1.2}).
If $(a,a') \in \varepsilon$,
so there is a path
$a = a_1 , \ldots , a_n =a'$
as in (\ref{eq9.2.9}),
then $h \lhd a_j$ (resp. $a_j \lhd h$, $(h\sslash a_j)$, $(a_j\sslash h)$)
is a path, showing $(h \lhd a, h \lhd a')$
 (resp. $(a \lhd h , a' \lhd h)$,
$(h\sslash a, h\sslash a')$, $(a\sslash h, a'\sslash h)$) is in $\varepsilon$.
This follows from the \vspace{3mm}\\
\textbf{\large Basic identities 9.2.3}
\begin{equation*}
\label{eq3.2.12}
\begin{array}{clcl}
  & h \lhd (d \lhd c\sslash b) = ((h \lhd d) \lhd c)\sslash b & , &
   h \lhd (d\sslash (b \lhd c)) = (h \lhd d)\sslash( b \lhd c)\\
\mbox{resp.}& ((d \lhd c)\sslash b) \lhd h = (d \lhd (c \lhd \tilde{h}))\sslash \tilde{b}
&, & (d\sslash (b \lhd c))\lhd h = (d \lhd \tilde{h})\sslash (\tilde{b} \lhd \tilde{c})\\
& h\sslash ((d \lhd c)\sslash b) = (h \lhd b)\sslash (d \lhd c) & ,& h\sslash (d\sslash (b \lhd c)) =
  ((h \lhd b)\lhd c)\sslash d\\
& ((d \lhd c)\sslash b)\sslash h = (d \lhd c)\sslash ( h \lhd b) &,&
(d\sslash (b \lhd c))\sslash h = d\sslash ((h \lhd b)\lhd c)
\end{array}
\end{equation*}
\end{proof}
\vspace{10pt}

It follows that
a functorial ideal $\sqa$ is stable,
 $\sqa = ZE(\sqa)$,
if and only if for all $b,d,c,\bar{c}$ as in (\ref{eq9.2.9}),
\begin{equation}
\begin{array}{cc}
 & (d \lhd c)\sslash b \in \sqa \Leftrightarrow (d \lhd \bar{c})\sslash b \in \sqa\\
\mbox{and} & d\sslash (b \lhd c) \in \sqa \Leftrightarrow d\sslash (b \lhd \bar{c}) \in \sqa
\end{array}
\end{equation}

\setcounter{equation}{0}
\section{Operations on functorial ideals}

The intersection of functorial ideals is a functorial ideal,
\begin{equation*}
\sqa_i \in fun \cdot il(A) \Rightarrow \bigcap\limits_i \sqa_i \
\in fun \cdot il(A)
\end{equation*}

Given a collection $B = \{ B_X \subseteq A_X \}$,
the functorial ideal generated by $B$ will be denoted by $\{B\}_A$,
it is the intersection of all
functorial
ideals containing $B$.
If $B$ satisfies $B \lhd A \subseteq B$,
it can be described explicitly as
$$ \{B\}_{A,X} = $$
$$= \left\{ a \in A_X , a=(d \lhd c)\sslash b\
\mbox{or} \ a=d\sslash (b \lhd c)\
\mbox{for some}\  d \in A_Y,
              b \in A_f ,
\right. $$
$$ f \in Set_{\bullet} (Z,X),
 c \in A_g ,\; g \in Set_{\bullet}(Z,Y),$$
 $$\mbox{or resp.} \ g \in Set_{\bullet}(Y,Z) $$
\begin{equation}
\begin{split}
\label{eq9.3.1}
\left. \mbox{and with} \ c= (c^{(z)}) , c^{(z)} \in B_{g^{-1}(z)} \,
           \mbox{for all} \, z \in Y \,  (\mbox{resp.} \, z \in Z) \right\}
\end{split}
\end{equation}

It is clear that the set described in (\ref{eq9.3.1}), contains $B_X$,
and is contained in $\{B\}_{A,X}$,
and we only have to check that it is a functorial ideal - this follows from the \emph{basic identities 9.2.3}.

In particular, for $\sqa_i \in fun \cdot il(A)$,
we can take $B=\bigcup\limits_i \sqa_i$, and we obtain the smallest
functorial
ideal
containing all the $\sqa_i$'s.
\begin{equation}
\sum\limits_i \sqa_i = \left\{ \bigcup\limits_i \sqa_i\right\}_A
\end{equation}
Thus $fun \cdot il(A)$ is a complete lattice, with minimal element the zero ideal
$0 = \{0_X\}$,
and maximal element the unit ideal $\{1\}_A = \{A_X\}$.

Note that for arbitrary $B$ we can similarly describe the
functorial ideal it generates as

\begin{equation}
\begin{split}
\label{eq9.3.3}
\{B\}_{A,X} =& \\
=& \left\{a \in A_X , a= (d\lhd c \lhd e)\sslash b
           \,\mbox{or} \, a = d\sslash (b \lhd c \lhd e) , \mbox{with} \,
            c^{(z)}\in B \, \mbox{for all }\, z \right\}
\end{split}
\end{equation}

Note that for a subset $B \subseteq A_{[1]}$, we have
\begin{equation}
\begin{split}
\{B\}_{A,X} = \left\{ a \in A_X , a= d\sslash (b \lhd c) , d \in A_Y ,
              b \in A_f ,
      \right.
\\
\left.
f \in Set_{\bullet}(Y,X),
 c = (c^{(y)}) \in A_{id_Y} =
               (A_{[1]})^Y, c^{(y)}\in B \cup B^t
               \right\}
\end{split}
\end{equation}
We do not need the $e$'s in (\ref{eq9.3.3})
because we have commutativity $c \lhd e = e \lhd \tilde{c}$,
c.f. (\ref{eq_up1}),
and we do not need the two shapes of (\ref{eq9.3.3}) because the
$(A_{[1]})^Y$-action is self-adjoint
$d\sslash (b \lhd c) = (d \lhd c^t )\sslash b$,
c.f. (\ref{eq_up3}).

\setcounter{equation}{0}
\section{Homogeneous ideals}
\defin{9.4.1} A functorial ideal $\sqa \in fun \cdot il(A)$ is called \emph{homogeneous}
if it is generated by $\sqa_{[1]}$.\vspace{3mm}\\
The subset $\sqa_{[1]} \subseteq A_{[1]}$ satisfies for all $X \in \FF$,
\begin{equation}
\label{homideal}
A_X \sslash \left( A_X \lhd (\sqa_{[1]})^X \right) \subseteq \sqa_{[1]}
\end{equation}
Conversely, if a subset $\sqa_{[1]} \subseteq A_{[1]}$ satisfies
(\ref{homideal}), then $\sqa_{[1]}^t = \sqa_{[1]}$ (because for $a
\in \sqa_{[1]}$, $a^t = 1\sslash (1 \lhd a) \in \sqa_{[1]}$), and
$\sqa_{[1]} \lhd A_{[1]} = \sqa_{[1]}$ (because for $a \in
\sqa_{[1]}$, $b \in A_{[1]}$, $a \lhd b = b\sslash (1 \lhd a^t )
\in \sqa_{[1]}$). Moreover, the functorial ideal $\sqb$ generated
by $\sqa_{[1]}$, has
\begin{equation}
\sqb_X = \bigcup\limits_{f \in Set_{\bullet}(Y,X)}
      A_Y \sslash \left( A_f \lhd (\sqa_{[1]})^Y\right)
\end{equation}
and in particular $\sqb_{[1]} = \sqa_{[1]}$.
Thus we identify the set of homogeneous
functorial
ideals with the collection of
subsets $\sqa_{[1]} \subseteq A_{[1]}$ satisfying (\ref{homideal}),
and we denote this set by $[1]\mbox{-}il(A)$,
\begin{equation}
\label{eq9.4.3}
[1]\mbox{-}il(A) = \{ \sqa \subseteq A_{[1]} , A_X\sslash (A_X \lhd (\sqa)^X)
       \subseteq \sqa \}
\end{equation}
The set $[1]\mbox{-}il(A)$ is a complete lattice, with minimal element $\{0\}$,
maximal element $\{1\}_A = A_{[1]}$.
For $\sqa_i \in [1]\mbox{-}il(A)$, we have
$\bigcap\limits_i\sqa_i \in [1]\mbox{-}il(A)$, and
$\sum\limits_i \sqa_i \in [1]\mbox{-}il(A)$,
where
\begin{equation}
\begin{split}
\label{eq9.4.4}
\sum\limits_i \sqa_i =& \\
=&  \left\{ a \in A_{[1]}, a = b\sslash (d \lhd c) ,
              b,d \in A_X , c=(c^{(x)}) \in (\bigcup\limits_i \sqa_i)^X
           \subseteq A_{id_X} , X \in \FF \right\}
\end{split}
\end{equation}
Note that the homogeneous functorial ideal generated by elements
$a_i \in A_{[1]}$, $i \in I$,
can be described as
\begin{equation}
\begin{split}
\label{eq9.4.5}
\{a_i\}_A =& \\
=& \left\{ a \in A_{[1]} , a= b\sslash (d \lhd c) ,
b,d \in A_X , c=(c^{(x)}) \in A_{id_X}, c^{(x)}=a_{i(x)}\right.\\
 & \left.\mbox{or} \, c^{(x)}=  a^t_{i(x)}
\mbox{for} \, x \in X \right\}
\end{split}
\end{equation}
We have also the operation of multiplication of homogeneous ideals.
For $\sqa_1 , \sqa_2 \in [1]\mbox{-}il(A)$,
we let $\sqa_1 \cdot \sqa_2$ denote the homogeneous ideal
generated by the product
$\sqa_1 \lhd \sqa_2 = \{ a_1 \lhd a_2, a_i \in \sqa_i \}$.

Thus
\begin{equation}
\begin{split}
\sqa_1 \cdot \sqa_2 =& \\
=& \left\{a \in A_{[1]} , a=b\sslash (d \lhd c) , b,d \in A_X ,
c=(c^{(x)}) \in (\sqa_1 \lhd \sqa_2)^X \subseteq A_{id_X} \right\}
\end{split}
\end{equation}
This operation is associative, we have for
$\sqa_1, \sqa_2, \sqa_3 \in [1]\mbox{-}il(A)$,
\begin{equation}
(\sqa_1 \cdot \sqa_2) \cdot  \sqa_3 =
 \sqa_1 \cdot \sqa_2 \cdot \sqa_3 =
\sqa_1\cdot ( \sqa_2 \cdot \sqa_3)
\end{equation}
with \\
\begin{equation*}
\begin{split}
\sqa_1 \cdot \sqa_2 \cdot \sqa_3 =& \\
=& \left\{ a \in A_{[1]} , a =b\sslash (d \lhd c) , b,d \in A_X ,
c=(c^{(x)}) \in (\sqa_1 \lhd \sqa_2 \lhd \sqa_3)^X \subseteq A_{id_X} \right\}
\end{split}
\end{equation*}
(use $(b\sslash (d \lhd a_1 \lhd a_2))\lhd a_3 =
b\sslash ( d \lhd a_1 \lhd a_2 \lhd \tilde{a}_3^t)$ and
$a_1 \lhd (b\sslash (d \lhd a_2 \lhd a_3)) =
b\sslash (d \circ a_2 \circ a_3 \circ \tilde{a}_1^t)$).

\vspace{5pt}

The multiplication of homogeneous ideals is clearly commutative,
$\sqa_1 \lhd \sqa_2 = \sqa_2 \lhd \sqa_1$,
has unit element $\{1\}_A= A_{[1]}$,
$\sqa \cdot \{1\}= \sqa$,
and has zero element $\{0\}$, $\sqa \cdot \{0\} = \{0\}$.
Thus
$[1]\mbox{-}il(A) \in Mon$.

\vspace{10pt}

For a homomorphism $\varphi \in \cG\cR_C (A,B)$,
and for $\sqb \in fun \cdot il(B)$,
(resp. $\sqb \in [1]\mbox{-}il(B)$),
its inverse image $\varphi^{*} (\sqb)_X = \varphi_X^{-1}(\sqb_X)$,
(resp.
$\varphi^*_{[1]}(\sqb) = \varphi^{-1}_{[1]}(\sqb) \subseteq A_{[1]}$) is clearly a (resp. homogeneous) functorial
ideal of $A$.
For $\sqa \in fun \cdot il(A)$,
(resp. $[1]\mbox{-}il(A)$)
let $\varphi_*(\sqa) \subseteq B$ denote the (homogeneous)
functorial ideal generated by the image $\varphi(\sqa)$ (resp.
$\varphi_{[1]}(\sqa)$).
We have  Galois correspondences
\begin{equation}
[1]\mbox{-}il(A)
\begin{array}{c}
\stackrel{\textstyle \varphi^*_{[1]}}{\longleftarrow}\\
\stackrel{\textstyle\longrightarrow}{\varphi_{*[1]}}
\end{array}
[1]\mbox{-}il(B)
\quad
\mbox{and} \quad
fun \cdot il(A)
\begin{array}{c}
\stackrel{\textstyle \varphi^*}{\longleftarrow}\\
\stackrel{\textstyle \longrightarrow}{\varphi_*}
\end{array}
fun \cdot il(B)
\end{equation}
The maps $\varphi^*$, $\varphi_*$ are monotone, and satisfy
\begin{equation}
\label{eq3.5.9}
\sqa \subseteq \varphi^* \varphi_* (\sqa) ,
\quad
 \varphi_* \varphi^* (\sqb) \subseteq \sqb
\end{equation}
It follows that we have
\begin{equation}
\varphi_*(\sqa) = \varphi_*\varphi^*\varphi_*(\sqa) \quad , \quad
\varphi^*(\sqb) =  \varphi^*\varphi_* \varphi^* (\sqb)
\end{equation}
and $\varphi^*$, $\varphi_*$ induce inverse bijections,
\begin{equation}
\begin{split}
\left\{ \sqa \in [1]\mbox{-}il(A) ,\sqa = \varphi^*\varphi_* (a) \right\} &= \\
=\left\{ \varphi^* (\sqb) , \sqb \in [1]\mbox{-}il(B) \right\} & \stackrel{\sim}{\longleftrightarrow}
\left\{ \varphi_* (\sqa) , \sqa \in[1]\mbox{-}il(A) \right\}=\\
= & \left\{ \sqb \in [1]\mbox{-}il(B) , \sqb = \varphi_*\varphi^* (\sqb)\right\}
\end{split}
\end{equation}
and similarly with $fun \cdot il(A)$ and $fun \cdot il(B)$.

\vspace{10pt}

\defin{9.4.2}
\setcounter{equation}{11}
For an equivalence ideal $\varepsilon \in eq(A)$,
and for a functorial ideal or a homogeneous ideal $\sqa$,
we say $\sqa$ is $\varepsilon$-\emph{stable} if for all
$(a,a') \in \varepsilon:$
$$a \in \sqa \Leftrightarrow a' \in \sqa .$$
We denote by $fun \cdot il(A)^{\varepsilon}$
(resp. $[1]\mbox{-}il(A)^{\varepsilon}$)
the set of $\varepsilon$-stable (homogeneous) ideals.

Letting $\pi_{\varepsilon} : A \twoheadrightarrow A/\varepsilon$
denote the canonical homomorphism, we have bijections
\begin{equation}
\begin{array}{c}
fun\cdot il(A)^{\varepsilon} \stackrel{ \sim}{\longleftrightarrow}
fun\cdot il(A/\varepsilon) \ \mbox{and} \
{[1]}\mbox{-}il(A)^{\varepsilon} \stackrel{\sim}{\longleftrightarrow}
{[1]}\mbox{-}il(A/\varepsilon)\\
\sqa \mapsto \pi_{\varepsilon}(\sqa)\\
\pi^{-1}_{\varepsilon}(\sqb) \mapsfrom \sqb
\end{array}
\end{equation}
\defin{9.4.3}
\setcounter{equation}{12}
For a (resp. homogeneous) functorial ideal $\sqa$,
we say $\sqa$ is \emph{stable} if it is $E(\sqa)$-stable.
We denote by $E\mbox{-}fun\cdot il(A)$, (resp. $E[1]\mbox{-}il(A)$), the set of stable (homogeneous) functorial ideals.
Note that by the explicit description of $E(\sqa)$ in
Claim 9.2.2,
a subset $\sqa \subseteq A_{[1]}$ is a stable homogeneous ideal if and only if
$$
\mbox{for} \; X,Y \in \FF , b, d \in A_{X \oplus Y},
   c, \bar{c} \in (A_{[1]})^{X \oplus Y},
$$
\[
\mbox{with} \;
c^{(x)}= \bar{c}^{(x)}\;\;
\mbox{for} \; x \in X
\mbox{and} \; c^{(y)} , \bar{c}^{(y)} \in \sqa\;\;
\mbox{for} \; y \in Y,
\]
\begin{equation}
 \mbox{have}:  b\sslash (d \lhd c) \in \sqa \Leftrightarrow
b\sslash (d \lhd \bar{c}) \in \sqa
\end{equation}

(taking $X = [0]$, $\bar{c}^{(y)}\equiv 0$,
we see that this condition includes $\sqa$ being a homogeneous ideal).

\section{Ideals and symmetric ideals}
\setcounter{equation}{0}
\defin{9.5.1}
For $A \in \cG\cR_C$, a subset $\sqa \subseteq A_{[1]}$ will be called \emph{ideal}
if for all $X \in \FF$,
$b,d \in A_X$, $c = (c^{(x)}) \in (\sqa)^X \subseteq A_{id_X}$,
we have
\begin{equation}
(b \lhd c)\sslash d \in \sqa
\end{equation}
We denote by $il(A)$ the set of ideals of $A$.

Comparing with the definition of homogeneous ideals (\ref{homideal}),
we see that the homogeneous ideals are precisely the self-adjoint ideals
\begin{equation}
[1]\mbox{-}il(A) = \{ \sqa \in il(A) , \sqa = \sqa^t \}
\end{equation}
The set $il(A)$ is a complete lattice,
with minimal element $(0)$,
maximal element $(1)=A_{[1]}$.
For $\sqa_i \in il(A)$, we have
$\bigcap\limits_i \sqa_i \in il(A)$,
and $\sum\limits_i \sqa_i \in il(A)$,
where (c.f. (\ref{eq9.4.4}))
\begin{equation}
\sum\limits_i \sqa_i =
\left\{ a \in A_{[1]}, a = (b \lhd c)\sslash d , b,d \in A_X ,
c = (c^{(x)}) \in (\bigcup\limits_i a_i )^X \subseteq A_{id_X} \right\}
\end{equation}
Note that the ideal generated by elements $a_i \in A_{[1]}$,
$i \in I$, can be described as (c.f. (\ref{eq9.4.5}))
\begin{equation}
(a_i)_A=
\left\{ a \in A_{[1]}, a = (b \lhd c)\sslash d , b,d \in A_X ,
c = (c^{(x)}) \in \{a_i\}^X \subseteq A_{id_X} \right\}
\end{equation}
In particular, for $a \in A_{[1]}$ the \emph{principal}
ideal generated by $a$ is just
\begin{equation}
(a)_A = a \lhd A_{[1]}
\end{equation}
Indeed, if $c^{(x)}= a \lhd e_x$, $e_x \in A_{[1]}$ for
$x \in X$,
than for $b,d \in A_X$,
\begin{equation}
(b \lhd c)\sslash d  = a \lhd ((b \lhd e)\sslash d) \in a \lhd A_{[1]}
\end{equation}
(while the homogeneous ideal generated by $a$ is the ideal
generated by $a$ and $a^t$,
c.f. (\ref{eq9.4.5}) for $\{a\}_A$).

\defin{9.5.2}
An ideal $\sqa \in il(A)$
will be called \emph{symmetric} if it is generated by
its subset of self-adjoint elements
\begin{equation}
\sqa^+ = \sqa \cap A^+_{[1]} =
\{a \in \sqa \ , \ a = a^t\}
\end{equation}
Such an ideal is clearly homogeneous.
We write $il^t(A)$ for the collection of symmetric ideals.
It is again a complete lattice.
We have multiplication of ideals, for
$\sqa_1 , \sqa_2 \in il(A)$,
\[
\sqa_1 \cdot \sqa_2 =
\]
\begin{equation}
=\left\{ a \in A_{[1]},
a = (b \lhd c)\sslash d, \, b,d \in A_X ,
c = (c^{(x)}) \in (\sqa_1 \lhd \sqa_2)^X \subseteq A_{id_X}
\right\}
\end{equation}
It is associative,
\begin{equation}
(\sqa_1 \cdot \sqa_2) \cdot \sqa_3 =
\sqa_1 \cdot \sqa_2 \cdot \sqa_3 =
\sqa_1 \cdot (\sqa_2 \cdot \sqa_3 )
\end{equation}
\[  \hspace{-1in} \mbox{with} \;
\sqa_1 \cdot \sqa_2 \cdot \sqa_3 =
\left\{ a \in A_{[1]},
a = (b \lhd c)\sslash d,  b,d \in A_X ,\right.\]
\[\hspace{2in} \left.
c = (c^{(x)}) \in (\sqa_1 \circ \sqa_2 \circ \sqa_3)^X \subseteq A_{id_X}
\right\}\]

(Use now:
$((b \lhd a_1 \lhd a_2 )\sslash d )\lhd a_3 =
(b \lhd a_1 \lhd a_2 \lhd \tilde{a}_3) \sslash d$,
and
$a_1 \lhd ((b \lhd a_2 \lhd a_3)\sslash d)  =
 (b \lhd \tilde{a}_1 \lhd a_2 \lhd a_3) \sslash d $).

It is clearly commutative, $\sqa_1 \cdot \sqa_2 = \sqa_2 \cdot \sqa_1$;
has unit $(1)$,
$\sqa \cdot (1) = \sqa$;
has zero $(0)$, $\sqa \cdot (0) = (0)$.

Thus, $il(A)$ is a monoid with involution, $[1]\mbox{\text{-}}il(A)$
is its submonoid of element fixed by the involution,
and $il^t(A)$ is a submonoid of $[1]\mbox{\text{-}}il(A)$.

We can also divide ideals. For $\sqa_0 , \sqa_1 \in il(A)$,
we let
\begin{equation}
(\sqa_0 : \sqa_1 ) =
\{ c \in A_{[1]} , c \lhd \sqa_1 \subseteq \sqa_0 \}
\end{equation}
This is an ideal,
$(\sqa_0 : \sqa_1 ) \in il(A)$.
Indeed, for $c^{(x)} \in (\sqa_0 : \sqa_1 )$,
and for any $b,d \in A_X$,
and any $a_1 \in \sqa_1$,
we have
\[((b \lhd (c^{(x)}))\sslash d) \lhd a_1 =
(b\lhd (c^{(x)} \lhd a_1 )) \sslash d \in \sqa_0
\]
and so $(b \lhd (c^{(x)}))\sslash d  \in (\sqa_0 : \sqa_1 )$.

\vspace{10pt}

For elements $m_1 , m_2 \in A_X$, we have their \emph{annihilator}
\begin{equation}\label{9.5.11}
ann_A (m_1, m_2 ) =
\{a \in A_{[1]} , a \lhd m_1 = a \lhd m_2 \}
\end{equation}
This is an ideal,
$ann_A(m_1, m_2 ) \in il(A)$.
Indeed for $c^{(y)} \in ann_A (m_1, m_2 )$,
and for any $b,d \in A_Y$ we have
\[ ((b \lhd (c^{(y)}))\sslash d) \lhd m_1 =
(b \lhd (c^{(y)}\lhd m_1 ))\sslash  \tilde{d}=
(b \lhd (c^{(y)}\lhd m_2))\sslash \tilde{d}=
((b \lhd (c^{(y)}))\sslash d) \lhd m_2
\]
and so
$(b \lhd (c^{(y)}))\sslash d\in ann_A (m_1, m_2)$

The ideal generated by the symmetric elements
$ann_A(m_1, m_2) \cap A^+_{[1]}$
is a subideal:
\begin{equation}\label{9.5.12}
ann_A^t(m_1, m_2) \subseteq ann_A(m_1, m_2),
\end{equation}
it is clearly a symmetric ideal.

\vspace{10pt}

Let $\varphi \in \cG\cR_C(A,B)$.
For an ideal $\sqb \in il(B)$,
its inverse image $\varphi^*(\sqb) = \varphi^{-1}_{[1]}(\sqb) \subseteq A_{[1]}$
is clearly an ideal of $A$.
For $\sqb \in il^t(B)$,
its inverse image $\varphi^*(\sqb)$ is the ideal of $A$ generated by
\begin{equation}
\{a=a^t \in A^+_{[1]} \ , \ \varphi_{[1]}(a) \in \sqb \}
\end{equation}
For an ideal $\sqa \in il(A)$ (resp. a symmetric ideal $\sqa \in il^t(A)$),
we let $\varphi_*(\sqa) \subseteq B_{[1]}$ denote the
ideal generated by the image $\varphi_{[1]}(\sqa)$,
$\varphi_*(\sqa) =
\lim\limits_{\stackrel{\textstyle \rightarrow}{X}}
(B_X \lhd ( \varphi_{[1]}(\sqa))^X \sslash B_X )$, (resp. by
$\varphi_{[1]}(\sqa^+)$).

We have a Galois correspondence
\begin{equation}\label{galois_ideals2}
il(A)
\begin{array}{c}
\stackrel{\textstyle \varphi^*}{ \longleftarrow}\\
\stackrel{\textstyle \longrightarrow}{\varphi_*}
\end{array}
il(B)
\ \ \mbox{and} \  \
il^t(A)
\begin{array}{c}
\stackrel{\textstyle \varphi^*}{ \longleftarrow}\\
\stackrel{\textstyle \longrightarrow}{\varphi_*}
\end{array}
il^t(B)
\end{equation}
The maps $\varphi^*$, $\varphi_*$ are monotone, and satisfy
\begin{equation}
\sqa \subseteq \varphi^* \varphi_* (\sqa)
\quad, \quad
\varphi_*\varphi^*(\sqb) \subseteq \sqb
\end{equation}
It follows that we have
\begin{equation}
\varphi_*(\sqa) = \varphi_*\varphi^*\varphi_*(\sqa)
\quad ,\quad
\varphi^*(\sqb)=\varphi^*\varphi_*\varphi^*(\sqb)
\end{equation}
and $\varphi^*$, $\varphi_*$ induce inverse bijections,
$$
\hspace{-4cm}
\{\sqa \in il(A) , \sqa=\varphi^*\varphi_* (\sqa) \}=
$$
\begin{equation}=
\{\varphi^*(\sqb) , \sqb \in il(B) \}
\stackrel{\sim}{\longleftrightarrow}
\{\varphi_*(\sqa) , \sqa \in il(A)\} =
\end{equation}
\[\hspace{4cm}
=
\{\sqb\in il(B) , \sqb = \varphi_* \varphi^*\varphi_* (\sqb) \}
\]
and similarly with symmetric ideals
$il^t(A)$ and $il^t(B)$.

\vspace{20pt}

\noindent In summary,
for a general $A \in \cG\cR_C$, we have defined the sets
\begin{equation}\label{fig3_5_14}
\begin{tikzcd}
il(A)&\mbox{[1]-}il(A)\arrow[hook]{l}\arrow[hook]{r}& fun\cdot il(A)\\
il^t(A)\arrow[hook]{ur}\arrow[hook]{u}&E\mbox{-}[1]\mbox{-}il(A)\arrow[hook]{r}\arrow[hook]{u}&E\mbox{-}fun\cdot il(A)\arrow[hook]{u} \\
E\mbox{-}il^t(A)\arrow[hook]{u}\arrow[hook]{ur}&&
\end{tikzcd}
\end{equation}


For a self adjoint $A$ we have equality
$il(A) =[1]\mbox{-}il(A)=il^t(A)$.
It is easy to check that for $A=\cG (B)$, $B$ a commutative ring,
all the inclusions in (\ref{fig3_5_14})
are equalities, and are identified with the set of ideals of $B$.
For $A = \cG(B)$, $B$ a commutative ring with involution
$( \ )^t: B \rightarrow B$,
$il(A)$ is just the set of ideals of $B$,
(on which we have involution $\sqb \mapsto \sqb^t$); the set
\begin{equation}
\begin{array}{ccccc}
[1]\mbox{-}il(A) &=& fun\cdot il (A)& = &eq (A)\\
||& & ||& & || \\
E\mbox{-}[1]\mbox{-}il(A)& =& E\mbox{-}fun\cdot il (A)& =& Z \cdot eq(A)
\end{array}
\end{equation}
is the set of ideals $\sqb$ of $B$ fixed by the involution
$\sqb = \sqb^t$
(so that $B/\sqb$ has involution);
and finally,
$E\mbox{-}il^t(A) =il^t(A)$
are the ideals
$\sqb \subseteq B$ that are generated by their symmetric elements
$\sqb^+ = \sqb \cap B^+$
(with
$B^+ = \{b \in B \ , \ b^t = b \} \subseteq B$
the subring of symmetric elements),
which in turn correspond bijectively with the ideals of $B^+$.

\chapter{Primes and Spectra}

\section{Maximal ideals and primes}

We say that an equivalence ideal $\varepsilon \in eq(A)$
is \emph{proper} if $(1,0) \not\in \varepsilon$,
or equivalently
$\varepsilon_X
\subsetneqq
 A_X \times A_X$
for some/all $X \in \FF$,
or equivalently $A/\varepsilon  \neq 0$.
We say that a functorial ideal, or an ideal,
$\sqa$ is \emph{proper} if $1 \not\in \sqa$,
or equivalently
$\sqa_{[1]} \subsetneqq A_{[1]}$.
Since a union of a chain of proper ideals is again a proper ideal,
an application of Zorn's lemma gives

\prop{10.1.1} For $A \in \cG\cR_C$, there exists maximal
proper ideal.

We let $Max(A) \subseteq il(A)$ denote the set of maximal ideals.

\defin{10.1.2} A (proper) ideal
$\sqp \in il(A)$ is called \emph{prime} if
${\cal S}_{\sqp} = A_{[1]}\setminus \sqp$ is closed with respect to multiplication,
i.e. if for all $a,b \in A_{[1]}$,
\begin{equation}
a \lhd b \in \sqp \quad \mbox{implies} \quad  a \in \sqp
\quad \mbox{or} \quad b \in \sqp
\end{equation}
We let $spec(A) \subseteq il(A) $ denote the set of primes of $A$.

\prop{10.1.3}
$$Max(A) \subseteq spec(A).$$

\begin{proof}
Let $\sqp \in Max(A)$,
and take any elements $a, a' \in A_{[1]} \setminus \sqp$.
Since $\sqp$ is maximal, the ideals $(\sqp, a)_A$
and $(\sqp, a')_A$ are the unit ideal.
We have therefore $1=(b \lhd c)\sslash d$, and $1=(b' \lhd c')\sslash d'$,
with $b,d \in A_X$, $b' ,d' \in A_{X'}$,
$c \in (\sqp \cup \{a\})^X$,
$c' \in (\sqp \cup \{a'\})^{X'}$.

Thus we have
\begin{equation}
1=1\lhd 1= ((b \lhd c)\sslash d)\lhd ((b' \lhd c')\sslash d')=
(b \lhd c \lhd \widetilde{b'} \lhd \widetilde{c'})\sslash (d' \lhd \widetilde{d})=
(b \lhd \widetilde{b'} \lhd \widetilde{c} \lhd \widetilde{c'})\sslash (d' \lhd \widetilde{d})
\end{equation}
But $\widetilde{c} \lhd \widetilde{c'} \in ((\sqp \cup \{a\}) \lhd (\sqp \cup \{a'\} ))^{X\prod X'}
\subseteq (\sqp \cup \{a\lhd a' \})^{X \prod X'}$, \\
and so
$1 \in (\sqp , a \lhd a')_A$,
hence $a \lhd a' \not\in \sqp$.
\end{proof}

Similarly, a union of a chain of proper symmetric ideals
is a proper symmetric ideal,
and the set of maximal symmetric ideals
$Max^t(A) \subseteq il^t(A)$ is non-empty.

\defin{10.1.4} A (proper) symmetric ideal
$\sqp \in il^t(A)$
is called \emph{symmetric prime} if
${\cal S}^+_{\sqp} = A^+_{[1]}\setminus \sqp$
is closed with respect to multiplication,
i.e. if for all $a = a^t$, $b = b^t \in A_{[1]}$,
\begin{equation}
a \lhd b \in \sqp \ \ \mbox{implies} \ \ a \in \sqp \ \mbox{or} \
b \in \sqp
\end{equation}
Now a maximal symmetric ideal $\sqp \in Max^t(A)$
is a symmetric prime.
We let $spec^t(A)$ denote the set of symmetric primes:
\begin{equation}
Max^t(A) \subseteq spec^t(A) \subseteq il^t(A)
\end{equation}
Note that for a prime $\sqp \in spec(A)$, the ideal $A \cdot
\sqp^+$, generated by the symmetric elements of $\sqp$, is a
symmetric prime, and we have a canonical map
\begin{equation}
spec(A) \twoheadrightarrow spec^t(A) \ ,  \
\sqp \mapsto A \cdot \sqp^+
\end{equation}

\section{The Zariski topology}
\label{sec4.2}
\defin{10.2.1} The \emph{closed sets} in $spec(A)$ are the set of the form
\begin{equation}
V(\sqa) = \{ \sqp \in spec(A) , \sqp \supseteq \sqa \},
\end{equation}
with $\sqa \subseteq A_{[1]}$,
which we may take to be an ideal $\sqa \in il(A)$. \vspace{2mm}\\

We have

\begin{equation}\label{zariski}
\begin{array}{ccl}
& (i) & V(\sum\limits_i \sqa_i) = \bigcap\limits_i V(\sqa_i),\\
 &(ii) & V(\sqa \cdot \sqa') = V(\sqa) \cup V(\sqa'),\\
 &(iii) & V(0) = spec(A) \quad , \quad V(1) = \emptyset
\end{array}
\end{equation}

This shows the sets $V(\sqa)$ define a topology on $spec(A)$,
the \emph{Zariski topology}.

\vspace{10pt}

The closed sets in $spec^t(A)$ are similarly given by
\begin{equation}
V^t(\sqa) = \{\sqp \in spec^t(A) \ , \ \sqp \supseteq \sqa \}
\end{equation}
with $\sqa \in il^t(A)$;
these satisfy (\ref{zariski}),
and define the Zariski topology on $spec^t(A)$.

\vspace{10pt}

For a subset $C \subseteq spec(A)$, we have the ideal,
\begin{equation}
I(C) = \bigcap\limits_{\sqp \in C} \sqp
\end{equation}
For a subset $C \subseteq spec^t(A)$,
we have the symmetric ideal $I^t(C)$ generated by
$\bigcap\limits_{\sqp \in C} \sqp^+$.

We have  Galois correspondences,
\begin{equation}
il(A)
\begin{array}{c}
\stackrel{\textstyle V}{\longrightarrow}\\
\stackrel{\textstyle \longleftarrow}{I}
\end{array}
\{C \subseteq spec(A)\}
\ \ \mbox{and} \ \ il^t(A)
\begin{array}{c}
\stackrel{\textstyle V^t}{\longrightarrow}\\
\stackrel{\textstyle \longleftarrow}{I^t}
\end{array}
\{ C \subseteq spec^t(A) \}
\end{equation}
The maps $V$, $I$  (resp. $V^t$, $I^t$), are monotone
\begin{equation}
\begin{array}{l}
\sqa_1 \subseteq \sqa_2 \Rightarrow V(\sqa_1) \supseteq V(\sqa_2)\\
C_1 \subseteq C_2 \Rightarrow I(C_1) \supseteq I(C_2)
\end{array}
\end{equation}
and we have
\begin{equation}
\sqa \subseteq IV(\sqa) \quad , \quad
C \subseteq VI(C)
\end{equation}
It follows that we have
\begin{equation}
V(\sqa) = VIV(\sqa) \quad \mbox{and} \quad
I(C) = IVI(C)
\end{equation}
and the maps $V$, $I$ induce inverse bijections
between radical ideals and closed subsets of $spec(A)$.
\begin{equation}
\label{eq4.2.8}
\{\sqa \in il(A) , \sqa = IV(\sqa)\}=
\hspace{2in}
\end{equation}
\[=\{I(C) , C \subseteq spec(A)\}
\stackrel{\textstyle \sim}{\longleftrightarrow}
\{ C \subseteq spec(A) , C=VI(C)\} =
\]
\[\hspace{2in}=\{ V(\sqa) , \sqa \in il(A) \}
\]

Similarly, we have a bijection between the radical symmetric ideals
$\sqa = I^t V^t(\sqa)$
 and the closed subset of
$spec^t(A)$.

\vspace{10pt}

\lem{10.2.2} For $\sqa \in il(A)$, we have
\begin{equation}
IV(\sqa) =
\{ a \in A_{[1]}, a^n \in \sqa \, \mbox{for some} \,
n > 0 \}
\stackrel{def}{=}
\sqrt{\sqa}
\end{equation}
\begin{proof}
If $a \in \sqrt{\sqa}$, say $a^n \in \sqa$,
then for all $\sqp \supseteq \sqa$,
$a \in \sqp $,
and so $\sqrt{\sqa} \subseteq \bigcap\limits_{\sqa \subseteq \sqp}
\sqp = IV(\sqa)$.

Assume $a \not\in \sqrt{\sqa}$, so $a^n \not\in \sqa$ for all $n$.
An application of Zorn's lemma gives that there exists a
maximal element $\sqp$ in the set
\begin{equation}
\label{eq4.2.10}
\{ \sqb \in il(A), \, \sqb \supseteq \sqa , \,
a^n \not\in \sqb \, \mbox{for all} \, n\}
\end{equation}
We claim $\sqp$ is prime. If $x, x' \in A_{[1]}\setminus \sqp$,
then the ideals $(\sqp,x)_A$,
$(\sqp, x')_A$ properly contain $\sqp$,
and by maximality of $\sqp$ in the set
(\ref{eq4.2.10}),
we must have $a^n \in (\sqp, x)_A$,
$a^{n'} \in (\sqp, x')_A$,
for some $n$, $n'$.
We get
\[
a^{n +n'}= a^n \lhd a^{n'}=
((b \lhd c)\sslash d) \lhd ((b' \lhd c')\sslash d' )=
(b \lhd \widetilde{b'}\lhd \widetilde{c} \lhd \widetilde{c'})\sslash (
 d' \lhd \widetilde{d})
\]
with $b, d \in A_X$, $b', d' \in A_{X'}$,
$c \in (\sqp \cup \{x \})^X$,
$c' \in (\sqp \cup \{x' \})^{X'}$.

But $\widetilde{c}\lhd \widetilde{c'} \in \left(
(\sqp \cup \{x\})\lhd (\sqp \cup \{x'\})
\right)^{X \prod X'} \subseteq
\left( \sqp \cup \{x \lhd x' \}\right)^{X \prod X}$,
and since $a^{n+n'} \not\in \sqp$,
we must have $x\lhd x' \not\in \sqp$;
and $\sqp$ is indeed prime.
Now $\sqp \supseteq \sqa$, and $a \not\in \sqp$,
so $a \not\in \bigcap\limits_{\sqa \subseteq \sqp} \sqp =
IV(\sqa)$.
\end{proof}
Similarly, for a symmetric ideal
$\sqa \in il^t(A)$, we have $I^+V^+(\sqa) = \sqrt{\sqa}^+$,
where $\sqrt{\sqa}^+$ is the ideal of $A$ generated by the set
$\{a=a^t \in A_{[1]}^+ \ , \ a^n \in \sqa \ \mbox{for some} \ n > 0 \}$.

\lem{10.2.3} For a subset $C \subseteq spec(A)$,
$VI(C) = \bar{C}$ the closure of $C$.\\

\begin{proof}
We have $C \subseteq VI(C)$,
and $VI(C)$ is closed. If $C \subseteq V(\sqa)$,
where we may assume $\sqa = \sqrt{\sqa}$,
then $VI(C) \subseteq VIV(\sqa) = V(\sqa)$, and so
\begin{equation}
VI(C) = \bigcap\limits_{C \subseteq
V(\sqa)} V(\sqa) = \bar{C}
\end{equation}
\end{proof}

Similarly, for $C \subseteq spec^t(A)$,
$V^+I^+(C) = \overline{C}$ the closure of $C$. \\
We can restate (\ref{eq4.2.8}),
\begin{equation}
\label{eq4.2.12}
\{\sqa \in il^{(t)}(A), \sqa = \sqrt{\sqa}^{(+)} \}
\stackrel{\textstyle \sim}{\longleftrightarrow}
\{C \subseteq spec^{(t)}(A) , C=\bar{C} \}
\end{equation}

\setcounter{equation}{0}
\section{Basic open sets}
A basis for the open sets of $spec(A)$ is given by the
\emph{basic open sets},
these are defined for $a \in A_{[1]}$ by
\begin{equation}
D_a = spec(A) \setminus V(a) =
\{ \sqp \in spec(A) , a \not\in \sqp \}
\end{equation}

A basis for the open sets of $spec^t(A)$ is given by
\begin{equation}
D^+_a = spec^t(A) \setminus V^+(a)
\end{equation}
with symmetric $a = a^t \in A_{[1]}^+$

We have,
\begin{equation}
\begin{array}{ll}
D_{a_1} \cap D_{a_2} = D_{a_1 \lhd a_2} &
D^+_{a_1} \cap D^+_{a_2} = D^+_{a_1 \lhd a_2}\\
D_1 = spec(A) \quad  , \quad D_0 = \emptyset &
D^+_1 = spec^t(A) \ , \ D^+_0= \emptyset
\end{array}
\end{equation}

That every open set is the union of basic open sets, is shown by
\begin{equation}
spec(A) \setminus V(\sqa) = \bigcup\limits_{a \in \sqa} D_a
\end{equation}
and
\begin{equation}
spec^t(A) \setminus V^+(\sqa) = \bigcup\limits_{a \in \sqa^+} D^+_a \ , \
\mbox{the union over} \ \sqa^+ = \sqa \cap A_{[1]}^+.
\end{equation}
\vspace{10pt}

Note that we have,
\begin{equation}\label{10.3.6}
D_a = spec(A) \Leftrightarrow a \lhd A_{[1]} = \{a \}_A = (1)
\end{equation}
\[
\Leftrightarrow  \, \mbox{there exists a (unique)} \,
a^{-1} \in A_{[1]} , a \lhd a^{-1} =1
\]
We say that such $a$ is \emph{invertible}, and we let $A^*$
denote the set of invertible elements.
Note that $A^*$ is an abelian group (with a non-trivial involution
for $A$ non self-adjoint), and
$A \mapsto A^*$ is a functor $\cG\cR_C \rightarrow Ab$ (= abelian groups).
Similarly, for $a = a^t$ symmetric, $D^+_a = spec^t(A)$
if and only if $a$ is invertible.

\vspace{10pt}

Note that we have,
\begin{equation}
D_a = \emptyset \Leftrightarrow
a \in \bigcap\limits_{\sqp \in spec(A)} \sqp = \sqrt{0}
\end{equation}
\[
\Leftrightarrow
\, \mbox{there exists} \, n > 0 \, \mbox{with} \, a^n =0
\]
We say that such $a$ is \emph{nilpotent}.
Similarly, for $a = a^t$, $D^+_a = \emptyset$ if and only if $a$ is nilpotent.

\lem{10.3.1} Let $\sqa = \sqrt{\sqa}
\in il(A)$
be a radical ideal. Then
\begin{equation}
V(\sqa) \, \mbox{is irreducible} \,
\Leftrightarrow
\sqa \, \mbox{is prime}
\end{equation}

\begin{proof}
$(\Leftarrow)$: If $\sqa$ is prime,
$V(\sqa)= VI$  $\{\sqa\} = \overline{\{\sqa\}}$
is the closure of a point,
hence irreducible.

$(\Rightarrow)$: For any $a \in A_{[1]}$, we have
$$
\label{eq4.3.7}
V(\sqa) \cap D_a \neq \emptyset
\Leftrightarrow \exists \sqp \in spec(A) ,
\sqp \supseteq \sqa , \sqp \not\ni a
$$
\[
\Leftrightarrow
a \not\in \bigcap\limits_{\sqa \subseteq \sqp} \sqp =
\sqrt{\sqa} = \sqa
\]
Hence for any basic open sets $D_a$, $D_b$, $a,b \in A_{[1]}$,
we have
\[
V(\sqa) \cap D_a \neq \emptyset
\, \mbox{and}\,
V(\sqa) \cap D_b \neq \emptyset
\Leftrightarrow
a \not\in \sqa \, \mbox{and} \, b \not\in\sqa
\]
If $V(\sqa)$ is irreducible this implies
\[
\emptyset \neq V(\sqa) \cap D_a \cap D_b =
V(\sqa) \cap D_{a \lhd b}
\Leftrightarrow
a \lhd b \not\in \sqa
\]
\end{proof}
Similarly, for a radical symmetric ideal $\sqa = \sqrt{\sqa}^+$,
the set $V^+(\sqa)$ is irreducible if and only if
$\sqa$ is a symmetric prime.

\vspace{10pt}

Thus the bijection (\ref{eq4.2.12}) induces a bijection
\begin{equation}
\begin{array}{rl}
spec(A) \stackrel{\textstyle \sim}{\longleftrightarrow}&
\big\{ C \subseteq spec(A) , C= \bar{C} \,
\mbox{closed and irreducible} \big\} \\
\sqp \mapsto & V(\sqp) = \bar{\{ \sqp \}}
\end{array}
\end{equation}
and
\begin{equation}\label{10.3.10}
spec^t(A)  \stackrel{\textstyle \sim}{\longleftrightarrow}
\big\{C \subseteq spec^t(A) \ , \
 C= \overline{C} \ \mbox{closed and irreducible} \big\}
\end{equation}
and the spaces $spec(A)$
and $spec^t(A)$
 are Sober spaces (or Zariski spaces): \\
every closed irreducible subset has a unique generic point.

\prop{10.3.2} For $a \in A_{[1]}$,
(Respectively, $a=a^t \in A_{[1]}^+$)
the basic open set $D_a$ (resp. $D_a^+$) is compact.
\setcounter{equation}{9}

In particular, $D_1 = spec(A)$ and $D_1^+ = spec^t(A)$ are compact.

\begin{proof}
We have to show that in every covering of $D_a$
by basic open sets $D_{g_i}$,
there is always a finite subcovering. We have
\begin{equation}
\label{eq4.3.10}
\begin{array}{ccl}
D_a \subseteq \bigcup\limits_{i}D_{g_i} &
\Leftrightarrow &
V(a) \supseteq \bigcap\limits_i V(g_i) = V(\sum\limits_i
     g_i \lhd A_{[1]})\\
&\Leftrightarrow &
\sqrt{\{a\}_A} = IV(a) \subseteq IV\left(
\sum\limits_i g_i \lhd A_{[1]}\right) =
\sqrt{\sum\limits_i g_i \lhd A_{[1]}}\\
&\Leftrightarrow &
\mbox{for some} \, n , a^n \in \sum\limits_i g_i \lhd A_{[1]}\\
&\Leftrightarrow &
\mbox{for some} \, n , X \in \FF, b, d, \in A_X ,
a^n = (b \lhd c)\sslash d ,\\
&& \hspace{2cm}
\mbox{with} \, c= (c^{(x)}) \in \left(
\{g_i \} \right)^X
\end{array}
\end{equation}
Thus $c^{(x)} = g_{i(x)}$,
and going backwards in the above equivalences we get \\
$D_a \subseteq \bigcup\limits_{x \in X} D_{g_{i(x)}}$,
 a finite subcovering.
\end{proof}

The canonical map
$\pi_A:spec(A) \twoheadrightarrow spec^t(A)$,
$\pi_A(\sqp) = A \cdot \sqp^+$, \\
is continuous:
\begin{equation}
\pi_A^{-1}\left(V^+(\sqa)\right)= V(\sqa) \ , \
\pi_A^{-1}(D_a^+)=D_a
\end{equation}

\section{Functoriality}

For a homomorphism $\varphi \in \cG\cR_C (A,B)$,
the pull-back of a (symmetric) prime is a
(symmetric) prime, and we have maps
\begin{equation}
\begin{array}{cl}
\varphi^* = spec (\varphi): &spec(B) \rightarrow spec(A)\\
& \sqq \mapsto \varphi^* (\sqq) = \varphi^{-1}_{[1]} (\sqq)
\end{array}
\end{equation}
and
\begin{equation}
\begin{array}{cl}
\varphi^* = spec^t (\varphi): &spec^t(B) \rightarrow spec^t(A)\\
& \sqq \mapsto \varphi^* (\sqq) = A \cdot \big( \varphi^{-1} (\sqq)
\cap A^+ \big)
\end{array}
 \end{equation}
The inverse image under $\varphi^*$ of a closed set is closed, we have
\begin{equation}
\varphi^{*-1} (V_A(\sqa)) =
\{\sqq \in spec(B) , \varphi^{-1}_{[1]}(\sqq) \supseteq \sqa \}=
\end{equation}
\[=
\{ \sqq \in spec(B) , \sqq \supseteq \varphi_{[1]}(\sqa) \}=
V_B(\varphi_{[1]}(\sqa))
\]
and similarly,
$ \varphi^{*-1}(V^+_A(\sqa)) =
V^+_B(\varphi_*(\sqa)) $
with
$\varphi_*(\sqa) = B \cdot \varphi (\sqa^+)$.

Also the inverse image under $\varphi^*$ of a basic open set is a basic open set, we have
\begin{equation}
\varphi^{*-1}(D_a) =
\{\sqq \in spec(B) ,  \varphi^{-1}_{[1]}(\sqq) \not\ni a\}=
\end{equation}
\[=
\{\sqq \in spec(B) , \varphi_{[1]}(a) \not\in \sqq\}=
D_{\varphi_{[1]}(a)}
\]
and similarly, $\varphi^{*-1}(D^+_a) = D^+_{\varphi_{[1]}(a)}$
for $a=a^t$.

Thus the maps $\varphi^*= spec(\varphi)$
and $spec^t(\varphi)$ are
 continuous,
and we see that $spec^{(t)}$ are contravariant functors from
$\cG\cR_C$ to the category $Top$, whose objects are (compact,
sober) topological spaces, and continuous maps,
\begin{equation}
spec \ , \ spec^t: (\cG\cR_C)^{op} \rightarrow Top
\end{equation}

\lem{10.4.1} For $\varphi \in \cG\cR_C(A,B)$,
and for $\sqb \in il(B)$, we have
\begin{equation}
V_A(\varphi_{[1]}^{-1}(\sqb)) =
\overline{\varphi^*(V_B(\sqb))}
\end{equation}
\begin{proof}
We may assume without loss of generality that
$\sqb = \sqrt{\sqb}$ is radical
(noting that
$\sqrt{\varphi^{-1}_{[1]}(\sqb)} =\varphi^{-1}_{[1]}(\sqrt{\sqb})$). \\
Put $\sqa = I\varphi^*(V(\sqb))$,
so that $V(\sqa) = \overline{\varphi^*(V(\sqb))}$ by lemma 10.2.3.

We have for any $a \in A_{[1]}$,
\begin{equation}
\begin{array}{cl}
a \in \sqa&\Leftrightarrow
a \in \sqp , \, \mbox{for every prime} \, \sqp \in \varphi^*(V(\sqb))\\
&\Leftrightarrow a \in \varphi^* (\sqq) = \varphi^{-1}_{[1]}(\sqq) ,
  \, \mbox{for every prime} \, \sqq \in V(\sqb)\\
&\Leftrightarrow
\varphi_{[1]}(a) \in \bigcap\limits_{\sqb \subseteq \sqq} \sqq =
       \sqrt{\sqb} = \sqb\\
&\Leftrightarrow a \in \varphi_{[1]}^{-1}(\sqb)
\end{array}
\end{equation}
Thus $\sqa = \varphi_{[1]}^{-1}(\sqb)$, and the lemma is proved.
\end{proof}

\chapter{Localization and sheaves}
\section{Localization}
For  $A \in \cG \cR_C$, and $s \in A_{[1]}$,
in the generalized ring
$A[1/s]$ obtained from $A$ by adding an inverse to $s$,
we have also an inverse to $s^t$,
$(s^t)^{-1} = (s^{-1})^t$,
and therefore also an inverse to the symmetric element
$s \lhd s^t$,
so $A[1/s] = A[1/s\lhd s^t]$.

We shall therefore localize only with respect to symmetric
elements!.

For $A \in \cG\cR_C$, a subset
$S \subseteq A_{[1]}^{+}$ is called \emph{multiplicative}
if
\begin{equation}
1 \in S , \ \mbox{and} \
S \lhd S \subseteq S.
\end{equation}
For such $S \subseteq A_{[1]}^+$, and for $X \in \FF$,
we let $(S^{-1}A)_X = (A_X \times S) /\hspace{-5pt}\approx$
denote the equivalence classes of $A_X \times S$ with respect to the equivalence relation defined by
\begin{equation}
(a_1 , s_1 ) \approx (a_2, s_2) \ \mbox{if and only if} \
s \lhd s_2 \lhd a_1 = s \lhd s_1 \lhd a_2
\ \mbox{for some} \ s \in S
\end{equation}
We write $a/s$ for the equivalence class $(a,s) /\hspace{-5pt}\approx$.
Note that by taking "common denominator"
we can write any element
$a = (a^{(y)}/s_y) \in (S^{-1}A)_f = \prod\limits_{y \in Y}
(S^{-1}A)_{f^{-1}(y)}$,
in the form
\begin{equation}
a = (\bar{a}^{(y)}/s), \left( \mbox{take} \ s = \prod s_y ,\
\bar{a}^{(y)}= (\prod\limits_{y' \neq y}s_{y'}) \lhd a^{(y)})\right)
\end{equation}
For $f \in Set_{\bullet}(X,Y)$,
$g \in Set_{\bullet}(Y,Z)$,
we have well-defined operations of multiplication and contraction,
independent of the choice of representatives,
\begin{equation}
\lhd: (S^{-1}A)_g \times (S^{-1}A)_f \longrightarrow
       (S^{-1}A)_{g \circ f} , \quad
 (a/s_1) \lhd (b/s_2) = (a \lhd b) / (s_1 \lhd s_2)
\end{equation}
\begin{equation}
(\, \sslash \, ):(S^{-1}A)_{g \circ f } \times (S^{-1}A)_f \longrightarrow
       (S^{-1}A)_{g} , \quad
(a/s_1) \sslash (b/s_2) = (a \sslash b) / (s_1 \lhd s_2)
\end{equation}

It is straightforward to check that these operations satisfy the
axioms of a commutative generalized ring. The canonical
homomorphism
\begin{equation}
\phi:A \rightarrow S^{-1}A , \quad \phi(a) = a/1
\end{equation}
satisfies the universal property
\begin{equation}
\begin{array}{rcl}
\cG\cR_C (S^{-1}A,B)& = &
\{\varphi \in \cG\cR_C (A,B) , \varphi(S) \subseteq B^*\}\\
\widetilde{\varphi} & \mapsto & \widetilde{\varphi} \circ \phi\\
\varphi(s)^{-1} \lhd \varphi (a) = \widetilde{\varphi} (a/s)
&\mapsfrom& \varphi
\end{array}
\end{equation}

\exmpl{11.1.1} For $s \in A_{[1]}^+$,
take $S = \{ s^n \ , \ n \geq 0 \}$.
We write $A_{s}$ for $S^{-1}A$,
and $\phi_{s} \in \cG\cR_C (A,A_s)$ satisfy
\begin{equation}
\cG\cR_C (A_{s},B) = \{ \varphi \in \cG\cR_C(A,B) , \varphi (s) \in B^* \}
\end{equation}

\vspace{10pt}

\exmpl{11.1.2}
\setcounter{equation}{8}
For $\sqp  \in spec^t(A)$,
take $S_{\smallsqp} = A^+_{[1]}\setminus \sqp$.
We write $A_{\smallsqp}$ for $S^{-1}_{\smallsqp} A$,
and $\phi_{\smallsqp} \in \cG\cR_C (A,A_{\smallsqp})$ satisfy
\begin{equation}
\cG\cR_C (A_{\smallsqp},B) =
\{ \varphi \in \cG\cR_C(A,B) ,
    \varphi (A_{[1]}^+\setminus \sqp) \subseteq B^* \}
\end{equation}

\section{Localization and ideals}
For an ideal $\sqa \in il(A)$, we let
\begin{equation}
S^{-1}\sqa = \{ a/s \in (S^{-1}A)_{[1]} , s \in S , a \in \sqa \}
\end{equation}
By using common denominator,
we see that $S^{-1}\sqa$ is an ideal of $S^{-1}A$, \\
$S^{-1}\sqa \in il(S^{-1}A)$:

For $b/s_1, d /s_2  \in (S^{-1}A)_X$, and for
$a_x/s_x \in S^{-1}\sqa$, $x \in X$, we have
$$
(b/s_1 \lhd (a_x/s_x))\sslash (d/s_2)=
((b \lhd (s'_x \lhd a_x))\sslash d)/(s_1 \lhd s_2^t \lhd
\prod\limits_{x \in X}s_x)
$$
\[
 \mbox{with} \ s'_x =\prod\limits_{x' \neq x} s_{x'}, \mbox{and this is in }S^{-1}\sqa \mbox{ since }
\sqa \mbox{ is an ideal.}
\]

If $\sqa = A \cdot \sqa^+$ is symmetric,
$S^{-1}\sqa$ is symmetric.

We have, therefore, the Galois correspondence
\begin{equation}
il^{(t)}(A)
\begin{array}{c}
\stackrel{\textstyle S^{-1}}{\longrightarrow}\\
\stackrel{\textstyle \longleftarrow}{\phi^*}
\end{array}
il^{(t)}(S^{-1}A)
\end{equation}
For $\sqb \in il(S^{-1}A)$, we have
\begin{equation}
S^{-1}\phi^* \sqb = \sqb
\end{equation}
Indeed, for an element $a/s \in \sqb$,
we have $a/1 \in \sqb$, or $a \in \phi^*\sqb$,
and $a/s \in S^{-1}\phi^*\sqb$;
hence $\sqb \subseteq S^{-1} \phi^* \sqb$,
and the reverse inclusion is clear. \\
We have immediately from the definitions,
for $\sqa \in il(A)$,
\begin{equation}
\phi^* S^{-1}\sqa =
\{ a \in A_{[1]}, \ \mbox{there exists} \ s \in S
\ \mbox{with} \ s \lhd a \in \sqa \} =
\bigcup\limits_{s \in S} (\sqa:s)
\end{equation}
In particular,
\begin{equation}
S^{-1}\sqa = (1)
\Leftrightarrow
\sqa \cap S \neq \emptyset
\end{equation}
We say that $\sqa \in il^{(t)}(A)$ is $S$-\emph{saturated}
if $\phi^* S^{-1}\sqa = \sqa$, that is if
\begin{equation}
\mbox{for all}\ s \in S,
a \in A_{[1]}^+: s \lhd a \in \sqa
\Rightarrow a \in \sqa
\end{equation}
We get that $S^{-1}$ and $\phi^*$ induce inverse bijections,
\begin{equation}
\{ \sqa \in il^{(t)}(A) , \quad \sqa \ \mbox{is} \
S\mbox{-saturated} \}
\stackrel{\textstyle \sim}{\longleftrightarrow}
il^{(t)}(S^{-1}A)
\end{equation}
For an $S$-saturated (symmetric) ideal
$\sqa \in il^{(t)}(A)$,
let $\pi_{\sqa}: A \twoheadrightarrow A/\sqa = A/E(\sqa)$
be the canonical homomorphism, and let
$\bar{S} = \pi_{\sqa}(S) \subseteq (A/\sqa)_{[1]}^+$,
then we have canonical isomorphism
\begin{equation}
\label{eq5.2.10}
\bar{S}^{-1}(A/\sqa) \cong S^{-1}A/S^{-1}\sqa
\end{equation}
Note that for a (symmetric) prime $\sqp \in spec^{(t)}(A)$,
$\sqp$ is $S$-saturated if and only if $\sqp \cap S = \emptyset$,
and
in this case $S^{-1}\sqp$ is a (symmetric) prime,
$S^{-1}\sqp \in spec^{(t)}(S^{-1}A)$.
Note that for a (symmetric) prime $\sqq \in spec^{(t)}(S^{-1}A)$,
$\phi^*(\sqq)$ is always an $S$-saturated (symmetric) prime.
We get the bijections:
\begin{equation}
\label{eq5.2.11}
\{\sqp \in spec^{(t)}(A) , \sqp \cap S = \emptyset \}
\stackrel{\textstyle \sim}{\longleftrightarrow}
spec^{(t)}(S^{-1}A)
\end{equation}
These are homeomorphisms for the Zariski topologies.

\vspace{10pt}

For $s \in A^+_{[1]}$, the homeomorphism (\ref{eq5.2.11})
gives for Example 11.1.1,
\begin{equation}
\phi_{s}^*: spec^t(A_{s})
\stackrel{\textstyle \sim}{\longrightarrow}
D_s^t \subseteq spec^t(A)
\end{equation}

For a symmetric prime $\sqp \in spec^t(A)$, the homeomorphism (\ref{eq5.2.11})
with $S = S_{\smallsqp}$,
of Example 11.1.2, reads
\begin{equation}
\phi_{\smallsqp}^*: spec^t(A_{\smallsqp})
\stackrel{\textstyle \sim}{\longrightarrow}
\{ \sqq \in spec^t(A) , \sqq \subseteq \sqp \}
\end{equation}
The generalized ring $A_{\smallsqp}$
is a \emph{local-generalized-ring}
in the sense that it has a unique maximal symmetric ideal
$m_{\smallsqp} = S_{\smallsqp}^{-1} \sqp$.
The \emph{residue field at $\sqp$} is defined by
\begin{equation}
\mathbb{F}_{\smallsqp} = A_{\smallsqp}/m_{\smallsqp} =
A_{\smallsqp}/E(m_{\smallsqp}).
\end{equation}
We have the canonical homomorphism
$\pi_{\smallsqp}: A \sur A/\sqp$,
and putting $\bar{S}_{\smallsqp} = \pi_{\smallsqp}
(S_{\smallsqp})$,
we have (\ref{eq5.2.10})
\begin{equation}
\mathbb{F}_{\smallsqp} = \bar{S}_{\smallsqp}^{-1}(A/\sqp)
\end{equation}
The square diagram
\begin{equation}
\label{eq5.2.17}
\end{equation}
\begin{diagram}
A & \rTo^{\phi_{\smallsqp}}& A_{\smallsqp}\\
\dOnto^{\pi_{\smallsqp}}& & \dOnto_{\pi_{{\smallsqm}_{\smallsqp}}}\\
A/\sqp& \rTo^{\phi_0}& {\mathbb F}_{\smallsqp}
\end{diagram}

is cartesian,
\begin{equation}
\mathbb{F}_{\smallsqp} =
(A/\sqp) \bigotimes\limits_{A} A_{\smallsqp}
\end{equation}
\[\cG\cR_C (\mathbb{F}_{\smallsqp}, B) =
\left\{\varphi \in \cG\cR_C(A,B) , \
       \varphi(\sqp) \equiv 0 , \
       \varphi(A_{[1]}^+\setminus \sqp) \subseteq B^*
\right\}
\]

\vspace{10pt}

For a homomorphism of generalized rings
$\varphi \in \cG\cR_C (A,B)$, and for \\
$\sqq \in spec^t(B)$
with
$\sqp = \varphi^*(\sqq) \in spec^t(A)$,
the square diagram (\ref{eq5.2.17})
is functorial,
and we have a commutative cube diagram

\begin{equation}
\begin{diagram}
 & & A &\rTo^{\hspace{0.5cm}\phi_{\smallsqp}} & & & A_{\smallsqp}\\
 &\ldTo^{\varphi} &\dTo& &   & \ldTo^{\varphi_{\smallsqp}}&\dTo\\
B&\rTo^{\hspace{2cm}\phi_{\smallsqq}}& & & B_{\smallsqq}& & \\
\dTo &              & & &\dTo     & & \\
     &              &A/\sqp& \rTo& & &{\mathbb F}_{\smallsqp}\\
  &\ldTo &                 &   &  & \ldTo& \\
B/\sqq& \rTo &             &  &{\mathbb F}_{\smallsqq}   & &
\end{diagram}
\end{equation}

Note that the homomorphism $\varphi_{\smallsqp}
\in \cG\cR_C(A_{\smallsqp},
 B_{\smallsqq} )$
is a \emph{local-homomorphism}
in the sense that
\begin{equation}
m_{\smallsqp}= \varphi^{*}_{\smallsqp}(m_{\smallsqq}) ,
\ \mbox{or equivalently} \
\varphi_{\smallsqp}(m_{\smallsqp}) \subseteq m_{\smallsqq}
\end{equation}

\defin{11.2.1} We let $\cL\cG\cR$
denote the subcategory of $\cG\cR_C$
with objects the local generalized rings,
and with maps
\begin{equation}
\cL\cG\cR(A,B) = \{ \varphi \in \cG\cR_C(A,B), \varphi^*(m_B) = m_A\}
\end{equation}

\setcounter{equation}{0}

\section{The structure sheaf $\cO_A$}

\defin{11.3.1}
For  $A \in \cG\cR_C$,
$U \subseteq spec^t(A)$ open,
$X \in \Fbb$,
we denote by $\cO_A(U)_X$
the set of \emph{sections}
\[f:U \rightarrow \coprod\limits_{\smallsqp \in U} (A_{\smallsqp})_X , \quad
f(\smallsqp) \in (A_{\smallsqp})_X
\]
such that $f$ is \emph{locally a fraction}:

\noindent for all  $\sqp \in U$,
there exists open $U_{\smallsqp} \subseteq U$,
$\sqp \in U_{\smallsqp}$,

and there exist $a \in A_X$, $s \in A^+_{[1]}\setminus
     \bigcup\limits_{\smallsqq \in U_{\smallsqp}}\sqq$,

such that for all $\sqq \in U_{\smallsqp}$:
\begin{equation}\label{loc}
f(\sqq) \equiv a/s \in (A_{\smallsqq})_X.
\end{equation}

Note that $\cO_A(U)$ is a commutative generalized ring,
and for $U' \subseteq U$ restriction gives a homomorphism of generalized rings
\begin{equation*}
\cO_A(U) \rightarrow \cO_A(U') , \quad
f \mapsto f|_{U'}
\end{equation*}
Thus $\cO_A$ is a pre-sheaf of generalized rings over $spec^t(A)$,
and by the local nature of the condition (\ref{loc})
it is clear that it is a \emph{sheaf of generalized rings},
i.e. for $X \in \Fbb$,
$U \mapsto \cO_A(U)_X$ is a sheaf of sets.
It is also clear that the \emph{stalks} are given by
\begin{equation}
\cO_{A,\smallsqp}= \lim\limits_{\stackrel{\longrightarrow}{\smallsqp \in U}}
\cO_A(U)
\stackrel{\textstyle \sim}{\longrightarrow} A_{\smallsqp}
\end{equation}
\[\left(f/\approx\right) \mapsto f(\sqp)
\]

\thm{11.3.2} For $s \in A_{[1]}^+$,
we have a canonical isomorphism
\begin{equation}
\Psi : A_{s} \stackrel{\textstyle \sim}{\longrightarrow}
\cO_A(D_{s}^+) , \quad
\Psi(a/s^n) =
\{f(\sqp) \equiv a/s^n \}
\end{equation}
In particular for $s=1$,
\[
A \stackrel{\textstyle \sim}{\rightarrow} \cO_A(spec^t(A))
\]

\begin{proof} The  map $\Psi$ which takes
$a/s^n \in A_{s}$ to the constant section $f$
with $f(\sqp) \equiv a/s^n$ for all $\sqp \in D_{s}^+$,
is clearly well-defined, and is a homomorphism of generalized rings.

$\Psi$ \underline{is injective}: Assume $\Psi(a_1/s^{n_1})= \Psi(a_2 / s^{n_2})$,
and let

$\sqa = ann_A^t(s^{n_2}\lhd a_1 ,  s^{n_1}\lhd a_2) \in
il^t(A)$,
cf. (\ref{9.5.11}- \ref{9.5.12}  ).
We have,

$$
\hspace{-2cm}
 a_1/s^{n_1} = a_2/s^{n_2} \ \mbox{in} \ A_{\smallsqp}
    \ \mbox{for all} \ \sqp \in D_{s}^+
$$
\begin{equation}
\begin{array}{l}
\Rightarrow s_{\smallsqp} \lhd s^{n_2} \lhd a_1 =
            s_{\smallsqp} \lhd s^{n_1} \lhd a_2
             \ \mbox{with} \ s_{\smallsqp} \in A_{[1]}^+\setminus \sqp
             \ \mbox{for} \ \sqp \in D_{s}^+\\
\Rightarrow \sqa \not\subseteq \sqp \ \mbox{for} \
            \sqp \in D_{s}^+ \\
\Rightarrow V^+(\sqa) \cap D^+_{s} = \emptyset \\
\Rightarrow V^+(\sqa) \subseteq V^+(s) \\
\Rightarrow s \in I^+V^+(\sqa) = \sqrt{\sqa}^+ \\
\Rightarrow s^n \in \sqa \ \mbox{for some} \ n\\
\Rightarrow s^{n + n_2} \lhd a_1 = s^{n+n_1} \lhd a_2 \\
\Rightarrow a_1/s^{n_1} = a_2/s^{n_2} \ \mbox{in} \ A_{s}
\end{array}
\end{equation}
$\Psi$ \underline{is surjective:} Fix
  $f \in \cO_A(D_{s}^+)_X$.
Since $D_{s}^+$ is compact (Proposition 10.3.2),
we can cover $D^+_{s}$ by a finite collection of basic open sets,
$D_{s}^+ = D^+_{g_1}\cup \ldots \cup D^+_{g_N}$,
$g_i^t = g_i$,
such that on $D^+_{g_i}$ the section $f$ is constant,
\[
f(\sqp) = a_i/s_i \ \mbox{for} \
          \sqp \in D^+_{g_i},
\]
We have $V^+(s_i) \subseteq V^+(g_i)$, hence
$g_i \in I^+V^+(s_i) = \sqrt{s_i}^+$, hence for some $n_i$,
and some $c_i \in A_{[1]}$,
$g_i^{n_i}= c_i \lhd s_i$.
Thus our section $f$ is given on $D^+_{g_i}$ by
$a_i/s_i = c_i \lhd a_i/g_i^{n_i}$.
Noting that $D^+_{g_i}= D^+_{g_i^{n_i}}$,
we may replace $g_i^{n_i}$ by $g_i$,
and replace $c_i \lhd a_i$ by $a_i$,
and we have
\[f(\sqp) = a_i/g_i \ \mbox{for} \
          \sqp \in D^+_{g_i},
\]
On the set $D^+_{g_i\lhd g_j} = D^+_{g_i}\cap D^+_{g_j}$, $i \neq j$,
our section $f$ is given by both $a_i/g_i$ and $a_j/g_j$.
By the injectivity of $\Psi$, we have
\[
a_i/g_i = a_j/g_j \ \mbox{in} \ A_{g_i\lhd g_j}
\]
Thus for some $n$ we have
\[
(g_i\lhd g_j)^n \lhd g_j \lhd a_i = (g_i\lhd g_j)^n\lhd g_i \lhd a_j
\]
By finiteness we may assume one $n$ works for all $i,j \leq N$.\\
Replacing $g_i^n\circ a_i$ by $a_i$,
and replacing $g_i^{n+1}$ by $g_i$,
we may assume $f \equiv a_i/g_i$ on $D^+_{g_i}$,
and
\begin{equation}
g_j \lhd a_i = g_i \lhd a_j \ \mbox{for all} \ i,j
\end{equation}
We have $D^+_{s} \subseteq \bigcup D^+_{g_i}$,
hence by (\ref{10.3.10})
we have
\begin{equation}
s^M = (b\lhd c)\sslash d
\end{equation}
with $b,d \in A_Y$,
$c =(c^{(y)}) \in \left(A_{[1]}\right)^Y$ with
$c^{(y)} = g_{i(y)}$,
$$i(y):Y \rightarrow \{1, \ldots, N\}$$.
Define $a \in A_X$ by
\begin{equation}
a = (b \lhd e)\sslash \tilde{d}, \ \mbox{with} \
    \tilde{d} \in A_{\pi_X} = (A_Y)^X ,
    \tilde{d}^{(x)} \equiv d ,
\end{equation}
\[
    e \in A_{\pi_Y} = (A_X)^Y , e^{(y)} = a_{i(y)}
\]

We have for $j = 1, \ldots, N$
\begin{equation}
\begin{array}{ll}
g_j \lhd a &= g_j \lhd ((b \lhd e) \sslash \tilde{d})=
               (b \lhd g_j\lhd e) \sslash \tilde{d}=
               (b \lhd (g_j \lhd a_{i(y)}))\sslash \tilde{d}\\
  & = (b \lhd (g_{i(y)} \lhd a_j ))\sslash \tilde{d} =
   (b \lhd c \lhd \tilde{a}_j) \sslash \tilde{d} =
((b \lhd c) \sslash d) \lhd a_j = s^M \lhd a_j
\end{array}
\end{equation}
Thus we have in $A_{s}$,
$a_j/g_j = a/s^M$ for all $j$,
and our section $f$ is constant
$f=\Psi(a/s^M)$, and
$\Psi$ is surjective.
\end{proof}
\chapter{Schemes}
\setcounter{equation}{0}
\section{Locally generalized ringed spaces}

\defin{12.1.1}
For a topological space $\cX$, we let $\cG\cR_C/\cX$
denote the category
of sheaves of generalized rings over $\cX$.
Its objects are pre-sheaves $\cO$ of
commutative generalized rings,
i.e. functors $U \mapsto \cO(U): \cC_{\cX} \rightarrow \cG\cR_C$,
(with $\cC_{\cX}$ the category of open sets of $\cX$,
with $\cC_{\cX}(U,U')=\{j^U_{U'}\}$
for $U' \subseteq U$,
otherwise
$\cC_{\cX}(U,U') = \emptyset$),
such that for all $X \in \Fbb$,
$U \mapsto \cO(U)_X$
is a sheaf.
The maps $\cG\cR_C/\cX(\cO,\cO')$
are natural transformations of functors $\varphi = \{\varphi(U)\}$,
$\varphi(U) \in \cG\cR_C(\cO(U), \cO'(U))$.

\vspace{10pt}

\defin{12.1.2}
We denote by $\cG\cR\cS$
the category of (commutative) generalized ringed spaces.
Its objects are pairs $(\cX, \cO_{\cX})$,
with $\cX \in Top$,
and $\cO_{\cX} \in \cG\cR_C/\cX$.
The maps $f \in \cG\cR\cS (\cX,\cY)$
are pairs of a continuous function $f \in Top(\cX, \cY)$,
and a map of sheaves of generalized rings over $\cY$,
$f^{\sharp} \in \cG\cR_C/\cY (\cO_{\cY}, f_*\cO_{\cX})$;
explicitly, for all open subsets $U \subseteq \cY$,
we have a homomorphism of generalized rings
\begin{equation}
f^{\sharp}_U = \{f^{\sharp}_{U,X}\} \in \cG\cR_C (\cO_{\cY}(U),
\cO_{\cX}(f^{-1}U))
\end{equation}
and these homomorphisms are compatible with restrictions: for $U'
\subseteq U \subseteq \cY$ open, and for $a \in \cO_{\cY}(U)_X$,
we have $f^{\sharp}_{U,X}(a) |_{f^{-1}(U')} = f^{\sharp}_{U',X}
(a|_{U'})$ in $\cO_{\cX} (f^{-1} U')_X$.

\rem{12.1.3} For a continuous map
$f \in Top(\cX, \cY)$,
we have a pair of adjoint functors
\begin{equation}
\cG\cR_C/\cX
\begin{array}{c}
\stackrel{f^*}{\curvearrowleft}\\
\stackrel{\begin{turn}{180}
$\curvearrowleft$
\end{turn}}{f_*}
\end{array}
\cG\cR_C/\cY
\end{equation}
For sheaves of generalized rings $\cO_{\cX} \in \cG\cR_C/\cX$,
$\cO_{\cY} \in \cG\cR_C/\cY$,
we have
\begin{equation}
f_*\cO_{\cX}(U)= \cO_{\cX}(f^{-1}U) \ , \ U \subseteq \cY \
\mbox{open;}
\end{equation}
\[
f^*\cO_{\cY} (U)_X = \
\mbox{sheaf associated to the pre-sheaf} \]
\begin{equation}
U \mapsto
\lim\limits_{\begin{array}{c}
\stackrel{\textstyle \longrightarrow}{V \subseteq \cY \ \mbox{open}}\\
f(U) \subseteq V
\end{array}}
\cO_{\cY}(V)_X;
\end{equation}
and we have adjunction,
\begin{equation}
\cG\cR_C/\cY (\cO_{\cY} , f_*\cO_{\cX}) =
\cG\cR_C/\cX(f^*\cO_{\cY} ,\cO_{\cX})
\end{equation}

\rem{12.1.4}
For a map of generalized ringed spaces
$f \in {\cal GRS(X,Y)}$,
and for a point $x \in \cX$,
we get the induced homomorphism on \emph{stalks} \\
$f^{\sharp}_x \in \cG\cR(\cO_{\cY,f(x)}, \cO_{\cX,x})$,
via

\begin{flushleft}
\begin{equation}\label{12.1.6}
\displaystyle{f^{\sharp}_x: \cO_{\cY,f(x)}=
\lim\limits_{\begin{array}{c}
\stackrel{\textstyle \longrightarrow}{V \subseteq \cY \ \mbox{open}}\\
f(x) \in V
\end{array}}
\hspace{-1cm}
\cO_{\cY}(V)
\xrightarrow{\lim\limits_{\rightarrow} f_V^{\sharp}}
\lim\limits_{\begin{array}{c}
\stackrel{\textstyle \longrightarrow}{V \subseteq \cY \ \mbox{open}}\\
x \in f^{-1}V
\end{array}}
\hspace{-1cm}
\cO_{\cX}(f^{-1}V)
\rightarrow}
\end{equation}
\end{flushleft}

\begin{flushright}
$\displaystyle{\rightarrow
\lim\limits_{\begin{array}{c}
\stackrel{\textstyle \longrightarrow}{U \subseteq \cX \ \mbox{open}}\\
x \in U
\end{array}}
\hspace{-1cm}
\cO_{\cX}(U)=\cO_{\cX,x}}$
\end{flushright}

\defin{12.1.5}
We let $\cal LGRS \subseteq GRS$
denote the subcategory of $\cal GRS$
of \emph{locally generalized ringed spaces}.
Its objects are the objects $(\cX, \cO_{\cX})\in \cG\cR\cS$
such that for all points $x \in \cX$ the stalk
$\cO_{\cX,x} \in \cL\cG\cR$
is a \emph{local}  generalized ring,
i.e. has a unique maximal symmetric ideal $m_{X,x}$.
The maps
$f \in \cL\cG\cR\cS (\cX, \cY)$
are the maps
$(f, f^{\sharp})\in \cG\cR\cS(\cX, \cY)$,
such that for all points $x \in \cX$,
the induced homomorphism on stalks (\ref{12.1.6}) is a \emph{local}
homomorphism,
$f_x^{\sharp} \in \cL\cG\cR (\cO_{\cY,f(x)},
\cO_{\cX,x})$,
$f_x^{\sharp}(m_{Y,f(x)}) \subseteq m_{X,x}$.

\thm{12.1.6}
The functor of global sections
\[
\Gamma: \cL\cG\cR\cS \rightarrow (\cG\cR_C)^{op} , \quad
\Gamma(\cX, \cO_{\cX})=\cO_{\cX}(\cX) , \]
\begin{equation}
\Gamma(f,f^{\sharp}) = f_{\cY}^{\sharp} \ \mbox{for} \
f \in \cL\cG\cR\cS (\cX, \cY)
\end{equation}
and the spectra functor
\[
spec^t:(\cG\cR_C)^{op} \rightarrow {\cal LGRS} ,\quad
spec^t(A) = (spec^t(A), \cO_A) ,\]
\begin{equation}
spec^t(\varphi) = \varphi^* \ \mbox{for} \
\varphi \in \cG\cR(A,B)
\end{equation}
are an adjoint pair:
\[
{\cal LGRS}(\cX , spec^t(A)) =
\cG\cR_C (A,\cO_{\cX}(\cX) ) \]
\begin{equation}\label{adjpair}
\mbox{functorially in} \
\cX \in \cL\cG\cR\cS  , \quad
A \in \cG\cR_C
\end{equation}

\begin{proof}
For a point $x\in \cX$ we have the canonical homomorphism of
taking the stalk at $x$ of a global section, $\phi_x \in
\cG\cR(\cO_{\cX}(\cX) , \cO_{\cX,x})$. Since $\cO_{\cX,x}$ is
local with a unique maximal symmetric ideal $m_{\cX,x}$, we get by
pullback a symmetric prime $\sqp_x = \phi_x^*(m_{\cX,x}) \in
spec^t(\cO_{\cX}(\cX))$. Thus we have a canonical map
\begin{equation}
\sqp:\cX \rightarrow  spec^t(\cO_{\cX}(\cX)) ,\quad
x \mapsto \sqp_x
\end{equation}

The map $\sqp$ is continuous:
For a global section
$g = g^t \in \cO_{\cX}(\cX)^+_{[1]}$,
we have the basic open set
$D_g^+ \subseteq spec^t(\cO_{\cX}(\cX))$,
and
\begin{equation}
\sqp^{-1}(D^+_g) =
\{x \in \cX , \ \sqp_x \in D^+_g \}=
\{x \in \cX , \phi_x(g) \not\in m_{\cX,x}\}
\end{equation}

This set is open in $\cX$, because if $\phi_x(g) \not\in m_{\cX,x}$
we have in $\cO_{\cX,x}$ some $v_x$
with $v_x \lhd \phi_x(g) =1$,
hence there is an open set $U \subseteq \cX$,
with $x \in U$,
and an element $v \in\cO_{\cX}(U)^+_{[1]}$ with $v \lhd g|_{U}=1$,
and for all $x'\in U$,
$v_{x'}\circ \phi_{x'}(g) =1$,
and
$\phi_{x'}(g) \not\in m_{\cX,x'}$.
This shows $\sqp$ is continuous.
The uniqueness of the inverse $v_x= \phi_x(g)^{-1}$ for
$x \in \sqp^{-1}(D_g^+)$
shows we have a well defined inverse
$v = (g|_{\smallsqp^{-1}(D_g)})^{-1} \in \cO_{\cX}(\sqp^{-1}(D_g^+))^+_{[1]}$.
Thus we have a homomorphism of generalized rings
\begin{equation}
\label{eq6.1.8}
\sqp^{\sharp}_{D_g}: \cO_{\cX}(\cX)_g
\rightarrow
\cO_{\cX}(\sqp^{-1}(D_g^+)) , \quad
a/g^n \mapsto v^n \lhd \left(a|_{\smallsqp^{-1}(D_g^+)}\right)
\end{equation}
The collection of homomorphisms
$\{ \sqp^{\sharp}_{D_g} , g \in \cO_{\cX}(\cX)_{[1]}^+ \}$,
are compatible with restrictions,
and the sheaf property gives homomorphisms
\begin{equation}
\sqp_U^{\sharp} \in \cG\cR \left(\cO_{spec \cO_{\cX}(\cX)}(U) ,
\cO_{\cX}(\sqp^{-1}(U))\right).
\end{equation}
Thus we have a map of generalized ringed spaces
\begin{equation}
\sqp=  (\sqp, \sqp^{\sharp}) \in \cG\cR\cS (\cX ,
spec(\cO_{\cX}(\cX))).
\end{equation}
For a point $x \in \cX$, we can take the direct limit of
$\sqp^{\sharp}_{D_g}$, over all global sections
$g \in \cO_{\cX}(\cX)^+_{[1]}$
with
$\phi_x(g) \not\in m_{\cX,x}$,
and we get a local homomorphism
$\sqp_x^{\sharp} \in \cL\cG\cR
(\cO_{\cX}(\cX)_{\smallsqp_x}, \cO_{\cX,x})$.
This shows $\sqp$ is a map of \emph{locally}-ringed spaces,
$\sqp \in {\cal LGRS}(\cX, spec^t(\cO_{\cX}(\cX)))$.

Given a homomorphism of generalized rings
$\varphi \in \cG\cR_C(A, \cO_{\cX}(\cX))$,
we get the map in $\cal LGRS$
\begin{equation}
(spec^t \varphi)\circ \sqp: \cX \rightarrow spec^t (\cO_{\cX}(\cX))
\rightarrow spec^t(A)
\end{equation}
Given a map of locally ringed spaces
$f = (f, f^{\sharp}) \in {\cal LGRS}(\cX, spec^t(A))$,
we get a homomorphism in $\cG\cR_C$,
\begin{equation}
\Gamma (f) = f^{\sharp}_{spec(A)}:A = \cO_A
(spec^t(A)) \rightarrow \cO_{\cX}(\cX)
\end{equation}
These correspondences give the functorial bijection of
(\ref{adjpair}),
we need only show they are inverses of each other.
First for $\varphi \in \cG\cR(A, \cO_{\cX}(\cX))$,
we have
\begin{equation}
\Gamma(spec^t(\varphi)\circ \sqp ) =
\Gamma(\sqp) \circ \Gamma (spec^t(\varphi)) =
id_{\cO_{\cX}(\cX)} \circ \varphi = \varphi
\end{equation}
Fix a map $f=(f, f^{\sharp}) \in {\cal LGRS}(\cX, spec^{t}(A))$.
For a point $x \in \cX$,
we have a commutative square in $\cG\cR_C$
\begin{equation}
\begin{array}{ccc}
A=\cO_A(spec^t(A))& \stackrel{\textstyle \Gamma (f)}{\longrightarrow}&
\cO_{\cX}(\cX)\\
\phi_{f(x)} \downarrow & &\downarrow \phi_x\\
A_{f(x)}= \cO_{A,f(x)}&
 \stackrel{\textstyle f^{\sharp}_x}{\longrightarrow}&
\cO_{\cX,x}
\end{array}
\end{equation}
Since the homomorphism $f_x^{\sharp}$ is assumed to be local, we get
\begin{equation}
\Gamma (f)^*(\sqp_x) = \Gamma (f)^{*}(\phi_{x}^{*}(
  m_{\cX,x})) = \phi_{f(x)}^{*}
(f_x^{\sharp *} (m_{\cX,x}))=
\end{equation}
\[=
\phi_{f(x)}^{*}(m_{A_{f(x)}})=f(x)
\]
This shows that
$(spec^t\ \Gamma(f)) \circ \sqp = f$
as continuous maps.

For a symmetric element $s=s^t \in A_{[1]}$,
we have the commutative square in $\cG\cR_C$,
\begin{equation}
\begin{array}{ccl}
A=\cO_A(spec^t(A))&
\stackrel{\textstyle \Gamma(f)}{\longrightarrow}&
\cO_{\cX}(\cX)\\
\downarrow& &\downarrow\\
A_{s}=\cO_A(D^+_{s})&
\stackrel{\textstyle f_{D_{s}}^{\sharp}}{\longrightarrow}&
\cO_{\cX}(f^{-1}(D^+_{s}))=\cO_{\cX}(D^+_{f^{\sharp}(s)})
\end{array}
\end{equation}
Thus for $a/s^n \in A_{s}$,
we must have
\begin{equation}
f_{D_{s}^+}^{\sharp}(a/s^n)=
\left(\Gamma(f)(s^n)|_{f^{-1}(D^+_{s})}\right)^{-1}
\circ
\left(\Gamma(f)(a)\right)|_{f^{-1}(D^+_{s})}
\end{equation}
\[
=\sqp^{\sharp}_{D^+_{f^{\sharp}(s)}}\circ \Gamma(f) (a/s^n)
\]
This shows that
$f = (spec^t \ \Gamma(f) ) \circ \sqp$
also as maps of generalized-ringed spaces.
\end{proof}

\setcounter{equation}{0}
\section{Schemes}

We define the \emph{Grothendieck-generalized-schemes} to be the
objects of the full subcategory $\cG \cG \cS \subseteq \cL \cG \cR
\cS$, consisting of the $(\cX , \cO_{\cX})$'s which are locally
affine. Later we shall define the category of
\emph{generalized-schemes} $\cG \cS$ to be the pro-category of
$\cG \cG \cS$.

\defin{12.2.1}
An object $\cX =(\cX, \cO_{\cX}) \in \cL\cG\cR\cS$
will be called a
\emph{Grothendieck-generalized-scheme}
if it is locally isomorphic to
$spec(A)$'s:
there exists a covering of $\cX$ by open sets $U_i$,
$\cX = \bigcup\limits_{i}U_i$,
such that the canonical maps are isomorphisms
\begin{equation}
\sqp:(U_i, \cO_{\cX}|_{U_i})
\stackrel{\sim}{\longrightarrow}
spec^t (\cO_{\cX}(U_i))
\end{equation}
We let $\cG\cG\cS$ denote the full sub-category of
$\cL\cG\cR\cS$,
with objects the Grothendieck-generalized-schemes.

\vspace{10pt}

\noindent {\bf Open subschemes 12.2.2} \;
Note that for $\cX \in \cG\cG\cS$,
and for an open set $U \subseteq \cX$,
we have the \emph{open subscheme} of $\cX$
given by $(U, \cO_{\cX}|_{U})$.
That this is again a scheme,
$(U, \cO_{\cX}|_{U})\in \cG\cG\cS$,
follows from the existence of affine basis for the
Zariski topology on $spec^t(A)$,
$A \in \cG\cR_C$,
namely $(D^+_{s}, \cO_A|_{D^+_{s}}) \cong spec^t(A_{s})$
for $s \in A^+_{[1]}$.

\vspace{10pt}

\noindent {\bf Gluing schemes 12.2.3} The local nature of the
definition of Grothendieck-generalized-scheme implies that
$\cG\cG\cS$ admits gluing:

Given $\cX_i \in \cG\cG\cS$,
and open subsets $U_{ij} \subseteq \cX_i$,
and maps $\varphi_{ij} \in \cG\cG\cS (U_{ij}, U_{ij})$,
satisfying the consistency conditions
\begin{equation}
\begin{split}
(i)& \;\;\; U_{ii}=\cX_i, \text{and }\varphi_{ii}=id_{\cX_i}, \\
(ii)& \;\;\; \varphi_{ij}(U_{ij}\cap U_{ik})= U_{ji}\cap U_{jk},
\text{and }\varphi_{jk} \circ \varphi_{ij}= \varphi_{ik}
\text{on } U_{ij}\cap U_{ik},
\end{split}
\end{equation}

\vspace{3mm}

\noindent there exists $\cX \in\cG\cG\cS $, and maps
$\varphi_i \in \cG\cG\cS (\cX_i , \cX)$ such that
\begin{equation}
\begin{split}
(i)& \varphi_i \text{ is an isomorphism of } \cX_i \text{onto an open subset } \varphi_i(\cX_i) \subseteq \cX \\
(ii)& \cX = \bigcup\limits_i \varphi_i(\cX_i) \\
(iii)& \varphi_i(\cX_i) \cap \varphi_j(\cX_j) =
               \varphi_i(U_{ij}),
\text{and } \varphi_j\circ\varphi_{ij}=\varphi_i \text{on } U_{ij}.
\end{split}
\end{equation}

\vspace{10pt}

\noindent {\bf Ordinary Schemes 12.2.4}
For an ordinary scheme
$(\cX, \cO_{\cX})$,
with $\cO_{\cX}$
 a sheaf of commutative rings, there is a covering by open sets
$\cX= \bigcup\limits_i U_i$,
with $(U_i,\cO_{\cX}|_{U_i}) \cong spec(A_i)$,
the ordinary spectrum of the commutative ring
$A_i = \cO_{\cX}(U_i)$.
We then have Grothendieck-generalized schemes
$\cX_i = spec^t(\cG(A_i)) = (U_i, \cG(\cO_{\cX})|_{U_i})$.
These can be glued along $U_{ij}= U_i \cap U_j$,
to a Grothendieck-generalized scheme denoted by
$\cG(\cX) = (\cX, \cO_{\cG(\cX)} = \cG(\cO_{\cX}))$.
It is just the underlying topological space $\cX$
with the sheaf of generalized rings $\cG(\cO_{\cX})$
associated to the sheaf of commutative rings $\cO_{\cX}$
via the functor $\cG : Ring \rightarrow \cG\cR$,
(\ref{embedding}).
Denoting by $\cR\cS$ the category of
(ordinary, commutative) ringed spaces,
the functor
$\cG$ applied to a sheaf of commutative rings $\cO$,
gives a sheaf of commutative generalized rings $\cG(\cO)$,
and we have a functor $\cG:\cR\cS \rightarrow \cG\cR\cS$.
Denoting by $\cL\cR\cS$ (resp. by $\cS$) the
category of
locally-(commutative)-ringed
spaces (resp. the full subcategory of ordinary schemes),
the fact that $\cG$ is fully-faithful implies that we have
full-embeddings of categories:

\begin{equation}
\begin{diagram}
\mbox{\underline{ordinary}}& &\mbox{\underline{generalized}}\\
\cL \cR \cS & \rInto^{\cG} & {\cal LGRS}\\
\uInto & & \uInto\\
\cS & \rInto^{\cG}& \cG\cG \cS\\
\uInto^{spec}& & \uInto_{spec^t}\\
Ring^{op}& \rInto^{\cG}& (\cG\cR_C)^{op}
\end{diagram}
\end{equation}

The generalized Grothendieck scheme
$\cG(X)$, associated with an ordinary scheme $X \in S$,
has always a unique map:
$\cG(X) \rightarrow spec^t\cG ({\mathbb Z})$.

\vspace{20pt}

\begin{center}
\underline{Some Examples of Schemes over $\mathbb F$}
\end{center}

\exmpl{12.2.5} \emph{The (symmetric) affine line} over $\mathbb F$ is given by,
cf. (\ref{eq8.3.63}),
\begin{equation}
{\mathbb A}^1_+ = spec \, \Delta_+^{[1]} = spec \, {\mathbb F}
\{z^{\mathbb N}\}
\end{equation}
We have $F\{z^{\mathbb N}\}_{[1]} = z^{\mathbb N} \cup \{ 0 \}$;
$(0)$ is a prime, the generic point of ${\mathbb A}^1$;
and $(z)$ is a prime, the closed point of
${\mathbb A}_+^1 = \{ (0), (z) \}$.

Similarly, we have
\begin{equation}
\mathbb{A}^1=spec^t\Delta^{[1]}=spec^t\FF\{z^{\tiny \NN}\cdot (z^t)^{\tiny \NN}\}.
\end{equation}
The symmetric prime $(0)$ is the generic point of $\mathbb{A}^1$,
and the symmetric prime $(z\cdot z^t)$ is the closed point:
$\mathbb{A}^1=\{(0),(z\cdot z^t)\}$. The homomorphism
$\FF\{z^{\tiny \NN}\cdot (z^t)^{\tiny \NN}\}\sur \FF\{z^{\tiny
\NN}\},\;z,z^t\mapsto z$, gives the immersion
\begin{equation}
\mathbb{A}^1_+\inj \mathbb{A}^1
\end{equation}

\vspace{20pt}

\exmpl{12.2.6}
\emph{The (symmetric) multiplicative group} over $\mathbb F$
is given by
\begin{equation}
{\mathbb G}_m^+ = spec \,\FF\{z^{\mathbb Z}\}=
\{(0)\} \subseteq {\mathbb A}^1_+
\end{equation}
Note that for a commutative generalized ring $A$,
\begin{equation}
\cG\cR_C (\FF\{z^{\mathbb Z}\},A) =A^*\cap A^+_{[1]}=
\{a \in A_{[1]}^+, \ \mbox{there is} \ a^{-1} \in A_{[1]}^+ ,
a \circ a^{-1}=1 \}
\end{equation}
\noindent Similarly,
\begin{equation}
{\mathbb G}_m=spec^t\;\FF[z^{\tiny \ZZ}\cdot (z^t)^{\tiny \ZZ}]=\{(0)\}\subseteq \mathbb{A}^1,
\end{equation}
and \\
$\cG\cR_C(\FF[z^{\tiny \ZZ}\cdot (z^t)^{\tiny \ZZ}],A)=A^*$. We have an immersion
\begin{equation}
{\mathbb G}_m^+\inj {\mathbb G}_m
\end{equation}
\vspace{10pt}

\exmpl{12.2.7} \emph{The (symmetric) projective line} over $\mathbb F$
is obtained by gluing two (symmetric) affine lines along
${\mathbb G}_m^{(+)}$
\begin{equation}
{\mathbb P}^1_+ = spec  \, \FF\{z^{\mathbb N}\}
\coprod\limits_{spec  \, \FF\{z^{\mathbb Z}\}}
spec \, \FF\{(z^{-1})^{\mathbb N}\}=
\{ m_1 , m_0, m_{\infty}\}
\end{equation}
Respectfully, the full
\begin{equation}
\mathbb{P}^1=spec^t\; \FF\{z^{\tiny \NN}\cdot (z^t)^{\tiny \NN}\}\underset{spec^t\;\FF\{z^{\tiny \ZZ}\cdot (z^t)^{\tiny \ZZ}\}}{\coprod} spec^t\; \FF\{(z^{-1})^{\tiny \NN}\cdot ((z^t)^{-1})^{\tiny \NN}\}
\end{equation}
It has a generic point $m_1 = (0)$,
and two closed points $m_0 = (z)$, (resp. ($z\cdot z^t$)) $m_{\infty}= (z^{-1})$, (resp. $(z^{-1}\cdot (z^t)^{-1})$).
We have the immersion
\begin{equation}
\mathbb{P}_+^1\inj \mathbb{P}^1
\end{equation}

Interchanging $z$ and $z^{-1}$ we get an involutive automorphism

\begin{equation}\label{inversion}
\displaystyle{I: {\mathbb P}^1_+
\xrightarrow{\sim} {\mathbb P}^1_+}, \quad
I \circ I = id_{{\mathbb P}^1_+}
\end{equation}

interchanging $m_0$ and $m_{\infty}$.

\vspace{10pt}

Every rational number $f \in {\mathbb Q}^*$,
defines a geometric map
\begin{equation}
f_{\mathbb Z} \in \cG\cG\cS (spec \, \cG({\mathbb Z}), {\mathbb P^1_+})
\end{equation}

If $f = \pm 1$ this is given by the constant map
\begin{equation}
\mathbb F \{z^{\mathbb Z}\} \twoheadrightarrow {\mathbb F}\{\pm 1\}
\subseteq \cG (\mathbb Z) ,
z \mapsto f = \pm 1
\end{equation}
If $f \neq \pm 1$,
let $N_0 = \prod\limits_{\nu_p (f) >0} p$,
$N_{\infty} = \prod\limits_{\nu_p (f) <0} p$,
then
\begin{equation}
spec \, \cG (\ZZ) = spec \,  \cG (\ZZ [\frac{1}{N_0}])
\coprod\limits_{spec \, \cG(\mathbb Z [\frac{1}{N_0N_{\infty}}])}
spec \, \cG (\mathbb Z [\frac{1}{N_{\infty}}])
\end{equation}
and the geometric map $f_{\mathbb Z}$ is given by the spec-maps associated
to the homomorphisms:
\begin{equation}\label{geometricf}
\begin{array}{ccc}
\mathbb F \{z^{\mathbb N}\} & \longrightarrow &
       \cG(\mathbb Z [\frac{1}{N_{\infty}}])\\
| \bigcap & & | \bigcap\\
\mathbb F \{z^{\mathbb N}\} & \longrightarrow &
       \cG(\mathbb Z [\frac{1}{N_0 \cdot N_{\infty}}])\\
| \bigcup & & | \bigcup\\
\mathbb F \{(z^{-1})^{\mathbb N}\} & \longrightarrow &
       \cG(\mathbb Z [\frac{1}{N_{0}}])
\end{array}
\end{equation}
\[
\begin{array}{c}
z \mapsto f\\
z^{-1} \mapsto f^{-1}
\end{array}
\]

\section{Projective limits}
\setcounter{equation}{0}
\noindent The category of locally generalized ringed spaces
$\cal LGRS$
admits directed inverse limits. For a partially ordered set
$J$,
which is directed
(for $j_1, j_2 \in J$, have $j \in J$ with $j \geq j_1$,
$j \geq j_2$)
and for a functor $\cX:J\rightarrow \cL\cG\cR\cS$,
$J \ni j \mapsto \cX_j$,
$j_1 \geq j_2 \mapsto \pi_{j_2}^{j_1}\in \cL\cG\cR\cS (\cX_{j_1} ,
 \cX_{j_2})$,
we have the inverse limit
$\lim\limits_{\stackrel{\longleftarrow}{J}}\cX \in {\cal LGRS}$.
The underlying topological space of
$\lim\limits_{\stackrel{\longleftarrow}{J}}\cX$
is the inverse limit of the sets $\cX_j$,
with basis for the topology given by the sets
$\pi_j^{-1}(U_j)$,
with $U_j \subseteq \cX_j$
open, and where
$\pi_j : \lim\limits_{\stackrel{\longleftarrow}{j \in J}}
\cX_j \rightarrow \cX_j$
denote the projection.
The sheaf of generalized rings
$\cO_{\lim\limits_{\leftarrow}\cX}$ over
$\lim\limits_{\stackrel{\longleftarrow}{J}}\cX_j$,
is the sheaf associated to the pre-sheaf
$U \mapsto \lim\limits_{\stackrel{\longrightarrow}{J}}
\pi_j^* \cO_{\cX_j}(U)$.
For a point
$x =(x_j) \in \lim\limits_{\longleftarrow}\cX_j$,
the stalk $\cO_{\lim\limits_{\longleftarrow}\cX,x}$
is the direct limit of the local-generalized-rings
$\cO_{\cX_j, x_j}$,
and hence is local, and
$(\lim\limits_{\longleftarrow}\cX_j,
\cO_{\lim\limits_{\longleftarrow}\cX})
\in {\cal LGRS}$.

An alternative explicit description of the sections
$s \in \cO_{\lim\limits_{\longleftarrow}\cX}(U)$,
for $U \subseteq \lim\limits_{\longleftarrow}\cX_j$
open,
are as maps
\begin{equation}
s:U \rightarrow \coprod\limits_{x \in U}
    \cO_{\lim\limits_{\longleftarrow}\cX,x},
\ \mbox{with} \
s(x) \in \cO_{\lim\limits_{\longleftarrow}\cX,x}
\end{equation}
such that for all $x \in U$,
there exists $j \in J$, and open subset $U_j \subseteq \cX_j$,
with $x \in \pi_j^{-1}(U_j) \subseteq U$
and there exists a section $s_j \in \cO_{\cX_j}(U_j)$,
such that for all $y \in \pi_j^{-1}(U_j)$,
we have
$s(y) = \pi_j^{\sharp}(s_j) |_{y}$.

We have the universal property
\begin{equation}
\label{universal}
{\cal LGRS}(Z, \lim\limits_{\stackrel{\longleftarrow}{J}}\cX_j)=
\lim\limits_{\stackrel{\longleftarrow}{j \in J}} {\cal LGRS} (Z, \cX_j)
\end{equation}
Note that if
$\cX_j = spec^t(A_j)$ are affine generalized schemes, then the inverse limit
\begin{equation}
\lim\limits_{\stackrel{\longleftarrow}{J}}(spec^t(A_j)) =
spec^t(\lim\limits_{\stackrel{\longrightarrow}{J}}A_j)
\end{equation}
is the affine generalized scheme associated to
$\lim\limits_{\stackrel{\longrightarrow}{J}}A_j$
the direct limit of the $A_j$'s computed in $\cal GR$.
(Hence in $Set_0$, cf. (\ref{eq2.7.3})).

Note on the other hand that the category of Grothendieck-
generalized schemes $\cal GGR$ is not closed under directed
inverse limits (just as in the "classical" counterparts, the
category $\cal LRS$ of locally ringed spaces (resp. $Ring$) is
closed under directed inverse (resp. direct) limits , while the
category $\cS$ of schemes is not closed under directed inverse
limits). The  point is: for a point $x = (x_j) \in
\lim\limits_{\longleftarrow} \cX_j$, in the inverse limit of the
Grothendieck (generalized) schemes $\cX_j$, while each $x_j \in
\cX_j$ has an open affine neighborhood, $x_j \in spec^t A_j
\subseteq \cX_j$, there may \emph{not} be an open affine
neighborhood of $x$ in $\lim\limits_{\longleftarrow} \cX_j$.

\defin{12.3.1}
The category of generalized schemes $\cal GS$
is the category of pro-objects of the category of
Grothendieck-generalized schemes,
\begin{equation}
\cG\cS = pro\mbox{-}\cG\cG\cS.
\end{equation}

Thus the objects of $\cal GS$
are inverse systems
$\cX=(\{\cX_j\}_{j \in J}, \{\pi_{j_2}^{j_1}\}_{j_1 \geq j_2})$,
where $J$ is a directed partially ordered set,
$\cX_j \in {\cal GGS}$ for $j \in J$,
and $\pi_{j_2}^{j_1}\in {\cal GGS}(\cX_{j_1}, \cX_{j_2})$
for $j_1 \geq j_2$, $j_1, j_2 \in J$,
with $\pi_{j}^{j} = id_{\cX_j}$,
and $\pi_{j_3}^{j_2} \circ \pi_{j_2}^{j_1}
=\pi_{j_3}^{j_1}$ for $j_1 \geq j_2 \geq j_3$.
The maps from such an object to another object
$\cY = (\{\cY_i\}_{i \in I}, \{\pi_{i_2}^{i_1}\}_{i_1 \geq i_2})$
are given by
\begin{equation}
\cG\cS (\cX, \cY )=
\lim\limits_{\stackrel{\longleftarrow}{I}}
\lim\limits_{\stackrel{\longrightarrow}{J}}
\cG\cG\cS (\cX_j , \cY_i )
\end{equation}
i.e. the maps $\varphi \in \cG\cS(\cX, \cY)$
are a collection of maps
$\varphi_i^j \in \cG\cG\cS (\cX_j, \cY_i)$ defined for all
$i \in I$,
and for $j \geq \tau(i)$ sufficiently large
(depending on $i$),
these maps satisfy:
\begin{itemize}
\item[]
for all $i \in I$, and for $j_1 \geq j_2$ sufficiently large in $J$:
\begin{equation}
\varphi_i^{j_1} =\varphi_i^{j_2} \circ \pi_{j_2}^{j_1}
\end{equation}
\item[] for all $i_1 \geq i_2$ in $I$,
and for $j \in J$ sufficiently large:
\begin{equation}
\pi_{i_2}^{i_1} \circ \varphi_{i_1}^{j} = \varphi_{i_2}^{j}
\end{equation}
\end{itemize}
 The maps $\varphi = \{\varphi_i^j \}_{j \geq \tau(i)}$,
and $\tilde{\varphi} = \{\tilde{\varphi}_i^j\}_{j \geq
       \tilde{\tau}(i)}$,
are considered equivalent if
\begin{itemize}
\item[] for all $i \in I$, and for $j \in J$
sufficiently large:
\begin{equation}
\varphi_i^j = \tilde{\varphi}_i^j
\end{equation}
\end{itemize}
The composition of $\varphi = \{\varphi_i^j\}_{j \geq \tau(i)} \in
\cG\cS(\cX, \cY )$,
with $\psi = \{\psi^i_k \}_{i \geq \sigma(k)} \in
\cG\cS(\cY, \cZ)$,
is given by
$\psi \circ \varphi = \{ \psi^i_k \circ \varphi_i^j \}_{j \geq \tau(\sigma(k))} \in \cG\cS(\cX, \cZ)$.

There is a canonical map (which in general is not injective or surjective, but is so for"finitely- presented" $\{X_j\}$ and $\{Y_i\}$, see \cite{EGA}),
\begin{equation}
\label{cannonicalmap}
\lim\limits_{\stackrel{\longrightarrow}{J}} {\cal LGRS}
(\cX_j , \cY_i) \longrightarrow {\cal LGRS}
(\lim\limits_{\stackrel{\longleftarrow}{J}} \cX_j, \cY_i)
\end{equation}
By the universal property (\ref{universal}) we have bijection
\begin{equation}
\label{uni_bijection}
\lim\limits_{\stackrel{\longleftarrow}{I}} {\cal LGRS}
(\lim\limits_{\stackrel{\longleftarrow}{J}} \cX_j, \cY_i)=
{\cal LGRS}
(\lim\limits_{\stackrel{\longleftarrow}{J}} \cX_j ,
\lim\limits_{\stackrel{\longleftarrow}{I}} \cY_i)
\end{equation}
Composing (\ref{cannonicalmap}) and (\ref{uni_bijection})
we obtain a map
\begin{equation}
\cL: \lim\limits_{\stackrel{\longleftarrow}{I}}
     \lim\limits_{\stackrel{\longrightarrow}{J}}
{\cal LGRS} (\cX_j , \cY_i) \longrightarrow {\cal LGRS}
(\lim\limits_{\stackrel{\longleftarrow}{J}} \cX_j ,
\lim\limits_{\stackrel{\longleftarrow}{I}} \cY_i)
\end{equation}
Thus we have a functor
\begin{equation}
\cL: {\cal GS} \longrightarrow {\cal LGRS} , \quad
\cL(\{\cX_j\}_{j \in J}) =
\lim\limits_{\stackrel{\longleftarrow}{J}} \cX_j
\end{equation}
We view the category $\cG\cG\cS$ as a full subcategory of
$\cG\cS$
(consisting of the objects
$\cX = \{ \cX_j\}_{j \in J}$,
with indexing set $J$ reduced to a singleton).

\section{The compactified $\overline{spec \mathbb{Z}}$}
\label{sec6.4} We denote by $\eta$ the \emph{real prime} of
$\mathbb{Q}$, so $|\ |_{\eta}: \QQ \rightarrow [0, \infty )$ is
the usual (non archimedean) absolute value, and we let
$\cO_{\eta}$ denote the associated generalized ring
(\ref{realprimes}), $\cO_{\eta} \subseteq \cG(\mathbb{Q})$. For a
square-free integer $N \geq 2$, we have the sub-generalized-ring
\begin{equation}
A_N = \cG(\mathbb{Z}[\frac{1}{N}]) \cap \cO_{\eta}
\subseteq \cG(\mathbb{Q})
\end{equation}
The localization of $A_N$ with respect to
$\frac{1}{N} \in A_{N,[1]}$
gives
$(A_N)_{\frac{1}{N}}= \cG(\mathbb{Z}[\frac{1}{N}])$, \\
so the inclusion $j_N: A_N \hookrightarrow \cG(\mathbb{Z}[\frac{1}{N}])$
gives the basic open set
\begin{equation}
j_N^*:spec(\mathbb{Z}[\frac{1}{N}])=
      spec^t \, \cG(\mathbb{Z}[\frac{1}{N}])
\stackrel{\textstyle \sim}{\longrightarrow}
D^+_{\frac{1}{N}} \subseteq spec^t(A_N)
\end{equation}
The inclusion
$i_N:A_N \hookrightarrow \cO_{\eta}$,
gives the real prime $\eta_N \in spec^t \,(A_N)$,
\begin{equation}
\eta_N = i^*_N (\sqm_{\eta}), \quad
(\eta_N)_X = \{a = (a_x) \in (\ZZ[\frac{1}{N}])^X,
||a||^2 = \sum\limits_{x \in X}|a_x|^2 < 1 \}
\end{equation}
Note that $\eta_N$ is the unique maximal ideal of $A_N$,
and $A_N$ is a local generalized ring.
Let $\cX_N$ denote the Grothendieck generalized scheme obtained by gluing
$spec^t \, (A_N)$ with $spec^t \,  \cG(\mathbb{Z})$ along the common (basic) open set
$spec^t(\cG(\mathbb{Z}[\frac{1}{N}]))$.
The open sets of $\cX_N$ are the open sets
$U_{M} = spec(\mathbb{Z}[\frac{1}{M}]) \subseteq spec(\mathbb{Z})$,
(and $\cO_{\cX_N}(U_M)= \cG(\mathbb{Z}[\frac{1}{M}])$),
as well as the sets $\{\eta_N\} \cup U_M $,
\emph{with} $M$ \emph{dividing} $N$
(and $\cO_{\cX_N}(\{\eta_N\} \cup U_M)=A_M$, $M|N$).

For $ N_2$ dividing $N_1$,
we have a map $\pi_{N_2}^{N_1}\in \cG\cG\cS (\cX_{N_1}, \cX_{N_2})$
induced by the inclusions \\
$A_{N_2} \hookrightarrow A_{N_1}$,
and $\cG(\mathbb{Z}[\frac{1}{N_2}]) \hookrightarrow
     \cG(\mathbb{Z}[\frac{1}{N_1}])$.\\
Note that $\pi_{N_2}^{N_1}$ is a bijection on points,
and that moreover,
\begin{equation}
(\pi_{N_2}^{N_1})_* \cO_{\cX_{N_1}}= \cO_{\cX_{N_2}} \; \text{and } (\pi_{N_2}^{N_1})^{\sharp} \text{is the identity map of }
\cO_{\cX_{N_2}}
\end{equation}
But there are more open sets in $\cX_{N_1}$
then there are in $\cX_{N_2}$! .

The compactified
$\overline{spec \,  \mathbb{Z}}$
is the object of $\cG\cS = pro \cG\cG\cS$ given by
$(\{\cX_N\}, \{\pi_{N_2}^{N_1}\}_{N_2|N_1})$,
\begin{equation}\label{spec_z}
\overline{spec \, \mathbb{Z}} =
\left\{
\cX_N = spec(A_N) \coprod\limits_{spec \, \cG(\mathbb{Z}[\frac{1}{N}])}
        spec \, \cG(\mathbb{Z})
\right\}_{N \geq 2 \ \mbox{square free}}
\end{equation}
Note that the associated locally-generalized-ring space
\begin{equation}
\cX= \cL(\overline{spec \mathbb{Z}})=
\lim\limits_{\stackrel{\longleftarrow}{N}} \cX_N \in {\cal LGRS}
\end{equation}
has underlying topological space
$\cX = \{\eta\} \coprod spec(\mathbb{Z})$, with open sets \\
$U_M = spec(\mathbb{Z}[\frac{1}{M}])$
(and $\cO_{\cX}(U_M) = \cG(\mathbb{Z}[\frac{1}{M}])$),
as well as the sets $\{\eta\} \coprod U_M$,
\emph{with no restrictions on} $M$,
and $\cO_{\cX}(\{\eta\} \coprod U_M) = A_M$
for $M \geq 2$,
while the global sections are
$\cO_{\cX}(\cX)= {\mathbb F}\{\pm 1\}$.

The stalks of $\cO_{\cX}$ are given by
\begin{equation}
\begin{array}{lll}
\cO_{\cX, p} = \cG(\mathbb{Z}_{(p)}) , & p \in spec(\mathbb{Z}) , &
\mathbb{Z}_{(p)}= \{ \frac{m}{n} \in \mathbb{Q} , p \not\mid n\},\\
\cO_{\cX, \eta} = \cO_{\eta} & &
\end{array}
\end{equation}

\vspace{10pt}

Similarly for a number field $K$,
with ring of integers $\cO_K$, and with real
and complex
 primes $\eta_i$,
$i = 1, \ldots , \gamma= \gamma_{\mathbb{R}} + \gamma_{\mathbb{C}}$,
we have the sub-generalized-ring of $\cG(K)$ given by
\begin{equation}
A_{N,i}= \cG(\cO_K[\frac{1}{N}]) \cap \cO_{K,\eta_i}
\subseteq \cG(K)
\end{equation}
Let $\cX_N$ be the Grothendieck generalized scheme obtained by gluing
$\{spec (A_{N,i})\}_{i \leq \gamma }$ and
$\{spec (\cG(\cO_K))\}$ along the common (basic) open set
$spec \left(\cG(\cO_K[\frac{1}{N}])\right)$. \\
For $N_2|N_1$,
we have $\pi_{N_2}^{N_1} \in \cG\cG\cR(\cX_{N_1}, \cX_{N_2})$
induced by the inclusions $A_{N_2,i} \hookrightarrow A_{N_1,i}$.
We get the compactified $\overline{spec(\cO_K)}$,
it is the object of $\cG\cS$ given by the $\cX_N$'s and
$\pi_{N_2}^{N_1}$'s.\\
The space
\begin{equation}
\cX_K = \cL(\overline{spec(\cO_K)})=
\lim\limits_{\stackrel{\longleftarrow}{N}}\cX_N \in {\cal LGRS}
\end{equation}
has for points the set
$spec(\cO_K) \coprod \{\eta_i\}_{i \leq \gamma}$,
and for open subsets the sets $U \coprod \{\eta_i\}_{i \in I}$,
$U \subseteq spec(\cO_K)$ open,
$I \subseteq \{1, \ldots ,\gamma\}$, where
\begin{equation}
\cO_{\cX_K}(U \coprod \{\eta_i\}_{i \in I}) =
\bigcap\limits_{p \in U}\cG(\cO_{K,p})\cap
\bigcap\limits_{i \in I} \cO_{K,\eta_i}
\end{equation}
In particular, the global sections are
\begin{equation}
\cO_{\cX_K}(\cX_K) = \bigcap\limits_{p \in spec\cO_K}
\cG(\cO_{K,p}) \cap \bigcap\limits_{i \leq \gamma} \cO_{K,\eta_i}=
\mathbb{F}\{\mu_K\}
\end{equation}
with $\mu_K \subseteq \cO^*_K$
the group of roots of unity in $\cO^*_K$.

\vspace{20pt}

Returning for simplicity to the rational
$\cX = \overline{spec {\mathbb Z}}$
case of (\ref{spec_z}),
every rational number $f \in \mathbb Q^*$,
defines a geometric map
\begin{equation}
\underline{f} \in \cG\cS (\cX , \mathbb P^1_+)
\end{equation}
i.e. a collection of maps
$\underline{f}_N \in \cG\cG\cS (\cX_N , \mathbb P^1_+)$,
for $N$ divisible by $N_0 \cdot N_{\infty}$, \\
$N_0 = \prod\limits_{\nu_p(f) >0} p$,
$N_{\infty} = \prod\limits_{\nu_p(f) <0} p$,
with $\underline{f}_N \circ \pi_N^M = \underline{f}_M$. \\
For $f = \pm 1$ it is the constant map given by
\begin{equation}
\mathbb F\{z^{\mathbb Z}\} \twoheadrightarrow \mathbb F \{\pm 1\} =
\cG(\mathbb Z) \cap A_N , \ \mbox{for any} \ N
\end{equation}
\[z \mapsto f = \pm 1
\]
For $f \neq \pm 1$,
 we may assume $|f|_{\eta} < 1$, by the commutativity of
the diagram
(with $I$ the inversion (\ref{inversion})),

\noindent
\begin{equation}
\begin{diagram}
 & & \mathbb P^1_+\\
& \ruTo^{\underline{f}}& \\
\cX & &\uTo^{\wr}_{I}\\
  &\rdTo_{\underline{f^{-1}}} &\dTo\\
& & \mathbb P^1_+
\end{diagram}
\end{equation}

Thus for $N$ divisible by $N_0 \cdot N_{\infty}$ we have $f \in
A_N$, and the map $\underline{f}_N$ is given by $\underline{f}_N =
f_{\mathbb Z} \coprod f_{\eta , N}$, with $f_{\eta , N} = spec
f^{\sharp}_{\eta , N}$,
$f^{\sharp}_{\eta , N} \in \cG\cR ({\mathbb F}\{z^{\mathbb N}\}, A_N)$,\\
the unique homomorphism with $f^{\sharp}_{\eta , N}(z) =f$,
and $f_{\mathbb Z}$ is as in (\ref{geometricf}).


\begin{equation}
\begin{array}{ccclcccc}
&\cX_N= & spec \, \cG (\ZZ )
\coprod\limits_{spec \, \cG (\mathbb Z [\frac{1}{N}])} spec \ A_N&
\cG (\mathbb Z [\frac{1}{N}])\supseteq  A_N  \ni f & \\
\vspace{5pt}\\
&\stackrel{\displaystyle{\downarrow}}{\mid} \underline{f}_N= & \;\;\;
\stackrel{\displaystyle{\downarrow}}{\mid} f_{\small \ZZ} \coprod f_{\eta , N}
&\stackrel{\displaystyle{\uparrow}}{\mid} \hspace{0.3cm}
f^{\sharp}_{\eta, N}\hspace{0.6cm}
\stackrel{\displaystyle{\uparrow}}{\mid}\hspace{0.4cm}
\stackrel{\displaystyle{\uparrow}}{\perp}& \\
\vspace{5pt}\\
&\mathbb P^1_+ = & spec \, \FF \{(z^{-1})^{\small \NN}\}
\coprod\limits_{spec \, \FF \{z^{\ZZ}\}}
spec \, \FF \{z^{\NN}\}&
\FF \{z^{\mathbb Z}\} \supseteq \FF \{z^{\mathbb N}\} \ni z&
\end{array}
\end{equation}

Similarly for a number field $K$,
every element $f \in K^*$ defines a geometric map
\begin{equation}
\underline{f} \in \cG\cS (\overline{spec \, \cO_K}, \mathbb P^1_+)
\end{equation}

\setcounter{equation}{0}
\chapter{Products}
\section{Tensor product}\label{13.1}

The category $\cG\cR_C$ of
commutative generalized rings has tensor-products, i.e. fibred sums:
Given homomorphisms $\varphi^j \in \cG\cR_C (A,B^j)$
$j=0,1$, there exists $B^0 \bigotimes\limits_{A}B^1 \in \cG\cR_C$,
and homomorphisms
$\psi^j \in \cG\cR_C (B^j,B^0 \bigotimes\limits_{A}B^1)$,
such that\\
$\psi^0 \circ \varphi^0 = \psi^1 \circ \varphi^1$,
and for any $C \in \cG\cR_C$,
$$
\cG\cR_C (B^0 \bigotimes\limits_{A}B^1, C) =
\cG\cR_C(B^0,C) \prod\limits_{\cG\cR_C(A,C)} \cG\cR_C (B^1,C)
$$
So given homomorphisms $f^j \in \cG\cR_C (B^j,C)$
with $f^0 \circ \varphi^0 = f^1 \circ \varphi^1$,
there exists a unique homomorphism
$f^0 \otimes f^1 \in \cG\cR_C (B^0 \bigotimes\limits_{A}B^1 ,C)$,
such that $(f^0 \otimes f^1) \circ \psi^j = f^j$.

The construction of $B^0 \bigotimes\limits_{A}B^1$ goes as follows.
First for a finite set\\
$\{b_1^0, \ldots , b_n^0,b_1^1, \ldots , b_m^1 \}$, where $b_i^j
\in B^j_{X_i^j}$, we have the free commutative generalized ring on
the sets $\{X_1^0, \ldots , X_n^0,X_1^1, \ldots , X_m^1 \}$
(\ref{8.3.54}), and we write $\underline{b}_1^0, \ldots ,
\underline{b}_n^0, \underline{b}_1^1, \ldots ,\underline{b}_m^1 $
for its canonical generators. Taking the direct limit over such
finite subsets, cf. \S \ref{Limits}, we have the free generalized
ring $\Delta$ with generators $\underline{b}$, with $b \in B^0_X$
or $b \in B^1_X$, and any $X \in \FF$. We divide $\Delta$ by the
equivalence ideal $\varepsilon_A$ generated by
\begin{equation}
\begin{array}{ccl}
\underline{b} \lhd \underline{b'} \sim \underline{b \lhd b'}&,&
b,b' \in B^j \ , \ j=0,1;\\
\underline{b}\sslash \underline{b'} \sim \underline{b\sslash b'}&,&
b,b' \in B^j \ , \ j=0,1;\\
\underline{1^j} \sim 1 &,&
\mbox{where} \ 1^j \in B^j_{[1]} \ \mbox{is the unit};\\
\underline{\varphi^0(a) } \sim \underline{\varphi^1(a)}&,&
\mbox{for} \ a \in A
\end{array}
\end{equation}
The quotient generalized ring
$\Delta/\varepsilon_A$ is the tensor product
$B^0 \bigotimes\limits_{A}B^1$,\\
the homomorphism $\psi^j$ is given by
$\psi^j(b) = \underline{b} \  mod \ \varepsilon_A$,
$b \in B^j$.

Note that every element of $(B^0 \bigotimes\limits_{A}B^1)_X$
can be expressed (non-uniquely) as
\begin{equation}\label{13.1.2}
(\underline{a},\underline{b}) = (\underline{a}_1 \lhd
\underline{a}_2 \lhd \cdots \lhd \underline{a}_n)\sslash
(\underline{b}_1 \lhd \cdots \lhd \underline{b}_m)\ mod \ \varepsilon_A
\end{equation}
with $a_i \in B_{f_i}^{i(mod 2)}$,
$b_j \in B_{g_j}^{j(mod 2)}$,
and $f_1 \circ \cdots \circ f_n = c_X \circ g_1 \circ \cdots \circ g_m$
(where $c_X \in Set_{\bullet}(X,[1])$ is the canonical map,
$c_X(x) =1$ for all $x \in X$).
These elements  are multiplied and contracted by the formulas of multiplication
\ref{eq8.2.7} and contraction \ref{eq8.2.8}.

\exmpl{13.1.1}
For monoids $M_0$, $M_1$, $N$,
and homomorphisms $\psi^i \in Mon (N, M_i)$,
$i=0,1$, we have (by adjunction (\ref{mon_adj})),
\begin{equation}
\mathbb{F}\{M_0\} \bigotimes\limits_{\mathbb{F}\{N\}}
\mathbb{F}\{M_1\} = \mathbb{F}\{M_0 \bigotimes\limits_{N} M_1\}
\end{equation}
where $M_0 \bigotimes\limits_{N} M_1$
is the fibered sum in the category $Mon$.
The monoid $M_0 \bigotimes\limits_{N} M_1$ is given by elements
$m_0 \otimes m_1$, $m_i \in M_i$,
with relations
\begin{equation}\label{eq13.1.4}
m_0 \otimes 0 = 0 \otimes 0 = 0 \otimes m_1 \ , \ m_i \in M_i
\end{equation}
and
\begin{equation}
m_0\cdot \psi^0(n) \otimes m_1 =
m_0  \otimes \psi^1(n) \cdot m_1 \ , \ n \in N
\end{equation}

\exmpl{13.1.2}
For a commutative ring $B$,
let $B^!$ denote the underlying multiplicative monoid of $B$
(i.e. forget addition),
and let $\mathbb{F}\{B^!\}$ denote the associated generalized ring, cf. \S \ref{s8.3.5}
From the identity map $B^!
\xrightarrow{=}\cG(B)_{[1]}$,
we obtain by adjunction (\ref{mon_adj})
the canonical injective homomorphism
$J_B \in \cG\cR(\mathbb{F}\{B^!\} , \cG(B))$.
The unique homomorphism of rigs
$\mathbb{N} \rightarrow B$,
gives the unique homomorphism of generalized rings
$I_B \in \cG\cR (\cG(\mathbb{N}), \cG(B) )$.
We get a canonical homomorphism of generalized rings,

\begin{equation}\label{Psi}
\displaystyle{\Psi_B = I_B \otimes J_B \in \cG\cR_C\left(\cG(\mathbb{N})
\bigotimes\limits_{\mathbb{F}} \mathbb{F} \{B^!\} , \cG(B)
\right)}
\end{equation}
The homomorphism $\Psi_B$ is always surjective
(as follows from (\ref{8.3.19})).\\
For any monoid $B$, the
elements of the
generalized ring
${\cal N}^B = \cG(\mathbb{N}) \bigotimes\limits_{\mathbb{F}} \mathbb{F}\{B\}$,
can be described as
in (\ref{13.1.2}), but we can move the elements of $\mathbb F \{B\}$
to the right (using (\ref{eq_up1})),
and we can take the elements of $\cG(\mathbbm N)$
to be the generators ${\mathbbm 1}_Z$;
thus we can write every element of ${ \cN}_{X}^{B}$
as $\left({\mathbbm 1}_{\widetilde{X}} \circ \mu, {\mathbbm 1}_{\pi}\right)$,
with $\pi \in Set_{\bullet}(\widetilde{X}, X)$,
and $\mu \in (B)^{\widetilde{X}}$,

\begin{equation}
\displaystyle{{\cal N}_X^B = \left\{
(\pi: \widetilde{X} \rightarrow X, \mu: \widetilde{X} \rightarrow B) \right\}/\approx}
\end{equation}

The elements of ${\cal N}_X^B$ are (isomorphism classes of) sets over
$X \prod B$,
where the equivalence relation $\approx$
is invariant by isomorphisms, i.e.

$(\pi:\widetilde{X} \rightarrow X, \mu: \widetilde{X} \rightarrow B)
\approx
(\pi':\widetilde{X'} \rightarrow X,
\mu': \widetilde{X'} \rightarrow B)$ if there is a bijection
$\sigma:\widetilde{X} \rightarrow \widetilde{X}'$,
$\pi=\pi' \circ \sigma$,
$\mu= \mu' \circ \sigma$,

and by zero, i.e.

\begin{equation}
\displaystyle{(\widetilde{X} , \pi , \mu)
\approx (\widetilde{X} \setminus \{x\} ,
\pi|_{\widetilde{X} \setminus \{x\}},
\mu|_{\widetilde{X} \setminus \{x\}})}\text{if } \mu(x)=0.
\end{equation}

For $f \in Set_{\bullet}(X,Y)$, and for
$(\widetilde{X} , \mu ) \in {\cal N}_X^B$,
$(Z, \lambda) \in {\cal N}_f^B$,
we have the contraction, cf. \S 8.2.8,

\begin{equation}
\displaystyle{\left( (\widetilde{X}, \mu)\sslash (Z, \lambda) \right) =
\left( \widetilde{X} \prod\limits_{X} Z,
(\mu\sslash \lambda)\right)}
\end{equation}
\[
(\mu\sslash \lambda)(x,z) = \mu(x)\cdot \lambda(z)
\]

For $(\widetilde{Y}, \mu ) \in {\cal N}_Y^B$
we have the multiplication, cf. \S 8.2.7.

\begin{equation}
\displaystyle{(\widetilde{Y}, \mu) \lhd
(Z, \lambda) =
(\widetilde{Y} \prod\limits_Y Z, \mu \lhd \lambda)}
\end{equation}
\[\mu \lhd \lambda (y,z) = \mu(y) \cdot \lambda(z)
\]
For a commutative rig $B$, the canonical homomorphism
(\ref{Psi})
$\Psi_B \in \cG\cR ({\cal N}^{B!}, \cG(B))$ is given in this description as

\begin{equation}
\displaystyle{\left(\Psi_B (\widetilde{X}, \mu )\right)_x = \sum\limits_{\stackrel{\textstyle{\tilde{x} \in \widetilde{X} }}{\pi(\tilde{x})=x}}\mu(\tilde{x})}
\end{equation}

To get such a surjective homomorphism we can use any multiplicative submonoid
$B_0 \subseteq B^!$
such that $\mathbb{N}\{B_0\}=B$.
For example, for $B=\mathbb{Z}$ the integers,
we can take $B_0 = \{0, \pm 1 \}$,
and we get a surjective homomorphism

\begin{equation}
\displaystyle{\Psi \in \cG\cR(\cG(\mathbb{N}) \bigotimes\limits_{\mathbb{F}} \mathbb{F} \{\pm 1\} , \cG(\mathbb{Z}))}
\end{equation}

\addtocounter{subsection}{2}
\subsection*{13.1.3 \; generators and relations for $\cG(B)$, $B$ commutative ring.}

We have a surjective homomorphism
\begin{equation}
\begin{split}
\Phi:\Delta^{[2]}\sur \cG(\NN) \\
\delta=\delta^{[2]}\mapsto (1,1)\in \cG(\NN)_{[2]} \\
\Phi_X(F_1,\{\bar{F}_x\},\sigma)/\approx=(\hash\partial \bar{F}_x)_{x\in X}
\end{split}
\end{equation}

\thm{13.1.4}
The equivalence- ideal $\KER(\Phi)=\Delta^{[2]}\underset{\cG(\NN)}{\prod} \Delta^{[2]}$ is generated by
\begin{equation}\label{rel13.1.14}
\begin{matrix}
 \textbf{Zero}: &  \delta\lhd 1_1=1, \hspace{10mm} 1_1\in \FF_{[2],[1]},1_1(1)=1\in [2] \\
 \textbf{Ass:} &  \delta\underset{1}{\lhd} \delta= \delta\underset{2}{\lhd} \delta\\
\textbf{Comm:} & \delta\lhd 1_{\tiny \left(\begin{matrix}0 &1 \\ 1 & 0 \end{matrix}\right)}=\delta \\
\end{matrix}
\end{equation}

We get a surjective homomorphism
\begin{equation}
\Phi_{\tiny \ZZ}: \FF\{\pm 1\}\underset{\FF}{\otimes} \Delta^{[2]}\sur \cG(\ZZ)
\end{equation}
\thm{13.1.5}
The equivalence ideal $\KER(\Phi_{\ZZ})$ is generated by the relations (\ref{rel13.1.14}), and the relation
\begin{equation}
\textbf{Cancelation:}\hspace{16mm} (\delta\lhd (-1)^i_{i=1,2})\sslash \delta =0 \hfar\hfar
\end{equation}
For a commutative ring $B$ we get a surjective homomorphism
\begin{equation}
\Phi_B: \FF\{B^!\} \underset{\FF}{\otimes} \Delta^{[2]}\sur \cG(B)
\end{equation}
\thm{13.1.6} The equivalence- ideal $\KER(\Psi_B)$ is generated by the relations (\ref{rel13.1.14}), and the relations for $b_1,b_2\in B$

\begin{equation}
\textbf{($b_1$,$b_2$):}\hspace{26mm}(\delta\lhd (b_i)_{i=1,2})\sslash \delta=(b_1+b_2)\in \FF\{B^!\}_{[1]}\equiv B. \hfar\hfar
\end{equation}
The proofs of theorems 13.1.4-6 are the same as the proof given in Theorem 2.10,1. \vspace{3mm}\\

Every element $G\in (\FF\{B^!\}\underset{\FF}{\otimes}
\Delta^{[2]})_X$ can be represented (after moving the elements of
$\FF\{B^!\}$ to the boundary using (\ref{eq_up1})) as $G=(G_1
;\{\bar{G}_x\}_{x\in X}; \sigma ; \mu )$ $G_1,\bar{G}_x$ are
$\{1,2\}=[2]$- labelled binary trees
\begin{equation}
\begin{split}
\sigma:\partial G_1\iso \underset{x\in X}{\partial \bar{G}_x}\\
\mu:\partial G_1\rightarrow B
\end{split}
\end{equation}
and
\begin{equation}
\Phi_B(G)=\left(\underset{z\in \partial\bar{G}_x}{\Sigma} \mu(\sigma^{-1}(z))\right)_{x\in X}
\end{equation}
Note that the associated graph $\tilde{G}\in Graph_{[1],X}$ from $X$ to $[1]$, obtained by going from $X$ up the trees $\bar{G}_x$, and than via $\sigma^{-1}$, down the tree $G_1$,
\begin{equation}
\tilde{G}_0\equiv G_1\underset{\partial G_1}{\coprod}(\underset{x\in X}{\coprod}\bar{G}_x)\equiv G_1\circ
(\underset{x\in X}{\coprod}\bar{G}_x)
\end{equation}
is already in the ("left")$\circ$("right") form of the proof of theorem 2.10.1.

\section{The arithmetical plane $\cG(\mathbb{N}) \bigotimes_{\small \FF}\cG(\mathbb{N})
$}

We next give a description of the arithmetical plane
$\cG(\mathbb{N}) \bigotimes\limits_{\mathbb{F}}\cG(\mathbb{N})$.

An \emph{oriented-tree} is a (rooted) tree $F$ together with a map
\begin{equation}
\varepsilon_F: F \setminus \partial F \rightarrow \{0,1 \}
\end{equation}
It is \emph{1-reduced} if $\nu(a) \neq 1$
for all $a \in F$. \\
If for some $a \in F$,
$S_F^{-1}(a) = \{a'\}$, we obtain by 1-reduction the tree
\begin{equation}
\label{eq13.2.1}
 1_a(F) = F \setminus \{ a \}
\end{equation}
with $S_{F'}(a') = S_F(a)$.

For every oriented tree $F$ there is a unique 1-reduced tree
$F_{1-\mbox{red}}$;
it is obtained from $F$ by a finite sequence of 1-reductions.

The oriented tree $F$ is \emph{$\lhd$- reduced} if for all
$a \in F \setminus (\partial F \coprod \{0_F\})$,
$\varepsilon(a) \neq \varepsilon(S(a))$. \\
If for some $a  \in F \setminus (\partial F \coprod \{0_F\})$,
$\varepsilon(a) = \varepsilon(S(a))$,
we obtain by \emph{$\lhd$- reduction} the tree
\begin{equation}
\label{eq13.2.2}
O_a(F) = F \setminus \{a\}
\end{equation}
with $S_{O_a(F)}(a') = S_F(a)$
if $S_F(a')=a$.

For every oriented tree $F$ there is a unique $\lhd$-reduced tree
$F_{\lhd-\mbox{red}}$;
it is obtained from $F$ by a finite sequence of $\lhd$-reductions.
For a $\lhd$-reduced oriented tree $F$, the orientation
$\varepsilon_F$ is
completely determined by its value at the root
$\varepsilon_F(0_F)$,
since $\varepsilon_F(x) \equiv \varepsilon_F(0_F)+ ht(x) (mod 2 )$.
Thus we view $\lhd$-reduced oriented trees $F$ as ordinary trees together with an orientation of the root
$\varepsilon_F = \varepsilon_F(0_F) \in \{0, 1 \}$.

Note that the operations of
 1-reduction and
$\lhd$-reduction do not alter the boundary of a tree.

We let $\approx$ denote the equivalence relation on oriented trees generated by
$1$-reductions and $\lhd$-reductions.
We let $[F]$ denote the equivalence class
of the oriented tree $F$.
Thus $[F]=[F']$ if and only if there exist
$F=F_0, F_1, \ldots , F_l = F'$,
such that for $j=1, \ldots , l$,
the pair $\{F_j , F_{j-1} \}$
is related by $1$-reduction,
or $\lhd$-reduction;
it follows that there is a canonical identification of the
boundaries:
$\partial F = \partial F'$.

\vspace{10pt}

For a finite set $X \in \Fbb$,
let $\Upsilon_X$ denote the collection of isomorphism classes of data
\begin{equation}
\Upsilon_X = \{F=([F_1]; [\bar{F}_x], x \in X ; \sigma_F) \}/\cong
\end{equation}
where $F_1$, $\bar{F}_x$
are oriented trees taken modulo $\approx$-equivalence,
and $\sigma_F$ is a bijection
$\sigma_F: \partial F_1 \xrightarrow{\sim} \coprod\limits_{x \in X}
\partial \bar{F}_x$
and the data is taken up to isomorphism and consistent-commutativity.
Thus explicitly, the data $F$ is equivalent to the data $F'$,
if and only if there exists
$F=F^0, F^1, \ldots , F^l=F'$
such that for $j=1, \ldots ,l$
the pair $\{F^j, F^{j-1} \} = \{ G,G' \}$
is related by either:

\vspace{10pt}

\noindent {\bf Isomorphism:}
have isomorphism
$\tau_1: G_1 \xrightarrow{\sim} G'_1,
\tau_x: \bar{G}_x \xrightarrow{\sim} \bar{G}'_x, x \in X$
such that
$\sigma_{G'} \circ \tau_1(b) = \tau_x \circ \sigma_G(b)$
for
$b \in \partial G_1$,
$\sigma_G(b) \in \partial \bar{G}_x$.

\vspace{10pt}

\noindent {\bf $1$-reduction:}
have $G' = 1_aG$,
for some
$a \in G_1 \coprod {\displaystyle \coprod\limits_{x \in X} } \bar{G}_x$
with $\nu(a) =1$,
cf. (\ref{eq13.2.1}).

\vspace{10pt}

\noindent {\bf $\lhd$-reduction:}
have $G' =O_a G$,
for some \\$a \in \left(G_1 \setminus (\partial G_1 \coprod \{ 0\} )
\right)\coprod \displaystyle{\coprod\limits_{x \in X}}
\bar{G}_x \setminus (\partial \bar{G}_x \coprod \{0 \})$
with $\varepsilon (a) = \varepsilon (S(a))$, cf. (\ref{eq13.2.2}).

\vspace{10pt}

\noindent {\bf Consistent-commutativity:} $\{G,G'\}$ of the form (\ref{eq8.3.48}) or (\ref{eq8.3.50}).\\
The operations of multiplication (\ref{eq8.3.44}),
and of contraction (\ref{eq8.3.46}),
induce well defined operations on equivalent classes of data,
and make $\Upsilon$ into a
commutative generalized ring.
It is straightforward to check that
\begin{equation}
\begin{array}{cl}
F \lhd G & \cong F \lhd (1_a G) \cong F \lhd (O_a G)\\
          &\cong (1_a F) \lhd G \cong (O_a F) \lhd G \\
F\sslash G & \cong (F \sslash 1_a G) \cong (F \sslash O_a G)\\
      &\cong (1_a F\sslash  G) \cong (O_a F\sslash G)
\end{array}
\end{equation}

whenever the operations $1_a$, $O_a$ are relevant, and that
$\Upsilon$ satisfies the axioms of a commutative generalized ring.

Note that for $\varepsilon =0, 1$,
we have the elements
\begin{equation}
\delta_X^{\varepsilon} = \left( [X \coprod \{ 0 \}];
      [0_x] , x \in X ; \sigma \right) \in \Upsilon_X
\end{equation}
where
$X \coprod \{ 0 \}$ is the oriented tree with
$\varepsilon (0) = \varepsilon$,
$S(x) =0$ for $x \in X$,
and $\sigma :X \xrightarrow{\sim}
    \coprod\limits_{x \in X} \{ 0_x \}$ is the natural bijection
$\sigma (x) = 0_x$.

For $f \in Set_{\bullet}(X,Y)$,
and $(\delta_{f}^{\varepsilon})^{(y)} =
      \delta_{f^{-1}(y)}^{\varepsilon} $,
$y \in Y$,
we have via $\lhd$-reduction
\begin{equation}
\label{eq13.2.7}
\delta_Y^{\varepsilon} \lhd \delta_f^{\varepsilon}
\cong
\delta_{D(f)}^{\varepsilon}
\end{equation}
we also have by $1$-reduction
\begin{equation}
\label{eq13.2.8}
\delta_{[1]}^{\varepsilon} =
([\{0\} \coprod \{1\}]; [0_1];\sigma)\cong
([0_1];[0_1];id)=
1 \in \Upsilon_{[1]}
\end{equation}
Thus we get homomorphisms,
$\Psi^{\varepsilon} \in \cG\cR(\cG(\mathbb{N}), \Upsilon)$
with
$\Psi^{\varepsilon}(\mathbbm{1}_X) = \delta_X^{\varepsilon}$.
It is clear that $\Upsilon$ is generated by the
$\delta_X^{\varepsilon}$, and the only relations they satisfy
are (\ref{eq13.2.7}), (\ref{eq13.2.8}), and consistent commutativity.
It follows that $\Upsilon$ is the sum of $\cG(\mathbb{N})$
with itself in the category of commutative generalized rings:
for any $A \in \cG\cR_C$,
\begin{equation}
\begin{array}{ccl}
\cG\cR_C(\cG(\mathbbm{N}), A) \times \cG\cR_C(\cG(\mathbb{N}), A)&
\stackrel{\sim}{\leftrightarrow}&
\cG\cR_C(\Upsilon, A)\\
(\varphi \circ \psi^0 , \varphi \circ \psi^1)&
\mapsfrom&
\varphi\\
\varphi^0, \varphi^1 & \mapsto &
\varphi^0 \otimes \varphi^1 (\delta_X^{\varepsilon}) :=
\varphi^{\epsilon}(\mathbbm{1}_X)
\end{array}
\end{equation}

\vspace{10pt}

The diagonal homomorphism
\begin{equation}
\nabla \in \cG\cR_C(\Upsilon, \cG(\mathbb{N}))
\end{equation}
is determined by
$\nabla_X (\delta_X^{\varepsilon}) = \mathbbm{1}_X$,
and is given explicitly by
\begin{equation}
\label{eq13.2.11}
\nabla_X ([F_1]; [\bar{F}_x]_{x \in X}; \sigma )  =
\left( \sharp \partial \overline{F}_{x}\right)_{x \in X}
\end{equation}
The homomorphism $\nabla$ is surjective,
but it is not injective.

\vspace{10pt}

For a  monoid $B$,
the tensor product
$\Upsilon \bigotimes\limits_{\mathbb{F}} \mathbb{F}\{B\}$
can be described as isomorphism classes of data
\begin{equation}
(\Upsilon \bigotimes\limits_{\mathbb{F}}\mathbb{F}\{B\})_X
:=
\{F=([F_1]; \{[\bar{F}_x]\}_{x \in X}; \sigma_F; \mu_F)\}/\cong
\end{equation}
Here the data $([F_1]; \{[\bar{F}_x]\}_{x \in X}; \sigma_F)$
is the data for $\Upsilon_X$,
and $\mu_F$
is a map\\ $\mu_F: \partial F_1 \rightarrow B$,
and isomorphisms are required to preserve the $B$-valued maps,
and the zero law holds in the form:
\begin{equation}
\mu_F(b) = 0 , \ \underline{\sigma_F}(b) = x_0 \Rightarrow
F \approx ([F_1 \setminus \{b\}]; \{[\bar{F}_x ]\}_{x \neq x_0}
\cup \{[\bar{F}_{x_0}] \setminus \{\sigma _F(b)\})
\end{equation}
The operations of multiplication and contraction are the
given ones on the $\Upsilon$-part of the data
(i.e. given by (\ref{eq8.3.44}) and (\ref{eq8.3.46})),
and are given on the $B$-valued maps by
(using the notations of (\ref{eq8.3.45}) and (\ref{eq8.3.47})):
\begin{equation}
\begin{array}{ll}
\mu_{G \lhd F}(b,a) = \mu_G(b) \cdot
\mu_{F_{\underline{\tau}(b)}}(a)&
, b \in \partial G_1 , a \in \partial F_{\underline{\tau}(b)}\\
\mu_{G\sslash F}(b,a) = \mu_G(b) \cdot
\mu_{F_{f \circ \underline{\tau}(b)}}(\sigma^{-1}a)&
, b \in \partial G_1 , a \in \partial \bar{F}_{\underline{\tau}(b)}
\end{array}
\end{equation}

For commutative rings $B_0$, $B_1$,
taking $B = B^!_0 \otimes B_1^!$
(the sum in $Mon$, cf. (\ref{eq13.1.4})),
we get the generalized ring
\begin{equation}
\Upsilon \bigotimes\limits_{\mathbb{F}} \mathbb{F}\{B_0^! \otimes
B_1^!\}= \cG(\mathbb{N}) \bigotimes\limits_{\mathbb{F}}
\mathbb{F}\{B_0^!\} \bigotimes\limits_{\mathbb{F}}
\cG(\mathbb{N}) \bigotimes\limits_{\mathbb{F}}
\mathbb{F}\{B_1^!\}
\end{equation}
which maps surjectively onto
$\cG(B_0)\bigotimes\limits_{\mathbb{F}} \cG(B_1)$.

For the integers $\mathbb{Z}$, taking $B=\{0, \pm 1\}$,
we get the generalized ring
\begin{equation}
\Upsilon \bigotimes\limits_{\mathbb{F}} \mathbb{F}\{\pm 1\} =
(\cG(\mathbb{N}) \bigotimes\limits_{\mathbb{F}} \mathbb{F}\{\pm 1\})
\bigotimes\limits_{\mathbb{F}\{\pm 1\}}
(\cG(\mathbb{N}) \otimes \mathbb{F}\{\pm 1\})
\end{equation}
with a surjective homomorphism
\begin{equation}
\pi:\Upsilon \bigotimes\limits_{\mathbb{F}} \mathbb{F} \{\pm 1\}
\twoheadrightarrow
\cG(\mathbb Z) \bigotimes\limits_{\mathbb{F} \{\pm 1\}} \cG(\mathbb Z).
\end{equation}
Note that $KER(\pi)=E(\mathfrak{a})$ is the equivalence ideal
generated by the homogeneous ideal $\mathfrak{a}$ generated by the
elements giving the "cancellations "on the left and right
$\cG(\ZZ)$'s:
\begin{equation}
x_{\epsilon}=(F_{\epsilon}; F_{\epsilon} ; \sigma=id \; ;\mu=id)\in (\Upsilon\otimes_{\FF}\FF\{\pm 1\})_{[1]},\;\;\; \epsilon=0,1
\end{equation}
with the reduced oriented tree $F_{\epsilon}=\{0\}\amalg\{\pm 1\}, \mathcal{S}(\pm 1)=0,\epsilon(0)=\epsilon $.

\section{Products of Grothendieck-Generalized-schemes}
\label{sec13.3}

The category $\cal{GGS}$
has fibred products:

Given maps $f^j \in \cG\cG\cS(X^j, Y)$,
there exists $X^0 \Pi_YX^1 \in \cG\cG\cS$,
and maps $\pi_j \in \cG\cG\cS(X^0 \Pi_YX^1, X^j)$,
with $f^0 \circ \pi_0 = f^1\circ\pi_1$,
and for any $g^j \in \cG\cG\cS(Z,X^j)$,
with $f^0 \circ g^0 = f^1\circ g^1$,
there exists a unique map
$g^0 \pi g^1 \in \cG\cG\cS(Z,X^0\Pi_YX^1)$,
such that $\pi_j \circ (g^0 \pi g^1) = g^j$, $j=0,1$.

Writing $Y=\bigcup\limits_i spec^t (A_i)$,
$(f^j)^{-1}(spec^t(A_i))= \bigcup\limits_k spec^t (B^j_{i,k})$,
the fibred product $X^0 \Pi_Y X^1$
is obtained by gluing
$spec^t (B^0_{i,k_0} \bigotimes\limits_{A_i} B^1_{i,k_1})$.
See the construction of fibred product of ordinary schemes,
e.g. \cite[Theorem 3.3, p. 87]{Hart}.

\section{Products of Generalized-schemes}
\label{sec13.4}

The category $\cal{GS}$
has fibred products.
This is an immediate corollary of
\ref{sec13.3}. Given maps
\begin{equation*}
\varphi = \{\varphi_i^j\}_{j \geq \sigma(i)} \in
\cG\cS ( \{X_j\}_{j \in J}, \{Y_i\}_{i \in I}),
\end{equation*}
and
\begin{equation}
\varphi' = \{\varphi_i^{j'}\}_{j' \geq \sigma'(i)} \in
\cG\cS ( \{X'_{j'}\}_{j' \in J'}, \{Y_i\}_{i \in I})
\end{equation}
the fibred product of $\varphi$
and $\varphi'$ in $\cG\cS$ is given by the inverse system
$\{X_j \Pi_{Y_i} X'_{j'}\}$,
the indexing set is
\begin{equation}
\{(j,j',i) \in J \times J' \times I\  | \  j \geq \sigma(i) , \ j' \geq \sigma'(i) \}
\end{equation}

\section{The Arithmetical plane: $\mathbb{X} = \overline{spec \mathbb{Z}}
\prod\limits_{\mathbb{F}\{\pm 1\}}
\overline{spec \mathbb{Z}} $}

This is a special case of \S \ref{sec13.4}:
The (compactified) arithmetical plane $\mathbb{X}$
is given by the inverse system
$\{ \cX_N \prod\limits_{spec \mathbb{F}\{\pm 1\}} \cX_M \}$,
with indexing set\\
$\{ (N,M) \in \mathbb{N} \times \mathbb{M} |
N, M \ \mbox{square-free} \}$
and with
\begin{equation}
\cX_N = spec \, \cG(\mathbb{Z})
\coprod\limits_{spec \cG( \mathbb{Z} [\frac{1}{N}])}
spec (\cG( \mathbb{Z} [\frac{1}{N}]) \cap \cO_{\eta})
\end{equation}
as in (\ref{spec_z}).
This generalized scheme $\mathbb{X}$ contains the
affine open dense subset,
\begin{equation}
spec \, \cG(\mathbb{Z}) \prod\limits_{\mathbb{F}\{\pm 1\}}
spec \cG(\mathbb{Z}) =
spec^t (\cG(\mathbb{Z}) \bigotimes\limits_{\mathbb{F} \{\pm 1\}}
\cG(\mathbb{Z}))
\end{equation}
e.g.  basis for neighborhoods of $(p, \eta)$
is given by
\begin{equation}
spec^t  \ \left[\cG(\mathbb{Z}[\frac{1}{N}])
\bigotimes\limits_{\mathbb{F} \{\pm 1\}}
(\cG(\mathbb{Z}[\frac{1}{M}])\cap \cO_{\eta})\right]
\end{equation}
where $p$ does not divide $N$, and $M$ is arbitrary.

\vspace{10pt}

Similarly, for any number field $K$ we
have the compactified surface
\begin{equation}
\overline{spec \cO_K} \prod\limits_{spec \mathbb{F}\{\mu_K\}}
\overline{spec \cO_K}
\end{equation}
It contains the affine open dense subset
$spec^t (\cG(\cO_K) \bigotimes\limits_{\mathbb{F} \{\mu_K\}}
\cG(\cO_K))
$.

\setcounter{equation}{0}
\chapter{Modules and differentials}
\section{\Amd}\label{14.1}
\defin{14.1.1}
Let $A\in \cG\cR$. An $\Amd$ is a functor $M\in (Ab)^{\FF}$, with $M_{[0]}=\{0\}$, and operations:
\begin{equation}
\begin{split}
\textbf{multiplication:} \text{ for } z_0\in Z\in \FF, \\
\underline{\;\;}\underset{z_0}{\lhd}\underline{\;\;} : M_Z\times A_X\rightarrow M_{Z\underset{z_0}{\lhd} X}
\end{split}
\end{equation}
\begin{equation}
\begin{split}
\textbf{contraction:} \text{ for } X\subseteq Y\in \FF, \\
\underline{\;\;}\sslash \underline{\;\;} :  M_Y\times A_X\rightarrow M_{Y\slash X}
\end{split}
\end{equation}
These are assumed to satisfy the Disjointness Axioms $I,II,III,$ and so for $f\in \Set(Y,Z)$, (with $A_f=\underset{z\in Z}{\prod} A_{f^{-1}(z)}$), we have "multiple"
\begin{equation}
\begin{split}
&\textbf{multiplication:} \\
&\underline{\;\;}\lhd\underline{\;\;} : M_Z\times A_f\rightarrow M_{Y}
\end{split}
\end{equation}
\begin{equation}
\begin{split}
&\textbf{contraction:} \\
&\underline{\;\;}\sslash \underline{\;\;} :  M_Y\times A_f\rightarrow M_Z
\end{split}
\end{equation}
We further assume these operations satisfy the following axioms:\vspace{3mm}\\
\textbf{Unit and Functoriality:}
\begin{equation} \text{ for } f\in \FF_{Z,Y}, m\in M_Y, m\sslash 1_f=m\lhd 1_{f^t}=f_M(m)
\end{equation}
\textbf{Homomorphism:}
\begin{equation}
\begin{split}
(m_1+m_2)\lhd a=(m_1\lhd a)+(m_2\lhd a) \\
 (m_1+m_2)\sslash a=(m_1\sslash a)+(m_2\sslash a)
\end{split}
\end{equation}
\textbf{Associativity:}
\begin{equation}
(m\lhd a_0)\lhd a_1=m\lhd(a_0\lhd a_1)
\end{equation}
\textbf{Left Adjunction:}
\begin{equation}
(m\sslash a_0)\sslash a_1=m\sslash (a_1\lhd a_0)
\end{equation}
\textbf{Right Adjunction:}
\begin{equation}
m\sslash (b\sslash a) = (m\lhd a) \sslash b
\end{equation} \vspace{3mm}\\
\textbf{Left Linearity:}
\begin{equation}
m\lhd (b\sslash a) = (m\lhd b)\sslash a
\end{equation}
We are only interested in the \emph{commutative-} $\Amd$s that further satisfy  \vspace{3mm} \\
\textbf{Right Linearity:}
\begin{equation}
(m\sslash a)\lhd b=(m\lhd f^*b)\sslash g^*a\;\;\text{ for }\;\;\; a\in A_f,b\in A_g .
\end{equation}
Note that for $M$ an $\FF$- module, i.e. a functor $M\in (Ab)^{\FF}$ with $M_{[0]}=\{0\}$, and for $X\subseteq Y\in \FF$, $M_X$ is a subgroup of $M_Y$, a direct summand, and the projection $M_Y\sur M_X$ is denoted by $m\mapsto m|_X$. \\
In particular we have the "matrix - coefficient" map
\begin{equation}
\begin{split}
&J_X:M_X\rightarrow (M_{[1]})^X \\
&J_X(m)=(m|_{\{x\}})_{x\in X} \hspace{10mm} (\text{identifying }\;\;\; M_{[1]}=M_{\{x\}}\subseteq M_X ).
\end{split}
\end{equation}
If these maps are injective for all $X\in \FF$ we say $M$ is a "matrix"- $\Amd$.
\defin{14.1.2} A homomorphism of $\Amd$s $\varphi:M\rightarrow M'$ is a natural transformation of functors that commutes with the $A$- action. \\
Thus we have an abelian category $\Amod$. It is complete and co- complete: all (co-) limits can be taken pointwise. It has enough projectives and injectives: the evaluation functor at $X\in \FF$, $i^X:\Amod\rightarrow Ab$, $i^XM=M_X$, has a left (resp. right) adjoint $i_!^X$ (resp.$i^X_*$). \\
In particular, we have the \vspace{5mm} \\
\textbf{\large 14.1.3 Free $\Amod$ of degree $X\in \FF$:}
\begin{equation}
\begin{split}
&A^X=i_i^X\ZZ, \hspace{10mm} \text{generated by } \delta_X\in A_X^X, \text{ and} \\
&\Amod(A^X,M)\equiv M_X \text{ via } \varphi\mapsto \varphi_X(\delta_X).
\end{split}
\end{equation}
The elements of degree $Z\in \FF$, $m\in A^X_Z$, are linear combinations
\begin{equation}
m=\sum_{i=1}^k m_i\cdot (\delta_X\lhd a_i)\sslash b_i \;\;\;\; \text{ with } m_i\in \ZZ, a_i\in A_{f_i:X_i\rightarrow X}, b_i\in A_{g_i:X_i\rightarrow Z}
\end{equation}
and for $f\in \Set (Y,Z)$, (resp. $f\in \Set (Z,Y)$), and $a\in A_f$, the action of $a$ on such an element $m$ is given by
\begin{equation}
\begin{split}
&m\lhd a= \sum_{i=1}^k m_i (\delta_X\lhd (a_i\lhd g_i^* a))\sslash f^* b_i \\
\text{ resp. } &m\sslash a =\sum_{i=1}^k m_i (\delta_X\lhd a_i)\sslash (a\lhd b_i).
\end{split}
\end{equation}
The elements $(\delta_X\lhd a)\sslash b$ are subjected to the axioms of an $\Amod$, and we have
\begin{equation}
\begin{matrix}
(\delta_X\lhd (a\sslash c))\sslash b = (\delta_X\lhd a)\sslash (b\lhd c)&&\text{ for } \;\; \xymatrix@C-1pc@R-1pc{Y\ar[rd]^c\ar[dd]^a&& \\&W\ar[dl]\ar[dr]^b& \\ X&&Z } \\
\text{and }\; (\delta_X\lhd (a\lhd c))\sslash b = (\delta_X\lhd a)\sslash (b\sslash c) &&\text{ for } \;\; \xymatrix@C-1pc@R-1pc{Y\ar[d]_c\ar[ddr]^b& \\ W\ar[d]_a\ar[rd]& \\X &Z }
\end{matrix}
\end{equation}

\exmpl{14.1.4} For $A=\cG(R)$, $R$ a commutative rig (so every $a=(a_x)\in A_X$, can be written as $a=\mathbbm{1}_X\lhd <a_x>$, with $<a_x>\in (A_{[1]})^X=A_{id_X}$ ), we have
\begin{equation}
(A^X)_Z:= \ZZ\cdot A_{X\otimes Z}\slash \ZZ\cdot  0_{X\otimes Z} = \ZZ \cdot (R^{X\otimes Z})\slash \ZZ\cdot 0_{X\otimes Z},
\end{equation}
the free abelian group on $X$ by $Z$ matrices over $R$ (modulo $\ZZ\cdot 0_{X\otimes Y}$, because $A_{[0]}^X=\{0\}$).
For $a=(a_{x,z})\in R^{X\otimes Z}$, we have the generators
\begin{equation}
\delta_X\cdot (a):= (\delta_X\underset{x\in X}{\lhd} (a_{x,z}))\sslash (\mathbbm{1}_X)_{z\in Z}=(\delta_X\underset{x\in X}{\lhd} \mathbbm{1}_Z) \sslash (a_{x,z})
\end{equation}
and the $\cG(R)$- action on these generators is the diagonal action at each $x\in X$. \vspace{3mm}\\
\textbf{\large Localization 14.1.5} \vspace{3mm}\\
For a multiplicative set $S\subseteq A_{[1]}^+$, and $M\in \Amod$, we have the localization $S^{-1}M\in S^{-1}\Amod$ described as in ($\ref{localization26}\text{-}\ref{localization27}$). In particular, we have the localizations $M_{\mathfrak{p}}=S^{-1}_{\mathfrak{p}}M\in A_{\mathfrak{p}}\text{-mod}, \mathfrak{p}\in \Spec^t(A)$, and $M_s\in A_s\text{-mod}, s=s^t\in A_{[1]}^+$.
\defin{14.1.6} For $(X,\mathcal{O}_X)\in \cG\cR\cS$, an $\mathcal{O}_X\textit{-module}$ is a functor (where $C_X$ is the category of open subsets of $X$ and inclusions),
\begin{equation}
\begin{split}
M:&C_X^{op}\times \FF\rightarrow Ab \\
&U,Z\mapsto M(U)_Z
\end{split}
\end{equation}
such that for $U\subseteq X$ open, $M(U)=\{M(U)_Z\}\in \mathcal{O}_X(U)\text{-mod}$; for $U\subseteq U'\subseteq X$ open, the homomorphisms $M(U')_Z\rightarrow M(U)_Z,m\mapsto m|_U$ is compatible with the operations: $(m\lhd a)|_u=m|_u\lhd a|_u, (m\sslash a)|_u=m|_u\sslash a|_u$ ; and for fixed $Z\in \FF:\;  U\mapsto M(U)_Z$ is a sheaf. \\
A homomorphism of $\mathcal{O}_X\text{-modules}$ $\varphi:M\rightarrow M'$, is a natural transformation of functors $\varphi(U)_Z:M(U)_Z\rightarrow M'(U)_Z$, such that for $U\subseteq X$ open, $\{\varphi(U)_Z\}_{Z\in \FF}\in \mathcal{O}_X(U)\text{-mod}\bigg (M(U),M'(U) \bigg )$, and for $Z\in \FF, \{\varphi(U)_Z\}_{U\subseteq X}\in (Ab)^{C_X^{op}}(M_Z,M'_Z)$. \\
Thus we have an abelian category: $\mathcal{O}_X\text{-mod}$.  \\
For $A\in \cG\cR_C, M\in \Amod$, we have $\tilde{M}\in \mathcal{O}_A\text{-mod}$, defined as in Definition $7.3.2$. \\
It has stalks at $\mathfrak{p}\in \Spec^t A$ given by $(\tilde{M})_{\mathfrak{p}}=M_{\mathfrak{p}}$, cf. ($\ref{spectracondition}$), and it has global sections over a basic open set $D_s^+\subseteq \Spec^t A$ given by $\tilde{M}(D_s^+)=M_s, s=s^t\in A_{[1]}^+$ cf. Theorem $11.3.2$ (The commutativity of $M$ is essential !) .\\
For $X\in \cG\cG\cS$, a Grothendieck- generalized- scheme, we have the full subcategory of "quasi-coherent" $\mathcal{O}_X\text{-modules}$,  $q.c. \mathcal{O}_X\text{-mod}\subseteq  \mathcal{O}_X\text{-mod}$. Its objects are the $\mathcal{O}_X\text{-modules}$ satisfying the equivalent conditions of theorem $7.3.3$, and for $X=\Spec^t A$ affine, localization gives an equivalence
\begin{equation}
\begin{split}
\Amod &\iso q.c.\mathcal{O}_A\text{-mod}\subseteq \mathcal{O}_A\text{-mod} \\
M &\mapsto \tilde{M}
\end{split}
\end{equation}
\textbf{\large Restriction and extension of scalars 14.1.7} \vspace{3mm}\\
For $\varphi\in \cG\cR_C(B,A)$, we have the adjoint functors (using geometric notations):
\begin{equation}
\xymatrix{\Amod \ar@<1ex>@{->}[r]^{\varphi_*}& B\text{-}mod\ar@<1ex>@{->}[l]^{\varphi^*} }
\end{equation}
The right adjoint takes $N\in \Amod$ to $\varphi_* N=N$ with $B$-action given using $\varphi: \;\; n\lhd b=n\lhd \varphi(b), \; n\sslash b=n\sslash \varphi(b)$. \\
The left adjoint takes $M\in B\text{-mod}$ into the $A\text{-module}$ $\varphi^*M=M^A$, whose elements in degree $Z\in \FF, m\in (M^A)_Z$, can be written as sums
\begin{equation}
m=\sum_{i=1}^k ([m_i]\lhd a_i)\sslash a_i' \text{ with } m_i\in M_{X_i}, a_i\in A_{f_i:Y_i\rightarrow X_i}, \; a_i'\in A_{g_i:Y_i\rightarrow Z}
\end{equation}
and $a\in A_f, f\in \Set(Y,Z)$, (resp. $f\in \Set(Z,Y)$ ), acting via
\begin{equation}
\begin{split}
&m\lhd a = \sum_{i=1}^k ([m_i]\lhd (a_i \lhd g_i^*a))\sslash f^* a_i'  \\
\text{resp. } & m\sslash a=\sum_{i=1}^k ([m_i]\lhd a_i)\sslash (a\lhd a_i')
\end{split}
\end{equation}
and we have the relations:
\begin{equation}
\begin{split}
([m+m']\lhd a)\sslash a' &=([m]\lhd a)\sslash a'+([m']\lhd a)\sslash a' \\
([m\lhd b]\lhd a)\sslash a' &=([m]\lhd (\varphi(b)\lhd a))\sslash a' \\
([m\sslash b]\lhd a)\sslash a'&=([m]\lhd h^*(a))\sslash (a'\lhd f^*\varphi (b)), \;\; b\in B_h, a\in A_f .
\end{split}
\end{equation}
\section{Derivations and differentials }
\defin{14.2.1} For $A\in \cG\cR_C, M\in \Amod$, we define the even and odd infinitisimal extensions $A\prod^{\pm} M$, an abelian group object of $\cG\cR\slash A$, by
\begin{equation}
\FF\ni X\mapsto (A\prod M)_X=A_X\prod M_X
\end{equation}
and for $f\in \Set(Y,Z), a=(a^{(z)})\in A_f, m=(m^{(z)})\in M_f:= \underset{z\in Z}{\prod} M_{f^{-1}(z)}$, we have \vspace{3mm} \\
\textbf{multiplication:} for $a_Z\in A_Z,m_Z\in M_Z $,
\begin{equation}
(a_Z,m_Z)\lhd (a,m):= (a_Z\lhd a, m_Z\lhd a+\sum_{z\in Z} m^{(z)}\underset{y\in f^{-1}(z)}{\lhd} (a_Z|_z))
\end{equation}
\textbf{$\pm$-contraction:} for $a_Y\in A_Y ,m_Y\in M_Y $,
\begin{equation}
(a_Y,m_Y)\sslash (a,m):= (a_Y\sslash a, m_Y\sslash a\pm\sum_{z\in Z} m^{(z)}\sslash (a_Y|_{f^{-1}(z)}))
\end{equation}
We have,
\begin{equation}
\begin{matrix}
\emph{projection:}& \pi\in \cG\cR (A\prod M,A),\;\;  \pi (a,m)=a, \\
\emph{addition:} & \mu\in \cG\cR\slash A \bigg ( (A\prod M) \underset{A}{\prod} (A\prod M), A\prod M    \bigg ), \mu((a,m),(a,m'))=(a,m+m') \\
\emph{unit:} & \epsilon \in \cG\cR\slash A(A,A\prod M), \epsilon(a)=(a,0),  \\
\emph{antipode:} & S\in \cG\cR\slash A(A\prod M,A\prod M), S(a,m)=(a,-m).
\end{matrix}
\end{equation}
(Note: When $A\in \cG\cR_C$, $M$ a (commutative) $\Amd$, the generalized rings $A\prod^{\pm} M$ need \emph{not} be commutative).
\defin{14.2.2} For $\varphi \in \cG\cR(C,A), M\in \Amod$, an even/odd $C$ -\emph{linear derivations from $A$ to $M$}, is a collection of maps $\delta=\{\delta_X:A_X\rightarrow M_X\}_{X\in \FF}$ satisfying:
\begin{equation}
\begin{matrix}
(\star)\; \emph{Leibnitz: } \text{ for } f\in \Set(Y,Z),a_Z\in A_Z, a_Y\in A_Y, a_f\in A_f, \\
\delta(a_Z\lhd a^f)= \delta(a_Z)\lhd a^f+\underset{z\in Z}{\sum}\delta(a^f_z)\underset{y\in f^{-1}z}{\lhd}(a_Z|_z) \\
\delta(a_Y\sslash a^f)= \delta(a_Y)\sslash a^f\pm \underset{z\in Z}{\sum}\delta(a^f_Z)\sslash (a_Y|_{f^{-1}(z)}) \\
(\star\star)\; C\emph{-linear: }\text{ for }c\in C, \; \delta{\varphi(c)}\equiv 0.  \hspace{40mm}
\end{matrix}
\end{equation}
Note that, since $\FF\subseteq C$, we get $\delta(a\lhd 1_{f^t})=\delta(a)\lhd 1_{f^t}$, and $\delta$ is always a natural transformation. \\
We denote by $\Der^{\pm}_C(A,M)$ the collection of even/odd $C$ -linear derivations $\delta:A\rightarrow M$.
These are functors: $\Amod\rightarrow Ab, M\mapsto \Der^{\pm}_C(A,M)$, represented by \vspace{3mm}\\
\textbf{\large The module of even/odd differentials $\Omega^{\pm}(A\slash C)\in \Amod:$ 14.2.3} \vspace{3mm}\\
\begin{equation}
\begin{split}
\Der^{\pm}_C(A,M) & \equiv \Amod(\Omega^{\pm}(A\slash C),M) \\
\varphi\circ d^{\pm}_{A\slash C} &\mapsfrom \varphi
\end{split}
\end{equation}
with the universal even/odd derivation $d^{\pm}_{A\slash C}:A\rightarrow \Omega^{\pm}(A\slash C)$. \\
For $Z\in \FF$, the elements of $\Omega^{\pm}(A\slash C)_Z$ are sums of the form
\begin{equation}
\begin{split}
&\sum_{i=1}^k m_i\cdot (d^{\pm}(a_i)\lhd a_i')\sslash a_i'' \\
&m_i\in \ZZ, a_i\in A_{W_i}, a_i'\in A_{f_i:Y_i\rightarrow W_i}, a_i''\in A_{g_i:Y_i\rightarrow Z}
\end{split}
\end{equation}
subjected to the $\pm$ Leibnitz and $C$ -linearity relations. \vspace{3mm}\\
\textbf{\large 14.2.4 \; Example} \vspace{3mm}\\
For $A=C[\delta_W]=C\underset{\FF}{\otimes} \Delta^W$, the free \emph{commutative} generalized ring over $C$ generated by $\delta_W\in A_W$, we have
\begin{equation}
\Omega^{\pm}(C[\delta_W]\slash C)=\text{ free } C[\delta_W] \text{-module generated by }d^{\pm}(\delta_W) \text{ in degree } W,
\end{equation}
modulo "$\pm$ almost - linearity" - the derived commutativity relation (i.e. the $C[\delta_W]$ -sub-module generated by the difference of $d^{\pm}$ applied to $(a\sslash \delta_W)\lhd \delta_W$ and to $(a \underset{W}{\lhd}\delta_W)\sslash_W (\delta_W), a\in C[\delta_W]_Y, W\subseteq Y$). \vspace{5mm}\\
For $\varphi \in \cG\cR_C(C,A), B\in C\setminus \cG\cR\slash A$, i.e. $\varphi=\pi\circ \epsilon, \epsilon\in \cG\cR(C,B), \pi\in \cG\cR(B,A)$, and for $M\in \Amod$ we have the natural identifications
\begin{equation}
\begin{split}
 C\setminus \cG\cR\slash A(B,A\Pi^{\pm} M) & \equiv \Der^{\pm}_C(B,M)\equiv B\text{-mod}(\Omega^{\pm}(B\slash C),M)\equiv \Amod(\Omega^{\pm}(B\slash C)^A,M) \\
\pi\prod \delta & \mapsfrom \delta, \hspace{10mm} \varphi\circ d_{B\slash C}\mapsfrom \varphi
\end{split}
\end{equation}
We obtain, \vspace{3mm}\\
\textbf{\large  The adjunctions\; 14.2.5} \vspace{3mm}\\
\begin{equation}\label{Der_adj}
\xymatrix{\Omega^{\pm}(\slfrac{B}{C})^A&\Amod \ar@/^/@{|->}[d]& M \ar@{|->}[d] \\
 B\ar@{|->}[u]&C \diagdown \cG\cR \diagup A\ar@/^/@{|->}[u]& A\Pi^{\pm} M}
\end{equation}
Thus the even and odd differentials satisfy the ($\cG\cR$ -analogues of) all the properties (0 to 5) of \S 7.7. \vspace{3mm} \\
\textbf{\large 14.2.6 \; Example} \vspace{3mm}\\
In particular we have the $\cG(\NN)$, (resp. $\cG(\ZZ)$) -modules $\Omega^{\NN}_{\pm}=\Omega^{\pm}(\cG(\NN)\slash \FF)$, (resp. $\Omega_{\pm}^{\ZZ}=\Omega^{\pm}(\cG(\ZZ)\slash \FF\{\pm 1\})$), with the universal even/odd derivation in degree $X\in \FF$
\begin{equation}
\begin{split}
&d^{\pm}_X:\cG(\NN)_X=\NN^X\rightarrow (\Omega^{\NN}_{\pm})_X \\
\text{resp. }\;\;\; &d^{\pm}_X:\cG(\ZZ)_X=\ZZ^X\rightarrow (\Omega^{\ZZ}_{\pm})_X
\end{split}
\end{equation}
The module $\Omega^{\NN}_{\pm}$ (resp. $\Omega^{\ZZ}_{\pm}$) is obtained from the free $\cG(\NN)$ (resp. $\cG(\ZZ)$) module of degree $[2]$, with generator $d^{\pm}(1,1)$, modulo the derived relations (\ref{almostlinear}- \ref{GR_relations}). Thus $\Omega_X^{\NN}$ (resp. $\Omega_X^{\ZZ}$) is the free abelian group with generators $\left\{\begin{matrix} (a_x) \\ (a_x') \end{matrix}\right\}, (a_x),(a_x')\in \NN^X$ (resp. $\ZZ^X$),
modulo the relations, $a,a',a''\in \NN^X$ (resp. $\ZZ^X$): \\
\begin{equation}\label{almostlinear}
\begin{matrix}
\emph{$\pm$ almost linearity:} & \left\{\begin{matrix} x\cdot a \\ x\cdot a' \end{matrix}\right\}+\left\{\begin{matrix} y\cdot a\\ y\cdot a' \end{matrix}\right\}\pm \left\{\begin{matrix} x\cdot a \\ y\cdot a \end{matrix}\right\}\pm \left\{\begin{matrix} x\cdot a' \\ y\cdot a' \end{matrix}\right\}= \\
&\left\{\begin{matrix} (x+y)\cdot a\\ (x+y)\cdot a' \end{matrix}\right\}\pm \left\{\begin{matrix} x\cdot (a+a') \\ y\cdot (a+a') \end{matrix}\right\}, \; (x,y\in \NN \; \text{resp. }\ZZ).
\end{matrix}\tag{14.2.12$\pm$}
\end{equation}
(apply $d^{\pm}$ to the identity
\begin{equation*}
((x,y)\sslash (1,1))\lhd (1,1)= ((x,y)\underset{\{1,2\}}{\lhd}(1,1))\underset{\{1,2\}}{\sslash}(1,1)
\end{equation*}
or schematically
\begin{equation}
\begin{tikzpicture}[x=0.8cm, y=0.4cm]
           \vertex(lu) at (-1,1) [label=above:$x$]{};
     \vertex (ld) at (-1,-1) [label=above:$y$] {};
           \vertex (lm) at (0,0) {};
           \vertex (lru) at (1,1) [label=above:$a$] {};
    \vertex (lrd) at (1,-1) [label=above:$a'$]{};
    \node(m) at (2,0){$=$};
    \vertex (ru) at (3,1) [label=above:$x$] {};
    \vertex (rd) at (3,-1) [label=above:$y$]{};
           \vertex (rmu) at (4,1.5){};
           \vertex (rmmu) at (4,0.5) {};
    \vertex (rmmd) at (4,-0.5) {};
    \vertex (rmd) at (4,-1.5)  {};
    \vertex (rru) at (5,1) [label=above:$a$] {};
    \vertex (rrd) at (5,-1) [label=above:$a'$]{};

\path
        (lm) edge (lu)
        (lm) edge (ld)
        (lru) edge (lm)
        (lrd) edge (lm)

        (rmu) edge (ru)
        (rmmd) edge (ru)
        (rmmu) edge (rd)
        (rmd) edge (rd)

        (rru) edge (rmu)
        (rru) edge (rmmu)
        (rrd) edge (rmmd)
        (rrd) edge (rmd)
        ;

\end{tikzpicture}
\end{equation}
).
\setcounter{equation}{12}
\begin{equation}\label{GR_relations}
\begin{matrix}
\emph{cocycle:}\hspace{3mm}&\left\{\begin{matrix} a+a' \\ a'' \end{matrix}\right\}+\left\{\begin{matrix} a \\ a' \end{matrix}\right\}=\left\{\begin{matrix} a \\ a'+a'' \end{matrix}\right\}+\left\{\begin{matrix} a' \\ a'' \end{matrix}\right\}\\
(\scriptsize \text{apply } d^{\pm} \text{to the identity }& (1,1)\lhd\{1,(1,1)\}=(1,1)\lhd \{(1,1),1\} \\
\scriptsize \text{or schematically}& \\
 & \begin{tikzpicture}[x=0.8cm, y=0.4cm]
    \vertex (l1) at (-1,0.5) {};
    \vertex (l2) at (0,1) {};
    \vertex (l3) at (0,0) {};
    \vertex (l4) at (1,1.5) [label=right:$a$] {};
    \vertex (l5) at (1,0.5) [label=right:$a'$]{};
    \vertex (l6) at (1,-0.5) [label=right:$a''$]{};
    \node(m) at (2,0){$=$};
    \vertex (r1) at (3,0.5){};
    \vertex (r2) at (4,1){};
    \vertex (r3) at (4,0) {};
    \vertex (r4) at (5,1.5) [label=right:$a$] {};
    \vertex (r5) at (5,0.5) [label=right:$a'$]{};
    \vertex (r6) at (5,-0.5) [label=right:$a''$]{};

\path
        (l1) edge (l2)
        (l1) edge (l3)
        (l2) edge (l4)
        (l2) edge (l5)
        (l3) edge (l6)

        (r1) edge (r2)
        (r1) edge (r3)
        (r2) edge (r4)
        (r3) edge (r5)
        (r3) edge (r6)

        ;

\end{tikzpicture} \\
&  ).\\
\emph{normalized:} & \left\{\begin{matrix} a \\ 0 \end{matrix}\right\}=0=\left\{\begin{matrix} 0 \\ a' \end{matrix}\right\}\\
\scriptsize (\text{apply } d^{\pm} \text{ to the identity }& (1,1)\sslash (1,0)=1 )\\
\emph{symmetric:} & \left\{\begin{matrix} a' \\ a \end{matrix}\right\}=\left\{\begin{matrix} a\\ a' \end{matrix}\right\}\\
\scriptsize (\text{apply } d^{\pm} \text{ to the identity }& (1,1)\lhd 1_{\left( \tiny \begin{matrix}0& 1 \\ 1& 0  \end{matrix}\right)}=(1,1) )\\
\text{(resp. and cancellation:} &2\cdot \left\{\begin{matrix} a \\ -a\end{matrix}\right\}=0 \text{)} \\
( \scriptsize \text{apply } d^{\pm} \text{ to the identity }& ((1,1)\lhd \{+1,-1\})\sslash (1,1)=0)
\end{matrix}
\end{equation}
The $-1$- almost- linearity in (\ref{almostlinear}) is a consequence of the other relations (\ref{GR_relations}), (applying the cocycle relation to the term in the square brackets) we have:
\begin{equation}\label{calc}
\begin{split}
&\left[\left\{\begin{matrix} (x+y)a\\ (x+y)a' \end{matrix}\right\}\right]-\left\{\begin{matrix} x\cdot (a+a') \\ y\cdot (a+a') \end{matrix}\right\}\\
 &=-\left\{\begin{matrix} xa \\ ya\end{matrix}\right\}+\left[\left\{\begin{matrix} xa \\ xa'+y(a+a') \end{matrix}\right\}\right]+
\left\{\begin{matrix} (x+y)a' \\ ya \end{matrix}\right\}-\left\{\begin{matrix} x\cdot (a+a') \\ y\cdot (a+a') \end{matrix}\right\}\\
&=-\left\{\begin{matrix} xa\\ ya \end{matrix}\right\}- \left\{\begin{matrix} y(a+a') \\ xa' \end{matrix}\right\}+\left\{\begin{matrix} x(a+a' )\\ y(a+a') \end{matrix}\right\}+\left\{\begin{matrix} xa \\ xa' \end{matrix}\right\}+\left[\left\{\begin{matrix} (x+y)a' \\ ya \end{matrix}\right\}\right]- \left\{\begin{matrix} x(a+a') \\ y(a+a') \end{matrix}\right\}\\
&=-\left\{\begin{matrix} xa\\ ya \end{matrix}\right\}- \left\{\begin{matrix} y(a+a') \\ xa' \end{matrix}\right\}+\left\{\begin{matrix} xa\\ xa' \end{matrix}\right\}-\left\{\begin{matrix} xa' \\ ya' \end{matrix}\right\}+\left\{\begin{matrix} xa' \\ y(a+a') \end{matrix}\right\}+ \left\{\begin{matrix} ya \\ ya' \end{matrix}\right\}\\
&= -\left\{\begin{matrix} xa \\ ya \end{matrix}\right\}+\left\{\begin{matrix} xa \\ xa' \end{matrix}\right\}-\left\{\begin{matrix} xa' \\ ya' \end{matrix} \right\}+\left\{\begin{matrix} ya \\ ya' \end{matrix}\right\}.
\end{split}\end{equation}
Thus only the relations (\ref{GR_relations}) are involved in $\Omega_{-}^{\NN}, \Omega_{-}^{\ZZ}$. The universal odd derivation $d_X^{-}:\NN^X\rightarrow (\Omega_{-}^{\NN})_X$, is trivial for $X=[1]$, but is non-trivial for $|X|>1$. Letting $\mathcal{N}_X^{\NN}$ (resp. $\mathcal{N}_X^{\ZZ}$) denote the free abelian group with generators $[a], a\in \NN^X$ (resp. $\ZZ^X$), modulo the relation $[k\cdot a]=k\cdot [a], k\in\NN$ (resp. $\ZZ$), we have the exact sequence
\begin{equation}
\xymatrix{(\Omega_{-}^{\NN})_X\ar[r]^{\partial}\ar@{^{(}->}[d]&\mathcal{N}_X^{\NN}\ar[r]^{\pi}\ar@{^{(}->}[d]&\cG(\ZZ)_X\ar@{=}[d]\ar[r]&0 \\
(\Omega_{-}^{\ZZ})_X\ar[r]^{\partial}&\mathcal{N}_X^{\ZZ}\ar[r]^{\pi}&\cG(\ZZ)_X\ar[r]&0}
\end{equation}
with $\partial \left\{ \begin{matrix} a\\ a' \end{matrix} \right\}=[a+a']-[a]-[a'], \pi(\sum_i k_i[a_i])=\sum_i k_i\cdot a_i$. \vspace{3mm}\\
The $+1$- almost- linearity, is thus equivalent modulo 2 torsion (i.e. after taking $\underset{\ZZ}{\otimes} \ZZ[\frac{1}{2}]$), to left and right \emph{linearity}: $a,a'\in \NN^X$ (resp. $\ZZ^X$), $x,y\in \NN$ (resp. $\ZZ$)
\begin{equation}\label{rightlinear}
right:  \left\{\begin{matrix} x(a+a')\\ y(a+a')\end{matrix}\right\}= \left\{\begin{matrix} xa\\ ya \end{matrix}\right\}+
 \left\{\begin{matrix} xa'\\ ya'\end{matrix}\right\}\hspace{25mm}
\end{equation}
\begin{equation}\label{leftlinear}
left: \left\{\begin{matrix} (x+y)a\\ (x+y)a' \end{matrix}\right\}=\left\{\begin{matrix} xa\\ xa'\end{matrix}\right\}+
\left\{\begin{matrix} ya\\ ya'\end{matrix}\right\}\iff \left\{\begin{matrix} ka\\ ka'\end{matrix}\right\}
=k\left\{\begin{matrix} a'\\ a'\end{matrix}\right\}
\end{equation}
Conversely, left and right lineary (\ref{rightlinear}-\ref{leftlinear}) imply $+$ and $-$ almost linearity (\ref{almostlinear}). For $X=[1]$, the right linearity conditions (\ref{rightlinear}) and the left one (\ref{leftlinear}) are one and the same ! \\
Thus letting $\bar{\Omega}_X^{\NN}$ (resp. $\bar{\Omega}_X^{\ZZ}$) denote the free abelian group with generators $\left\{\begin{matrix} a \\ a' \end{matrix} \right\},a,a'\in \NN^X$ (resp. $\ZZ^X$), modulo the relations (\ref{GR_relations}) and (\ref{rightlinear}-\ref{leftlinear}), we get the even and odd derivation
\begin{equation}
\begin{split}
d_{[X]}^{\pm}:\NN^X\rightarrow \bar{\Omega}_{[X]}^{\NN} \\
d^{\pm}_{[X]}:\ZZ^X\rightarrow \bar{\Omega}_{[X]}^{\ZZ}
\end{split}
\end{equation}
In particular taking $X=[1]$, we get the even differential
\begin{equation}
\begin{split}
d_{[1]}^+:\NN\rightarrow \bar{\Omega}_{[1]}^{\NN} \\
\text{resp. } d^+_{[1]}:\ZZ\rightarrow \bar{\Omega}_{[1]}^{\ZZ},
\end{split}
\end{equation}
It satisfies:
\begin{equation}
d^{+}_{[1]}(0)=d^{+}_{[1]}(1)=0
\end{equation}
For $n\geq 2$,
\begin{equation}
\begin{split}
&d^+_{[1]}(n)= 2\cdot (\left\{\begin{matrix} n-1 \\ 1\end{matrix}\right\}+\left\{\begin{matrix} n-2 \\ 1\end{matrix}\right\}+\dots +\left\{\begin{matrix} 1 \\1 \end{matrix}\right\}) \\
\text{resp. }& d^+_{[1]}(-n)= 2\cdot (\left\{\begin{matrix} 1-n \\ -1\end{matrix}\right\} + \left\{\begin{matrix} 2-n \\ -1\end{matrix}\right\}+\dots +\left\{\begin{matrix} -1\\ -1 \end{matrix}\right\})=-d^+_{[1]}(n).
\end{split}
\end{equation}
(use induction on $n$, and apply $d_{[1]}^+$ to the identity $(n+1)=((1,1)\lhd \{n,1\})\sslash (1,1)$ for the induction step)
\begin{equation}
\emph{Leibnitz} \; d_{[1]}^+(n\cdot m)=m\cdot d_{[1]}^+(n)+n\cdot d^+_{[1]}(m)
\end{equation}
It follows that
\begin{equation}
d^+_{[1]}(q^n)=n\cdot q^{n-1}\cdot d^+_{[1]}(q)
\end{equation}
and hence
\begin{equation}
d^+_{[1]}(n)= \sum_p v_p(n)\frac{n}{p} d^+_{[1]}(p)
\end{equation}
After extending scalars to $\QQ$ this is equivalent to
\begin{equation}
\frac{d^+_{[1]}(n)}{n}=\sum_p v_p(n)\frac{d^+_{[1]}(p)}{p}
\end{equation}
This is the arithmetical analogue of the formula for $f=f(z)\in \CC(z)$,
\begin{equation*}
\frac{df}{f}=d \log f=\sum_{\alpha\in \CC} v_{\alpha}(f)\frac{d(z-\alpha)}{z-\alpha}=\sum_{\alpha\in \CC}\frac{v_{\alpha}(f)}{z-\alpha}dz
\end{equation*}
The derivation $d^+_{[1]}$ is not additive, but we do have the identity
\begin{equation}
d_{[1]}^+(n_1+n_2)=d_{[1]}(n_1)+d^+_{[1]}(n_2)+2\left\{\begin{matrix}n_1 \\ n_2\end{matrix}\right\}.
\end{equation}
(apply $d^+_{[1]}$ to the identity $(n_1+n_2)=((1,1)\lhd (n_i))\sslash (1,1)$ and use Leibnitz). \\
To see that $\bar{\Omega}_{[1]}^{\NN}$ (resp. $\bar{\Omega}_{[1]}^{\ZZ}$) is non- trivial, note that for each prime $p$ we have a homomorphism $\varphi_p$ from it onto $\ZZ$, given on the generators $\left\{\begin{matrix} a \\ a'\end{matrix}\right\}$ by
\begin{equation}
\begin{split}
\varphi_p(\left\{\begin{matrix}a \\ a' \end{matrix}\right\})= v_p(a+a')\cdot \frac{a+a'}{p}-v_p(a)\cdot \frac{a}{p}
-v_p(a')\cdot \frac{a'}{p} \\
(v_p(p^n\cdot a)=n\; \text{for } p\nmid a, \text{ the $p$- adic valuation}).
\end{split}
\end{equation}
(Indeed, $\bar{\Omega}_{[1]}^{\NN}$ is the free abelian group generated by $\left\{\begin{matrix}1\\p-1\end{matrix} \right\}$, $p$ prime).\\
The generalized ring $H=\cG(\ZZ)\prod \bar{\Omega}^{\ZZ}$, is an
(absolute) abelian group via the maps
\begin{equation}
\begin{split}
H_Z\times H_Z&\rightarrow H_Z \\
(a_1,m_1)+(a_2,m_2)&:=(a_1+a_2,\left\{\begin{matrix} a_1\\ a_2\end{matrix} \right\}+m_1+m_2)
\end{split}
\end{equation}
and we get an exact sequence
\begin{equation}
0\rightarrow \bar{\Omega}(\ZZ)\rightarrow H \rightarrow \cG(\ZZ)\rightarrow 0.
\end{equation}

\let\cleardoublepage\clearpage
\numberwithin{equation}{chapter}
\begin{appendices}
\chapter{Beta integrals and the local factors of zeta}
\setcounter{equation}{0}
We shall concentrate on the case of the rational numbers $\QQ$. We denote by $p$ the close points of $\varprojlim(\bar{\Spec \ZZ})$, that is the finite primes (denoted by "$p\neq \eta $"), and the real prime (denoted by "$p=\eta$"); when we want to emphasise that a formula holds for \textbf{all} primes, finite or real, we write "$p\geq \eta$". For each $p\geq \eta$ we have the completion $\QQ_{p}$, the $p$- adic numbers for $p\neq \eta$, $\QQ_{\eta}=\RR$ the reals, and we have the (maximal- compact) sub- generalized ring $\mathcal{G}(\ZZ_{p})\subseteq \mathcal{G}(\QQ_{p})$, with $\mathcal{G}(\ZZ_{\eta})=\ZZ_{\eta}$ the real prime (cf. \S \ref{generealprimes}). We let for $p\geq \eta$, 
\begin{equation}
S_{p}^n=\{(x_1,...,x_n)\in \mathcal{G}(\ZZ_{p})_{[n]}, \;\;\; |x_1,\dots ,x_n|_p=1\}
\end{equation}
where,
\begin{equation}
|x_1,\dots,x_n|_{p}=
\begin{cases}
  \Max\{|x_1|_{p},\dots, |x_n|_{p}\}& p\neq \eta \\
   (|x_1|^2+\dots +|x_n|^2)^{1/2} & p= \eta
  \end{cases}
\end{equation}
(cf. \S \ref{2.5}). Thus for $p=\eta$ we have the $(n-1)$ dimensional sphere (while for $p\neq \eta$, we get an open subset of $n$ dimensional space). Note that for all $p\geq \eta$ (with $\ZZ_{\eta}^*=\{\pm 1\}$):
\begin{equation}
S_{p}^n\slash \ZZ_{p}^* \cong \mathbb{P}^{n-1}(\ZZ_{p})\cong \mathbb{P}^{n-1}(\QQ_{p})
\end{equation}
For $p\neq \eta$, 
\begin{equation}
S_{p}^n=\varprojlim_k S^n(\ZZ\slash p^k), \;\;\;\; \mathbb{P}^{n-1}(\ZZ_{p})=\varprojlim_k \mathbb{P}^{n-1}(\ZZ\slash p^k).
\end{equation}
For all $p\geq \eta$, $S_p^n$ is a homogenuous space of the compact group $GL_n(\ZZ_p)$ (where $GL_n(\ZZ_{\eta})=O(n)$ the orthogonal group). \\
For all $p\geq \eta$, we denote by $\sigma_p^n$ the \textbf{unique} $GL_n(\ZZ_p)$- invariant probability measure on $S_p^n$. We write the local factors of the zeta function as 
\begin{equation}
\zeta_p(s):=\begin{cases}
(1-p^{-s})^{-1}& p\neq \eta \\
2^{\frac{s}{2}}\cdot \Gamma(\frac{s}{2})& p=\eta
\end{cases}
\end{equation}   
(so the global zeta, $\zeta_{\mathbb{A}}(s)=\prod_{p\geq \eta} \zeta_p(s), \; \Re (s)>1$, satisfies the functional equation $\zeta_{\mathbb{A}}(s)=(2\pi)^{s-\frac{1}{2}}\zeta_{\mathbb{A}}(1-s)$). We shall eventually recover these local factors $\zeta_p(s), p\geq \eta$, in a uniform way. We denote for $p\geq \eta$, the Beta function
\begin{equation}
\beta_p(\alpha_1,...,\alpha_n)=\frac{\zeta_p(\alpha_1)\cdots \zeta_p(\alpha_n)}{\zeta_p(\alpha_1+\cdots+\alpha_n)}=
\begin{cases}
\frac{1-p^{-(\alpha_1+\dots+\alpha_n)}}{(1-p^{-\alpha_1})\cdots (1-p^{-\alpha_n})}& p\neq \eta \\[2em]
\frac{\Gamma(\frac{\alpha_1}{2})\cdots \Gamma(\frac{\alpha_n}{2})}{\Gamma(\frac{\alpha_1+\dots+\alpha_n}{2})} & p=\eta
\end{cases}
\end{equation}
and the normalized Beta function
\begin{equation}
B_p(\alpha_1,\dots,\alpha_n)=\frac{\beta_p(\alpha_1,\dots,\alpha_n)}{\beta_p(1,\dots,1)}=\frac{\zeta_p(n)}{\zeta_p(1)^n}\cdot \frac{\zeta_p(\alpha_1)\dots \zeta_p(\alpha_n)}{\zeta_p(\alpha_1+\cdots +\alpha_n)} .
\end{equation}
We have for all $p\geq \eta$ the following Beta- integral:
\begin{equation}\label{beint}
\textbf{Beta-$\int$:}\hspace{10mm}\int_{S_p^n} |x_1|_p^{\alpha_1-1}\cdots |x_n|_p^{\alpha_n-1} \sigma_p^n(dx)=  B_p(\alpha_1,\dots,\alpha_n), \;\; \Re(\alpha_i)>0.
\end{equation}
here $x=(x_1,\dots,x_n)\in S_p^n$. \\
(This can be verified directly for $p\neq \eta$, and for $p=\eta$, or can be obtained as appropiate limits of a $q$- analogue. cf. \cite{H08}.) \\
Note that for real $\alpha_i>0$, we get a probability measure, the "Beta- measure" 
\begin{equation}
\mu_p^{(\alpha_1\dots \alpha_n)}(dx)=|x_1|_p^{\alpha_1-1}\cdots |x_n|_p^{\alpha_n-1}\frac{\sigma_p(dx)}{B_p(\alpha_1,\dots, \alpha_n)}
\end{equation}
on $S_p^n$, it is $\ZZ_p^*$- invariant, hence can be viewed as a probability measure on $\mathbb{P}^{n-1}(\QQ_p)$ (hence for $p\neq \eta$, we get a Markov chain on $\coprod_k \mathbb{P}^{n-1}(\ZZ\slash p^k)$, cf. \cite{H08} for the $q$- analogue, which also gives the real analogue). \\
On the other hand we have for all $p\geq \eta$ , and all $y=(y_1,\dots, y_n)\in \QQ_p^n$:  
\begin{equation}\label{sslashint}\begin{split}
\mathbf{\sslash\text{-}\int:}  \hspace{10mm}&\int_{S_p^n} |x_1y_1+\dots x_n y_x|_p^{s-1} \sigma_p^n(dx)=\int_{S_p^n} |x\sslash y|_p^{s-1} \sigma_p^n(dx) \\ &=\frac{\zeta_p(n)}{\zeta_p(1)}\frac{\zeta_p(s)}{\zeta_p(n-1+s)}|y|_p^{s-1}. 
\end{split}\end{equation}
Note that the operations of multiplications and contraction are very natural for the Beta measures: we have for $N=n_1+\dots + n_k,\; p\geq \eta$, a surjection
\begin{equation}\begin{split}
\lhd: &S_p^k\times S_p^{n_1}\times \cdots \times S_p^{n_k}\sur S_p^N \\
&(t,x^{(1)},\dots, x^{(k)})\mapsto t\lhd (x^{(i)})  
\end{split}\end{equation}
and the measure $\sigma^N$ is obtained as the image of the measures $\sigma^{n_j}$ and the Beta measure $\mu_p^{(n_1,\dots , n_k)}$,  
\begin{equation}\label{multformula}\begin{split}
&\textbf{multiplication formula:}  \\
&\int_{S_p^N} f(x)\sigma_p^N(dx) =\\ 
&\int_{S_p^k}\sigma_p^k(dt)\cdot \frac{|t_1|^{n_1-1}\dots |t_k|^{n_k-1}}{B_p(n_1,\dots,n_k)}\int_{S_p^{n_1}} \sigma_p^{n_1}(dx^{(1)}) \dots \int_{S_p^{n_k}} \sigma_p^{n_k} (dx^{(n)})\hsm f(t\lhd (x^{(i)}))
\end{split}\end{equation}
Writing $x=(x_1,\dots, x_k), y=(y_1,\dots, y_k)\in \QQ_p^N $, with $x_j,y_j\in \QQ_p^{n_j}$, we have for $p\geq \eta, \; N=n_1+\dots+n_k $ : 
\begin{equation}\label{Contint}\begin{split}
\textbf{Contraction $\int$ :}\hspace{10mm}  \int_{S_p^N} |x_1\sslash y_1|_p^{\alpha_1-1}\dots |x_k\sslash y_k|_p^{\alpha_k-1}\sigma_N(dx)= \\
\frac{\zeta_p(N)\zeta_p(\alpha_1)\dots \zeta_p(\alpha_k)}{\zeta_p(1)^k\cdot \zeta_p(N-k+(\alpha_1+\dots+\alpha_k))}|y_1|_p^{\alpha_1-1}\dots |y_k|_p^{\alpha_k-1}
\end{split}\end{equation}
Note that ($\ref{Contint}$) for $k=1$ is ($\ref{sslashint}$) and ($\ref{Contint}$) for $n_1=\dots=n_N=1$, and $y_j\equiv 1$ is ($\ref{beint}$), and conversly, ($\ref{Contint}$) follows form ($\ref{beint}$), ($\ref{sslashint}$) and the multiplication formula, ($\ref{multformula}$) . \\  
Taking the vectors $y=\mathbbm{1}_N=(1,1,\dots,1)\in \QQ_p^N$, we get a probability measure $\varphi_p^N=\sigma_p^N\sslash \mathbbm{1}_N$ on $\QQ_p$ for all $p\geq \eta$: 
\begin{equation}
\int_{\QQ_p} f(x) \varphi_p^N(dx)= \int_{S_p^N} f(x_1+\dots + x_N) \sigma_p^N(dx)
\end{equation}
We have: 
\begin{equation}
\varphi_p^N(dx)=\bigg [\frac{1-p^{1-N}}{1-p^{-N}}\phi_{\tiny \ZZ_p}(x)+\frac{p^{1-N}}{1-p^{-N}}\phi_{\tiny \ZZ_p^*}(x)\bigg ]dx. \;\; p\neq \eta,
\end{equation}
here $dx$ is the additive Haar measure on $\ZZ_p, dx(\ZZ_p)=1$, and $\phi_{\tiny \ZZ_p}$ (resp $\phi_{\tiny \ZZ_p^*}$) is the characteristic function of $\ZZ_p$ (resp. $\ZZ_p^*$)  ; 
\begin{equation}
\varphi_{\tiny \eta}^N(dx)=\frac{\Gamma(\frac{N}{2})}{\sqrt{\pi N}\cdot \Gamma(\frac{N-1}{2})}\cdot (1-\frac{|x|^2}{N})^{\frac{N-1}{2}-1}dx, \;\; x\in [-\sqrt{N},\sqrt{N}], 
\end{equation} 
with the usual Haar measure $dx$ on $\RR$, $dx([0,1])=1$. \\
When we take the limit $N\rightarrow \infty$ we obtain the measures 
\begin{equation}\begin{split}
&\varphi_p^{\infty}(dx) =\phi_{\tiny \ZZ_p} (x) dx, \;\; p\neq \eta, \\
&\varphi_\eta^{\infty}(dx)=e^{\tiny -\frac{|x|^2}{2}}\frac{dx}{\sqrt{2\pi}}. \;\; p=\eta,\;\;  dx([0,1])=1,
\end{split}\end{equation}
and for all $p\geq \eta$ we have from (\ref{sslashint}), with $y_i=1$, in the limit $n\rightarrow \infty$: 
\begin{equation}
\frac{\zeta_p(s)}{\zeta_p(1)}=\lim_{N\rightarrow \infty} \int_{\QQ_p} |x|_p^{s-1}\varphi_p^N(dx)=\lim_{N\rightarrow \infty} \int_{S_p^N} |x_1+\dots+x_N|_p^{s-1}\sigma_p^N(dx).
\end{equation}
\end{appendices}


\end{document}